\documentclass[a4paper, 11pt]{report}

\usepackage[T1]{fontenc}
\usepackage[francais]{babel}
\usepackage[latin1]{inputenc}  
\usepackage{amsfonts}
\usepackage{amsmath}
\usepackage{amssymb}
\usepackage{bbm}
\usepackage{tikz}
\usepackage{pgflibrarysnakes}
\usepackage[all]{xy}
\usepackage{graphicx}
\usepackage{cancel}
\usepackage{multirow}
\usepackage{url}

\begin{document}

\begin{titlepage}
\begin{center}
\vspace{-4cm}\includegraphics[height=1cm]{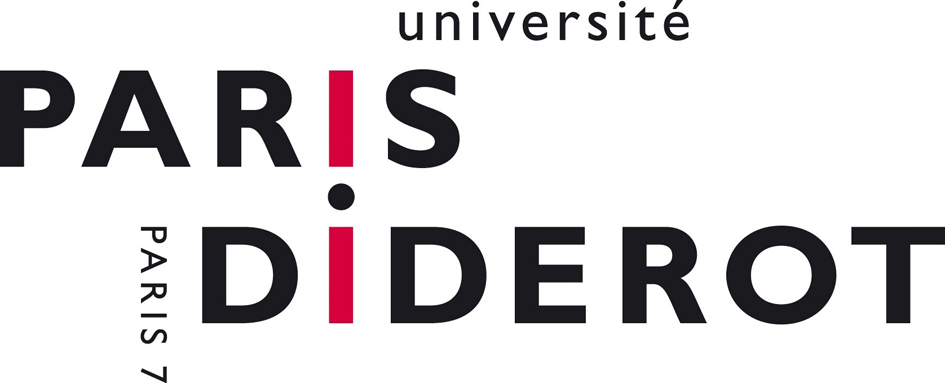} \hspace{7cm} \includegraphics[height=1.5cm]{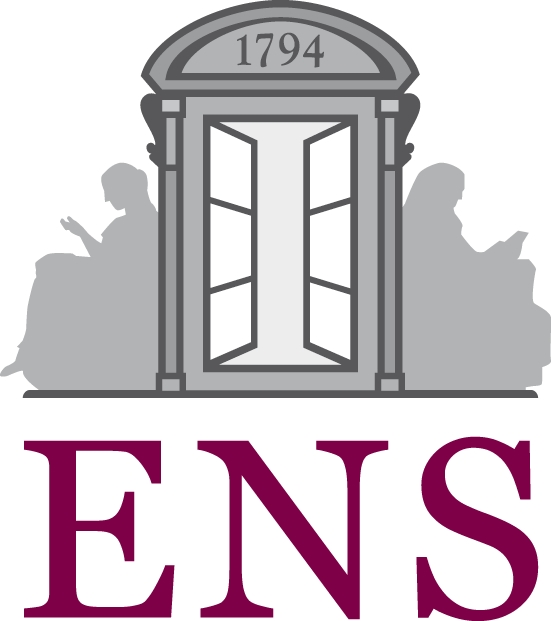} \\[2cm]
\textsc{\large UNIVERSIT\' E PARIS DIDEROT (PARIS 7)} \\[0.5cm]
\textsc{\large Sciences mathématiques de Paris Centre (ED386)} \\[1cm]

\textsc{\LARGE TH\`ESE DE DOCTORAT} \\[0.5cm]
\textsf{\large Discipline : Mathématiques} \\[1.5cm]

présentée par \\[0.5cm]
\LARGE Tony \textsc{Ly} \\[1cm]

% Title
\hrule ~\\[0.2cm]
{\LARGE \bfseries Représentations modulo $p$ de $\mathrm{GL}(m,D)$, $D$
  algèbre à division sur un corps local} \\[0.4cm]

\hrule ~\\[1cm]
\normalsize dirigée par Marie-France \textsc{Vignéras}
\end{center}

\vfill

% Bottom of the page
\small
Soutenue le 22 novembre 2013 devant le jury composé de: \\
\begin{tabular}{lllll}
$\qquad$ & & & & \\
& M. & Ioan \textsc{Badulescu} & Université Montpellier II & Examinateur \\
& M. & Laurent \textsc{Berger} & E.N.S. de Lyon & Examinateur \\
& M. & Elmar \textsc{Grosse-Klönne} & Humboldt-Universität zu Berlin & Examinateur \\
& M. & Guy \textsc{Henniart} & Université Paris-Sud 11 & Rapporteur \\
& $\textrm{M}^{\textrm{me}}$ & Ariane \textsc{Mézard} & Université Pierre et Marie Curie 6 & Examinatrice \\
& M. & Vincent \textsc{Sécherre} & Université de Versailles Saint-Quentin & Examinateur \\
& $\textrm{M}^{\textrm{me}}$ & Marie-France \textsc{Vignéras} & Université Paris Diderot 7 & Directrice de thèse 
\end{tabular} \\

\vspace{0.8cm}

Rapporteur absent lors de la soutenance: \\
\begin{tabular}{lllll}
& \phantom{$\textrm{M}^{\textrm{me}}$} & & \\ 
$\qquad$ & M. & Florian \textsc{Herzig} $\qquad \quad \ \, $ & University of Toronto & $\phantom{,}$
\end{tabular}

\end{titlepage}

\newtheorem{lem}{Lemme}[section]
\newtheorem{ax}[lem]{Axiome}
\newtheorem{defi}[lem]{Définition}
\newtheorem{prop}[lem]{Proposition}
\newtheorem{theo}[lem]{Théorème}
\newtheorem{cor}[lem]{Corollaire}
\newtheorem{por}[lem]{Porisme}
\newtheorem{ex}[lem]{Exemple}
\newtheorem{meth}[lem]{Méthode}
\newtheorem{err}{Erreur}
\newtheorem{conj}{Conjecture}
\newtheorem{ques}[conj]{Question}
\newtheorem{hypo}{Hypothèse}
\newtheorem{hypo*}{Hypothèse}

\newtheorem{Alem}{Lemma}[section]
\newtheorem{Aax}[Alem]{Axiom}
\newtheorem{Adefi}[Alem]{Definition}
\newtheorem{Aprop}[Alem]{Proposition}
\newtheorem{Atheo}[Alem]{Theorem}
\newtheorem{Acor}[Alem]{Corollary}
\newtheorem{Apor}[Alem]{Porism}
\newtheorem{Aex}[Alem]{Example}
\newtheorem{Ameth}[Alem]{Method}
\newtheorem{Aerr}{Error}

\newcommand{\A}{\mathcal{A}}
\newcommand{\Aa}{\mathbf{A}}
\newcommand{\AAa}{\mathbb{A}}
\newcommand{\AAA}{\mathfrak{A}}
\newcommand{\aaa}{\mathfrak{a}}
\newcommand{\Ab}{\mathfrak{Ab}}
\newcommand{\ab}{\mathrm{ab}}
\newcommand{\Ad}{\mathrm{Ad}}
\newcommand{\adm}{\mathrm{adm}}
\newcommand{\al}{\mathrm{al}}
\newcommand{\Ann}{\mathrm{Ann} \,}
\newcommand{\Ass}{\mathrm{Ass}}
\newcommand{\Aut}{\mathrm{Aut}}
\newcommand{\AV}{\mathfrak{Av}}
\newcommand{\B}{\mathcal{B}}
\newcommand{\bbb}{\mathfrak{b}}
\newcommand{\BE}{B_{(E)}}
\newcommand{\C}{\mathbb{C}}
\newcommand{\Cc}{\mathcal{C}}
\newcommand{\Card}{\mathrm{Card}}
\newcommand{\cd}{\mathrm{cd} \,}
\newcommand{\ch}{\mathrm{ch}}
\newcommand{\Char}{\mathrm{char}}
\newcommand{\Coker}{\mathrm{Coker}}
\newcommand{\cyc}{\mathrm{cyc}}
\newcommand{\D}{\mathrm{D}}
\newcommand{\Dd}{\mathcal{D}}
\newcommand{\DD}{\mathbf{D}}
\newcommand{\dd}{\mathrm{d}}
\newcommand{\Dlim}{\mathop{D}\limits_\to}
\newcommand{\Diag}{\mathrm{Diag}}
\newcommand{\Dir}{\mathrm{Dir}}
\newcommand{\Div}{\mathrm{div}}
\newcommand{\DR}{\mathrm{DR}}
\newcommand{\E}{\mathcal{E}}
\newcommand{\Ee}{\mathbf{E}}
\newcommand{\End}{\mathrm{End}}
\newcommand{\et}{\mathrm{et}}
\newcommand{\etd}{\mathrm{etd}}
\newcommand{\Ext}{\mathrm{Ext}}
\newcommand{\Extc}{\mathcal{E} \! xt}
\newcommand{\F}{\mathbb{F}}
\newcommand{\Ff}{\mathcal{F}}
\newcommand{\Faisc}{\mathfrak{Faisceaux}}
\newcommand{\Fct}{\mathrm{Fct}}
\newcommand{\fd}{\mathrm{fd}}
\newcommand{\Fil}{\mathrm{Fil}}
\newcommand{\Fin}{\mathrm{Fin}}
\newcommand{\Fix}{\mathrm{Fix}}
\newcommand{\Fp}{\mathbb{F}_p}
\newcommand{\Fpbar}{\overline{\mathbb{F}}_p}
\newcommand{\Frac}{\mathrm{Frac}}
\newcommand{\G}{\mathcal{G}}
\newcommand{\Ga}{\mathbb{G}_a}
\newcommand{\GE}{G_{(E)}}
\newcommand{\Gm}{\mathbb{G}_m}
\newcommand{\Gal}{\mathrm{Gal}}
\newcommand{\gen}{\mathrm{gen}}
\newcommand{\GL}{\mathrm{GL}}
\newcommand{\gld}{\mathrm{gld} \,}
\newcommand{\Gr}{\mathrm{Gr}}
\newcommand{\gr}{\mathrm{gr} \,}
\newcommand{\grr}{\mathrm{gr}}
\newcommand{\Hh}{\mathcal{H}}
\newcommand{\hh}{\mathfrak{h}}
\newcommand{\Hom}{\mathrm{Hom}}
\newcommand{\I}{\mathcal{I}}
\newcommand{\It}{\widetilde{I}}
\newcommand{\id}{\mathbbm{1}}
\newcommand{\Id}{\mathrm{id}}
\newcommand{\IE}{I_{(E)}}
\newcommand{\im}{\mathrm{Im}}
\newcommand{\ind}{\mathrm{ind}}
\newcommand{\Ind}{\mathrm{Ind}}
\newcommand{\J}{\mathcal{J}}
\newcommand{\Jac}{\mathrm{Jac}}
\newcommand{\K}{\mathcal{K}}
\newcommand{\Kbar}{\overline{K}}
\newcommand{\Kt}{\widetilde{K}}
\newcommand{\KE}{K_{(E)}}
\newcommand{\Ker}{\mathrm{Ker}}
\newcommand{\Ll}{\mathcal{L}}
\newcommand{\LL}{\mathbb{L}}
\newcommand{\LA}{\mathcal{LA}}
\newcommand{\lgld}{\mathrm{lgld} \,}
\newcommand{\Lie}{\mathrm{Lie}}
\newcommand{\Long}{\mathrm{long}}
\newcommand{\M}{\mathfrak{M}}
\newcommand{\Mm}{\mathcal{M}}
\newcommand{\m}{\mathfrak{m}}
\newcommand{\mapsfrom}{\mathrel{\reflectbox{\ensuremath{\mapsto}}}}
\newcommand{\Mc}{\mathrm{Mc}}
\newcommand{\mdt}{\textsl{MdT}}
\newcommand{\Modd}{\mathfrak{Mod}}
\newcommand{\Mod}[1]{\mathfrak{Mod} \, \Hh_{\overline{\F}_p}(G,#1)}
\newcommand{\Mpel}[1]{\mathcal{M}^\mathrm{PEL}_{#1}}
\newcommand{\N}{\mathbb{N}}
\newcommand{\Nn}{\mathcal{N}}
\newcommand{\NN}{\mathfrak{N}}
\newcommand{\nr}{\mathrm{nr}}
\newcommand{\Nrm}{\mathrm{Nrm}}
\newcommand{\ns}{\mathrm{ns}}
\newcommand{\Oo}{\mathcal{O}}
\newcommand{\OO}{\mathrm{O}}
\newcommand{\op}{\mathrm{op}}
\newcommand{\ord}{\mathrm{ord}}
\newcommand{\Ord}{\mathrm{Ord}}
\newcommand{\Pp}{\mathcal{P}}
\newcommand{\PP}{\mathfrak{P}}
\newcommand{\PPp}{\mathbf{P}}
\newcommand{\PPP}{\mathbb{P}}
\newcommand{\pp}{\mathfrak{p}}
\newcommand{\Par}{\mathrm{Par}}
\newcommand{\pc}{\mathrm{pc}}
\newcommand{\pd}{\mathrm{pd}}
\newcommand{\PGL}{\mathrm{PGL}}
\newcommand{\pr}{\mathrm{pr}}
\newcommand{\Phired}{\Phi_{\mathrm{red}}}
\newcommand{\Phibarred}{\overline{\Phi}_{\mathrm{red}}}
\newcommand{\Pic}{\mathrm{Pic}}
\newcommand{\PSL}{\mathrm{PSL}}
\newcommand{\Q}{\mathbb{Q}}
\newcommand{\Qq}{\mathcal{Q}}
\newcommand{\q}{\mathfrak{q}}
\newcommand{\Qp}{\mathbb{Q}_p}
\newcommand{\R}{\mathbb{R}}
\newcommand{\Rr}{\mathcal{R}}
\newcommand{\Rad}{\mathrm{Rad}}
\newcommand{\ram}{\mathrm{ram}}
\newcommand{\red}{\mathrm{red}}
\newcommand{\Rep}{\mathfrak{Rep}_{\overline{\F}_p} \, G}
\newcommand{\Res}{\mathrm{Res}}
\newcommand{\rg}{\mathrm{rg}}
\newcommand{\rgld}{\mathrm{rgld} \,}
\newcommand{\RHom}{\mathbb{R}\mathrm{Hom}}
\newcommand{\Ss}{\mathfrak{S}}
\newcommand{\Sss}{\mathcal{S}}
\newcommand{\sep}{\mathrm{sep}}
\newcommand{\Sh}{\mathrm{Sh}}
\newcommand{\SL}{\mathrm{SL}}
\newcommand{\sm}{\mathrm{sm}}
\newcommand{\Spec}{\mathrm{Spec} \,}
\newcommand{\Smax}{\mathrm{Spec}_\mathrm{max} \,}
\newcommand{\sss}{\mathrm{ss}}
\newcommand{\St}{\mathrm{St}}
\newcommand{\Stab}{\mathrm{Stab}}
\newcommand{\Sub}{\mathrm{Sub}}
\newcommand{\Supp}{\mathrm{Supp}}
\newcommand{\Sym}{\mathrm{Sym}}
\newcommand{\T}{\mathcal{T}}
\newcommand{\Tgt}{\mathrm{Tgt}}
\newcommand{\tor}{\mathrm{tor}}
\newcommand{\Tor}{\mathrm{Tor}}
\newcommand{\tors}{\mathrm{tors}}
\newcommand{\Tot}{\mathrm{Tot}}
\newcommand{\tr}{\mathrm{tr}}
\newcommand{\Tr}{\mathrm{Tr}}
\newcommand{\uu}{\mathfrak{u}}
\newcommand{\U}{\mathcal{U}}
\newcommand{\UE}{U_{(E)}}
\newcommand{\un}{\mathrm{un}}
\newcommand{\V}{\mathcal{V}}
\newcommand{\VV}{\mathbf{V}}
\newcommand{\Vect}{\mathrm{Vect}}
\newcommand{\W}{\mathrm{W}}
\newcommand{\Ww}{\mathcal{W}}
\newcommand{\X}{\mathcal{X}}
\newcommand{\Z}{\mathbb{Z}}
\newcommand{\Zz}{\mathcal{Z}}
\newcommand{\zp}{\mathcal{Z}_p}
\newcommand{\Zp}{\mathbb{Z}_p}

\pagestyle{empty}
\setcounter{page}{0}

\newpage
$\phantom{,}$
\newpage

\pagestyle{plain}
{\LARGE \textbf{Remerciements}} \\
\vspace{1cm}

C'est un plaisir de commencer par remercier celle qui aura été pour moi une directrice de thèse idéale, Marie-France Vignéras. Elle a su me laisser la liberté de trouver mes marques et mes envies, tout en faisant partager sa motivation et son enthousiasme communicatif pour chacune de mes avancées. Enfin, je ne peux m'empêcher de souligner l'apport qu'auront constituées ses questions et critiques pour la préparation de ce manuscrit, ainsi que la générosité avec laquelle elle m'a partagé son expérience mathématique. \\
Je suis aussi infiniment reconnaissant envers Guy Henniart et Florian Herzig d'avoir accepté de relire ce document. Leurs remarques pertinentes ont contribué à l'amélioration de ce texte et je les en remercie: ils ont accompli cette tâche disgracieuse en portant une attention minutieuse à mon travail. \\
Je suis honoré que Ioan Badulescu, Laurent Berger, Elmar Grosse-Klönne, Ariane Mézard et Vincent Sécherre fassent partie de mon jury. Leur importante contribution au domaine n'est plus à démontrer et leur présence à ma soutenance est une grande gratification. \\ 

Je n'oublie pas que ma formation mathématique a débuté il y a maintenant de longues années déjà: une grande influence est restée de ce que m'ont transmis Nicolas Tosel, Rachel Ollivier et Alex Paulin, et ils en sont remerciés. \\
Une partie non négligeable des mathématiques que j'ai apprises a bénéficié du cadre propice de groupes de travail: à ce titre, je voudrais adresser un merci à Shenghao Sun et ses compères du \textit{Shimura varieties subseminar}, à Stefano Morra, Gabriel Dospinescu, Benjamin Schraen, Raphaël Beuzart et Arthur-César Le Bras. \\
Il m'a été très agréable de discuter avec ceux qui ont toujours porté un grand intérêt à mon travail: Yongquan Hu, Alberto Minguez et Benjamin Schraen font partie de ceux-là. Le mathématicien est un être social et c'est un plaisir de compter aussi Ramla Abdellatif, Alexis Bouthier, Eugen Hellman, Arthur Laurent, Jean-Baptiste Teyssier, Paul-James White parmi mes amis du métier. \\
Enfin je ne saurais trop apprécier le cadre agréable dans lequel j'ai pu mener à bien mon travail: l'ambiance des toits et la chaleur apportée par les membres du DMA (avec mention spéciale à la bibliothèque et au secrétariat) passés et présents n'y sont pas étrangères. \\

Pour finir, un petit mot mais un grand merci à ma famille et mes amis, à ceux que je ne nommerai pas mais à qui je pense aussi lorsque je me remémore ces années de thèse, qui ont su égayer mon c\oe ur et mon esprit de leur présence inconditionnelle.

\newpage 
$\phantom{,}$
\newpage
\tableofcontents
\newpage

\chapter*{Introduction}
\addcontentsline{toc}{chapter}{Introduction}

\section{Motivation historique}

Soient $F$ un corps local et $G$ un groupe réductif connexe sur $F$. Soit $R$ un corps de coefficients algébriquement clos, sur lequel on prendra nos représentations. \\ En $1967$, une lettre de Robert Langlands à André Weil (\cite{Lan67}) ouvre la voie à de nombreuses questions qui constituent un gigantesque domaine de recherche actuel avec l'étiquette \og programme de Langlands \fg. Ces questions couvrent en particulier les conjectures de Langlands locale dont la philosophie est la suivante: il existe une bijection naturelle (en un certain sens) 
\[ \left( \begin{array}{cc} \textrm{certaines morphismes} \\ \textrm{de } \Gal(\overline{F}/F) \\ \textrm{dans } {}^L G \end{array} \right) \to \left( \begin{array}{ccc} \textrm{certains paquets} \\ \textrm{de représentations} \\ \textrm{de } G(F) \textrm{ sur } R \end{array} \right) \] 
où ${}^L G$ est un groupe (le \og dual de Langlands \fg, défini pour la première fois dans la lettre mentionnée ci-dessus) associé à $G$ sur $R$. \\ Soit $p$ un nombre premier. Lorsque $F$ est une extension finie de $\Qp$, que $G$ est un groupe \og sympathique \fg\ et $R$ est $\C$, $\overline{\F}_\ell$ ou $\overline{\Q}_\ell$ pour $\ell \neq p$ un nombre premier, une série de travaux permettent de mettre la main sur l'existence d'une telle bijection. Citons entre autres \cite{HarTay01}, \cite{Hen00} et \cite{Vig01} pour le cas du groupe général linéaire $\GL(n,F)$.

\subsection{Le cas modulo $p$ et $p$-adique}

Cette thèse s'intéresse à des représentations modulo $p$, c'est-à-dire que l'on prend pour corps des coefficients $R$ une clôture algébrique $\overline{\F}_p$ fixée du corps fini $\Fp$. Dans ce cadre-là, très peu de choses sont connues et le gros des connaissances s'accumule autour du cas de $\GL(2,\Qp)$. La bijection prend alors la forme, pour $\F_{p^f} \subsetneq \Fpbar$ fini: 
\[\left( \begin{array}{ccc} \textrm{représentations semisimples} \\ \textrm{de dimension } 2 \\ \textrm{de } \Gal(\overline{\Q}_p/\Qp) \textrm{ sur } \F_{p^f} \end{array} \right) \to \left( \begin{array}{cccc} \textrm{certains paquets de} \\ \textrm{représentations lisses} \\ \textrm{irréductibles admissibles} \\ \textrm{de } \GL(2,\Qp) \textrm{ sur } \F_{p^f} \end{array} \right) .\] Le côté représentations galoisiennes est décrypté à travers la théorie des $(\varphi,\Gamma)$-modules de Fontaine (\cite{Fon91}). Une première forme de cette correspondance a d'abord été établie numériquement par Breuil (\cite{Bre03}), en terminant la classification de telles représentations de $\GL(2,\Qp)$. Forts de cette correspondance modulo $p$, les travaux de Breuil (\cite{Bre03b}, \cite{Bre04}), Berger (\cite{BerBre10}), Kisin (\cite{Kis10}), Emerton (\cite{Eme10}), Colmez (\cite{Col10}, \cite{Col10b}, \cite{Col10c}) et Paskunas (\cite{Pas10}) permettent d'établir une correspondance (fonctorielle) à coefficients $p$-adiques compatible à la correspondance modulo $p$ évoquée ci-dessus. 

\subsection{Les représentations modulo $p$ de groupes réductifs $p$-adiques}

Comme on vient de le voir, dans le cas de $\GL(2,\Qp)$, une étape importante dans l'établissement de la correspondance modulo $p$ a été la classification des représentations de $\GL(2,\Qp)$. On s'intéresse uniquement aux représentations irréductibles, lisses, admissibles. L'histoire débute avec la classification des représentations de Steinberg et des séries principales de $\GL(2,F)$ par Barthel et Livné (\cite{BarLiv95}, \cite{BarLiv94}), où $F$ est une extension finie de $\Qp$. Ils mettent de plus en lumière l'existence de représentations qui ne sont pas de cette nature (c'est-à-dire non sous-quotient d'une induite à partir d'un sous-groupe de Borel): les représentations dites supersingulières. La compréhension de ces représentations pour $\GL(2,\Qp)$ est due aux travaux de Breuil (\cite{Bre03}). On peut citer aussi que par la suite des travaux indépendants de Berger (\cite{Ber10}) et d'Emerton (\cite{Eme08}) ont établi l'irréductibilité des restrictions au Borel des représentations supersingulières de $\GL(2,\Qp)$. Le cas de $\GL(2,F)$ pour $F \neq \Qp$ est réellement différent et plus compliqué comme l'attestent les travaux de Paskunas (\cite{Pas04}), Breuil-Paskunas (\cite{BrePas12}), Hu (\cite{Hu10}), ou encore Schraen (\cite{Sch12}). \\
Récemment, les travaux de Grosse-Klönne (\cite{Grk09}) et de Herzig (\cite{Her11b}) ont permis de comprendre les sous-quotients d'induites paraboliques pour $\GL(n,F)$. Peu après, Abe a étendu cela au cas de tout groupe réductif connexe déployé sur $F$ (\cite{Abe11}). Aucune avancée n'est à noter pour la compréhension des supersingulières de $\GL(n,F)$ (y compris $F = \Qp$) pour $n \geq 3$. \\
Le premier cas de groupe non déployé a été étudié par Abdellatif avec l'étude des représentations de $U(2,1)(E/F)$ où $E/F$ est une extension séparable quadratique non ramifiée. Dans \cite{Abd11}, elle résout complètement le cas des sous-quotients d'induites paraboliques. Ensuite, Koziol et Xu utilisent la machinerie de Paskunas (\cite{Pas04}) pour exhiber des représentations supersingulières de $U(2,1)(E/F)$ lorsque $F$ est de corps résiduel $\Fp$ (\cite{KozXu12}). \\
Enfin, on mentionne que les représentations de $\mathrm{SL}(2,\Qp)$ et $U(1,1)(\Q_{p^2}/\Qp)$ sont complètement comprises: comme ces groupes sont en un sens \og proches \fg\ de $\GL(2,\Qp)$, Abdellatif (\cite{Abd11}) et Koziol (\cite{Koz12}) utilisent le cas de $\GL(2,\Qp)$ pour mettre la main sur toutes les représentations supersingulières du groupe qui les intéresse respectivement. \\

L'objet de notre travail a été de nous attaquer au premier cas de groupe non quasi-déployé, $U(2,1)(E/F)$ et $U(1,1)(\Q_{p^2}/\Qp)$ étant quasi-déployés. Nous nous intéressons dans un premier temps à $\GL(2,D)$ pour $D$ une algèbre à division de centre $F$, puis à $\GL(m,D)$ pour $m \geq 3$. Les différences pour l'étude des induites paraboliques sont déjà notables, et aucune approche des représentations supersingulières n'a pu être sérieusement tentée.

\section{Présentation des résultats}

\subsection{Notations}

Soit $F$ un corps local non archimédien localement compact, de caractéristique résiduelle $p$. Soit $D$ une algèbre à division de centre $F$. Son degré sur $F$ est un $d^2$ pour un certain entier $d \geq 1$. L'algèbre $D$ possède un sous-corps commutatif maximal $E$ tel que $E/F$ est une extension non ramifiée de degré $d$, unique \` a conjugaison dans $D$ pr\` es (voir \cite{Ser67}, Appendix). L'extension $D$ sur $E$ est alors totalement ramifiée au sens suivant: $D$ contient une uniformisante $\varpi_D$ de $D$ de sorte que $\varpi_F := \varpi_D^d$ soit une uniformisante de $F$ donc de $E$. On note $\Oo_D$ l'anneau des entiers de $D$, $k_D$ son corps résiduel et $D(1)$ le pro-$p$-groupe $1 + \varpi_D \Oo_D$. \\
Parce que $E$ déploie $D$, le diagramme commutatif suivant existe et définit $\det_0$:
\[ \xymatrix{ \GL(m,D \otimes_F E) \ar[r]^\sim & \GL(md,E) \ar[d]^{\det} \\ \GL(m,D) \ar@{^(->}[u] \ar[r]^{\det_0} & E^\times } .\] Par théorème de Skolem-Noether, l'action de $\sigma \in \Gal(E/F)$ sur $\GL(m,D) \to \GL(md,E)$ coïncide avec la conjugaison par un élément $c_\sigma \in \GL(md,E)$. En particulier, l'image de $\det_0$ est fixée point par point par $\Gal(E/F)$ et est donc incluse dans $F^\times$: on note $\det_G : \GL(m,D) \to F^\times$ le morphisme ainsi obtenu. Pour $m=1$, $\det_G$ devient le morphisme de groupes $\Nrm : D^\times \to F^\times$ appelé norme réduite. \\

Soient $m \geq 2$ un entier et $G =\GL(m,D)$. On note $B$ le sous-groupe de $G$ des matrices à coefficients triangulaires supérieurs: c'est un parabolique $F$-déployé minimal de $G$. \\

Toutes les représentations, sauf mention explicite du contraire, seront prises à coefficients dans $\Fpbar$.

\subsection{Le cas de $\GL(2,D)$}

On prend ici $m=2$. La principale différence avec les cas quasi-déployés déjà étudiés dans la littérature réside dans le fait que les représentations lisses irréductibles admissibles de $B$ ne sont pas toutes des caractères. En effet, dans le cas quasi-déployé, il existe un sous-groupe de Borel $B$ de $G$ défini sur $F$; et $B$ possède une décomposition de Levi $B=TU$ avec $U$ radical unipotent de $B$ et $T$ un tore (non nécessairement déployé). En confondant abusivement les groupes et leurs $F$-points, $U \cap K$ est un pro-$p$-groupe et toute représentation lisse irréductible admissible de $B$ sur $\Fpbar$ est triviale sur $U$, donc se factorise par $B \twoheadrightarrow B/U \simeq T$. Ceci explique pourquoi, dans le cas quasi-déployé, on ne travaille qu'avec des caractères de $B$. \\ 
Revenons-en à $G = \GL(2,D)$. Il n'est pas quasi-déployé parce que $D^\times$ ne l'est pas. A ce titre, commençons par dire un mot sur les représentations lisses irréductibles de $D^\times$. Comme $D(1)$ est un pro-$p$-groupe, toute représentation lisse irréductible de $D^\times$ se factorise par \[ D^\times \twoheadrightarrow D^\times /D(1) \simeq k_D^\times \rtimes \varpi_D^\Z .\] On remarque alors que $\varpi_D^{d \Z}$ est dans le centre de $D^\times$ et son action sur $k_D^\times$ dans le produit semi-direct précédent est triviale. Il en résulte que toute représentation lisse irréductible admissible de $D^\times$ est de dimension finie, égale à un diviseur de $d$. \\
Soit $\rho$ une représentation lisse irréductible admissible de $B$: elle se factorise par \[ B \twoheadrightarrow B/U \simeq D^\times \times D^\times \] où $U$ désigne le radical unipotent de $B$; et $\rho$ se décompose donc en $\rho_1 \otimes \rho_2$ où $\rho_1$ et $\rho_2$ sont deux représentations de dimensions finies de $D^\times$ comme expliqué ci-dessus. Le $\det_G$ d'un élément de $B$ se lit \[ \det\nolimits_G \left( \begin{array}{cc} a & b \\ 0 & d \end{array} \right) = \Nrm(a) \, \Nrm(d) .\] On dira qu'une représentation $\rho$ de $B$ se factorise par un caractère de $F^\times$ s'il existe un caractère $\rho_0 : F^\times \to \Fpbar^\times$ vérifiant $\rho = \rho_0 \circ \det_G$ sur $B$. \\
On commence par établir un critère d'irréductibilité pour l'induite parabolique $\Ind_B^G \, \rho$.

\begin{theo} \label{IntThm1}
Soit $\rho$ une représentation lisse irréductible admissible de $B$.
\begin{itemize}
\item[(i)] Si $\rho$ se factorise par un caractère $\rho_0$ de $F^\times$, alors $\Ind_B^G \, \rho$ est extension non scindée de deux représentations de $G$ irréductibles non isomorphes \[ 0 \to \rho_0 \circ \det\nolimits_G \to \Ind_B^G \, \rho \to \St \, \rho_0 \to 0 ,\] définissant par là la représentation de Steinberg $\St \, \rho_0$.
\item[(ii)] Si $\rho$ ne se factorise pas par $\det_G$, alors $\Ind_B^G \, \rho$ est irréductible.
\end{itemize}
\end{theo}
\textsf{Remarque :}
Les représentations lisses irréductibles qui apparaissent dans le théorème \ref{IntThm1} sont admissibles mais on ne s'y attachera pas dans un premier temps. En effet, ceci sera établi plus généralement au cours des parties \ref{IntSt} et \ref{IntGLm}. \\

On va présenter trois preuves distinctes du théorème \ref{IntThm1}.(ii). La première consiste à regarder la restriction à $B$ de la représentation à étudier; la deuxième (sous l'hypothèse supplémentaire que $p$ ne divise pas $2d$) consiste à étudier la simplicité de l'espace $(\Ind_B^G \, \rho)^{I(1)}$ des $I(1)$-invariants en tant que module à droite sur l'algèbre de Hecke du pro-$p$-Iwahori $\Hh(G,I(1)) = \End_G(\ind_{I(1)}^G \, \id)$, où $\ind_H^G$ désigne le foncteur d'induction lisse à support compact modulo $H$ pour $H$ un sous-groupe ouvert de $G$ et $\id$ la représentation triviale. La combinaison de ces deux méthodes (sans l'hypothèse supplémentaire $p \nmid 2d$) nous fournit aussi une preuve du théorème \ref{IntThm1}.(i). Quant à la dernière preuve présentée du théorème \ref{IntThm1}.(ii), elle repose sur un isomorphisme entre une telle représentation $\Ind_B^G \, \rho$ et un conoyau dans une représentation obtenue par induction compacte. \\

Précisons un peu la nature de ces derniers objets. Pour cela, notons $K$ le sous-groupe compact maximal $\GL(2,\Oo_D)$ de $G$. Soit $V$ une représentation lisse irréductible de $K$; elle est de dimension finie et on sait la décrire explicitement. L'induite compacte $\ind_K^G \, V$ est alors un module à gauche sur l'algèbre de Hecke $\Hh = \End_G(\ind_K^G \, V)$. On montre que cette dernière est commutative et les représentations évoquées précédemment seront des conoyaux de certains opérateurs de $\Hh$ dans $\ind_K^G \, V$: ce sont des $\ind_K^G \, V \otimes_{\Hh, \chi} \Fpbar$ pour des caractères $\chi : \Hh \to \Fpbar$ de $\Fpbar$-algèbres. \\
On dira qu'une représentation admissible est supersingulière si elle un quotient irréductible d'un $\ind_K^G \, V \otimes_{\Hh, \chi} \Fpbar$ pour un caractère $\chi$ s'annulant sur un certain opérateur explicite. Enfin, les comparaisons entre les induites paraboliques $\Ind_B^G \, \rho$ et ces $\ind_K^G \, V \otimes_{\Hh, \chi} \Fpbar$ fournissent le théorème de classification suivant. 

\begin{theo} \label{IntThm2}
Soit $\pi$ une représentation lisse irréductible admissible de $G$. Alors $\pi$ est dans l'un des cas suivants: 
\begin{itemize}
\item[(a)] un caractère $\rho_0 \circ \det_G$, où $\rho_0$ est un caractère de $F^\times$;
\item[(b)] une Steinberg $\St \, \rho_0$, où $\rho_0$ est un caractère de $F^\times$;
\item[(c)] une induite parabolique $\Ind_B^G \, \rho$, où $\rho$ est une représentation irréductible de $B$ qui ne se factorise pas par un caractère;
\item[(d)] une supersingulière.
\end{itemize}
\end{theo}
\textsf{Remarque :}
Ce sera bien un théorème de classification au sens où l'on prouvera qu'il n'y a pas d'entrelacement non trivial entre ces représentations.

\begin{cor}
Une représentation supersingulière de $G$ est supercuspidale, c'est-à-dire qu'elle n'est pas sous-quotient d'une induite parabolique; et réciproquement.
\end{cor}

\subsection{Les réprésentations de Steinberg généralisées} \label{IntSt}

Dans le but de généraliser les énoncés des théorèmes \ref{IntThm1} et \ref{IntThm2} à $\GL(m,D)$ pour $m \geq 3$, on a besoin de comprendre les \og briques de base \fg\ correspondant aux représentations de Steinberg dans le cas de $\GL(2,D)$; on les appellera donc représentations de Steinberg généralisées. Comme il n'est pas vraiment plus difficile de l'obtenir dans le cadre d'un groupe réductif connexe quelconque, c'est ce que l'on fera. Enfin, on souligne que les changements à apporter par rapport aux travaux de Grosse-Klönne et de Herzig qui ont traité le cas d'un groupe déployé sont minimes. \\ 

Soit $G$ le groupe des $F$-points d'un groupe réductif connexe défini sur $F$. Soit $P$ un sous-groupe parabolique standard de $G$, c'est-à-dire qui contient $B$, un parabolique déployé minimal fixé. L'espace des fonctions $P$-invariantes à gauche $\Ind_P^G \, \id$ contient l'espace $\Ind_Q^G \, \id$ pour tout sous-groupe parabolique $Q$ de $G$ contenant $P$ et on définit la représentation de Steinberg relative à $P$: \[ \St_P \id = \dfrac{\Ind_P^G \, \id}{\sum_{Q \gneq P} \Ind_Q^G \, \id} .\]
Soit $K$ un sous-groupe parahorique maximal spécial de $G$. Notons $K(1)$ le pro-$p$-radical de $K$. Pour tout sous-groupe $H$ de $G$, on note $\overline{H}$ le sous-groupe $(H \cap K)/(H \cap K(1))$ de $\overline{G} = K/K(1)$. Par exemple, lorsque $G$ est $\GL(m,D)$, $K(1)$ est le sous-groupe de $K$ composé des matrices congrues à l'identité modulo $(\varpi_D)$ et $\overline{G}$ est le groupe fini $\GL(m,k_D)$. On note $\red:K \to \overline{G}$ la réduction modulo $K(1)$ et on rappelle que le sous-groupe d'Iwahori standard de $G$ est $I = \red^{-1}(\overline{B})$. On commence par récolter des informations sur l'espace des $I$-invariants de $\St_P \id$ en le comparant à son analogue fini, l'espace des $\overline{B}$-invariants de $\overline{\St}_{\overline{P}} \id = \dfrac{\Ind_{\overline{P}}^{\overline{G}} \, \id}{\sum_{\overline{Q} \gneq \overline{P}} \Ind_{\overline{Q}}^{\overline{G}} \, \id}$.

\begin{prop}
$\phantom{,}$
\begin{itemize}
\item[(i)] Il existe  une injection $K$-équivariante $\overline{\St}_{\overline{P}} \id \hookrightarrow \St_P \id$ qui induit un isomorphisme \[ (\overline{\St}_{\overline{P}} \id)^{\overline{B}} \xrightarrow{\sim} (\St_P \id)^I \] de $\Fpbar$-espaces vectoriels.
\item[(ii)] Il existe une fonction $f_{\overline{P}} \in \Ind_{\overline{P}}^{\overline{G}} \id$ explicite telle que tout sous-$\End_{\overline{G}}(\Ind_{\overline{B}}^{\overline{G}} \, \id)$-module à droite non nul de $(\overline{\St}_{\overline{P}} \id)^{\overline{B}}$ contient la classe de $f_{\overline{P}}$ dans $\overline{\St}_{\overline{P}} \id$.
\end{itemize}
\end{prop}

La preuve de cette proposition est uniquement combinatoire et suit les arguments de \cite{Grk09} malgré la présence d'un système de racines éventuellement non réduit. On en déduit ensuite, par la machinerie de Herzig, le résultat voulu suivant.

\begin{theo}
$\phantom{,}$
\begin{itemize}
\item[(i)] Le $K$-socle\footnote{par définition, c'est la plus grande sous-$K$-représentation semi-simple} de $\St_P \id$ est irréductible et on peut l'identifier explicitement.
\item[(ii)] La représentation de Steinberg généralisée $\St_P \id$ est irréductible et admissible.
\end{itemize}
\end{theo}

\subsection{Le cas de $\GL(m,D)$} \label{IntGLm}

Maintenant que l'on dispose de toutes les informations voulues sur les Steinberg généralisées, on retourne au cadre $G = \GL(m,D)$ précédent. Le travail que l'on effectue, dans la lignée de celui de Herzig pour le cas déployé $\GL(m,F)$, est une généralisation des méthodes employées dans la troisième preuve présentée du théorème \ref{IntThm1} pour $\GL(2,D)$. \\

Soient $S$ un tore $F$-déployé maximal de $\GL(m,E) \leq \GL(m,D)$; c'est aussi un tore $F$-déployé maximal de $\GL(m,D)$. On note $X_*(S)$ le sous-groupe des cocaractères associés et $A$ le centralisateur de $S$ dans $\GL(m,D)$. Quitte à remplacer $S$ par l'un de ses conjugués, on peut supposer que $A$ est composé des matrices diagonales, et en particulier est inclus dans $B$; c'est donc ce que l'on fera par la suite. \\
Dans un premier temps, on doit comprendre un peu mieux l'algèbre de Hecke-Satake $\Hh(G,K,V) = \End_G(\ind_K^G \, V)$. Grâce au travail de Henniart-Vignéras (\cite{HenVig11}), on sait que c'est une $\Fpbar$-algèbre commutative qui s'injecte via la transformée de Satake modulo $p$ dans $\Hh \big( A,A \cap K , V^{U \cap K} \big)$ où $U$ est le radical unipotent de $B$ et $K$ le compact maximal $\GL(m,\Oo_D)$. D'autre part, pour inverser la transformée de Satake, dans le cas déployé et où $V$ est la représentation triviale de $K$, on dispose de la formule de Lusztig-Kato. On fera attention que ladite formule concerne des représentations à coefficients complexes. Ici, pour se faire, on va comparer les valeurs de transformées de Satake d'opérateurs de Hecke pour $\GL(m,D)$ avec celles d'opérateurs pour le groupe déployé $\GL(m,E) \leq \GL(m,D)$. Pour $H$ un sous-groupe de $G$, on note $H_{(E)}$ le sous-groupe $H \cap \GL(m,E)$ de $\GL(m,E)$ correspondant; en particulier, $G_{(E)}$ désigne $\GL(m,E)$ lui-même. Enfin, $W$ désignera le groupe de Weyl fini associé à $S$ (invariablement pour $G$ ou $G_{(E)}$), que l'on verra concr\` etement comme le sous-groupe de $G_{(E)} \leq G$ des matrices de permutation. \\
On prend ensuite les notations suivantes pour les transformées de Satake classiques, à coefficients complexes:
\[ \Sss_{(E)} : \Hh(G_{(E)},K_{(E)},\id) \xrightarrow{\sim} \C[X_*(S)]^W ,\] \[ \overline{\Sss} : \Hh(G,K,\id) \xrightarrow{\sim} \C[A/A \cap K]^W \xrightarrow{\sim} \C[X_*(S)]^W .\]
On remarque que l'isomorphisme $W$-équivariant de $\Fpbar$-algèbres \[ \C[A/A \cap K] \xrightarrow{\sim} \C[X_*(S)] \] utilisé pour définir $\overline{\Sss}$ provient du morphisme de groupes $W$-équivariant suivant:
 \[ \begin{array}{cccccc} \psi: & A/A \cap K & \xrightarrow{\sim} & S/S \cap K & \xrightarrow{\sim} & X_*(S) \\ & x & \mapsto & x^d & & \\ & & & \lambda(\varpi_F) & \mapsfrom & \lambda \end{array} \] 
Pour tout élément $\lambda$ de $X_*(S)$, on fixe $t_\lambda \in A$ un représentant de l'image réciproque $\psi^{-1}(\lambda)$ de $\lambda$ dans $A/A \cap K$. L'énoncé de comparaison prend alors la forme de la conjecture suivante, que l'on sait prouver en petite dimension comme on l'expliquera par la suite. On supposera la conjecture satisfaite dans toute la suite de l'introduction.

\begin{conj} \label{incipitconj}
Soit $\lambda$ un élément de $X_*(S)$. Alors on a \[ \overline{\Sss} \big( \id_{K t_\lambda K} \big) = \Sss_{(E)} \big( \id_{K_{(E)} \lambda(\varpi_F) K_{(E)}} \big) ,\] où $\id_{K t K}$ (respectivement $\id_{K_{(E)} t K_{(E)}}$) désigne l'opérateur de $\Hh(G,K,\id)$ (resp. $\Hh(G_{(E)},K_{(E)},\id)$) de support $KtK$ (resp. $K_{(E)} t K_{(E)}$) pour $t \in G$ (resp. $t \in G_{(E)}$).
\end{conj}

Cette identité nous permet de ramener l'étude de l'inversion de Satake modulo $p$ de $\GL(m,D)$ à un cas déployé lorsque $V = \id$ est la représentation triviale. \\
Dans le cas déployé de $\GL(m,E)$, on peut ramener les calculs voulus sur la transformée de Satake à ce cas $V=\id$. Ceci n'est malheureusement plus le cas pour $\GL(m,D)$, et la raison pour cela réside dans le fait suivant qui n'est plus valable pour $\GL(m,D)$: tout caractère de $K_{(E)}$ à valeurs dans $\Fpbar^\times$ se prolonge en un caractère de $G_{(E)}$. \\
Dans le cas de $\GL(m,D)$, on doit se contenter de l'avatar plus faible suivant.

\begin{prop}
Tout caractère $\xi$ de $K$ se prolonge en une application ensembliste du support de l'algèbre de Hecke $\Hh(G,K,\xi)$ à valeurs dans $\Fpbar^\times$. De plus, ce prolongement est multiplicatif pour $m \leq 2$ (en un sens convenable).
\end{prop}

Pour $m=3$, déjà, on ne peut plus garantir la qualité de multiplicativité à un tel prolongement. Et c'est la raison pour laquelle on doit rajouter l'hypothèse supplémentaire suivante, que l'on supposera donc satisfaite dans toute la suite de l'introduction.

\begin{hypo*} 
L'une des deux conditions suivantes est satisfaite:
\begin{itemize}
\item[(a)] le degré $d$ est un nombre premier, ou égal à $1$;
\item[(b)] l'entier $m$ est inférieur ou égal à $3$. 
\end{itemize}
\end{hypo*}

Dans le cadre de l'hypothèse, on va pouvoir utiliser le calcul de l'inversion de Satake sans coefficients (c'est-à-dire $V \! = \! \id$) et le prolongement ensembliste d'un caractère $\xi$ de $K$ pour proposer un argument de \og changement de poids \fg, essentiel dans la machinerie à la Herzig. Compte tenu de la remarque précédente, la présence du cas $m=3$ dans l'hypothèse peut être surprenant. Cependant, cela s'explique par le fait que l'on utilise alors uniquement un prolongement de caractère pour un Levi $\GL(2,D) \times \GL(1,D)$. \\

Aussi, les cas $m \leq 3$ sont des cas pour lesquels la conjecture \ref{incipitconj} est vérifiée. Remarquons ici que l'essence de cet énoncé réside dans la comparaison de cardinaux d'intersections de doubles classes de Cartan et de classes d'Iwasawa, pour $\GL(m,D)$ et $\GL(m,E)$ respectivement. 

\begin{prop}
Supposons $m \leq 3$. Alors la conjecture \ref{incipitconj} est vraie. 
\end{prop}

Avant d'énoncer le changement de poids précédemment évoqué, on a encore besoin de deux définitions. \\
D'abord, deux caractères $\xi$ et $\xi'$ de $A \cap K$ sont dits conjugués s'il existe un élément $a \in A$ vérifiant $\xi' = \xi(a \cdot a^{-1})$. Puis, les $K$-représentations irréductibles lisses $V$ (appelées \og poids \fg\ par Herzig) sont classifiées par le caractère $V^{U \cap K}$ de $A \cap K$ et le paramètre $M_V$ de régularité: $M_V$ est le Levi standard maximal tel que $M \cap K$ stabilise la droite $V^{U \cap K}$. Et $V$ est dit $M$-régulier pour tout Levi standard $M$ contenant $M_V$. Signalons par ailleurs que l'on a $M_V=G$ si et seulement si $V$ est de dimension $1$. \\
Soit $\chi: \Hh(G,K,V) \to \Fpbar$ un caractère de $\Fpbar$-algèbres. On définit le paramètre $M_\chi$ de régularité de $\chi$ comme le plus petit Levi standard $M$ tel que $\chi$ se factorise par la transformée de Satake partielle $\Hh(G,K,V) \to \Hh \big( M,M \cap K, V^{N \cap K} \big)$, où $N$ désigne le radical unipotent du parabolique standard de Levi $M$. \\

Soient $M = \prod\limits_{i=1}^r M_i$ un Levi standard de $G$ où chaque $M_i$ est de taille $m_i \geq 1$, et $P$ le parabolique standard associé. Soit $\sigma = \bigotimes\limits_{i=1}^r \sigma_i$ une représentation irréductible admissible de $M$, où chaque $\sigma_i$ est dans l'un des cas de figure suivant:
\begin{itemize}
\item[(a)] $\sigma_i$ est supersingulière avec $m_i > 1$;
\item[(b)] $\sigma_i$ est une représentation de $D^\times$ (ainsi $m_i = 1$) avec $\dim_{\Fpbar} \sigma_i >1$;
\item[(c)] $\sigma_i$ est une Steinberg généralisée $\big( \rho_i^0 \circ \det_G \big) \otimes \St_{Q_i} \Fpbar$ pour un parabolique standard $Q_i$ de $M_i$ et un caractère $\rho_i^0$ de $F^\times$.
\end{itemize}
Supposons de plus $\rho_i^0  \circ \Nrm \neq \rho_{i+1}^0 \circ \Nrm$ s'il existe deux blocs adjacents dans le cas (c). \\
On s'intéresse aux représentations irréductibles $V \subseteq (\Ind_{P^-}^G \, \sigma)|_K$. Pour $M_i$ dans le cas $(c)$, on note $M_i'$ le Levi standard de $Q_i$. Les tels $M_i'$ et $M_\chi$ engendrent un Levi standard que l'on notera $M_\chi'$. La proposition de changement de poids s'énonce alors comme suit.

\begin{prop} \label{IntProp4}
Soient $V$ une $K$-représentation lisse irréductible du $K$-socle de $\Ind_{P^-}^G \, \sigma$ et $\chi : \Hh(G,K,V) \to \Fpbar$ un caractère de $\Fpbar$-algèbres. On suppose l'existence d'un morphisme non nul \[ \ind_K^G \, V \otimes_{\Hh, \chi} \Fpbar \to \Ind_{P^-}^G \, \sigma \] de $G$-représentations. Il existe une $K$-représentation lisse irréductible $V'$ satisfaisant aux trois conditions suivantes:
\begin{itemize}
\item[(i)] $V^{U \cap K}$ et $(V')^{U \cap K}$ sont conjuguées en tant que $(A \cap K)$-représentations (de sorte que l'on peut identifier $\Hh(G,K,V)$ à $\Hh(G,K,V')$; on notera $\Hh$ chacune de ces deux algèbres);
\item[(ii)] on a un isomorphisme de $G$-représentations \[ \ind_K^G \, V \otimes_{\Hh,\chi} \Fpbar \xrightarrow{\sim} \ind_K^G \, V' \otimes_{\Hh,\chi} \Fpbar ;\]
\item[(iii)] $V'$ est $M_\chi'$-régulière.
\end{itemize}
\end{prop}

La preuve de cette proposition est présentée ici à travers des calculs sur des opérateurs de Hecke et leur transformée de Satake; c'est uniquement dans cette proposition que l'hypothèse précédente est nécessaire, et tous les autres arguments ne l'utilisent qu'à travers cette proposition. \\
Ajoutons qu'une $G$-représentation lisse irréductible $\pi$ est dite supersingulière si, pour toute sous-$K$-représentation irréductible $V$ de $\pi|_K$ et tout caractère $\chi : \Hh(G,K,V) \to \Fpbar$ de $\Fpbar$-algèbres induisant un morphisme de $G$-représentations non nul \[ \ind_K^G \, V \otimes_{\Hh,\chi} \Fpbar \to \pi ,\] le paramètre de régularité $M_\chi$ associé à $\chi$ est $G$. \\ 
Maintenant que l'on a à disposition la définition d'une représentation supersingulière, l'énoncé de changement de poids et le travail de Henniart et Vignéras (\cite{HenVig11b}) sur la comparaison entre induction parabolique et induction compacte, on peut en déduire un critère d'irréductibilité, qui est la généralisation directe du théorème \ref{IntThm1}.(ii).

\begin{theo} \label{IntThm4}
Soient $P$ et $\sigma$ comme avant la proposition \ref{IntProp4}. Alors $\Ind_P^G \, \sigma$ est une $G$-représentation irréductible admissible.
\end{theo}
\textsf{Remarque :} 
Comme pour le théorème \ref{IntThm1}.(i), on va donner les composants de Jordan-Hölder de $\Ind_P^G \, \sigma$ dans le cas de réductibilité (c'est-à-dire si on enlève l'hypothèse $\rho_i^0 \circ \Nrm \neq \rho_{i+1}^0 \circ \Nrm$). \\

Enfin, avec cela, on peut reprendre les calculs de Herzig avec le foncteur des parties ordinaires d'Emerton (étudié par Vignéras lorsque $F$ n'est pas de caractéristique nulle) et conclure quant à un énoncé de classification.

\begin{theo}
Soit $\pi$ une représentation irréductible admissible de $G$. Alors il existe un unique parabolique standard $P$ de Levi standard $M$ et une unique $M$-représentation irréductible admissible $\sigma$ de $M$ (à isomorphisme près) vérifiant les conditions du théorème \ref{IntThm4} tels que $\pi$ soit isomorphe à $\Ind_P^G \, \sigma$.
\end{theo}

Comme pour $\GL(2,D)$, on en déduit l'équivalence de notions suivante.

\begin{cor}
Une représentation supersingulière de $G$ est supercuspidale, et réciproquement.
\end{cor}

Concluons en disant que nos résultats sont de la même veine que dans le cas déployé $\GL(m,F)$ (correspondant donc à $d=1$) et les méthodes employées sont adaptées du cas déployé.
\chapter{Des représentations modulo $p$ de $\mathrm{GL}(2,D)$, $D$ algèbre à division sur un corps local} \label{Ly11b}

\section{Introduction}

L'histoire de l'étude des représentations lisses modulo $p$ de $\GL(2,D)$ commence avec Barthel-Livné en 1994-95 avec le cas déployé $D=F$ (voir \cite{BarLiv94}, \cite{BarLiv95}). On complète ici leur classification avec le cas non nécessairement déployé. Soient $G = \GL(2,D)$ et $B$ le parabolique minimal composé des matrices triangulaires supérieures inversibles. On commence par établir le critère d'irrréductibilité suivant.

\begin{theo} \label{mainthm}
Soient $\rho_1$ et $\rho_2$ deux représentations lisses irréductibles de dimensions finies de $D^\times$.
\begin{itemize}
\item[(i)] Supposons $\rho := \rho_1 = \rho_2$ de dimension $1$. Alors $\Ind_B^G \, \rho \otimes \rho$ est admissible et est extension non scindée de deux représentations irréductibles admissibles non isomorphes.
\item[(ii)] Dans tous les autres cas de figure, $\Ind_B^G \, \rho_1 \otimes \rho_2$ est irréductible admissible. 
\end{itemize}
\end{theo}

Un autre travail de l'auteur sur les représentations de Steinberg généralisées établit le \textit{(i)} dans un degré de généralité supérieur (voir chapitre \ref{Ly11}). Cependant une preuve pédestre sera présentée ici. \\
Pour le \textit{(ii)}, on présentera trois preuves distinctes. La première est notable pour sa simplicité et repose uniquement sur une étude de la restriction à $B$. La seconde, conditionnellement à une hypothèse sur $p$, s'appuie sur l'étude du module des invariants sous le pro-$p$-Iwahori. Quoique plus pénible, cette méthode a l'avantage d'exhiber des modules simples de Hecke-Iwahori, auxquels l'auteur consacrera un article ultérieur. Enfin, la dernière méthode présentée, plus proche des articles de Barthel-Livné (et par conséquent du travail récent de Herzig dans \cite{Her11b}) aboutit de plus au théorème de classification suivant.

\begin{theo} \label{mainthm2}
Les représentations lisses irréductibles admissibles de $G$ sur $\Fpbar$ sont les suivantes:
\begin{itemize}
\item[(a)] les caractères;
\item[(b)] les représentations de Steinberg;
\item[(c)] les induites paraboliques irréductibles;
\item[(d)] les supercuspidales.
\end{itemize}
De plus, les supercuspidales sont exactement les supersingulières.
\end{theo} 

On fait remarquer que ce résultat ne fait que débuter l'étude des représentations lisses admissibles modulo $p$ de $G$ puisqu'à l'heure actuelle les représentations supersingulières ne sont connues explicitement que dans le cas $D = F = \Qp$ grâce au travail de Barthel-Livné (voir \cite{BarLiv94}, \cite{BarLiv95}) et de Breuil (voir \cite{Bre03}).

\section{Notations et généralités}

Soient $p$ un nombre premier et $\Fpbar$ une clôture algébrique fixée de $\Fp$; tout corps fini de caractéristique $p$ sera vu comme un sous-corps de $\Fpbar$. \\

L'objet de cette section est de donner les notations et de présenter des faits généraux bien connus sur la théorie des représentations ou bien sur $D^\times$ ou  $\GL(2,D)$.

\subsection{Représentations et algèbres de Hecke} \label{parGeneralHecke}

On pourra se reporter au paragraphe 2 de \cite{BarLiv94}. Dans tout ce qui suit, toute représentation considérée sera lisse (et on oubliera souvent de le mentionner), à coefficients dans $\Fpbar$. Soient $G$ un groupe topologique et $H \leq G$ un sous-groupe fermé. On notera $\ind_H^G$ le foncteur d'induction compacte lisse et $\Ind_H^G$ celui d'induction lisse. L'action sur une induite se fera par translation à droite, à savoir $g. f : x \mapsto f(xg)$. On suppose $H$ ouvert. Soit $V$ une représentation de $H$. Pour tous $g \in G$ et $v \in V$, on définit $[g,v]$ comme étant la fonction de $\ind_H^G \, V$ à support dans $H g^{-1}$ et prenant pour valeur $v$ en $g^{-1}$. En particulier, si $h$ est un élément de $H$, on a $[gh,v] = [g,hv]$ et $g[1,v]=[g,v]$. \\ 
Pour toutes représentations $V_1$ et $V_2$ de $H$, on définit l'espace d'entrelacements \[ \Hh(G,H,V_1,V_2) := \Hom_G \big( \ind_H^G \, V_1, \ind_H^G \, V_2 \big) .\] On se permettra aussi de noter les algèbres de Hecke \[ \Hh(G,H,V) := \Hh(G,H,V,V), \quad \Hh(G,H) := \Hh(G,H,\id) ,\] où $\id$ désigne la représentation triviale de $H$. \\ 
Lorsque $H$ est ouvert dans $G$, par réciprocité de Frobenius, on a \[ \Hh(G,H,V_1,V_2) \simeq \Hom_H(V_1, \ind_H^G \, V_2) .\] 
Supposons que, pour tout $g \in G$, la double classe $HgH$ est union finie de classes à gauche (ou à droite), et que $V_1$ et $V_2$ sont finiment engendrées en tant que $H$-représentations. Alors on peut (voir \cite{HenVig11}, section 2.2) exhiber l'isomorphisme 
\[\begin{array}{ccc} \Hom_H(V_1, \ind_H^G \, V_2) & \simeq & \left\{ \begin{tabular}{r|l} \multirow{2}{*}{$G \xrightarrow{f} \Hom_{\Fpbar}(V_1,V_2)$} & $f(h_2 gh_1) = h_2 f(g) h_1$ \\ & $H \backslash \Supp \, f / H$ fini \end{tabular} \right\} \\ \big( v \mapsto (g \mapsto f(g)(v) )\big) & \mapsfrom & f \\ f & \mapsto & \big( g \mapsto (v_1 \mapsto f(v_1)(g)) \big) . \end{array} \] 
On va souvent assimiler les opérateurs de $\Hh(G,H,V_1,V_2)$ à des fonctions sur $G$ à travers cet isomorphisme. \\
Pour tous $V_1, V_2, V_3$ et tous $f_1 \in \Hh(G,H,V_1,V_2)$, $f_2 \in \Hh(G,H,V_2,V_3)$, on a la convolée $f_2 * f_1 \in \Hh(G,H,V_1,V_3)$ définie par 
\begin{equation} \label{defconvol} f_2 * f_1  : g \mapsto \sum_{x \in G/H} f_2(x) f_1(x^{-1} g) .\end{equation}
Aussi, quand on voit les opérateurs de $\Hh(G,H,V)$ comme fonctions sur $G$, la structure multiplicative est donnée par la convolution, qui est bien compatible à la multiplication par composition de la première définition (voir \cite{BarLiv94}, Proposition 5).

\subsection{Des notations pour $D$}

Soit $F$ un corps local non archimédien localement compact à corps résiduel fini $k_F$ de caractéristique $p$; on notera $\Oo_F$ son anneau de valuation, $\m_F$ son idéal maximal et $q = p^f$ le cardinal de $k_F$.\\
Soit $D$ une algèbre à division de centre $F$: elle est de degré $d^2$ sur $F$ pour un certain entier $d \geq 1$, et d'invariant de Brauer $a_D/d$ avec $a_D$ un entier de $[1,d[$ premier à $d$. Alors $D$ possède un sous-corps commutatif maximal $E$ tel que $E/F$ est une extension non ramifiée de degré $d$, unique \` a conjugaison dans $D$ pr\` es (voir \cite{Ser67}, Appendix); on note $\Oo_E$ son anneau d'entiers. On peut munir $D$ d'une valuation discrète  $v_D$; on fixe une uniformisante $\varpi$ de $D$ (et on normalise à $v_D(\varpi) = 1$) de sorte que $\varpi_F := \varpi^d$ est une uniformisante de $F$. On note $\Oo_D$ l'anneau de valuation de $D$, $\m_D = (\varpi)$ son idéal maximal et $k_D = \Oo_D/\m_D$ son corps résiduel (qui est de cardinal $q^d$). On note $x \mapsto \overline{x}$ l'application de réduction de $\Oo_D$ à $k_D$. L'action de l'uniformisante sur le corps résiduel est $\varpi \overline{x} \varpi^{-1} = \overline{x}^{q^{a_D}}$ pour tout $x \in \Oo_D$. On remarque que $k_D$ est aussi le corps résiduel de $E$. Ainsi, on note $[ \phantom{x}]$ l'application de Teichmüller $k_D^\times \to \Oo_E^\times$, que l'on prolonge en $k_D \to \Oo_E$ en envoyant $0$ sur $0$. \\ Pour plus de précisions on pourra aller voir le chapitre 17 de \cite{Pie82} ou le chapitre 14 de \cite{Rei75}. \\

Fixons un générateur $\mu$ de $k_D^\times$. Cela définit un isomorphisme $\{ \phantom{x} \} : k_D^\times \xrightarrow{\sim} \Z/(q^d-1) \Z$ envoyant $\mu$ sur $1$. \\
Pour deux entiers positifs $a$ et $b$, on notera $a \wedge b$ leur plus grand commun diviseur.

\subsection{Des notations pour $\GL(2,D)$}

Le groupe $\GL(2,D)$ est noté $G$. Notons $B = \left( \begin{array}{cc} D^\times & D \\ 0 & D^\times \end{array} \right)$ le sous-groupe des matrices triangulaires supérieures de $G$; c'est un parabolique minimal de $G$ sur $F$. On note $U$ son radical unipotent et $A$ son sous-groupe de Levi constitué des matrices diagonales. \\
Soit $K$ le compact maximal $\GL(2,\Oo_D)$ de $G$. Notons $K(1)$ le sous-groupe normal $1 + \varpi M_2(\Oo_D)$ de $K$, et $\overline{K}$ le groupe quotient $K/K(1)$, alors isomorphe au groupe fini $\GL(2,k_D)$. On définit $U_0 = U \cap K$. Si $\pi$ désigne la projection canonique $K \twoheadrightarrow \overline{K}$, le sous-groupe d'Iwahori de $K$ est $I = \pi^{-1}(\pi(B \cap K))$ et $I(1)$ désigne l'unique pro-$p$-Sylow de $I$, appelé pro-$p$-Iwahori. On remarque $I(1) = \pi^{-1}(\pi(U_0))$. \\

On note $W_e$ le sous-groupe de $G$ engendré par $\omega := \left( \begin{array}{cc} 0 & 1 \\ \varpi & 0 \end{array} \right)$ et $s := \left( \begin{array}{cc} 0 & 1 \\ 1 & 0 \end{array} \right)$. Si $T$ désigne le tore diagonal $\left( \begin{array}{cc} F^\times & 0 \\ 0 & F^\times \end{array} \right)$, qui est déployé maximal, $W_e$ est isomorphe au groupe de Weyl étendu $N_G(T)/(T \cap K)$. On note respectivement $\varphi_1$ et $\varphi_2$ les éléments particuliers $\left( \begin{array}{cc} \varpi & 0 \\ 0 & 1 \end{array} \right)$ et $\left( \begin{array}{cc} 1 & 0 \\ 0 & \varpi \end{array} \right)$ de $W_e$; et $\varpi$ désigne $\varphi_1 \varphi_2$. On note $A_\Lambda$ le sous-groupe de $A$ engendré par $\varphi_1$ et $\varphi_2$: c'est un système de représentants de $A/ A \cap K$. Enfin, pour $x , y \in D^\times$, on notera $\delta^x_y := \left( \begin{array}{cc} x & 0 \\ 0 & y \end{array} \right) \in A$. \\
On rappelle que $[\phantom{x}] : k_D \to \Oo_E$ désigne l'application de Teichmüller; on la prolongera en \[ \begin{array}{cccc} [\phantom{x}] : & k_D^\N & \to & \Oo_D \\ & (x_j)_{j} & \mapsto & \sum_j \varpi^j [x_j] \end{array} .\] De plus, pour tout $i \in \N$, on plongera naturellement $k_D^i$ dans $k_D^\N$ par l'application $(x_0, \dots, x_{i-1}) \mapsto (x_0, \dots, x_{i-1}, 0, 0, \dots )$. 
Pour $x \in D$, notons $n(x)$ la matrice de $U$ suivante : $\left( \begin{array}{cc} 1 & x \\ 0 & 1 \end{array} \right)$. On note aussi \[ \overline{n}(x) = sn(x) s \in U^- := sUs .\] 
% Un système de représentants de $G/(I \varpi^\Z)$ est \[ \{ g^\epsilon_x, g^\epsilon_x \omega \ | \ \epsilon \in \{0,1\} , \, x \in k_D^i, \, i \geq 0 \} .\] Un système de représentants de $G/(K \varpi^\Z)$ est donné par \[ \{ g^\epsilon_x \omega \ | \ \epsilon \in \{0,1\} , \, x \in k_D^i, \, i \geq 0 \} .\] (à reprouver si cela doit réapparaître!)

\subsection{Les \og déterminants \fg \, dans $\GL(2,D)$} \label{widetildedet}

Lorsque $D$ n'est pas commutatif (c'est-à-dire $d>1$), il n'y a pas de déterminant canonique $\GL(2,D) \to D^\times$. Cependant, on va tout de même présenter un morphisme de groupes $\det_G : \GL(2,D) \to F^\times$ qu'on conviendra d'appeler déterminant. \\
Parce que $E$ déploie $D$, on a le diagramme commutatif suivant entre algèbres de matrices, qui nous permet de définir $\det_0$:
\[ \xymatrix{ M_2(D) \otimes_F E \ar[r]^{\phantom{xx} \sim} & M_{2d}(E) \ar[d]^{\det} \\ M_2(D) \ar@{^{(}->}[u] \ar[r]^{\det_0} & E } .\]
L'isomorphisme de déploiement \[ \eta_D : M_2(D) \otimes_F E \xrightarrow{\sim} M_{2d}(E) \] n'est pas $\Gal(E/F)$-équivariant mais on sait par application du théorème de Skolem-Noether que, pour tout élément $\sigma \in \Gal(E/F)$, il existe une matrice inversible $c_\sigma \in \GL(2d,E)$ telle que pour tout $x \otimes 1 \in M_2(D) \otimes_F E$, on a : \[ \sigma \circ \eta_D (x \otimes 1) = c_\sigma \big( \eta_D(x \otimes 1) \big) c_\sigma^{-1} .\] En particulier, $\det_0(x)$ est fixe par $\Gal(E/F)$ pour tout $x \in M_2(D)$, et $\det_0$ est donc à valeurs dans $F$. On définit alors la restriction de $\det_0$ à $\GL(2,D)$, $\det_G : \GL(2,D) \to F^\times$, qui est un morphisme de groupes. Lorsque l'on fait la démarche précédente pour $\GL(1,D)$, $\det_G$ devient le morphisme de groupes $\Nrm: D^\times \to F^\times$ appelé norme réduite. \\
On remarquera que si $A \in \GL(2,D)$ est une matrice triangulaire supérieure de coefficients diagonaux $(a_1,a_2)$, alors on a 
\begin{equation} \label{detNrm} \det\nolimits_G(A) = \Nrm(a_1) \cdot \Nrm(a_2) .\end{equation}
Pour plus de détails, on pourra se reporter au \S 12 $\mathrm{n^o}$3 de \cite{Bou58II8}, au \S 4 de  \cite{Bou58II8}, ainsi qu'au paragraphe 16 de \cite{Pie82}. \\

Un autre morphisme dérivé du déterminant va aussi entrer en jeu lorsque l'on s'intéressera aux représentations de $\GL(2,\Oo_D)$. Il s'agit de la composée \[ \overline{\det} : \GL(2,\Oo_D) \to \GL(2,k_D) \xrightarrow{\det} k_D^\times ,\] où la première application est la projection canonique $K \to K/K(1)$, qui correspond concrètement à la réduction modulo $(\varpi)$ coefficient par coefficient. \\ On fait la remarque que si on note $N_{k_D/k_F}$ la norme de l'extension galoisienne $k_D/k_F$, alors on peut vérifier la compatibilité que résume le diagramme commutatif suivant:
\[\xymatrix{ \GL(2,\Oo_D) \ar[d]_{\det_G} \ar[r]^{\phantom{xxx} \overline{\det}} &  k_D^\times \ar[d]^{N_{k_D/k_F}} \\ \Oo_F^\times \ar[r]^{\overline{\cdot}} & k_F^\times } .\] Pour une preuve, on pourra se reporter à la Proposition 16.1 du chapitre \ref{Ly12}.

\subsection{Représentations irréductibles de $\GL(2,k_D)$} \label{repfini}

On rappelle ici les principales propriétés des représentations irréductibles de $\overline{K} = \GL(2,k_D)$. On se référera au paragraphe 1 de \cite{BarLiv94} pour les preuves. \\

On commence par exhiber quelques représentations de $\overline{K}$. On peut former les puissances symétriques $\Sym^r \, \Fpbar^2$ de la représentation standard de $\overline{K}$ sur $\Fpbar^2$, pour $r \geq 0$. L'espace $\Sym^r \, \Fpbar^2$ est celui des polynômes homogènes de degré $r$ en deux variables; il a donc dimension $(r+1)$ et admet \[ X^r, X^{r-1} Y, X^{r-2} Y^2, \dots, X Y^{r-1}, Y^r \] 
comme base sur $\Fpbar$. L'action de $\overline{K}$ sur cette base est donnée par: 
\[ \left( \begin{array}{cc} a & b \\ c & d \end{array} \right) X^{r-i} Y^i = (aX+cY)^{r-i} (bX+dY)^i .\]
Ensuite, on peut tordre cette action par le Frobenius absolu $\varphi_0 : x \mapsto x^p$ et ses puissances: l'action de $\overline{K}$ sur $\big( \Sym^r \, \Fpbar^2 \big)^{\varphi_0^j}$ est donnée par 
\[ \left( \begin{array}{cc} a & b \\ c & d \end{array} \right) X^{r-i} Y^i = (a^{p^j} X+c^{p^j} Y)^{r-i} (b^{p^j} X+d^{p^j} Y)^i .\]
Pour tout vecteur $\vec{r} = (r_0, \dots, r_{fd-1}) \in [0,p-1]^{fd}$, on forme la représentation 
\[ \Sym^{\vec{r}} \, \Fpbar^2 = \bigotimes_{j=0}^{fd-1} \big( \Sym^{r_j} \, \Fpbar^2 \big)^{\varphi_0^j} \]
de dimension $\prod_j (r_j+1)$ de $\overline{K}$. Enfin on peut tordre cette représentation par un caractère de $\overline{K}$. Il se trouve que ce procédé permet de construire toutes les représentations irréductibles de $\overline{K}$.

\begin{prop}
Les classes d'isomorphisme de représentations irréductibles de $\overline{K}$ sont en bijection via 
\[ (\vec{r},\chi) \mapsto V(\vec{r},\chi) = ( \chi \circ \det) \otimes \Sym^{\vec{r}} \, \Fpbar^2 \]
avec les paires $(\vec{r},\chi)$ où $\vec{r}$ est un élément de $[0,p-1]^{fd}$ et $\chi$ est un caractère $k_D^\times \to \Fpbar^\times$.
\end{prop}

Il est important de décrire l'espace d'invariants ou de coinvariants de ces représentations par 
\[ \overline{U} = (U \cap K)/(U \cap K(1)) , \quad \overline{U}^- = (U^- \cap K)/(U^- \cap K(1)) .\]
De même, on note 
\[ \overline{B} = (B \cap K)/(B \cap K(1)) , \quad \overline{B}^- = (B^- \cap K)/(B^- \cap K(1)) .\]
Enfin, pour $\vec{s} = (s_0, \dots, s_{fd-1}) \in [0,p-1]^{fd}$ avec $s_i \leq r_i$ pour tout $i$, on note 
\[ X^{\vec{s}} = X^{s_0} \otimes X^{s_1} \otimes \dots \otimes X^{s_{fd-1}} \in \Sym^{\vec{r}} \, \Fpbar^2 ,\]
et de même pour $Y^{\vec{s}}$; et pour $x \in k_D^\times$, on écrit 
\begin{equation} \label{caractkD} x^{\vec{s}} = x^{s_0}(x^p)^{s_1} \dots (x^{p^{fd-1}})^{s_{fd-1}} = x^{s_0 + p s_1 + \dots + p^{fd-1} s_{fd-1}} .\end{equation}

\begin{lem} \label{coinvV}
Soient $\vec{r} \in [0,p-1]^{fd}$ et $\chi$ un caractère de $k_D^\times$.
\begin{itemize}
\item[(i)] L'espace de coinvariants $V(\vec{r},\chi)_{\overline{U}}$ est isomorphe à l'espace d'invariants $V(\vec{r},\chi)^{\overline{U}^-}$. Il est porté par la droite $\Fpbar Y^{\vec{r}}$, sur laquelle l'action de $\overline{B}^-$ est donnée par \[ \left( \begin{array}{cc} a & 0 \\ c & d \end{array} \right) Y^{\vec{r}} = \chi(ad) d^{\vec{r}} Y^{\vec{r}} .\]
\item[(ii)] L'espace d'invariants $V(\vec{r},\chi)^{\overline{U}}$ est isomorphe à l'espace de coinvariants $V(\vec{r},\chi)_{\overline{U}^-}$. Il est porté par la droite $\Fpbar X^{\vec{r}}$, sur laquelle l'action de $\overline{B}$ est donnée par \[ \left( \begin{array}{cc} a & b \\ 0 & d \end{array} \right) X^{\vec{r}} = \chi(ad) a^{\vec{r}} X^{\vec{r}} .\]
\end{itemize}
\end{lem}

Les $fd$-uplets $\vec{0} = (0, \dots, 0)$ et $\vec{p-1} = (p-1, \dots, p-1)$ de $[0,p-1]^{fd}$ jouent un rôle particulier. En effet, fixons un caractère $\chi$ de $k_D^\times$; les caractères $\left( \begin{array}{cc} a & 0 \\ c & d \end{array} \right) \mapsto \chi(ad) d^{\vec{r}}$ par lesquels $\overline{B}^-$ agit sur $V(\vec{r},\chi)_{\overline{U}}$ sont deux à deux distincts sauf pour $\vec{r} = \vec{0}$ et $\vec{r} = \vec{p-1}$. Aussi, suivant la dénomination de \cite{HenVig11b}, on dira que les $V(\vec{r},\chi)$ sont alors $\overline{B}^-$-réguliers pour $\vec{r} \neq \vec{0}$, et non $\overline{B}^-$-réguliers pour $\vec{r} = \vec{0}$. On remarque que, dans ce cas particulier de $\GL(2,k_D)$, une représentation irréductible de $\overline{K}$ est $\overline{B}^-$-régulière si et seulement si elle est $\overline{B}$-régulière: on dira donc tout simplement que $V(\vec{r},\chi)$ est \textit{régulière} pour $\vec{r} \neq \vec{0}$. Et c'est aussi équivalent au fait que $V$ est de dimension strictement supérieure à $1$. \\

Soit $V$ une représentation lisse irréductible de $K$. Comme $K(1)$ est un pro-$p$-groupe, l'espace $V^{K(1)}$ des $K(1)$-invariants de $V$ est non réduit à $0$ (voir \cite{BarLiv94}, Lemma 3). Et parce que $K(1)$ est normal dans $K$, $V^{K(1)}$ est une $K$-représentation qui est donc égale à $V$ par irréductibilité. Ainsi l'action de $K$ sur $V$ se factorise par $\overline{K} = K/K(1)$. Réciproquement, toute représentation irréductible de $\overline{K}$ définit une représentation lisse irréductible de $K$ par inflation. De ce fait, il n'y a pas lieu dans la suite de faire de distinction entre une $\overline{K}$-représentation irréductible et une $K$-représentation lisse irréductible.

\section{Représentations irréductibles de $D^\times$} \label{Dtimes}

Soit $V$ une représentation irréductible de dimension finie de $D^\times$ sur $\Fpbar$. Commençons par remarquer que $D(1):= 1+ \varpi \Oo_D$ est un pro-$p$-groupe, de sorte que $V^{D(1)}$ est non réduit à $0$ (voir \cite{BarLiv94}, Lemma 3). De plus, $D(1)$ est normalisé par $D^\times$ et donc $D^\times$ agit sur $V^{D(1)}$; par irréductibilité de $V$, l'action de $D^\times$ se factorise\footnote{dans le langage de la théorie complexe, on peut dire que toute représentation irréductible est de niveau $0$} par $D^\times/D(1)$, et $V$ s'identifie à une représentation de $D^\times/D(1)$. Par la suite, on se servira de ce fait pour confondre les représentations irréductibles de $D^\times$ et celles de $D^\times/D(1)$. \\
Rappelons que l'on a la décomposition suivante : 
\begin{equation} \label{decompunits} D^\times/D(1) \simeq k_D^\times \rtimes \varpi^{\Z} \simeq \Z/(q^d-1)\Z \rtimes \Z \end{equation} 
où l'action de $\varpi$ sur $k_D^\times$ est donnée par $x \mapsto \overline{\varpi [x] \varpi^{-1}} = x^{q_D}$ en notant $q_D := q^{a_D}$. \\

Soient $\sigma$ un caractère de $k_D^\times$ et $N_\sigma$ le stabilisateur de $\sigma$ dans $D^\times/D(1)$ pour l'action via conjugaison: $\delta \in D^\times/D(1)$ envoie $\sigma$ sur ${}^\delta \sigma : x \mapsto \sigma(\delta^{-1} x \delta)$. Ce stabilisateur $N_\sigma$ contient $k_D^\times \varpi^{d \Z}$ et est donc égal à $k_D^\times \rtimes \varpi^{d_0 \Z}$ pour un certain entier $d_0 \geq 1$ divisant $d$. Prolonger $\sigma$ en un caractère de $N_\sigma$ revient à choisir l'image de $\varpi^{d_0}$ et, pour tout $\eta \in \Fpbar^\times$, on note $\widetilde{\sigma}_\eta$ le caractère défini par: \[ \widetilde{\sigma}_\eta|_{k_D^\times} = \sigma \quad \textrm{et} \quad \widetilde{\sigma}_\eta(\varpi^{d_0})= \eta .\]
On forme alors l'induite \[ \rho(\sigma,\eta) := \Ind_{N_\sigma}^{D^\times/D(1)} \, \widetilde{\sigma}_\eta ,\] que l'on voit aussi bien comme représentation de $D^\times/D(1)$ ou comme représentation de $D^\times$.

\begin{prop} \label{propDtimes}
Soient $\sigma$ et $\tau$ deux caractères de $k_D^\times$. Soient $\eta, \nu \in \Fpbar^\times$.
\begin{itemize}
\item[(i)] La représentation $\rho(\sigma,\eta)$ est irréductible.
\item[(ii)] Les induites $\rho(\sigma,\eta)$ et $\rho(\tau,\nu)$ sont isomorphes si et seulement si les deux conditions suivantes sont satisfaites:
\begin{itemize}
\item[(a)] il existe $\delta \in D^\times/D(1)$ tel que l'on ait $\tau = {}^\delta \sigma$;
\item[(b)] on a l'égalité $\eta = \nu$.
\end{itemize}
\end{itemize}
\end{prop}
\textsf{Preuve :} 
Voir \cite{Ser98}, paragraphe 8.2. \hfill$\Box$\\

En comparaison avec $\GL(2,F)$, notre description précédente de $\rho(\sigma,\eta)$ est semblable à une induite parabolique. Il nous sera utile plus tard de la voir comme une induite compacte aussi. \\
Pour $\sigma$ un caractère de $\Oo_D^\times$, on remarque l'isomorphisme de $\Fpbar$-espaces vectoriels suivant pour décrire l'algèbre de Hecke $\Hh(D^\times,\Oo_D^\times,\sigma)$, où le second isomorphisme est la réciprocité de Frobenius
\begin{equation} \label{HeckeD} \begin{array}{ccccc} \Fpbar[\varpi^{d_0 \Z}] & \xrightarrow{\sim} & \Hom_{\Oo_D^\times} \big( \sigma, \ind_{\Oo_D^\times}^{D^\times} \, \sigma \big) & \xrightarrow{\sim} & \Hh(D^\times, \Oo_D^\times, \sigma) \\ \varpi^{d_0 i} & \mapsto & \big( 1 \mapsto [\varpi^{- d_0 i} , 1] \big) & & \end{array} .\end{equation}
Parce que $\Oo_D^\times$ est normal dans $D^\times$, les opérateurs de $\Hh(D^\times,\Oo_D^\times,\sigma)$ sont portées par des classes simples et il s'ensuit que (\ref{HeckeD}) est en fait un isomorphisme de $\Fpbar$-algèbres. En identifiant $\Hh(D^\times,\Oo_D^\times,\sigma)$ à $\Fpbar[\varpi^{d_0 \Z}]$ via (\ref{HeckeD}), pour tout $\eta \in \Fpbar^\times$, on définit un caractère de $\Fpbar$-algèbres de $\Hh(D^\times, \Oo_D^\times, \sigma)$ en posant $\chi_\eta(\varpi^{d_0 i}) = \eta^i$ pour tout $i \in \Z$.

\begin{prop} \label{Dtimescompact}
Soient $\sigma$ un caractère de $k_D^\times$, que l'on voit comme un caractère de $\Oo_D^\times$ par inflation, et $\eta \in \Fpbar^\times$. On a l'isomorphisme de représentations de $D^\times$ 
\[ \ind_{\Oo_D^\times}^{D^\times} \, \sigma \otimes_{\chi_\eta} \Fpbar \xrightarrow{\sim} \rho(\sigma,\eta) .\]
\end{prop}
\textsf{Preuve :} \\
L'inclusion $\sigma \subseteq \rho(\sigma,\eta)$ induit par réciprocité de Frobenius une surjection $D^\times$-équivariante $\ind_{\Oo_D^\times}^{D^\times} \, \sigma \to \rho(\sigma,\eta)$ (puisque $\rho(\sigma,\eta)$ est irréductible). De plus, cette dernière se factorise par \[ \ind_{\Oo_D^\times}^{D^\times} \, \sigma \otimes_{\chi_\eta} \Fpbar \to \rho(\sigma,\eta) ,\]
qui est un isomorphisme car les deux termes sont de dimension $d_0$ sur $\Fpbar$. \hfill$\Box$ \\

Décrivons la restriction d'une induite $\rho(\sigma,\eta)$ de dimension $d_0$ à $k_D^\times$. Parce que $k_D^\times$ est cyclique d'ordre premier à $p$ et que $\Fpbar$ est algébriquement clos, l'action de $k_D^\times$ sur $\rho(\sigma,\eta)$ est diagonalisable.

\begin{lem} \label{generatedfield}
Soient $z_0$ un vecteur propre de $k_D^\times$ dans $\rho(\sigma,\eta)$ et $\xi_0 \in \Fpbar^\times$ le scalaire tel que l'on ait $\mu.z_0 = \xi_0 z_0$. Alors $\xi_0$ appartient à $\F_{q^{d_0}}$ et l'engendre en tant qu'extension sur $\F_q$.
\end{lem}
\textsf{Preuve :} \\
Il s'agit de voir que $\xi_0$ est un élément de degré $d_0$ sur $\F_q$. Or, par (\ref{decompunits}), le fait que $a_D$ est premier à $d$ et la définition de $d_0$, le sous-groupe de $\Gal(\Fpbar/\F_q)$ fixant $\xi_0$ est topologiquement engendré par $x \mapsto x^{q^{d_0}}$. Le lemme est prouvé. \hfill$\Box$\\

On remarque que l'on oubliera souvent par la suite de le préciser mais on ne considère que des représentations irréductibles de dimension finie de $D^\times$. Aussi, à cause du lemme de Schur et de (\ref{decompunits}), toute représentation irréductible admissible de $D^\times$ est de dimension finie. Parmi ces représentations irréductibles de $D^\times$ que l'on vient de décrire, les caractères jouent un rôle particulier. Précisons un peu leur structure.

\begin{lem} \label{Ddim1}
Soit $\rho$ une représentation irréductible de $D^\times$. Alors $\rho$ est un caractère si et seulement si il existe un caractère $\rho_0 : F^\times \to \Fpbar^\times$ tel que l'on ait $\rho = \rho_0 \circ \Nrm$.
\end{lem}
\textsf{Preuve :} \\
Une représentation $\rho_0 \circ \Nrm$ est un caractère de $D^\times$. Montrons que tout caractère $\rho$ de $D^\times$ possède une factorisation $\rho = \rho_0 \circ \Nrm$. Lorsque $\rho$ est un caractère, pour tout $x$ de $\Oo_D^\times$, on a $\rho(\varpi x \varpi^{-1}) = \rho(x)$. Comme l'action par conjugaison de $\varpi$ sur $k_D^\times$ est donnée par $y \mapsto y^{q_D}$, on a alors la factorisation
\begin{equation} \label{rhoentier} \rho |_{\Oo_D^\times} : \Oo_D^\times \twoheadrightarrow k_D^\times \xrightarrow{N_{k_D/k_F}} k_F^\times \xrightarrow{\overline{\rho}_0} \Fpbar^\times ,\end{equation}
où la première flèche est la réduction canonique et $N_{k_D/k_F}$ est la norme de l'extension galoisienne $k_D/k_F$. Par la preuve de \cite{MatNak43}, Satz 1, (\ref{rhoentier}) se réécrit 
\[\xymatrix{ \rho |_{\Oo_D^\times} : \Oo_D^\times \ar[r]^{\phantom{xxx} \Nrm} & \Oo_F^\times \ar@{->>}[r] \ar@/^2pc/[rr]^{\rho_0} & k_F^\times \ar[r]^{\overline{\rho}_0} & \Fpbar^\times },\]
où $\Oo_F^\times \twoheadrightarrow k_F^\times$ est la réduction canonique. Il reste à prolonger $\rho_0$ à $F^\times$ en posant \[ \rho_0 \big( (-1)^d \varpi_F \big) = \rho_0 \big( \Nrm(\varpi) \big) = \rho(\varpi) .\] La preuve est terminée. \hfill$\Box$

\section{Généralités et $I(1)$-invariants} \label{parI(1)}

Soient $\rho_1$ et $\rho_2$ deux représentations lisses irréductibles de $D^\times$, de dimension respective $d_1$ et $d_2$. On commence par utiliser le paragraphe \ref{Dtimes} pour expliciter des bases pour les espaces sous-jacents à $\rho_1$ et $\rho_2$ dans lesquelles l'action de $k_D^\times \hookrightarrow D^\times/D(1)$ est diagonale (voir (\ref{decompunits})). \\
Soient $v_0$ un vecteur propre dans $\rho_1$ pour l'action de $k_D^\times$ et $\xi_1 \in \F_{q^{d_1}}^\times \subseteq \Fpbar^\times$ la valeur propre associée à $\mu$ (voir lemme \ref{generatedfield}). On définit $v_a = \varpi^{-a}.v_0$ pour $1 \leq a \leq d_1 - 1$, de sorte que $(v_0, \dots, v_{d_1-1})$ forme une base de $\rho_1$; et l'action de $\mu$ est diagonale dans cette base, de valeur $\big( \xi_1, \xi_1^{q_D}, \dots, \xi_1^{q_D^{d_1-1}} \big)$. \\
De la même manière, pour $\rho_2$, on a un vecteur $w_0$ tel que $\mu.w_0 = \xi_2 w_0$ pour un $\xi_2 \in \F_{q^{d_2}}^\times \subseteq \Fpbar^\times$. Si on note $w_b = \varpi^{-b}.w_0$ pour $1 \leq b \leq d_2-1$, $(w_0, \dots, w_{d_2-1})$ est une base de $\rho_2$ dans laquelle l'action de $\mu$ est diagonale, donnée par $\big( \xi_2, \xi_2^{q_D}, \dots, \xi_2^{q_D^{d_2-1}} \big)$.

\begin{lem} \label{torscaract}
Soient $\xi$ un caractère lisse de $F^\times$ et $\rho_1$, $\rho_2$ deux représentations irréductibles de $D^\times$. Alors on a l'isomorphisme de $G$-représentations \[ (\xi \circ \det\nolimits_G) \circ \Ind_B^G \, \rho_1 \otimes \rho_2 \xrightarrow{\sim} \Ind_B^G \big( ((\xi \circ \Nrm) \otimes \rho_1) \otimes ((\xi \circ \Nrm) \otimes \rho_2) \big) .\]
\end{lem}
\textsf{Preuve :} \\
Le morphisme de $G$-représentations \[ \begin{array}{ccc} (\xi \circ \det_G) \otimes \Ind_B^G \, \rho_1 \otimes \rho_2 & \to & \Ind_B^G \big( (\xi \circ \det_G) \otimes (\rho_1 \otimes \rho_2) \big) \\ a \otimes f & \mapsto & (g \mapsto a \, \xi \circ \det_G(g) \, f(g)) \end{array} \] admet naturellement \[ f \mapsto 1 \otimes (g \mapsto \xi^{-1} \circ \det\nolimits_G(g) \, f(g)) \] pour inverse. Le résultat suit de (\ref{detNrm}). \hfill$\Box$\\

L'action de $\varpi^{d_1}$ sur $\rho_1$ est scalaire par la proposition \ref{propDtimes}. Choisissons $\xi : F^\times \to \Fpbar^\times$ tel que $\xi(\Nrm \, \varpi)^{d_1} = \rho_1(\varpi^{d_1}) \in \Fpbar^\times$ et $\xi |_{\Oo_F^\times}$ trivial: par le lemme \ref{torscaract}, on est donc ramené au cas où l'action de $\varpi^{d_1}$ sur $\rho_1$ est triviale, et c'est ce qu'on supposera par la suite. \\
Soit $V = \Ind_B^G \, \rho_1 \otimes \rho_2$. On va tâcher dans la suite de cette section de préparer l'étude de l'irréductibilité de $V$ en tant que $G$-représentation en énonçant des généralités sur $V^{I(1)}$. \\

Comme $I(1)$ est un pro-$p$-groupe, $V^{I(1)}$ est non réduit à $0$ (voir \cite{BarLiv94}, Lemma 3). C'est donc un objet intéressant à considérer, d'autant plus que, d'après la réciprocité de Frobenius, 
\begin{equation} \label{moddroite}  V^{I(1)} \simeq \Hom_{I(1)}(\id,V) \simeq \Hom_G(\ind_{I(1)}^G \, \id, V) \end{equation} 
est naturellement muni d'une structure de module à droite sur $\Hh(G,I(1))$.
Explicitement, pour $g \in G$, l'opérateur $T_g$ de $\Hh(G,I(1))$ de support $I(1) g I(1)$ et valant $\Id$ en $g$ agit sur $V^{I(1)}$ par 
\begin{equation} \label{I(1)agit} v \mapsto v T_g = \sum_{\gamma \in I(1) \backslash I(1) g I(1)} \gamma^{-1} v = \sum_{u \in (I(1) \cap g^{-1} I(1) g) \backslash I(1)} u^{-1} g^{-1} v .\end{equation}
Un tel opérateur $T_g$ est caractérisé par son support et on dira parfois simplement que $T_g$ est la fonction caractéristique de $I(1) g I(1)$. \\
Regardons l'espace des $I(1)$-invariants de $V$: en utilisant successivement les décomposition d'Iwasawa et de Bruhat, on a \[ G = B K = B I(1) \coprod B I(1) s I(1) .\] Aussi, on remarque les identités \[ B I(1) s I(1) = B \big( U^- \cap I(1) \big) s I(1) = B s I(1) = B \omega I(1) ,\] de sorte que l'on a finalement 
\begin{equation} \label{suppF} G = B I(1) \coprod B \omega I(1) .\end{equation} 
On fait remarquer que l'on préfère écrire $B \omega I(1)$ plutôt que $BsI(1)$ car $\omega$ a le bon goût de normaliser $I(1)$. \\
Comme l'espace sous-jacent à $\rho_1 \otimes \rho_2$ est de dimension $d_1 d_2$, on en déduit que $V^{I(1)}$ est de dimension $e := 2 d_1 d_2$. Pour $0 \leq a \leq d_1 -1$ et $0 \leq b \leq d_2 -1$ et $i = 0,1$, soit $f^i_{a,b} \in V$ la fonction $I(1)$-invariante de support $B \omega^i I(1)$ et prenant la valeur $v_a \otimes w_b$ en $\omega^i$. La famille $\Ff := (f^i_{a,b})_{i,a,b}$ (ordonnée lexicographiquement) forme une base du $\Fpbar$-espace vectoriel $V^{I(1)}$. 

\begin{prop} \label{Indadmiss}
Soient $\rho_1$ et $\rho_2$ deux représentations irréductibles de $D^\times$. L'induite parabolique $\Ind_B^G \, \rho_1 \otimes \rho_2$ est admissible.
\end{prop}
\textsf{Preuve :} \\
Par la discussion précédente, $(\Ind_B^G \, \rho_1 \otimes \rho_2)^{I(1)}$ est de dimension finie, égale à $2 d_1 d_2$. Comme $I(1)$ est un pro-$p$-sous-groupe ouvert du groupe localement profini $G$, le résultat suit de \cite{Pas04}, Theorem 6.3.2. \hfill$\Box$

\begin{lem} \label{critirred}
Soient $\rho_1$ et $\rho_2$ deux représentations irréductibles de $D^\times$.
\begin{itemize}
\item[(i)] La représentation $\Ind_B^G \, \rho_1 \otimes \rho_2$ est générée en tant que $B$-représentation par ses $I(1)$-invariants.
\item[(ii)] Soit $\pi$ une représentation de $G$ générée par ses $I(1)$-invariants. Si $\pi^{I(1)}$ est simple en tant que module à droite sur $\Hh(G,I(1))$, alors $\pi$ est irréductible.
\end{itemize}
\end{lem}
\textsf{Preuve :} \\
Pour le $(ii)$, voir \cite{Vig04} Criterium 4.5. Le $(i)$ se trouve dans \cite{Vig08}, Proposition 9 pour une induite de caractères; cependant, la preuve reste valable pour une induite parabolique de représentation de dimension finie\footnote{et même en général: on a uniquement besoin de la surjectivité de $\begin{array}{ccc} \big( \Ind_B^G \, \rho \big)^{I(1)} & \to & \rho \times \rho \\ f & \mapsto & (f(1),f(\omega)) \end{array}$}. Réexpliquons l'argument, qui est d'autant plus simplifié que le groupe ici considéré est de rang relatif $1$. \\ 
Soit $f$ une fonction non nulle de $V = \Ind_B^G \, \rho_1 \otimes \rho_2$. Rappelons la décomposition 
\begin{equation} \label{decompBU} G = B \coprod B sU = B \coprod B \omega U ,\end{equation} 
dérivée de la décomposition de Bruhat. Parce que $B U^-$ est dense dans $G$, si la valeur de $f$ en $1$ n'est pas nulle, il existe $k_0 \in \Z$ tel que $B \overline{n}(\varpi^{k_0} \Oo_D)$ est dans le support de $f$. Aussi, on a \[ B \overline{n}(\varpi^{k_0} \Oo_D) = B \varphi_2^{k_0-1} I(1) \varphi_2^{1-k_0} = B I(1) \varphi_2^{1-k_0} .\] Il existe alors des $\lambda^0_{a,b} \in \Fpbar$ tels que la fonction 
\begin{equation} \label{defg} g:= f- \varphi_2^{k_0-1} . \Big( \sum_{a,b} \lambda^0_{a,b} f^0_{a,b} \Big) \end{equation} 
prenne la valeur $0$ en $1$; par la décomposition (\ref{decompBU}), elle est donc à support dans $B \omega U$. \\
Soit $k_1 \geq 0$ un entier tel que $n(\varpi^{k_1} \Oo_D) \subseteq U$ soit un sous-groupe du stabilisateur de $g$ dans $G$: un tel $k_1$ existe par lissité de $V$. Ecrivons le support de $g$ sous la forme 
\[ \coprod\nolimits_{i \in \I_g} B \omega n(a_i + \varpi^{k_1} \Oo_D) ,\] où $\I_g$ est fini et les $a_i \in D$ sont distincts modulo $(\varpi^{k_1})$; on remarque que $g$ est constante sur chaque $\omega n(a_i + \varpi^{k_1} \Oo_D)$ par choix de $k_1$. Aussi, on a pour tout $i \in \I_g$: 
\begin{equation} \label{suppf1} B \omega n(a_i + \varpi^{k_1} \Oo_D) = B \omega \varphi_1^{k_1} I(1) \varphi_1^{-k_1}n(a_i) = B \omega I(1) \varphi_1^{- k_1} n(a_i) ;\end{equation} 
et cet ensemble est le support de la fonction $n(-a_i) \varphi_1^{k_1} . f_{a,b}^1$. Parce que les supports (\ref{suppf1}) sont distincts, il existe des scalaires $\lambda^i_{a,b} \in \Fpbar$ tels que la fonction
\begin{equation} \label{defg1} g_1 := g- \sum_{i \in \I_g} n(-a_i) \varphi_1^{k_1} . \Big( \sum_{a,b} \lambda^i_{a,b} f^1_{a,b} \Big) \end{equation} 
s'annule sur les $B \omega n(a_i+\varpi^{k_1} \Oo_D)$ pour tout $i \in \I_g$: la fonction $g_1$ est identiquement nulle. En observant (\ref{defg}) et (\ref{defg1}), le résultat est prouvé. \hfill$\Box$\\

Le pro-$p$-Iwahori $I(1)$ vérifie bien le fait que toute double classe $I(1) g I(1)$ est union finie de classes simples à droite. De ce fait, par le paragraphe \ref{parGeneralHecke}, on peut voir $\Hh(G,I(1))$ comme la $\Fpbar$-algèbre des fonctions sur $G$ à valeurs dans $\Fpbar$ et qui sont bi-$I(1)$-invariantes et à support compact; et c'est ce que l'on fera par la suite. Ainsi $\Hh(G,I(1))$ est engendrée vectoriellement par les fonctions caractéristiques des doubles classes $I(1) g I(1)$ pour $g \in G$, que l'on notera $\id_{I(1) g I(1)}$. Précisons un peu les générateurs de $\Hh(G,I(1))$ en tant que $\Fpbar$-algèbre. \\
Notons $C$ le sous-groupe de $G$ composé des matrices monomiales (c'est-à-dire ayant un seul coefficient non nul par ligne et par colonne). La décomposition de Bruhat pour la BN-paire $(I,C)$ nous donne 
\[ G = \coprod_{a,b \in \Z , \, i \in \{0,1\}} I \left( \begin{array}{cc} \varpi^a & 0 \\ 0 & \varpi^b \end{array} \right) s^i I = \coprod_{a,b \in \Z , \, i \in \{0,1\}} I(1) \left( \begin{array}{cc} \varpi^a & 0 \\ 0 & \varpi^b \end{array} \right) s^i I .\] 
On peut réécrire cela en 
\begin{equation} \label{BruhatI(1)} G = \coprod_{b,i,j,x,y} I(1) \omega^i \varphi_2^b \omega^j \delta^x_y I(1) ,\end{equation} 
où les indices parcourent $b \geq 0$, $i \in \Z$, $j \in \{0,1\}$ et $x,y \in \Oo_D^\times$. Parce que $\omega$ et $A \cap K$ normalisent $I(1)$, on a les égalités:
\begin{equation} \label{double=simple} I(1) \omega I(1) = \omega I(1) = I(1) \omega, \quad I(1) \delta^x_y I(1) = \delta^x_y I(1) = I(1) \delta^x_y \end{equation}
pour tous $x, y \in \Oo_D^\times$. On a alors 
\begin{equation} \label{doubleI(1)prod} I(1) \omega^i \varphi_2^b \omega^j \delta^x_y I(1) = (I(1) \omega I(1))^i (I(1) \varphi_2^b I(1)) (I(1) \omega I(1))^j (I(1) \delta^x_y I(1)) .\end{equation}
Il en découle l'énoncé suivant.

\begin{lem} \label{generateHeckeI}
La $\Fpbar$-algèbre de Hecke-Iwahori $\Hh(G,I(1))$ est de type fini, engendrée par les opérateurs $\id_{\omega^{\pm 1} I(1)}$, $\id_{I(1)sI(1)}$ et les $\id_{\delta^{[x]}_{[y]} I(1)}$ pour $x, y \in k_D^\times$.
\end{lem}
\textsf{Preuve :} \\
Grâce à (\ref{BruhatI(1)}), on sait que $\Hh(G,I(1))$ est engendrée vectoriellement par les fonctions caractéristiques des doubles classes (\ref{doubleI(1)prod}). Aussi, les égalités (\ref{doubleI(1)prod}) et (\ref{double=simple}) impliquent que $\Hh(G,I(1))$ est engéndrée en tant qu'algèbre par $\id_{\omega^{\pm 1} I(1)}$, les $\id_{I(1) \varphi_2^b I(1)}$ pour $b \geq 0$ et les $\id_{\delta^x_y I(1)}$ pour $x, y \in \Oo_D^\times$. En effet, on voit par exemple que l'on a, par définition (\ref{defconvol}) pour $g \in G$ : 
\begin{equation} \label{simplprod} \id_{I(1) \omega^i I(1)} * \id_{I(1) \varphi_2^b I(1)} (g) = \id_{\omega^i I(1)} (\omega^i) \id_{I(1) \varphi_2^b I(1)} (\omega^{-i} g) = \id_{I(1) \omega^i \varphi_2^b I(1)} (g)  ;\end{equation}
et les autres produits voulus s'obtiennent de manière identique. De plus, en utilisant successivement l'identité 
\[ I(1) \varphi_2 I(1) = \coprod_{x \in k_D} I(1) \varphi_2 n([x]) \]
et parce que l'on a l'égalité d'ensembles (pour $b \geq 1$)
\[ \{ \varphi_2^b n([x]) \ | \ x \in k_D^b \} = \{ \varphi_2 n([x_0]) \varphi_2 n([x_1]) \varphi_2 \dots \varphi_2 n([x_{b-1}]) \ | \ x_0, \dots, x_{b-1} \in k_D \} ,\]
on en déduit l'égalité de doubles classes \[ I(1) \varphi_2^b I(1) = \coprod_{x \in k_D^b} I(1) \varphi_2^b n([x]) = \big( I(1) \varphi_2 I(1) \big)^b .\] 
De cette manière, le support de la $b$-ème convolée \[ \big( \id_{I(1) \varphi_2 I(1)} \big)^{* b} := \id_{I(1) \varphi_2 I(1)} * \id_{I(1) \varphi_2 I(1)} * \cdots * \id_{I(1) \varphi_2 I(1)} \] est inclus dans $I(1) \varphi_2^b I(1)$. On établit enfin par récurrence sur $b \geq 1$ que sa valeur en $\varphi_2^b$ est $1$. Pour $b \geq 1$, c'est la définition; supposons donc l'affirmation vraie pour $b \geq 1$ et montrons-la au rang $b+1$: on a \[ \big( \id_{I(1) \varphi_2 I(1)} \big)^{* (b+1)} (\varphi_2^{b+1}) = \id_{I(1) \varphi_2^b I(1)} * \id_{I(1) \varphi_2 I(1)}(\varphi_2^{b+1}) \] par hypothèse de récurrence. La définition (\ref{defconvol}) donne ensuite 
\begin{align*} \big( \id_{I(1) \varphi_2 I(1)} \big)^{*(b+1)} & = \sum_{g \in I(1) \backslash I(1) \varphi_2 I(1)} (\varphi_2^{b+1} g^{-1}) \id_{I(1) \varphi_2 I(1)}(g) \\ & = \big| \{ x \in k_D \ | \ \varphi_2^{b+1} n([x]) \varphi_2^{-1} \in I(1) \varphi_2^b I(1) \} \big| .\end{align*} 
Or $\varphi_2^{b+1} n([x]) \varphi_2^{-1} = \varphi_2^b n([x] \varpi^{-1})$ appartient à $I(1) \varphi_2^b I(1)$ si et seulement si $x$ est nul. La quantité $\big( \id_{I(1) \varphi_2 I(1)} \big)^{* (b+1)} (\varphi_2^{b+1})$ est donc bien égale à $1$ et l'affirmation est prouvée par récurrence. En particulier, on sait à présent que $\id_{\omega^{\pm 1} I(1)}$, $\id_{I(1) \varphi_2 I(1)}$ et les $\id_{\delta^x_y I(1)}$ pour $x, y \in \Oo_D^\times$ génèrent l'algèbre $\Hh(G,I(1))$. Il reste à remarquer $\delta^x_y I(1) = \delta^{[\overline{x}]}_{[\overline{y}]} I(1)$ pour se limiter aux $\id_{\delta^{[x]}_{[y]} I(1)}$; et à la manière de (\ref{simplprod}), grâce (\ref{double=simple}) et à $\varphi_2 = \omega s$, on a l'identité \[ \id_{I(1) \varphi_2 I(1)} = \id_{\omega I(1)} * \id_{I(1) s I(1)} .\] Le résultat est prouvé. \hfill$\Box$\\

Remarquons que l'on aurait très bien pu garder $\id_{I(1) \varphi_2 I(1)}$ au lieu de $\id_{I(1) s I(1)}$ parmi les générateurs exhibés de $\Hh(G,I(1))$ mais il se trouve que la littérature a tendance à favoriser le second (car il fait sens dans un contexte de groupes finis aussi). \\

Pour $n \geq 1$ et $\nu \in \Fpbar^\times$, soit $C_n (\nu) := \left( \begin{array}{cccc} 0 & 1 & & \\ & \ddots & \ddots & \\ & & \ddots & 1 \\ \nu^{-1} & & & 0 \end{array}\right)$ la matrice cyclique et $1_n$ la matrice identité de $M_n(\Fpbar)$. On pose alors \[ H_\omega(\nu) := \left( \begin{array}{cc} 0 & 1_{d_1 d_2} \\ C_{d_1}(1) \otimes C_{d_2}(\nu) & 0 \end{array}\right) ,\] qui est une matrice de $M_e(\Fpbar)$. Pour $x,y \in \Oo_D^\times$, on note $\Delta(x,y) \in M_{d_1 d_2}(\Fpbar)$ la matrice diagonale de coefficients $\big( \xi_1^{\{\overline{x}^{q_D^{i_1}}\}} \xi_2^{\{\overline{y}^{q_D^{i_2}}\}} \big)_i$ pour $0 \leq i < d_1 d_2$ de division euclidienne $i = i_1 d_2 + i_2$, $0 \leq i_2 < d_2$. On pose alors 
\[ H_\delta(x,y) := \left( \begin{array}{cc} \Delta(x,y) & 0 \\ 0 & \Delta(y,x^{q_D}) \end{array} \right) \in M_e(\Fpbar) .\] Lorsque l'on a $d_1 = d_2$, on note $S_k \in M_{d_1 d_2}(\Fpbar)$, pour $0 \leq k < d_1$ un entier, la matrice définie par \[ S_k (v_a \otimes w_b) = \delta_{a+k,b+1} \, v_a \otimes w_b ,\] où $\delta$ désigne le symbole de Kronecker, valant $1$ si et seulement si $a+k$ est égal à $b+1$ modulo $d_1$. On définit alors, pour $\nu \in \Fpbar^\times$ et $\varepsilon \in \{ \pm 1\}$, \[ H_s(\nu, \varepsilon) := \begin{cases} \left( \begin{array}{cc} 0 & 0 \\ 1_{d_1} \otimes C_{d_2}(\nu) & \varepsilon S_k \end{array} \right) & \textrm{si } d_1 = d_2 \textrm{ et } \rho_1|_{\Oo_D^\times} = \rho_2(\varpi^k \cdot \varpi^{-k})|_{\Oo_D^\times}, \\ \left( \begin{array}{cc} 0 & 0 \\ 1_{d_1} \otimes C_{d_2}(\nu) & 0 \end{array} \right)  & \textrm{sinon.} \end{cases} \] 
On fait la remarque qu'un tel entier $0 \leq k < d_1$ est uniquement déterminé (s'il a lieu d'être). On rappelle aussi que l'on a posé $\lambda = \rho_2(\varpi^{d_2})^{-1} \in \Fpbar^\times$ au début du paragraphe \ref{parI(1)}. De plus, la valeur de $- \rho_1(-1)$ est un élément de $\{ \pm 1 \}$, que l'on va noter $\tau$. \\
Enfin, pour $x, y \in \Oo_D^\times$, on note de la manière suivante les opérateurs de $\Hh(G,I(1))$: \[ \Delta^x_y = \id_{(\delta^x_y)^{-1} I(1)}, \quad T^{(-1)} = \id_{\omega^{-1} I(1)}, \quad S = \id_{I(1) s I(1)} .\] Par le lemme \ref{generateHeckeI}, les $\Delta^x_y$, $(T^{(-1)})^{\pm 1}$ et $S$ engendrent $\Hh(G,I(1))$. On va maintenant calculer l'action de ces opérateurs sur $V^{I(1)}$.

\begin{lem} \label{I(1)actions}
Soient $\rho_1, \rho_2, V, \lambda, \tau$ comme précédemment. Soient $x$ et $y$ des éléments de $\Oo_D^\times$. Dans la base $\Ff$, les actions des opérateurs $\Delta^x_y$, $T^{(-1)}$ et $S$ sur $V^{I(1)}$ sont respectivement données\footnote{les images des vecteurs de base sont données par les colonnes des matrices respectives bien que les actions soient à droite} par $H_\delta(x,y)$, $H_\omega(\lambda)$ et $H_s(\lambda,\tau)$.
\end{lem}
\textsf{Preuve :} \\
Comme on a la décomposition (\ref{suppF}), les fonctions images des $f^i_{a,b}$ par chacun de ces opérateurs sont uniquement déterminées par leur valeur en $1$ et en $\omega$. Les actions de $\Delta^x_y$ et $T^{(-1)}$ sont faciles à calculer puisque ces opérateurs sont de support une simple classe (voir (\ref{double=simple})). Par exemple, en se rappelant la formule d'action (\ref{I(1)agit}) et $\omega^2 = \varpi$, on a :
\[ \big( f_{a,b}^0 . T^{(-1)} \big)(1) = f_{a,b}^0(\omega) = 0 ,\] 
\[ \big( f_{a,b}^0 . T^{(-1)} \big)(\omega) = f_{a,b}^0(\varpi) = (\rho_1(\varpi) \otimes \rho_2(\varpi)) \, v_a \otimes w_b ,\]
\[ \big( f^1_{a,b} . T^{(-1)} \big)(1) = f^1_{a,b}(\omega) = v_a \otimes w_b ,\]
\[ \big( f^1_{a,b} . T^{(-1)} \big)(\omega) = f^1_{a,b}(\varpi) = 0 .\]
Pour $\Delta^x_y$, cela résulte de la discussion au début du paragraphe \ref{parI(1)} et de l'identité $\omega \delta^x_y \omega^{-1} = \delta^y_{\varpi x \varpi^{-1}}$. Ainsi on a 
\[ \big( f^0_{a,b} . \Delta^x_y \big)(1) = f^0_{a,b}(\delta^x_y) = (\rho_1(x) \otimes \rho_2(y)) \, v_a \otimes w_b = \xi_1^{\{ \overline{x}^{q_D^a}\}} \xi_2^{\{ \overline{y}^{q_D^b}\}} \, v_a \otimes w_b ,\]
\[ \big( f^0_{a,b} . \Delta^x_y \big)(\omega) = f^0_{a,b} (\delta^y_{\varpi x \varpi^{-1}} \omega) =0 ,\]
\[ \big( f^1_{a,b} . \Delta^x_y \big)(1) = f^1_{a,b}(\delta^x_y) = 0 ,\] 
\[ \big( f^1_{a,b} . \Delta^x_y \big)(\omega) = f^1_{a,b}(\delta^y_{\varpi x \varpi^{-1}} \omega) = \xi_1^{\{ \overline{y}^{q_D^a}\}} \xi_2^{\{ \overline{x}^{q_D^{b+1}}\}} \, v_a \otimes w_b .\]
Pour $S$ par contre, le calcul est moins immédiat. On écrit la décomposition \[ I(1) s I(1) = \coprod_{x \in k_D} I(1) s n([-x]) ,\] ainsi que les relations $n([x]) s = \left( \begin{array}{cc} 1 & [x] \varpi^{-1} \\ 0 & \varpi^{-1} \end{array} \right) \omega$ et \[ \omega n([x]) s = \begin{cases} \left( \begin{array}{cc} -[x]^{-1} & \varpi^{-1} \\ 0 & \varpi [x] \varpi^{-1} \end{array} \right) \omega n([x]^{-1})& \textrm{si } x \neq 0, \\ \varphi_2 & \textrm{si } x=0 . \end{cases} \]
Parce que $|k_D|$ est divisible par $p$, on a : 
\[ \big( f^0_{a,b}.S \big)(1) = \sum_{x \in k_D} f^0_{a,b}(n([x]) s) =0 ,\] 
\[ \big( f^0_{a,b}.S \big)(\omega) = \sum_{x \in k_D} f^0_{a,b}(\omega n([x]) s) = (1 \otimes \rho_2(\varpi)) \, v_a \otimes w_b ,\] 
\[ \big( f_{a,b}^1 .S \big)(1) = \sum_{x \in k_D} f_{a,b}^1(n([x])s) = \sum_{x \in k_D} (1 \otimes \rho_2(\varpi)^{-1}) \, v_a \otimes w_b = 0 .\] 
Enfin, on calcule
\begin{equation} \label{f1S} \big( f_{a,b}^1 .S \big)(\omega) = \sum_{x \in k_D} f_{a,b}^1(\omega n([x])s) = \rho_1(-1) \sum_{x \in k_D^\times} (\rho_1([x])^{-1} \otimes \rho_2(\varpi [x] \varpi^{-1})) \, v_a \otimes w_b .\end{equation}
Dans cette dernière égalité, on a mis $\rho_1(-1)$ en facteur puisque, $-1$ étant central dans $D^\times$, $\rho_1(-1)$ agit par un scalaire. On commence par remarquer que l'action de $\Oo_D^\times$ à travers $\rho_1^{-1} \otimes \rho_2(\varpi \cdot \varpi^{-1})$ sur l'espace $| \rho_1 \otimes \rho_2 |$ est diagonale et qu'il n'existe pas d'entier $0 \leq k < d_2$ avec $\rho_1 |_{\Oo_D^\times} = \rho_2 (\varpi^k \cdot \varpi^{-k})|_{\Oo_D^\times}$ si $d_1$ et $d_2$ sont distincts (puisque les stabilisateurs dans $D^\times$ des caractères de $\Oo_D^\times$ les composant sont distincts). \\
Supposons dans un premier temps $\rho_1 |_{\Oo_D^\times} = \rho_2(\varpi^k \cdot \varpi^{-k})|_{\Oo_D^\times}$ pour un certain entier $0 \leq k < d_2$; c'est équivalent à la condition $\xi_2^{q_D^k} = \xi_1$ (en particulier $d_1 = d_2$). Si on a $b+1 = a+k$ modulo $d_1$, pour tout $x \in k_D^\times$, on a 
\begin{align*} (\rho_1([x])^{-1} \otimes \rho_2(\varpi [x] \varpi^{-1})) \, v_a \otimes w_b & = \xi_1^{-\{x^{q_D^a}\}} \xi_2^{\{ x^{q_D^{b+1}}\}} \, v_a \otimes w_b \\ & = \xi_1^{-\{x^{q_D^a}\}} \xi_2^{\{ x^{q_D^{a+k}}\}} \, v_a \otimes w_b = v_a \otimes w_b .\end{align*}
Alors (\ref{f1S}) est une somme de $(q^d-1)$ termes identiques et on a \[ \big( f^1_{a,b} . S \big)(\omega) = \tau \, v_a \otimes w_b .\]
Si on a $b+1 \neq a+k$ modulo $d_1$, pour tout générateur $y$ de $k_D^\times$, l'action de $\rho_1([y])^{-1} \otimes \rho_2(\varpi[y]\varpi^{-1})$ sur $v_a \otimes w_b$ est donnée par un scalaire différent de $1$. Mais alors (\ref{f1S}) vaut 
\begin{equation} \label{f1Snul} \sum_{j=1}^{q^d-1} \big( \rho_1([y])^{-1} \otimes \rho_2(\varpi[y]\varpi^{-1}) \big)^j \, v_a \otimes w_b = 0 .\end{equation}
Supposons maintenant $\rho_1 |_{\Oo_D^\times} \neq \rho_2(\varpi^k \cdot \varpi^{-k})|_{\Oo_D^\times})$ pour tout entier $0 \leq k < d_2$. Alors, pour tout générateur $y$ de $k_D^\times$, $\rho_1([y])^{-1} \otimes \rho_2(\varpi[y]\varpi^{-1})$ agit sur $v_a \otimes w_b$ par un scalaire différent de $1$. Comme en (\ref{f1Snul}), pour tous $a$ et $b$, on a alors $(f^1_{a,b}.S)(\omega)=0$. \hfill$\Box$

\section{Représentations de Steinberg}

Soit $\rho$ une représentation irréductible de dimension $1$ (autrement dit un caractère) de $D^\times$. Par le lemme \ref{Ddim1}, il existe un caractère lisse $\rho_0$ de $F^\times$ tel que $\rho$ se factorise en $\rho = \rho_0 \circ \Nrm$. Grâce à (\ref{detNrm}), le morphisme de groupes $\det_G$ nous permet de définir une injection 
\begin{equation} \label{csteinclus} \rho_0 \circ \det\nolimits_G \hookrightarrow \Ind_B^G \big( \rho_0 \circ \det\nolimits_G \big) = \Ind_B^G \, \rho \otimes \rho .\end{equation}
Définissons la représentation de Steinberg $\St_B \rho_0$ de $G$ comme la représentation quotient de $\Ind_B^G \, \rho \otimes \rho$ par $\rho_0 \circ \det_G$. Par le lemme \ref{torscaract} et parce que la tensorisation par un caractère de $G$ est exacte, on a: 
\begin{equation} \label{torsSt} \St_B \rho_0 \simeq (\rho_0 \circ \det\nolimits_G) \otimes \St_B \id .\end{equation}
Cette identité aurait pu servir pour définir $\St_B \rho_0$ de manière équivalente, et elle permet de ramener l'étude au cas du caractère trivial.

\begin{lem} \label{uniqueconstantes}
La représentation $\Ind_B^G \, \id$ contient une unique représentation propre et non nulle, la représentation triviale $\id$.
\end{lem}
\textsf{Preuve :} \\
Notons $V := \Ind_B^G \, \id$. Le fait que $V$ contienne $\id$ vient de (\ref{csteinclus}); venons-en à l'unicité. Soit $W$ une sous-représentation non nulle et propre de $V$. Parce que $I(1)$ est un pro-$p$-groupe, $W^{I(1)}$ n'est pas réduit à $0$. Et comme $V$ est engendrée par ses $I(1)$-invariants (lemme \ref{critirred}), qui forment un espace de dimension $2$, $W^{I(1)}$ est en fait une droite. Aussi, le lemme \ref{I(1)actions} nous donne explicitement l'action de l'algèbre de Hecke du pro-$p$-Iwahori sur $V^{I(1)}$. On cherche ainsi une droite de $V^{I(1)}$ stable sous l'action (à gauche) des trois matrices: 
\[ \left( \begin{array}{cc} 1 & 0 \\ 0 & 1 \end{array} \right), \quad \left( \begin{array}{cc} 0 & 1 \\ 1 & 0 \end{array} \right), \quad \left( \begin{array}{cc} 0 & 0 \\ 1 & -1 \end{array} \right) .\]
La seule telle droite est $\Fpbar(f^0_{0,0} + f^1_{0,0}) = W^{I(1)}$. Mais cette dernière engendre précisément $W = \id$, et cela fournit l'unicité voulue. \hfill$\Box$\\

Pour prouver l'irréductibilité de $\St_B \id$, on va exhiber un sous-espace des fonctions de $U$ dans $\Fpbar$ qui lui est isomorphe en tant que $B$-représentation. On dit qu'une fonction $f$ sur $U$ à valeurs dans un $\Fpbar$-espace vectoriel $L$ est lisse s'il existe un sous-groupe ouvert $U'$ de $U$ tel que $f(u \cdot)$ est constante sur $U'$ pour tout $u \in U$. On remarque que $t := \varphi_1^d$ contracte strictement $U$ et on note $U_n = t^n U_0 t^{-n}$ pour tout $n \in \Z$; de cette manière, les $U_n$ forment une base de voisinages ouverts de $U$. En particulier, une fonction $f$ sur $U$ à valeurs dans $L$ est lisse si et seulement si il existe $n \in \Z$ tel que $f(u \cdot)$ est constante sur $U_n$ pour tout $u \in U$. Soit $C_c^\infty(U)$ le $\Fpbar$-espace vectoriel des fonctions lisses à support compact sur $U$, à valeurs dans $\Fpbar$. On note \[ \iota_U : U \xrightarrow{\sim} B \backslash B \omega U \] l'isomorphisme topologique induit par $u \mapsto \omega u$. On fait de $C_c^\infty(U)$ une $B$-représentation en faisant agir $b \in B$ par:
\begin{equation} \label{Ccagit} f \mapsto \big( u \mapsto f(\iota_U^{-1}( \iota_U(u) b) ) \big) .\end{equation}
En particulier, on remarque que $U$ agit sur $C_c^\infty(U)$ par translation à droite et $a \in A$ agit par conjugaison: $f \mapsto f(a^{-1} \cdot a)$.

\begin{prop} \label{Cc(U)}
$\phantom{,}$
\begin{itemize}
\item[(i)] La $B$-représentation $C_c^\infty(U)$ est isomorphe à $\St_B \id$ et est irréductible en tant que représentation de $B$. La $G$-repr\' esentation $\St_B \id$ est alors irr\' eductible.
\item[(ii)] La suite exacte de $G$-représentations 
\begin{equation} \label{defStB} 0 \to \id \to \Ind_B^G \, \id \to \St_B \id \to 0 \end{equation}
est non scindée.
\end{itemize}
\end{prop}
\textsf{Remarque :}
La preuve du \textit{(i)} prouve même que la restriction de $C_c^\infty(U)$ à $t^{\Z} U$ est irréductible. \\
\textsf{Preuve :} \\
On va prouver l'irréductibilité dans \textit{(i)} à la manière du Théorème 5 de \cite{Vig08}. Soient $W$ un sous-espace non nul de $C_c^\infty(U)$ stable par $B$ et $f$ un élément non nul de $W$. Parce que $f$ est à support compact et que la famille $(U_m)_{m \in \Z}$ est décroissante et recouvre $U$, il existe $m \in \Z$ tel que le support de $f$ soit inclus dans $U_m$. Mais alors, le sous-espace $W_m$ de $W$ constitué des fonctions à support contenu dans $U_m$ n'est pas réduit à $0$. Parce que $U_m$ est un pro-$p$-groupe, $W_m^{U_m}$ contient aussi un vecteur non nul: à multiplication par un scalaire près, c'est la fonction caractéristique $\id_{U_m}$ de $U_m$. A fortiori $W$ contient $\id_{U_m}$, et donc les $ut^n. \id_{U_m} = \id_{U_{m-n} u}$ pour tous $n \in \Z$ et $u \in U$. Enfin, parce que toute fonction $h$ de $C_c^\infty(U)$ est lisse à support compact, il existe $k \in \Z$ tel que $h$ soit combinaison linéaire finie des $\id_{U_k u}$ pour $u \in U$. On a donc $W = C_c^\infty(U)$ et $C_c^\infty(U)$ est irréductible en tant que $B$-représentation. \\
La décomposition de Bruhat
\begin{equation} \label{wBruhat} G= B \coprod B \omega U ,\end{equation}
le fait que la suite exacte (\ref{defStB}) de représentations de $G$ (donc de $B$) définisse $\St_B \id$ et l'action (\ref{Ccagit}) montrent que la flèche
\begin{equation} \label{modelSt} \begin{array}{ccc} \St_B \id & \to & C_c^\infty(U) \\ f & \mapsto & (u \mapsto f(\omega u )) \end{array} \end{equation}
est un isomorphisme de $B$-représentations. Ainsi $\St_B \id$ est $B$-irréductible, et donc $G$-irréductible. La suite exacte (\ref{defStB}) de $G$-représentations est non scindée à cause du lemme \ref{uniqueconstantes}. \hfill$\Box$\\

Le même énoncé est valable après torsion par un caractère de $G$.

\begin{cor} \label{St2}
Soient $\rho_0$ un caractère de $F^\times$ et $\rho = \rho_0 \circ \Nrm$. La suite exacte de $G$-représentations \[ 0 \to \rho_0 \circ \det\nolimits_G \to \Ind_B^G \, \rho \otimes \rho \to \St_B \rho_0 \to 0 \] est non scindée et fait de $\Ind_B^G \, \rho \otimes \rho$ une extension non triviale entre deux représentations irréductibles.
\end{cor}

Bien que l'on sache que $\Ind_B^G \, \rho \otimes \rho$ est admissible (voir proposition \ref{Indadmiss}), l'admissibilité du quotient $\St_B \rho_0$ n'est pas automatique. Lorsque $F$ est un corps $p$-adique, on peut invoquer \cite{Vig07}. Dans le cas général des représentations de Steinberg, on peut se reporter au chapitre \ref{Ly11}. Ici, on peut donner une preuve élémentaire.

\begin{prop} \label{StI(1)}
Soit $\rho_0$ un caractère de $F^\times$.
\begin{itemize}
\item[(i)] L'espace des $I(1)$-invariants $\big( \St_B \rho_0 \big)^{I(1)}$ est de dimension $1$.
\item[(ii)] La $G$-représentation $\St_B \rho_0$ est admissible.
\end{itemize}
\end{prop}
\textsf{Preuve :} \\
Parce que $I(1)$ est un pro-$p$-sous-groupe ouvert du groupe localement profini $G$, \textit{(i)} implique directement l'admissibilité de $\St_B \rho_0$ par \cite{Pas04}, Theorem 6.3.2. Comme on a l'isomorphisme (\ref{torsSt}) de $G$-représentations et que $I(1)$ agit trivialement sur $\rho_0 \circ \det_G$, les espaces $\big( \St_B \rho_0 \big)^{I(1)}$ et $\big( \St_B \id \big)^{I(1)}$ sont de même dimension. On va donc supposer $\rho_0 = \id$ par la suite. \\
En appliquant le foncteur exact à gauche des $I(1)$-invariants à la suite exacte courte
\[ 0 \to \id \to \Ind_B^G \, \id \to \St_B \id \to 0 ,\]
on obtient la suite exacte
\begin{equation} \label{seqI(1)inv} 0 \to \id \to \big( \Ind_B^G \, \id \big)^{I(1)} \xrightarrow{\pr} \big( \St_B \id \big)^{I(1)} .\end{equation}
L'espace $\big( \Ind_B^G \, \id \big)^{I(1)}$ est de dimension $2$ par la discussion avant la proposition \ref{Indadmiss} et on veut montrer que (\ref{seqI(1)inv}) se complète en une suite exacte courte, ce qui donnera un isomorphisme de $\Fpbar$-espaces vectoriels
\[ \id \oplus \big( \St_B \id \big)^{I(1)} \xrightarrow{\sim} \big( \Ind_B^G \, \id \big)^{I(1)} \simeq \Fpbar^2 \] 
et le \textit{(i)} comme voulu. \\
Explicitement, prenons un élément $f \in \big( \St_B \id \big)^{I(1)}$ et montrons que $f$ se relève en un élément de $\big( \Ind_B^G \, \id \big)^{I(1)}$ pour assurer la surjectivité de $\pr$. On prend un relèvement de $f$ en une fonction $\widetilde{f} \in \Ind_B^G \, \id$, et on va montrer que $\widetilde{f}$ est invariante par $I(1)$. Parce que $f$ est invariante par $I(1)$, pour tout $i \in I(1)$, il existe une constante $\lambda_i \in \Fpbar$ telle que l'on a $\widetilde{f}( \cdot \, i) = \widetilde{f} + \lambda_i$. Prenons $a,b,c,d \in \Oo_D$; le petit calcul
\[ \left( \begin{array}{cc} 1 & 0 \\ \varpi c (1 + \varpi a)^{-1} & 1 \end{array} \right) \left( \begin{array}{cc} 1 + \varpi a & b \\ 0 & 1 + \varpi \big(d -c (1 + \varpi a)^{-1} b\big)\end{array} \right) = \left( \begin{array}{cc} 1 + \varpi a & b \\ \varpi c & 1 + \varpi d \end{array} \right) \]
montre que l'on peut écrire tout $i \in I(1)$ sous la forme $i = i^- i^+$ avec $i^- \in U^- \cap I(1)$ et $i^+ \in B \cap I(1)$. On écrit ensuite
\[ \lambda_i = \widetilde{f}(\omega i ) - \widetilde{f}( \omega) = \widetilde{f} ( \omega i^- \omega^{-1} \omega i^+) - \widetilde{f} (\omega) .\]
Parce que $\omega i^- \omega^{-1}$ est dans $B$ et que $\widetilde{f}$ est $B$-invariante à gauche, on obtient
\begin{equation} \label{lambdai} \lambda_i = \widetilde{f}(\omega i^+) - \widetilde{f} (\omega) = \lambda_{i^+} .\end{equation}
Maintenant on sait, par définition de $\lambda_{i^+}$ et par $B$-invariance à gauche:
\begin{equation} \label{lambdai+} \lambda_{i^+} = \widetilde{f} (i^+) - \widetilde{f}(1) = 0 .\end{equation}
La coordination de (\ref{lambdai}) et (\ref{lambdai+}) implique que $\lambda_i$ est nulle pour tout $i \in I(1)$. Il s'ensuit $\widetilde{f} \in \big( \Ind_B^G \, \id \big)^{I(1)}$ et la preuve est terminée. \hfill$\Box$

\section{Séries principales et restriction à $B$} 

On commence par s'intéresser aux représentations irréductibles de dimension finie de $G$.

\begin{lem} \label{factordet}
Toute représentation lisse de dimension finie de $G$ se factorise par $\det_G$.
\end{lem}
\textsf{Preuve :} \\
Soit $\rho$ une telle représentation, et soient $(e_1, \dots, e_n)$ une base de l'espace $|\rho|$ sous-jacent à $\rho$. Parce que $\rho$ est lisse et que $G$ possède une base de voisinages de l'unité composée de sous-groupes compacts ouverts (les $1 + \varpi^k M_2(\Oo_D)$ pour $k \geq 1$), pour tout $1 \leq i \leq n$, il existe un sous-groupe compact ouvert $H_i$ de $G$ stabilisant $e_i$. Parce que $(e_1, \dots, e_n)$ est une base de $|\rho|$, en notant $H = H_1 \cap \dots \cap H_n$, la restriction de $\rho$ à $H$ est triviale. Mais alors, le noyau de $\rho$ contient $H \cap U$ et $H \cap U^-$; et parce que $\ker \, \rho$ est normal dans $G$, en conjuguant par $\varphi_1^{\Z}$, il contient $U$ et $U^-$. Or, par \cite{Pie82}, Lemma 16.5.a, et \cite{MatNak43}, Satz 1, le sous-groupe de $G$ engendré par $U$ et $U^-$ est précisément le noyau de $\det_G$. Autrement dit, on a la factorisation 
\begin{equation} \label{factorrho} \xymatrix{ G \ar[rd]_{\det_G} \ar[r] \ar@/^2pc/[rr]^\rho & G/\ker(\det_G) \ar[r] \ar[d]^{\simeq} & \GL(|\rho|) \\ & F^\times & }\end{equation}
et le lemme est prouvé. \hfill$\Box$\\

De ce fait, en reprenant (\ref{factorrho}), comme $\Fpbar$ est algébriquement clos et $F^\times$ commutatif, si $\rho$ est irréductible alors $\rho$ est un caractère de $G$. \\

Soient $\rho_1$ et $\rho_2$ deux représentations irréductibles de $D^\times$, de dimension respective $d_1$ et $d_2$. On note $|\rho_1 \otimes \rho_2|$ l'espace sous-jacent à la $B$-représentation $\rho_1 \otimes \rho_2$. Comme précédemment avec la représentation de Steinberg, on va étudier $\Ind_B^G \, \rho_1 \otimes \rho_2$ en observant sa restriction à $B$. \\
Soit $C^\infty_c(U, \rho_1 \otimes \rho_2)$ le $\Fpbar$-espace vectoriel des fonctions lisses à support compact sur $U$, à valeurs dans $|\rho_1 \otimes \rho_2|$. On fait de $C_c^\infty(U,\rho_1 \otimes \rho_2)$ une $B$-représentation en faisant agir $b \in B$ par 
\begin{equation} \label{Ccoefagit} f \mapsto \big( u \mapsto \rho_1 \otimes \rho_2 (\omega ub \iota_U^{-1}(\iota_U(u) b)^{-1} \omega^{-1}) \, f(\iota_U^{-1}(\iota_U(u) b)) \big) .\end{equation}
On remarque que cette formule d'action fait de $C_c^\infty(U,\rho_1 \otimes \rho_2)$ une sous-$B$-représentation de $\Ind_B^G \, \rho_1 \otimes \rho_2$ et que l'on a la compatibilité $C_c^\infty(U,\id) = C_c^\infty(U)$. \\

On appelle \textit{série principale} une représentation $\Ind_B^G \, \rho_1 \otimes \rho_2$ de $G$ qui est irréductible. Conjointement au corollaire \ref{St2}, les séries principales sont exactement les $\Ind_B^G \, \rho_1 \otimes \rho_2$ avec $\rho_1, \rho_2$ satisfaisant aux hypothèses du théorème \ref{GL2preuve1}.(ii) ci-dessous.

\begin{theo} \label{GL2preuve1}
Soient $\rho_1, \rho_2$ deux représentations irréductibles de $D^\times$.
\begin{itemize}
\item[(i)] La $B$-représentation $C_c^\infty(U,\rho_1 \otimes \rho_2)$ est irréductible.
\item[(ii)] Supposons $\rho_1 \nsim \rho_2$ ou $\rho_1 \otimes \rho_2$ de dimension strictement supérieure à $1$. Alors $\Ind_B^G \, \rho_1 \otimes \rho_2$ est une représentation irréductible de $G$.
\end{itemize}
\end{theo}
\textsf{Preuve :} \\
On rappelle que $t$ désigne $\varphi_1^d$, qui appartient au centre de $A$, et que l'on a alors \[ (\rho_1 \otimes \rho_2) (t^{-1} \cdot t) = \rho_1 \otimes \rho_2 .\]
En tant que $t^\Z U$-représentation, $\rho_1 \otimes \rho_2$ est isomorphe à la somme directe de $d_1 d_2$ représentations triviales. En appliquant le foncteur exact\footnote{voir la Proposition 2.4 de \cite{BusHen06} appliqué à $H = \{1\}$ et à $G=U$ localement profini} $C^\infty_c(U, \cdot)$, on obtient alors la décomposition en somme directe de $t^\Z U$-représentations irréductibles (par la preuve de la proposition \ref{Cc(U)}.(i)): \[ C^\infty_c(U, \rho_1 \otimes \rho_2) \simeq C^\infty_c(U)^{d_1 d_2} .\] En particulier, l'ensemble des sous-$t^\Z U$-représentations de  $C^\infty_c(U, \rho_1 \otimes \rho_2)$ est en bijection avec l'ensemble des sous-espaces vectoriels de $|\rho_1 \otimes \rho_2|$ via $C_c^\infty(U, \cdot)$. \\
Soit $\pi$ une sous-$B$-représentation non nulle de $C^\infty_c(U, \rho_1 \otimes \rho_2)$. Par le fait précédent, il existe un sous-espace vectoriel non nul $W$ de $|\rho_1 \otimes \rho_2|$ vérifiant l'isomorphisme de $t^\Z U$-représentations 
\begin{equation} \label{isompi} \pi \simeq C^\infty_c(U,W) .\end{equation} 
Soient $f$ un élément non nul de $\pi$ et $u \in U$ tel que $f$ ne s'annule pas en $u$. Supposons $W \neq |\rho_1 \otimes \rho_2|$ et soit $v \in |\rho_1 \otimes \rho_2| \smallsetminus W$. Parce que $\rho_1 \otimes \rho_2$ est $A$-irréductible, il existe des $\lambda_i \in \overline{\mathbb{F}}_p$ et des $a_i \in A \subseteq B$ tels que l'on ait 
\[ \sum\nolimits_i \lambda_i \, \rho_1 \otimes \rho_2 (a_i) \, f(u) = v .\]
Parce que l'on a l'identité \[ a_i \omega u = \omega \big( (\omega^{-1} a_i \omega) u (\omega^{-1} a_i^{-1} \omega) \big) (\omega^{-1} a_i \omega) \] et que l'action de $B$ sur $C_c^\infty(U,\rho_1 \otimes \rho_2)$ est donnée par (\ref{Ccoefagit}), on a:
\[ \sum\nolimits_i \lambda_i \, \big( (\omega^{-1} a_i \omega) \big( u^{-1} (\omega^{-1} a_i \omega) u (\omega^{-1} a_i^{-1} \omega) \big)  .f \big)( u) = v .\] 
Ceci donne $v \in W$ à cause de (\ref{isompi}), ce qui est une contradiction. On a donc $W = |\rho_1 \otimes \rho_2|$ et le \textit{(i)} est prouvé. \\
Montrons maintenant \textit{(ii)}. Soit $\sigma$ une sous-$G$-représentation non nulle de $\Ind_B^G \, \rho_1 \otimes \rho_2$. On a la suite exacte de $B$-représentations
\[ 0 \to C_c^\infty(U, \rho_1 \otimes \rho_2) \to \Ind_B^G \, \rho_1 \otimes \rho_2 \xrightarrow{\pr} \rho_1 \otimes \rho_2 \to 0 \]
provenant de (\ref{wBruhat}). Supposons que $\sigma$ est d'intersection nulle avec $C_c^\infty(U, \rho_1 \otimes \rho_2)$. Alors $\sigma$ s'injecte dans $\rho_1 \otimes \rho_2$, est donc de dimension finie et se factorise par $\det_G$ d'après le lemme \ref{factordet}. Ceci implique alors que $\rho_1 = \rho_2$ est un caractère de $D^\times$ d'après (\ref{detNrm}) et le lemme \ref{Ddim1}. Cela contredit les hypothèses de l'énoncé: c'est donc que $\sigma$ intersecte $C_c^\infty(U,\rho_1 \otimes \rho_2)$ non trivialement. Par le \textit{(i)}, on a donc $C_c^\infty(U,\rho_1 \otimes \rho_2) \subseteq \sigma$. De plus, comme $C_c^\infty(U,\rho_1 \otimes \rho_2)$ n'est pas stable par $G$, $\sigma$ est d'image non nulle par $\pr$. Et comme $\rho_1 \otimes \rho_2$ est $B$-irréductible, $\pr(\sigma)$ est égal à $\rho_1 \otimes \rho_2$. Au final, $\sigma$ est égal à $\Ind_B^G \, \rho_1 \otimes \rho_2$ et \textit{(ii)} en résulte. \hfill$\Box$

\section{Séries principales et $I(1)$-invariants} \label{I(1)paraboliques}

Dans ce paragraphe, en faisant une petite hypothèse sur $p$, on va donner une seconde preuve de l'irréductibilité dans le théorème \ref{mainthm}.(ii). Pour ce faire, on va établir la simplicité du module $\big( \Ind_B^G \, \rho_1 \otimes \rho_2 \big)^{I(1)}$ de $\Fpbar$-dimension $e = 2 d_1 d_2$ sur l'algèbre de Hecke du pro-$p$-Iwahori. On rappelle que $\{ \phantom{x} \}$ est un isomorphisme $k_D^\times \xrightarrow{\sim} \Z/(q^d-1) \Z$ que l'on a fixé et on garde les notations du paragraphe \ref{parI(1)} (notamment $\xi_1$ et $\xi_2$). \\
Commençons par des lemmes calculatoires un peu techniques. 

\begin{lem} \label{Fourierfini}
Soient $g = d_1 \wedge d_2$ et $j$ un entier de $[0,g-1]$. Soit \[ F_j = \Vect_{\Fpbar} \left\{ \begin{array}{c|l} (x,y) \mapsto \xi_1^{\{x^{q_D^{k_1}}\}} \xi_2^{\{y^{q_D^{k_2}}\}} , & 0 \leq k_1 < d_1, 0 \leq k_2 < d_2, \\ (x,y) \mapsto \xi_1^{\{y^{q_D^{k_1}}\}} \xi_2^{\{x^{q_D^{k_2+1}}\}} & k_1 - k_2 = j \mod g \end{array}\right\} \] le sous-espace du $\Fpbar$-espace vectoriel des fonctions $(k_D^\times)^2 \to \Fpbar$.
\begin{itemize}
\item[(a)] Supposons $\xi_1$ et $\xi_2$ non Frobenius-conjugués\footnote{c'est-à-dire non dans la même orbite sous l'application $x \mapsto x^{q_D}$; on remarque que par le lemme \ref{generatedfield}, en particulier si $d_1 \neq d_2$, alors $\xi_1$ et $\xi_2$ sont non Frobenius-conjugués}. Alors $F_j$ est de dimension $e/g$.
\item[(b)] Supposons $\xi_1$ et $\xi_2$ Frobenius-conjugués. Alors il existe au plus un entier $j_0 \in [0,g-1]$ tel que $F_{j_0}$ est de dimension strictement inférieure à $e/g$.
\end{itemize}
\end{lem}
\textsf{Preuve :} \\
Les fonctions $(x,y) \mapsto \xi_1^{\{x^{q_D^{k_1}}\}} \xi_2^{\{y^{q_D^{k_2}}\}}$ et $(x,y) \mapsto \xi_1^{\{y^{q_D^{k_1}}\}} \xi_2^{\{x^{q_D^{k_2+1}}\}}$ sont chacunes des caractères de groupes $(k_D^\times)^2 \to \Fpbar^\times$. Par le théorème d'indépendance linéaire des caractères d'Artin, il suffit de montrer que ces caractères sont deux à deux distincts. \\
On prouve tout d'abord \textit{(a)} et on suppose donc $\xi_1$ et $\xi_2$ non Frobenius-conjugués. Il y a trois cas à envisager. \\
\underline{Cas 1} : $\xi_1^{\{x^{q_D^{k_1}}\}} \xi_2^{\{y^{q_D^{k_2}}\}} = \xi_1^{\{x^{q_D^{k_3}}\}} \xi_2^{\{y^{q_D^{k_4}}\}}$ pour tous $x, y \in k_D^\times$. \\
En spécialisant en $x= \mu$ et $y=1$, on obtient $\xi_1^{q_D^{k_1}} = \xi_1^{q_D^{k_3}}$. Par le lemme \ref{generatedfield}, et parce que $a_D$ est premier à $d$, donc à $d_1$, $\xi_1$ n'est fixe par aucune puissance $k$ de $x \mapsto x^{q_D}$ pour $1 \leq k \leq d_1-1$. De ce fait, $k_1$ et $k_3$ sont congrus modulo $d_1$ et donc égaux. De même, en spécialisant en $x=1$ et $y = \mu$, on $\xi_2^{q_D^{k_2}} = \xi_2^{q_D^{k_4}}$, et $k_2$ et $k_4$ sont alors congrus modulo $d_2$. Le Cas 1 est terminé. \\
\underline{Cas 2} : $\xi_1^{\{x^{q_D^{k_1}}\}} \xi_2^{\{y^{q_D^{k_2}}\}} = \xi_1^{\{y^{q_D^{k_3}}\}} \xi_2^{\{x^{q_D^{k_4+1}}\}}$ pour tous $x,y \in k_D^\times$. \\
En spécialisant en $x=1$ et $y=\mu$, on obtient $\xi_2^{q_D^{k_2}} = \xi_1^{q_D^{k_3}}$, ce qui contredit l'hypothèse que $\xi_1$ et $\xi_2$ sont non Frobenius-conjugués. Le Cas 2 ne se produit donc jamais. \\
\underline{Cas 3} : $\xi_1^{\{y^{q_D^{k_1}}\}} \xi_2^{\{x^{q_D^{k_2+1}}\}} = \xi_1^{\{y^{q_D^{k_3}}\}} \xi_2^{\{x^{q_D^{k_4+1}}\}}$ pour tous $x,y \in k_D^\times$.\\
Le raisonnement est identique au Cas 1. \\
Passons maintenant à la preuve de \textit{(b)}. Les Cas 1 et 3 n'utilisent pas l'hypothèse de Frobenius-conjugaison et sont donc encore valables. Réexaminons le cas restant. \\
\underline{Cas 2'} : $\xi_1^{\{x^{q_D^{k_1}}\}} \xi_2^{\{y^{q_D^{k_2}}\}} = \xi_1^{\{y^{q_D^{k_3}}\}} \xi_2^{\{x^{q_D^{k_4+1}}\}}$ pour tous $x,y \in k_D^\times$. \\
En spécialisant en $x=1$ et $y = \mu$, puis en $x=\mu$ et $y=1$, on obtient les identités
\begin{equation} \label{relxi} \xi_2^{q_D^{k_2}} = \xi_1^{q_D^{k_3}}, \quad \xi_1^{q_D^{k_1}} = \xi_2^{q_D^{k_4+1}}.\end{equation}
On les combine et on obtient \[ \xi_1 = \xi_2^{q_D^{k_4-k_1+1}} = \xi_1^{q_D^{k_3-k_2+k_4-k_1+1}} .\]
Par le lemme \ref{generatedfield}, cela implique
\begin{equation} \label{congru1} k_3 - k_2 + k_4 - k_1 + 1 = 0 \mod g .\end{equation}
On rappelle aussi que si on travaille sur un espace $F_j$ fixé on a:
\begin{equation} \label{congru2} k_1 - k_2 = j \mod g , \quad k_3 - k_4 = j \mod g .\end{equation}
En combinant (\ref{congru1}) et (\ref{congru2}), on a 
\begin{equation} \label{congru3} 2(k_2-k_4)+1 = 0 \mod g .\end{equation}
On remarque que cela nécessite que $g$ soit impair; supposons-le, puisque dans le cas contraire \textit{(b)} est prouvé. Comme $k_4-k_2$ appartient à $[1-g,g-1]$ (on rappelle que $d_1=d_2=g$ puisque $\xi_1$ et $\xi_2$ sont Frobenius-conjugués), (\ref{congru3}) nous dit en fait: \[ 2(k_4-k_2)+1 \in \{ -g,0,g\} .\] Le cas $2(k_4-k_2)+1=0$ est écarté par parité. Reste que dans les deux cas restants, on a \[ k_4-k_2 = \frac{g-1}{2} \mod g .\] En se servant de cela, de (\ref{congru2}) et de (\ref{relxi}), on a 
\begin{equation} \label{relxi2} \xi_1^{q_D^j} = \xi_2^{q_D^{\frac{g+1}{2}}} .\end{equation}
Comme $\xi_1$ et $\xi_2$ sont Frobenius-conjugués et engendrent tous deux $\F_{q^g}$, il existe un unique $j_0 \in [0,g-1]$ vérifiant (\ref{relxi2}). Et pour $j \neq j_0$, le Cas 2' est encore impossible et $F_j$ est de dimension $e/g$ comme convenu. \hfill$\Box$\\

On remarque que la preuve implique en fait qu'un $j_0$ comme en \textit{(ii)} existe si et seulement si $g$ est impair. Le cas contraire, tous les $F_j$ sont de dimension $e/g$. Cependant, on n'aura pas besoin d'utiliser cette remarque. Par contre, on va utiliser le fait que, sans aucune hypothèse sur $\xi_1$ et $\xi_2$, la preuve des Cas 1 et 3 précédents nous donne immédiatement l'énoncé suivant.

\begin{lem} \label{5.1bis}
Les sous-espaces \[ \Vect_{\Fpbar} \bigg\{  (x,y) \mapsto \xi_1^{\{x^{q_D^{k_1}}\}} \xi_2^{\{y^{q_D^{k_2}}\}} \ \bigg| \ 0 \leq k_1 < d_1, 0 \leq k_2 < d_2 \bigg\} ,\] \[ \Vect_{\Fpbar} \bigg\{ (x,y) \mapsto \xi_1^{\{y^{q_D^{k_1}}\}} \xi_2^{\{x^{q_D^{k_2+1}}\}} \ \bigg| \ 0 \leq k_1 < d_1, 0 \leq k_2 < d_2 \bigg\} \] du $\Fpbar$-espace vectoriel des fonctions $(k_D^\times)^2 \to \Fpbar$ sont tous deux de dimension $e/2 = d_1 d_2$.
\end{lem}

On donne un résultat qui constitue un premier pas vers la simplicité de $V^{I(1)}$ en tant que $\Hh(G,I(1))$-module (sous de bonnes conditions). On remarque que l'énoncé n'utilise pas l'action de $H_s(\lambda,\tau)$ et on va pouvoir l'utiliser par la suite pour des sous-espaces de $V^{I(1)}$ ne disposant pas d'une telle action, et avec $H_\delta$ au lieu de $H_\Delta$.

\begin{lem} \label{d1^d2=1}
Supposons que $p$ ne divise pas $e = 2 d_1 d_2$. Soient $W$ un $\Fpbar$-espace vectoriel de dimension $e$, de base 
\[\big(e^0_{0,0}, e^0_{0,1}, \dots, e^0_{0,d_2-1}, e^0_{1,0}, \dots, e^0_{d_1-1,d_2-1}, e^1_{0,0}, \dots, e^1_{0,d_2-1}, \dots, e^1_{d_1-1,d_2-1} \big)\]
 et $X$ un ensemble. Soient $(h^i_{a,b})$ pour $i=0,1$, $0 \leq a < d_1$, $0 \leq b < d_2$ une famille de fonctions $X \to \Fpbar$ et, pour tout $x \in X$, $H_\Delta(x)$ la matrice d'endomorphisme de $W$ diagonale telle que l'on ait $H_\Delta(x) e^i_{a,b} = h^i_{a,b} e^i_{a,b}$ pour tous $i,a,b$. Pour tout entier $j \in [0,g-1]$, on note \[ W_j = \Vect_{\Fpbar} \, \big\{ e^i_{a,b} \ \big| \ i=0,1 ; \, a-b = j \mod g \big\} ,\] sous-espace vectoriel de $W$. On suppose que la sous-famille $(h^i_{a,b})$ pour $i=0,1$ et $a-b =j$ modulo $g$ est libre, pour tout $j \in [0,g-1]$. Alors, pour tout $j$, $W_j$ est stable et est irréductible\footnote{c'est-à-dire qu'il n'existe pas de sous-espace propre non nul stable par les endomorphismes en question} sous l'action de $H_\omega(\lambda)$ et des $H_\Delta(x)$ pour tout $x \in X$.
\end{lem}
\textsf{Remarque :} 
Lorsque l'on a $p \nmid e$ et $d_1 \wedge d_2 = 1$, cela montre déjà la simplicité de $V^{I(1)}$. \\
\textsf{Preuve :} \\
Parce que les $(h^i_{a,b})$ pour $i=0,1$ et $a-b=j$ modulo $g$ forment une famille libre de fonctions, les seuls sous-espaces de $W_j$ qui sont stables par les $H_\Delta(x)$ sont ceux engendr\' es par des parties de $\big\{ e^i_{a,b} \ \big| \ i=0,1 ; \, a-b = j \mod g \big\}$. Soit $W'$ un tel sous-espace non r\' eduit \` a $\{0\}$: il contient un vecteur $x = e^i_{a,b}$. \\
Regardons maintenant l'action de $H_\omega(\lambda)$. Soit $\lambda_0$ un élément vérifiant $\lambda_0^{2d_2} = \lambda$ dans le corps algébriquement clos $\Fpbar$. L'image de $e^0_{j,0}$ par les itérés successifs de $\lambda_0 H_\omega(\lambda)$ sont: 
\[e^0_{j,0} \mapsto \lambda_0^{1-2d_2} e^1_{j-1,d_2-1} \mapsto \lambda_0^{2-2d_2} e^0_{j-1,d_2-1} \mapsto \lambda_0^{3-2d_2} e^1_{j-2,d_2-2} \mapsto \dots \] \begin{equation} \label{f0itere} \dots \mapsto \lambda_0^{-1} e^1_{j-d_2,0} \mapsto e^0_{j-d_2,0} \mapsto \lambda_0^{1-2d_2} e^1_{j-d_2-1,d_2-1} \mapsto \dots ,\end{equation}
où les indices $a,b$ dans $f^i_{a,b}$ sont respectivement vus modulo $d_1$ et $d_2$. Parce que $(1,1)$ engendre un sous-groupe d'ordre $e/g$ dans $\Z/d_1\Z \times \Z/d_2\Z$, les images dans (\ref{f0itere}) forment un cycle de longueur $e/g$ et parcourent toutes les droites $\Fpbar \, e^i_{a,b}$ pour $a-b =j$ modulo $g$. De ce fait, le sous-espace $W'$ de $W_j$ contenant $x = e^i_{a,b}$ ne peut \^ etre stable sous $H_\omega(\lambda)$ que s'il est tout $W_j$. Comme tout sous-espace de $W_j$ non r\' eduit \` a $\{0\}$ et stable par $H_\omega(\lambda)$ et les $H_\Delta(x)$ est $W_j$, cela prouve son irr\' eductibilit\' e. \hfill$\Box$

\begin{lem} \label{d1d2=1}
Supposons $d_1 d_2 = 1$ et $\rho_1$ non isomorphe à $\rho_2$. Alors l'espace $\big( \Ind_B^G \, \rho_1 \otimes \rho_2 \big)^{I(1)}$ est simple en tant que module à droite sur $\Hh(G,I(1))$.
\end{lem}
\textsf{Preuve :} \\
Parce que l'on a $d_1 = d_2 = 1$, $\xi_1$ et $\xi_2$ sont des éléments de $\F_q$, et ainsi ils sont fixés par $x \mapsto x^{q_D}$. On se place dans la base $(f^0_{0,0},f^1_{0,0})$ dans laquelle le lemme \ref{I(1)actions} nous donne les matrices respectives des actions de $T^{(-1)}$, de $\Delta^x_y$ et de $S$ suivantes:
\[ H_\omega(\lambda) = \left( \begin{array}{cc} 0 & 1 \\ \lambda^{-1} & 0 \end{array} \right) , \quad H_\delta(x,y) = \left( \begin{array}{cc} \xi_1^{\{\overline{x}\}} \xi_2^{\{\overline{y}\}} & 0 \\ 0 & \xi_1^{\{\overline{y}\}} \xi_2^{\{\overline{x}\}} \end{array} \right) , \quad H_s(\lambda,\tau) = \left( \begin{array}{cc} 0 & 0 \\ \lambda^{-1} & \tau \end{array} \right) .\]
On discute deux cas distincts. \\
\underline{Cas 1} : $\xi_1 \neq \xi_2$, c'est-à-dire $\rho_1 |_{\Oo_D^\times} \neq \rho_2 |_{\Oo_D^\times}$. \\
En prenant $\overline{x} = \mu$ et $\overline{y} =1$, on voit que $H_\delta(x,y)$ n'est pas une homothétie: les sous-espaces non triviaux et propres que $H_\delta(x,y)$ stabilise sont alors les droites $\Fpbar f^0_{0,0}$ et $\Fpbar f^1_{0,0}$. Or ces droites ne sont pas stabilisées par $H_\omega(\lambda)$. On a donc montré la simplicité voulue. \\
\underline{Cas 2} : $\xi_1 = \xi_2$, c'est-à-dire $\rho_1 |_{\Oo_D^\times} = \rho_2 |_{\Oo_D^\times}$. \\
Parce que l'on a $\tau^2 = 1$, les droites stables par $H_s(\lambda,\tau)$ sont $\Fpbar f^1_{0,0}$ et $\Fpbar(\lambda f^0_{0,0} -\tau f^1_{0,0})$. Parce que $\rho_1$ et $\rho_2$ ne sont pas isomorphes, on a nécessairement $\lambda \neq 1$. De ce fait, ni l'une ni l'autre de ces droites n'est stable par $H_\omega(\lambda)$. Là encore on a la simplicité voulue. \hfill$\Box$ 

\begin{theo} \label{irredgeneric}
Supposons l'une des deux conditions suivantes satisfaite :
\begin{itemize}
\item[(a)] $\rho_1$ et $\rho_2$ ne sont pas isomorphes; 
\item[(b)] $\rho_1 = \rho_2$ est de dimension $d_1 = d_2 > 1$.
\end{itemize}
Supposons de plus $p \nmid 2 d_1 d_2$. Alors l'espace $\big( \Ind_B^G \, \rho_1 \otimes \rho_2 \big)^{I(1)}$ est simple en tant que module à droite sur $\Hh(G,I(1))$.
\end{theo}

Avant que de prouver ce résultat, mentionnons la conséquence qu'était notre objectif, devenue immédiate après consultation du lemme \ref{critirred}. Comme annoncé précédemment, on va avoir besoin de faire une hypothèse sur $p$. Lorsque $\rho_1$ et $\rho_2$ parcourent l'ensemble des représentations irréductibles de $D^\times$, $p \nmid 2 d_1 d_2$ devient la condition suivante.

\begin{hypo} \label{p2d}
Le premier $p$ ne divise pas $2d$.
\end{hypo}

Faisons tout de suite la remarque qu'elle provient de la méthode employée dans ce paragraphe et l'irréductibilité des séries principales est vraie sans cette hypothèse (voir théorème \ref{GL2preuve1}).

\begin{cor} \label{GL2preuve2}
Supposons l'hypothèse \ref{p2d}. Soient $\rho_1$ et $\rho_2$ deux représentations irréductibles de $D^\times$ satisfaisant (a) ou (b) du théorème \ref{irredgeneric}. Alors $\Ind_B^G \, \rho_1 \otimes \rho_2$ est une représentation irréductible de $G$.
\end{cor}

\subsection{Preuve du théorème \ref{irredgeneric}}
Lorsque $d_1$ et $d_2$ sont égaux à $1$, par hypothèse $\rho_1$ et $\rho_2$ ne sont pas isomorphes. Le résultat est alors l'objet du lemme \ref{d1d2=1}. \\
Lorsque $d_1 d_2$ est strictement supérieur à $1$ mais que $d_1$ et $d_2$ sont premiers entre eux, c'est l'objet du lemme \ref{d1^d2=1} appliqué à $W= \big( \Ind_B^G \, \rho_1 \otimes \rho_2 \big)^{I(1)}$ et \[ \{ H_\Delta(x) \ | \ x \in X \} = \{ H_\delta(x,y) \ | \ x,y \in \Oo_D^\times \} ,\] ce qui est possible grâce au lemme \ref{Fourierfini}.(i). \\
On peut donc supposer à présent que le pgcd de $d_1$ et $d_2$, que l'on notera $g$, est strictement supérieur à $1$; et c'est ce qu'on fera par la suite. Notons $V = \Ind_B^G \, \rho_1 \otimes \rho_2$. Pour tout entier $j$ de $[0,g-1]$, on note\footnote{on verra parfois aussi l'indice $j$ de $W_j$ comme la classe qu'il représente dans $\Z/g\Z$} 
\begin{equation} \label{defWj} W_j = \Vect_{\Fpbar} \{ f^i_{a,b} \ | \ i=0,1 ; \ a-b = j \mod g \} ;\end{equation}
cela découpe une somme directe de $\Fpbar$-espaces vectoriels $V^{I(1)} = \bigoplus_j W_j$. Soient $W$ un sous-$\Hh(G,I(1))$-module non nul de $V^{I(1)}$ et $w$ un élément non nul de $W$. On écrit
\begin{equation} \label{wdecomp} w = w_0 + w_1 + \dots + w_{g-1} \in \bigoplus\nolimits_j W_j ;\end{equation}
et on note \[ \Gamma_w = \big| \{ j \in [0,g-1] \ | \ w_j \neq 0 \} \big| .\] C'est un entier strictement positif car $w$ est non nul. \\
On va montrer par récurrence sur $k \geq 1$ que tout élément non nul $w$ de $W$ avec $\Gamma_w = k$ engendre $V^{I(1)}$ en tant que $\Hh(G,I(1))$-module. On veut distinguer deux cas différents pour la preuve de cette récurrence. \\
On remarque que $W_j$ est stable par $T^{(-1)}$ et par les $\Delta^x_y$ pour tous $x,y \in \Oo_D^\times$. \\

\underline{Cas 1} : $\rho_1 |_{\Oo_D^\times} \nsim \rho_2 |_{\Oo_D^\times}$, c'est-à-dire $\xi_1$ et $\xi_2$ non Frobenius-conjugués. \\
Commençons la récurrence avec l'étape $k=1$. Dans ce cas-là, il existe un unique $i_0 \in [0,g-1]$ avec $w_{i_0} \neq 0$. Parce que $\xi_1$ et $\xi_2$ sont non Frobenius-conjugués, le lemme \ref{Fourierfini}.(i) nous dit que les coefficients diagonaux de $H_\delta(x,y)|_{W_{i_0}}$ forment une famille libre. On peut alors appliquer le lemme \ref{d1^d2=1} pour affirmer que $W_{i_0}$ est irréductible et il est donc inclus dans $W$. En particulier $W$ contient $f^0_{i_0,0}$; et son image par $S$ est un élément non nul de $W_{i_0+1}$. On se sert ensuite à nouveau des lemmes \ref{Fourierfini}.(i) et \ref{d1^d2=1} pour affirmer que $W_{i_0+1}$ est inclus dans $W$; puis l'image de $f^0_{i_0+1,0}$ par $S$ est non nulle et appartient à $W_{i_0+2}$, etc. Par une récurrence immédiate, on voit que $W$ contient $\bigoplus_j W_j = V^{I(1)}$ et on a terminé l'étape d'initiation. \\
Montrons maintenant que le rang $k \geq 1$ implique le rang $k+1$ dans la récurrence. On remarque que les cas $k>g$ ou $k+1>g$ n'ont pas de sens vu que l'on a $\Gamma_w \leq g$ pour tout $w \in W$; $k$ est ainsi inférieur à $g-1$. Soit $w$ un élément de $W$ avec $\Gamma_w = k+1$. Notons
\[ W^{(i)} = \Vect_{\Fpbar} \{ f^i_{a,b} \ | \ 0 \leq a < d_1, \ 0 \leq b < d_2 \} \] 
pour $i=0,1$. Supposons dans un premier temps $w$ inclus dans $W^{(0)}$ ou dans $W^{(1)}$, disons $W^{(0)}$, l'autre cas de figure étant identique. Comme $\Delta^{[x]}_{[y]}$ stabilise $W^{(0)}$ et $\Delta^{[x]}_{[y]} \big|_{W^{(0)}}$ est diagonalisable de valeurs propres deux à deux distinctes par le lemme \ref{5.1bis}, le lemme des noyaux nous affirme la décomposition
\[ W^{(0)} = \bigoplus\nolimits_j \ker \Big( \prod_{a-b = j} \Delta^{[x]}_{[y]} \big|_{W^{(0)}} - \xi_1^{\{x^{q_D^a}\}} \xi_2^{\{y^{q_D^b}\}} \Id_{W^{(0)}} \Big) .\]
En particulier, si $i_1$ est tel que l'on ait $w_{i_1} \neq 0$ (voir (\ref{wdecomp})), alors les fonctions 
\[ (x,y) \mapsto w'(x,y) = w \cdot \Big( \prod_{a-b = j} \Delta^{[x]}_{[y]} \big|_{W^{(0)}} - \xi_1^{\{x^{q_D^a}\}} \xi_2^{\{y^{q_D^b}\}} \Id_{W^{(0)}} \Big) \]
sont non identiquement nulles (car on a $\Gamma_w = k+1 \geq 2$) et vérifient $\Gamma_{w'(x,y)} < \Gamma_w$ pour tout $x,y \in k_D^\times$. Il reste à prendre $(x_0,y_0) \in (k_D^\times)^2$ avec $w'(x_0,y_0) \neq 0$ et appliquer l'hypothèse de récurrence à $w'(x_0,y_0)$. \\
Supposons maintenant $w$ inclus ni dans $W^{(0)}$, ni dans $W^{(1)}$. Alors l'image $w'$ de $w$ par $S$ est non nulle et vérifie $\Gamma_{w'} < \Gamma_w$. On applique l'hypothèse de récurrence à $w'$ pour conclure. \\

\underline{Cas 2} : $\rho_1 |_{\Oo_D^\times} \simeq \rho_2 |_{\Oo_D^\times}$, c'est-à-dire $\xi_1$ et $\xi_2$ Frobenius-conjugués. \\
Reprenons l'étape d'initiation $k=1$ dans ce cadre-là. Par le lemme \ref{Fourierfini}.(ii), il existe au plus un $j_0 \in [0,g-1]$ tel que l'on ne puisse pas appliquer le lemme \ref{d1^d2=1} à $W_{j_0}$. Supposons que ce $j_0$ existe bien (le cas contraire, on peut procéder exactement comme au Cas 1) et voyons comment contourner la difficulté. Si $i_0$ est distinct de $j_0$, alors par le raisonnement du Cas 1, on a déjà \[ W_{i_0} \oplus W_{i_0+1} \oplus \cdots \oplus W_{j_0-1} \subseteq W .\]
Or on remarque que la restriction de $S$ à $W^{(0)} \cap W_{j_0-1}$ induit un isomorphisme de $\Fpbar$-espaces vectoriels entre $W^{(0)} \cap W_{j_0-1}$ et $W^{(1)} \cap W_{j_0}$. De cette manière, $W^{(1)} \cap W_{j_0}$ est inclus dans $W$, et il reste à utiliser l'action de $T^{(-1)}$ pour conclure $W_{j_0} \subseteq W$. On poursuit ensuite pour montrer $W_{j_0+1}$ comme au Cas 1, etc. \\
Maintenant, si $i_0$ est égal à $j_0$, on regarde d'abord le noyau de $S$: c'est le sous-espace vectoriel de $V^{I(1)}$ engendré par 
\begin{equation} \label{baseKerS} \{f^1_{a,b} \ | \ a + k \neq b+1 \} \cup \{ \tau f^1_{a,b} - f^0_{a,b+1} \ | \ a+k=b+1 \neq d_2 \} \cup \{ \tau f^1_{d_1-k,d_2-1} - \lambda f^0_{d_1-k,0} \} .\end{equation}
Si $w$ est dans le noyau de $S$, comme on a $\Gamma_w = 1$, $w$ est dans l'espace engendré par $\{f^1_{a,b} \ | \ a + k \neq b+1 \}$. Mais alors $T^{(-1)}$ envoie $w$ sur un élément $w'$ non nul de $W^{(0)} \cap W_{j_0}$, et $S$ envoie $w'$ sur un élément non nul $w''$ de $W^{(1)} \cap W_{j_0+1}$. Reste à appliquer le cas $i_0 \neq j_0$ à $w''$. Si $w$ n'est pas dans le noyau de $S$, son image $w'$ par $S$ est dans $W^{(1)} \cap (W_{j_0} \oplus W_{j_0+1})$. Comme à la fin du Cas 1, on forme alors 
\[ w''(x,y) = w' \cdot \Big( \prod_{a-b = j_0} \Delta^{[x]}_{[y]} \big|_{W^{(1)}} - \xi_1^{\{y^{q_D^a}\}} \xi_2^{\{x^{q_D^{b+1}}\}} \Id_{W^{(0)}} \Big) \]
pour obtenir un $w''(x_1,y_1)$ non nul dans $W^{(1)} \cap W_{j_0+1}$ et conclure comme précédemment. \\
Montrons maintenant comment le rang $k \geq 1$ implique le rang $k+1$ de la récurrence. Prenons $w \in W$ avec $\Gamma_w = k+1$. Si $wS$ est non nul, on le note $w' \in W^{(1)}$ et on a $\Gamma_{w'} \leq k+2$. On utilise alors une fois de plus le lemme \ref{5.1bis} et une de ces images $w''$ par un polynôme en $\Delta^{[x]}_{[y]} \big|_{W^{(1)}}$ est non nulle et vérifie $\Gamma_{w''} \leq k$: on applique l'hypothèse de récurrence à $w''$. \\
Si $wS$ est nul, observons la base (\ref{baseKerS}) du noyau de $S$: ou bien $w$ est engendré par $\{f^1_{a,b} \ | \ a + k \neq b+1 \}$ et on prend, comme précédemment, une de ses images par un polynôme convenable en $\Delta^{[x]}_{[y]} \big|_{W^{(1)}}$. Ou bien $w$ possède une composante sur un $\tau f^1_{a,b} - f^0_{a,b+1}$ ou $\tau f^1_{d_1-k,d_2-1} - \lambda f^0_{d_1-k,0}$. Remarquons que ces vecteurs ne sont pas vecteurs propres pour les $\Delta^{[x]}_{[y]}$: si cela avait été le cas, on aurait \[ \xi_1^{\{y^{q_D^a}\}} \xi_2^{\{x^{q_D^{b+1}}\}} = \xi_1^{\{x^{q_D^a}\}} \xi_2^{\{y^{q_D^b}\}} \] pour tous $x,y \in k_D^\times$. Or cela impliquerait $\xi_1^{q_D^a} = \xi_2^{q_D^b}$, $\xi_2^{q_D^{b+1}} = \xi_1^{q_D^a}$, et donc $\xi_1^{q_D} = \xi_1$. Par le lemme \ref{generatedfield} et l'hypothèse $g>1$, c'est absurde. De ce fait, il existe $x_2, y_2 \in k_D^\times$ tels que $w \Delta^{[x_2]}_{[y_2]} = w'$ ne soit pas dans le noyau de $S$. Son image $w''$ par $S$ est alors non nulle, dans $W^{(1)}$, et on peut appliquer le procédé précédent pour faire diminuer $\Gamma_{w''}$. La récurrence est terminée. \hfill$\Box$ 

\section{Séries principales et induction compacte}

On tâche dans ce paragraphe d'établir une nouvelle preuve du théorème \ref{mainthm}.(ii) en utilisant une méthode à la Barthel-Livné-Herzig. Parce que ceci est précisément l'objet d'un autre travail de l'auteur (chapitre \ref{Ly12}) pour $\GL(m,D)$, on sera plus concis sur la présentation. On va rappeler la méthode en renvoyant au chapitre \ref{Ly12} pour les preuves, et on va surtout se concentrer sur l'argument de \og changement de poids \fg, pour lequel, dans la plupart des cas (voir l'hypothèse \ref{HypDiamond} ci-dessous) on peut se passer du gros du travail du chapitre \ref{Ly12} et effectuer un argument ad hoc. \\

Soient $\rho_1$ et $\rho_2$ deux représentations irréductibles de $D^\times$. On cherche à démontrer que $\Ind_B^G \, \rho_1 \otimes \rho_2$ est irréductible sous les conditions du théorème \ref{mainthm}.(ii). \\
Avant de commencer l'argument à proprement parler, on s'attache d'abord à un lemme facile mais utile au cours de la démonstration. On remarque que ce lemme est un corollaire facile des lemme \ref{factordet} et théorème \ref{GL2preuve1}, mais on veut s'en passer.

\begin{lem} \label{pasdedroite}
Soient $\rho_1$ et $\rho_2$ deux représentations irréductibles de $D^\times$ telles que l'une des deux conditions suivantes est satisfaite:
\begin{itemize}
\item[(a)] $\rho_1$ et $\rho_2$ ne sont pas isomorphes;
\item[(b)] $\rho_1 = \rho_2$ est de dimension $d_1 = d_2 >1$.
\end{itemize}
La représentation $\Ind_B^G \, \rho_1 \otimes \rho_2$ ne possède pas de sous-$G$-représentation de dimension $1$.
\end{lem}
\textsf{Preuve :} \\
Supposons l'existence d'une droite de $\Ind_B^G \, \rho_1 \otimes \rho_2$ stable par $G$, disons engendrée par une fonction $f: G \to | \rho_1 \otimes \rho_2|$ non nulle. Soit $g \in G$ un élément du support de $f$. Pour tout $b \in B$, on a 
\begin{equation} \label{bfg} ((g^{-1} bg).f)(g) = f(bg) = \rho_1 \otimes \rho_2 (b) \, f(g) .\end{equation}
Et (\ref{bfg}) est sur la droite $\Fpbar f(g)$ de $|\rho_1 \otimes \rho_2|$ car $\Fpbar f$ est une droite stable par $G$ dans $\Ind_B^G \, \rho_1 \otimes \rho_2$. Comme c'est valable pour tout $b \in B$ et que $\rho_1 \otimes \rho_2$ est irréductible, $\rho_1 \otimes \rho_2$ est de dimension $1$. Soit $a \in A$. Parce que l'action de $G$ sur $f$ est de donnée par un certain caractère $G \to \Fpbar^\times$, on a l'égalité 
\begin{equation} \label{conjf} (g^{-1} s a s g).f = (g^{-1} a g).f .\end{equation}
En regardant la valeur en $g$ et en utilisant (\ref{bfg}), (\ref{conjf}) devient 
\[ \rho_1 \otimes \rho_2(sas) \, f(g) = \rho_1 \otimes \rho_2(a) \, f(g) .\]
Parce que $f(g)$ est non nul, cela implique $\rho_1 = \rho_2$, ce qui, accumulé avec $d_1 = d_2 = 1$, finit par contredire l'hypothèse de l'énoncé. C'est donc que $f$ n'existe pas et le lemme est prouvé. \hfill$\Box$ \\

Soit $\pi$ une sous-$G$-représentation non nulle irréductible de $\Ind_B^G \, \rho_1 \otimes \rho_2$; en particulier, $\pi$ est admissible par la proposition \ref{Indadmiss}. On suppose que $\rho_1$ et $\rho_2$ vérifient l'hypothèse suivante.

\begin{hypo} \label{HypDiamond}
Si $\rho_1|_{\Oo_D^\times}$ et $\rho_2|_{\Oo_D^\times}$ sont isomorphes et ne sont pas des caractères, alors pour tout caractère $\chi_1 \otimes \chi_2 \subseteq (\rho_1 \otimes \rho_2)|_{B \cap K}$ avec $\chi_1 \neq \chi_2$, $\chi_2 \chi_1^{-1}$ correspond à $x \mapsto x^{\vec{s}}$ (voir (\ref{caractkD})) avec $\vec{s} = (s_0, s_1, \dots, s_{fd-1})$ tels que l'un au moins des $s_i$ n'appartient pas à $\{0,1,p-2,p-1\}$.
\end{hypo}
\textsf{Remarque :}
Le raisonnement général est valable sans elle mais elle est utile dans notre preuve du théorème \ref{chgtpds}. Dans le chapitre \ref{Ly12}, on prouve le théorème \ref{chgtpds} sans cette hypothèse et on recouvre complètement alors une autre preuve du théorème \ref{mainthm}.(ii). Cependant l'argument gagne beaucoup en complexité. \\
\textsf{Remarque :}
On voit que cela exclut d'office les premiers $p=2,3$. \\

Parce que $K(1)$ est un pro-$p$-groupe ouvert, $\pi^{K(1)}$ est non nul et de dimension finie. Comme $K(1)$ est normal dans $K$, $\pi^{K(1)}$ est une représentation de dimension finie de $K$, qui admet donc une sous-représentation irréductible, disons $V$. \\ 
On a l'alternative suivante: ou bien $\pi$ est de dimension finie et par le lemme \ref{factordet} et le fait que $\Fpbar$ est algébriquement clos, $\pi$ est un caractère de $G$, ce qui est exclu par le lemme \ref{pasdedroite}. Ou bien, par le théorème \ref{chgtpds} ci-dessous, on peut supposer que $V$ n'est pas un caractère (c'est-à-dire qu'elle est régulière au sens du paragraphe \ref{repfini}); c'est donc ce que l'on fera par la suite. \\
Comme l'algèbre de Hecke $\Hh(G,K,V)$ est commutative (par le paragraphe 2.10 de \cite{HenVig11}) et que
\[ \Hom_G(\ind_K^G \, V, \pi) \simeq \Hom_K(V,\pi) \] 
est non nul de dimension finie (par réciprocité de Frobenius et admissibilité de $\pi$), il existe un vecteur propre $\Phi$ dans $\Hom_G(\ind_K^G \, V,\pi)$ pour l'action de $\Hh(G,K,V)$. On note $\chi$ le caractère propre correspondant: c'est un caractère de $\Fpbar$-algèbres $\Hh(G,K,V) \to \Fpbar$. Et $\Phi$ se factorise donc par le morphisme de $G$-représentations
\begin{equation} \label{surjpi} \ind_K^G \, V \otimes_\chi \Fpbar \to \pi \end{equation}
que l'on notera encore $\Phi$. Parce que $\pi$ est irréductible, $\Phi$ est surjective. 
Observons que $\Phi$ provient par réciprocité de Frobenius d'un 
\[f \in \Hom_K(V,\pi) \subseteq \Hom_K \big( V, \Ind_B^G \, \rho_1 \otimes \rho_2 \big) \simeq \Hom_{A \cap K} \big( V_{U \cap K}, \rho_1 \otimes \rho_2 \big) ,\]
et on notera $f_1$ son image dans $\Hom_{A \cap K} \big( V_{U \cap K} , \rho_1 \otimes \rho_2 \big)$. 
Aussi, à cause des Lemmes 4.1 et 4.3 du chapitre \ref{Ly12}, $\chi$ se factorise à travers la transformée de Satake
\[ {}' \Sss_G : \Hh(G,K,V) \hookrightarrow \Hh \big( A,A \cap K, V_{U \cap K} \big) ;\]
on notera encore $\chi$ le caractère de $\Fpbar$-algèbres résultant $\Hh \big( A, A \cap K, V_{U \cap K} \big)$. \\
Parce que $V$ est régulière, par \cite{HenVig11b}, Theorem 1.2 et Corollary 1.3, on a l'isomorphisme de $G$-représentations 
\begin{equation} \label{comparInd} \iota_V : \ind_K^G \, V \otimes_\chi \Fpbar \xrightarrow{\sim} \Ind_B^G \big( \ind_{A \cap K}^A \, V_{U \cap K} \otimes_\chi \Fpbar \big) \end{equation}
caractérisé par 
\[ \iota_V \big( [1,v] \otimes 1 \big)(1) = [1, p_U(v)] \otimes 1 \]
pour tout $v \in V$ d'image $p_U(v)$ dans $V_{U \cap K}$. De plus, $f_1$ induit l'isomorphisme de $A$-représentations 
\[ \begin{array}{ccc} \ind_{A \cap K}^A \, V_{U \cap K} \otimes_\chi \Fpbar & \xrightarrow{\sim} & \rho_1 \otimes \rho_2 \\ \, [1,v] \otimes 1 & \mapsto & f_1(v) \end{array} \] 
comme on peut le voir à partir de la proposition \ref{Dtimescompact} (ou bien on peut se reporter au lemme \ref{parampos} et à la proposition \ref{cosocles}.(i) ultérieurs); on le notera encore $f_1$. En appliquant le foncteur exact $\Ind_B^G$, on obtient l'isomorphisme
\[ \Ind \, f_1 : \Ind_B^G \big( \ind_{A \cap K}^A \, V_{U \cap K} \otimes_\chi \Fpbar \big) \xrightarrow{\sim} \Ind_B^G \, \rho_1 \otimes \rho_2 \]
de $G$-représentations. Comme 
\[ (\Ind \, f_1) \circ \iota_V : \ind_K^G \, V \otimes_\chi \Fpbar \xrightarrow{\sim} \Ind_B^G \, \rho_1 \otimes \rho_2 \]
est déterminée par les valeurs 
\[ \big( (\Ind \, f_1) \circ \iota_V \big) ([1,v] \otimes 1) (1) = f_1 \circ p_U(v) = f(v)(1) \]
pour tout $v \in V$, on a $(\Ind \, f_1) \circ \iota_V = \Phi$. Il en résulte $\pi = \Ind_B^G \, \rho_1 \otimes \rho_2$ et $\Ind_B^G \, \rho_1 \otimes \rho_2$ est irréductible comme voulue. \\

On vient alors de démontrer l'énoncé suivant.

\begin{cor} \label{GL2preuve3}
Soient $\rho_1$ et $\rho_2$ deux représentations irréductibles de $D^\times$ telles que l'une des deux conditions suivantes est satisfaite:
\begin{itemize}
\item[(a)] $\rho_1$ et $\rho_2$ ne sont pas isomorphes;
\item[(b)] $\rho_1 = \rho_2$ est de dimension $d_1 = d_2 >1$.
\end{itemize}
Supposons de plus l'hypothèse \ref{HypDiamond}. Alors $\Ind_B^G \, \rho_1 \otimes \rho_2$ est une représentation irréductible de $G$.
\end{cor}

Enonçons maintenant l'énoncé de changement de poids qui nous a été utile.

\begin{theo} \label{chgtpds}
Soit $\pi$ une représentation irréductible non nulle de dimension infinie de $G$ qui s'injecte dans $\Ind_B^G \, \rho_1 \otimes \rho_2$ avec $\rho_1$ et $\rho_2$ satisfaisant l'hypothèse \ref{HypDiamond}. Alors $\pi$ contient une $K$-représentation irréductible régulière (c'est-à-dire qui n'est pas un caractère).
\end{theo}
\textsf{Remarque :}
Par le lemme \ref{factordet}, parce que $\Fpbar$ est algébriquement clos, il est équivalent de demander que $\pi$ ne soit pas un caractère de $G$. \\
\textsf{Remarque :} 
Au cours de la preuve, on verra que si $V$ est une $K$-représentation irréductible non régulière de $\pi|_K$, alors on peut choisir $V'$ régulière dans $\pi|_K$ telle que $(V')^{U \cap K}$ et $\varphi_2. V^{U \cap K}$ sont des $(A \cap K)$-représentations isomorphes, et sa classe d'isomorphisme est en fait uniquement déterminée.

\subsection{Preuve du théorème \ref{chgtpds}}

Soit $V$ une sous-$K$-représentation irréductible de $\pi$. Si $V$ est régulière, il n'y a rien à faire. \\
Supposons donc $V$ non régulière, c'est-à-dire égale à un caractère $\chi \circ \overline{\det}$ et portée sur une droite que l'on notera $\Fpbar v$. Si $\chi$ ne se factorise pas par la norme galoisienne $N_{k_D/k_F}$, alors $\Fpbar \varphi_2 v$ est une représentation de $K \cap \varphi_2 K \varphi_2^{-1} = I$ égale à \[ \chi \otimes \chi(\varpi^{-1} \cdot \varpi) \neq \chi \otimes \chi .\]
On a par suite une flèche non nulle 
\[ \Ind_I^K \, \chi \otimes \chi(\varpi^{-1} \cdot \varpi) \to \pi \subseteq \Ind_B^G \, \rho_1 \otimes \rho_2 = \Ind_{B \cap K}^K \, \rho_1 \otimes \rho_2 \]
de $K$-représentations. Par réciprocité de Frobenius, on a une flèche non nulle de $(B \cap K)$-représentations
\begin{equation} \label{ssrho12} \Ind_{\overline{B}}^{\overline{K}} \, \chi \otimes \chi(\varpi^{-1} [\cdot] \varpi) \to \rho_1 \otimes \rho_2 .\end{equation}
Le morphisme de (\ref{ssrho12}) se factorise en fait par le morphisme non nul de $(A \cap K)$-repr\' esentations 
\begin{equation} \label{ssrho12bis} \chi \otimes \chi(\varpi^{-1} [\cdot] \varpi) \oplus \big( \Ind_{\overline{A}}^{\overline{B}} \, \chi(\varpi^{-1} [\cdot] \varpi) \otimes \chi \big)_U \to \rho_1 \otimes \rho_2 .\end{equation}
En notant $\chi^\omega = \chi(\varpi^{-1} [\cdot] \varpi) \otimes \chi$, on cherche \` a d\' eterminer $\big( \Ind_{\overline{A}}^{\overline{B}} \, \chi^\omega \big)_U$. Une base de $\big( \Ind_{\overline{A}}^{\overline{B}} \, \chi^\omega \big)_U$ est donn\' ee par les $[u,v]$ pour $u \in \overline{U}$ et $v$ un \' el\' ement non nul fix\' e de la droite sous-jacente \` a $\chi^\omega$. Les $u.f - f$ pour $u \in \overline{U}$ et $f \in \Ind_{\overline{A}}^{\overline{B}} \, \chi^\omega$ sont inclus dans la sous-$\overline{B}$-repr\' esentation
\[ \big( \Ind_{\overline{A}}^{\overline{B}} \, \chi^\omega \big)_{(U)} = \big\{ \sum\nolimits_{u \in \overline{U}} \lambda_u [u,v] \ \big| \ \lambda_u \in \Fpbar, \, \sum \lambda_u = 0 \big\} ;\]
et ils l'engendrent bien comme on le voit en prenant $f = [1,v]$. Mais alors l'espace de coinvariants \[ \big( \Ind_{\overline{A}}^{\overline{B}} \, \chi^\omega \big)_{U} = \big( \Ind_{\overline{A}}^{\overline{B}} \, \chi^\omega \big) \big/ \big( \Ind_{\overline{A}}^{\overline{B}} \, \chi^\omega \big)_{(U)} \] est de dimension $1$, et il s'agit de d\' eterminer l'action de $A \cap K$ dessus. Et $a \in A \cap K$ agit par $\chi^\omega(\overline{a})$ sur l'image $\overline{[1,v]}$ de $[1,v]$ dans $\big( \Ind_{\overline{A}}^{\overline{B}} \, \chi^\omega \big)_{U}$; d'o\` u $\big( \Ind_{\overline{A}}^{\overline{B}} \, \chi^\omega \big)_{U} = \chi^\omega$ en tant que $(A \cap K)$-repr\' esentations. Et (\ref{ssrho12bis}) devient  
\[ \chi \otimes \chi(\varpi^{-1} \cdot \varpi) \oplus \chi(\varpi^{-1} \cdot \varpi) \otimes \chi \to \rho_1 \otimes \rho_2 .\]
Ainsi $\chi \otimes \chi(\varpi^{-1} \cdot \varpi)$ ou $\chi(\varpi^{-1} \cdot \varpi) \otimes \chi$ est un sous-caractère de $(\rho_1 \otimes \rho_2)|_{B \cap K}$, disons $\chi \otimes \chi(\varpi^{-1} \cdot \varpi)$. Par l'hypothèse \ref{HypDiamond}, il s'écrit $x \mapsto x^{\vec{s}}$ avec $\vec{s} = (s_0, s_1, \dots, s_{fd-1})$, l'un des $s_i$ n'appartenant pas à $\{0,1,p-2,p-1\}$. La Proposition 1.1 de \cite{Dia07} nous assure que dans ce cas-là, aucun des facteurs de Jordan-Hölder de $\Ind_I^K \, \chi \otimes \chi(\varpi^{-1} \cdot \varpi)$ n'est un caractère. La droite $\Fpbar \varphi_2 v$ engendre donc une $K$-représentation dans $\pi$, qui contient une sous-$K$-représentation irréductible régulière: la preuve est terminée. \\
Si $\chi$ se factorise par $N_{k_D/k_F}$, par la fin de la preuve de la Proposition 16.1 du chapitre \ref{Ly12}, il existe un caractère $\chi_0$ de $F^\times$ tel que $\chi \circ \overline{\det}$ soit la restriction à $K$ de $\chi_0 \circ \det_G$. En utilisant le lemme \ref{torscaract} pour tordre $\pi$ par $\chi_0^{-1} \circ \det_G$, on peut supposer que $\chi$ est le caractère trivial de $k_D^\times$. C'est donc ce que l'on fera par la suite. \\
La droite $\Fpbar v$ porte ainsi la représentation triviale de $K$. On regarde la droite $\Fpbar \varphi_2 v$: elle est stable par $K \cap \varphi_2 K \varphi_2^{-1} = I$, et on note $V_{1,0}$ la $K$-représentation qu'elle engendre. Parce que l'on a l'isomorphisme d'ensembles finis 
\[ I \backslash K \xrightarrow{\sim} \overline{B} \backslash \overline{K} = \big( (B \cap K) / (B \cap K(1)) \big) \big\backslash \overline{K} ,\]
et que $I$ est un sous-groupe compact ouvert de $K$, la réciprocité de Frobenius nous donne un morphisme non nul de représentations de $K$
\begin{equation} \label{generateV10} \Ind_I^K \, \id = \Ind^{\overline{K}}_{\overline{B}} \, \id \to V_{1,0} .\end{equation}
Par \cite{Pas04}, Lemma 3.1.7, on peut identifier $\Ind^{\overline{K}}_{\overline{B}} \, \id$: c'est la somme directe de la représentation triviale $\id$ de $\overline{K}$ et de la représentation de Steinberg finie $\overline{\St}$. Par \cite{CabEng04}, Definition 6.13, Theorem 6.10 et Theorem 6.12, $\overline{\St}$ est irréductible et de dimension $q^d$: c'est en fait la représentation que l'on a notée $\Sym^{\vec{p-1}} \, \Fpbar^2$, et qui est régulière. Le morphisme non nul (\ref{generateV10}) nous offre alors l'alternative suivante: ou bien $V_{1,0}$ contient $\overline{\St}$ et la preuve de la proposition est terminée, ou bien $V_{1,0}$ est encore la représentation triviale de $K$. \\
Dans ce second cas de figure, on va s'intéresser à $V_{2,0}$, la $K$-représentation engendrée par la $I$-représentation $\Fpbar \varphi_2^2 v$. De nouveau, par l'argument précédent, ou bien $V_{2,0}$ contient $\overline{\St}$, ou bien $V_{2,0}$ est $\id$. Dans le cas défavorable $V_{2,0} = \id$, on s'intéresse à $V_{3,0}$, etc. \\
S'il existe $i \geq 0$ (avec la convention $V_{0,0} = V$) tel que $V_{i,0}$ contienne $\overline{\St}$, la preuve est terminée. Supposons donc le cas contraire, c'est-à-dire que la $K$-représentation $V_{i,0}$ engendrée par $\Fpbar \varphi_2^i v$ est triviale pour tout $i \geq 0$. Mais alors, pour tout $j \geq 0$, la droite $\Fpbar \varpi^j \varphi_2^i v$ est stable par $K$ et est triviale: on la note $V_{i,j}$. On forme la sous-représentation de $\pi|_K$ 
\[ V_\infty = \sum_{i=0}^{2d-1} \sum_{j=0}^{d-1} V_{i,j} = \sum_{i=0}^{2d-1} \sum_{j=0}^{d-1} \Fpbar \varphi_2^i \varpi^j v ;\]
et on va montrer que $V_\infty$ est une représentation de $G$. Grâce à la décomposition de Cartan
\[ G = \coprod_{i \geq 0, \, j \in \Z} K \varphi_2^i \varpi^j K ,\]
il s'agit de voir que $V_\infty$ est stable par $\varpi^{\pm d}$ et $\varphi_2^{2d}$. Comme $\varpi^{d \Z}$ est un sous-groupe du centre de $G$, son action sur la représentation admissible irréductible $\pi$ est donnée par un caractère (par lemme de Schur): en particulier, elle stabilise $V_\infty$. Et parce que l'action de $s \in K$ est triviale sur chaque $V_{i,j}$ on a, pour tous $i, j \geq 0$: 
\[ \varphi_2^{2d} \varphi_2^i \varpi^j v = (\varphi_2 s \varphi_2 s \cdots \varphi_2 s) \varphi_2^i \varpi^j v = \varpi^d \varphi_2^i \varpi^j v .\]
Et comme l'action de $\varpi^d$ est scalaire, $\varphi_2^{2d}$ stabilise bien $V_\infty$ et $V_\infty$ est une $G$-représentation comme affirmé. Parce que $\pi$ est irréductible, $V_\infty$ est $\pi$, mais cela contredit le fait que $\pi$ est de dimension infinie. C'est donc que notre supposition était erronée: l'un des $V_{i,0}$ contient nécessairement $\overline{\St}$ et la preuve est terminée.

\section{Classification des représentations irréductibles admissibles}

Soit $\pi$ une représentation irréductible admissible de $G$. Par la discussion après le lemme \ref{pasdedroite}, il existe une $K$-représentation irréductible $V$ et un caractère $\chi : \Hh(G,K,V) \to \Fpbar$ de $\Fpbar$-algèbres tels que l'on ait une surjection 
\begin{equation} \label{generepi} \ind_K^G \, V \otimes_\chi \Fpbar \twoheadrightarrow \pi \end{equation}
de $G$-représentations. Si pour tous $V$ et $\chi$ satisfaisant (\ref{generepi}), $\chi$ ne se factorise pas par la transformée de Satake 
\[ {}' \Sss_G : \Hh(G,K,V) \to \Hh \big( A, A \cap K, V_{U \cap K} \big) ,\]
alors $\pi$ est dite \textit{supersingulière}. \\
On va tâcher de montrer que toute représentation $\pi$ qui n'est pas supersingulière a déjà été étudiée auparavant, c'est-à-dire qu'elle est un caractère, une représentation de Steinberg, ou bien une série principale. Dans ce paragraphe encore, on insiste sur le fait que l'on se passe complètement du travail de changement de poids du chapitre \ref{Ly12}. Et seule la proposition \ref{cosocles2}.(ii) fait appel au c\oe ur du travail du chapitre \ref{Ly12} avec le foncteur des parties ordinaires; tout le reste est \og élémentaire \fg. \\

Commençons pour cela par identifier les paramètres $(V,\chi)$ des séries principales comme en (\ref{comparInd}). Pour $\eta$ un caractère de $A \cap K$, on définit le sous-groupe suivant de $A_\Lambda$:
\[ \widetilde{Z}_A(\eta) = \{ a \in A_\Lambda \ | \ \forall x \in A \cap K, \, \eta(a^{-1}xa) = \eta(x) \} .\]
C'est un invariant important du caractère $\eta$, et cela transparaît sur $V$ aussi: pour toute représentation irréductible de $K$, on l'isomorphisme de $\Fpbar$-algèbres (voir le paragraphe 7.1 de \cite{HenVig11}): 
\begin{equation} \label{HeckeA} \begin{array}{rcl} \Hh \big( A,A \cap K, V_{U \cap K} \big) & \mathop{\leftarrow}\limits^{\sim} & \Fpbar \big[ \widetilde{Z}_A \big( V_{U \cap K} \big) \big] \\ \tau^{a_1}_{a_2} & \mapsfrom & \left( \begin{array}{cc} \varpi^{a_1} & 0 \\ 0 & \varpi^{a_2} \end{array} \right) \in \widetilde{Z}_A \big( V_{U \cap K} \big) \end{array} ,\end{equation}
où $\tau^{a_1}_{a_2}$ désigne l'opérateur de support $\left( \begin{array}{cc} \varpi^{a_1} & 0 \\ 0 & \varpi^{a_2} \end{array} \right) (A \cap K)$ et valant $\Id$ en $\left( \begin{array}{cc} \varpi^{a_1} & 0 \\ 0 & \varpi^{a_2} \end{array} \right)$.

\begin{lem} \label{parampos}
Soient $\rho_1$ et $\rho_2$ deux représentations irréductibles de $D^\times$. Soient $V$ une $K$-représentation irréductible et $\chi : \Hh(G,K,V) \to \Fpbar$ un caractère de $\Fpbar$-algèbres induisant un morphisme non nul de $G$-représentations 
\begin{equation} \label{morphnonnul} \ind_K^G \, V \otimes_\chi \Fpbar \to \Ind_B^G \, \rho_1 \otimes \rho_2 .\end{equation}
Alors $V_{U \cap K}$ est un caractère de $A \cap K$ qui est une sous-représentation de $\rho_1 \otimes \rho_2 |_{A \cap K}$, et $\chi$ se factorise à travers ${}' \Sss_G$ en un caractère de $\Hh \big( A, A \cap K, V_{U \cap K} \big)$ que l'on notera encore $\chi$. De plus, pour tous $a_1 , a_2 \in \Z$ tels que $\left( \begin{array}{cc} \varpi^{a_1} & 0 \\ 0 & \varpi^{a_2} \end{array} \right)$ soit un élément de $\widetilde{Z}_A \big( V_{U \cap K} \big)$, on a \[ \chi( \tau^{a_1}_{a_2} ) = \rho_1(\varpi^{a_1}) \rho_2(\varpi^{a_2}) \in \Fpbar^\times .\]
\end{lem}
\textsf{Preuve :} \\
L'existence d'un morphisme non nul (\ref{morphnonnul}) implique en particulier que l'espace $\Hom_G \big( \ind_K^G \, V, \Ind_B^G \, \rho_1 \otimes \rho_2 \big)$ est non réduit à $0$. En appliquant deux fois la réciprocité de Frobenius et en mettant à profit l'isomorphisme topologique $B \cap K \backslash K \xrightarrow{\sim} B \backslash G$, on a:
\[ \Hom_G \big( \ind_K^G \, V, \Ind_B^G \, \rho_1 \otimes \rho_2 \big) = \Hom_K \big( V, \Ind_{B \cap K}^K \, \rho_1 \otimes \rho_2 \big) \qquad \qquad \] 
\[ \qquad \qquad = \Hom_{B \cap K} (V, \rho_1 \otimes \rho_2) = \Hom_{A \cap K} \big( V_{U \cap K} , \rho_1 \otimes \rho_2 \big) .\]
De ce fait, $V_{U \cap K}$, qui est un caractère de $A \cap K$ par le lemme \ref{coinvV}.(i), s'injecte dans $\rho_1 \otimes \rho_2 |_{A \cap K}$. \\
Par le Lemme 4.6 du chapitre \ref{Ly12}, $\chi$ se factorise à travers ${}' \Sss_G$ et le caractère $\Hh \big( A, A \cap K, V_{U \cap K} \big) \to \Fpbar$ (que l'on note encore $\chi$) qui en résulte induit un morphisme non nul $\ind_{A \cap K}^A \, V_{U \cap K} \otimes_\chi \Fpbar \to \rho_1 \otimes \rho_2$. Lorsque $\left( \begin{array}{cc} \varpi^{a_1} & 0 \\ 0 & \varpi^{a_2} \end{array} \right)$ est un élément de $\widetilde{Z}_A \big( V_{U \cap K} \big)$, l'action de $\tau^{a_1}_{a_2}$ sur $\ind_{A \cap K}^A \, V_{U \cap K} \otimes_\chi \Fpbar$ est compatible avec celle de $\left( \begin{array}{cc} \varpi^{a_1} & 0 \\ 0 & \varpi^{a_2} \end{array} \right)$ sur $\rho_1 \otimes \rho_2$. Et cette dernière est scalaire, égale à $\rho_1(\varpi^{a_1}) \rho_2(\varpi^{a_2})$ par la proposition \ref{propDtimes}. \hfill$\Box$

\begin{prop} \label{cosocles}
Soient $V$ une $K$-représentation irréductible et $\chi$ un caractère de $\Hh \big( A, A \cap K, V_{U \cap K}\big)$. On note encore $\chi$ le caractère de $\Hh(G,K,V)$ que l'on obtient en le précomposant par ${}' \Sss_G$.
\begin{itemize}
\item[(i)] La $A$-représentation $\ind_{A \cap K}^A \, V_{U \cap K} \otimes_\chi \Fpbar$ est irréductible, de dimension finie et on la notera $\rho_1 \otimes \rho_2$. 
\item[(ii)] Supposons l'hypothèse \ref{HypDiamond}.
\begin{itemize}
\item[(a)] Supposons $V$ régulière. Alors on a un isomorphisme \[ \ind_K^G \, V \otimes_\chi \Fpbar \xrightarrow{\sim} \Ind_B^G \, \rho_1 \otimes \rho_2 \] de $G$-représentations. 
\item[(b)] Supposons $\widetilde{Z}_A \big( V_{U \cap K} \big) \neq A_\Lambda$ ou bien ($\widetilde{Z}_A \big( V_{U \cap K} \big) = A_\Lambda$ et $\chi(\tau^1_0) \neq \chi(\tau^0_1)$). Alors $\ind_K^G \, V \otimes_\chi \Fpbar$ est isomorphe à $\Ind_B^G \, \rho_1 \otimes \rho_2$.
\end{itemize}
\end{itemize}
\end{prop}
\textsf{Preuve :} \\
La $(A \cap K)$-représentation $V_{U \cap K}$ se décompose en $\sigma_1 \otimes \sigma_2$ où $\sigma_1, \sigma_2$ sont des caractères de $\Oo_D^\times$. Mais alors on a l'isomorphisme de groupes suivants, dans les notations du paragraphe \ref{Dtimes}:
\[ \left( \begin{array}{cc} N_{\sigma_1}/k_D^\times & 0 \\ 0 & N_{\sigma_2}/k_D^\times \end{array} \right) \simeq \widetilde{Z}_A \big( V_{U \cap K} \big) = \left( \begin{array}{cc} \varpi^{d_1 \Z} & 0 \\ 0 & \varpi^{d_2 \Z} \end{array} \right) \]
pour des entiers $d_1, d_2 \geq 1$ uniquement déterminés. Donc, par la proposition \ref{Dtimescompact}, $\ind_{A \cap K}^A \, V_{U \cap K} \otimes_\chi \Fpbar$ est isomorphe, en tant que $A$-représentation, à $\rho \big( \sigma_1, \chi(\tau^{d_1}_0) \big) \otimes \rho \big( \sigma_2 , \chi(\tau^0_{d_2})\big)$, qui est irréductible par la proposition \ref{propDtimes}. Le \textit{(i)} est prouvé. \\
Supposons $V$ régulière. Alors, par \cite{HenVig11b}, Theorem 1.2 et Corollary 1.3, on a l'isomorphisme de $G$-représentations 
\[ \ind_K^G \, V \otimes_\chi \Fpbar \xrightarrow{\sim} \Ind_B^G \big( \ind_{A \cap K}^A \, V_{U \cap K} \otimes_\chi \Fpbar \big) = \Ind_B^G \, \rho_1 \otimes \rho_2 .\]
Prouvons maintenant \textit{(ii).(b)}. La condition $\widetilde{Z}_A \big( V_{U \cap K} \big) \neq A_\Lambda$ est équivalente à ce que $\rho_1$ ou $\rho_2$ est de dimension strictement supérieure à $1$; et la condition ($\widetilde{Z}_A \big( V_{U \cap K} \big) = A_\Lambda$ et $\chi(\tau^1_0) \neq \chi(\tau^0_1)$) implique que $\rho_1$ et $\rho_2$ sont deux caractères de $D^\times$ non isomorphes. Dans les deux cas, on peut appliquer le théorème \ref{GL2preuve1} et $\Ind_B^G \, \rho_1 \otimes \rho_2$ est irréductible. Si $V$ est régulière, le résultat est trivial par \textit{(ii).(a)}. Supposons donc $V$ non régulière. \\
Par la coordination des lemmes \ref{parampos} et \ref{coinvV}, on peut prendre $V'$ une $K$-représentation irréductible régulière incluse dans $\Ind_B^G \, \rho_1 \otimes \rho_2$. En fait, on peut même supposer que $V'$ est l'unique $K$-représentation irréductible régulière telle que les $(A \cap K)$-représentations $V_{U \cap K}'$ et $\varphi_1. V_{U \cap K}$ sont isomorphes; c'est donc ce que l'on fera par la suite. En particulier, on a $\widetilde{Z}_A \big( V_{U \cap K}' \big) = \widetilde{Z}_A \big( V_{U \cap K} \big)$ et on peut identifier $\Hh \big( A, A \cap K, V_{U \cap K}' \big)$ et $\Hh \big( A, A \cap K, V_{U \cap K} \big)$ à travers (\ref{HeckeA}); on notera encore $\chi$ le caractère de $\Hh \big( A, A \cap K, V_{U \cap K}' \big)$ résultant de cette identification. Par la proposition \ref{chgtpds2} ci-dessous, on a alors l'isomorphisme de $G$-représentations \[ \ind_K^G \, V \otimes_\chi \Fpbar \xrightarrow{\sim} \ind_K^G \, V' \otimes_\chi \Fpbar \] et le résultat suit de \textit{(ii).(a)}. \hfill$\Box$

\begin{prop} \label{chgtpds2}
Soient $V$ un caractère de $K$ et $\chi$ un caractère de $\Hh(A,A \cap K,V)$ vérifiant $\chi(\tau^1_0) \neq \chi(\tau^0_1)$ si on a $\widetilde{Z}_A(V) = A_\Lambda$. Soit $V'$ l'unique $K$-représentation régulière telle que les $(A \cap K)$-représentations $V_{U \cap K}'$ et $\varphi_1. V_{U \cap K}$ sont isomorphes. On a un isomorphisme de $G$-représentations
\[ \ind_K^G \, V \otimes_\chi \Fpbar \xrightarrow{\sim} \ind_K^G \, V' \otimes_\chi \Fpbar .\]
\end{prop}
\textsf{Remarque :} 
Cet \' enonc\' e ne n\' ecessite pas l'hypoth\` ese \ref{HypDiamond} car on a le th\' eor\` eme \ref{GL2preuve1} \` a notre disposition (\` a la place du corollaire \ref{GL2preuve3}). Il en va de m\^ eme des r\' esultats qui suivent. \\
\textsf{Remarque :} 
La condition sur $\chi$ peut se réécrire, grâce à la Proposition 16.1 du chapitre \ref{Ly12}: $V$ ne se prolonge pas en un caractère de $G$ ou alors la restriction à $A$ du prolongement de $V$ à $G$ n'est pas isomorphe à $\ind_{A \cap K}^A \, V \otimes_{\chi} \Fpbar$. \\
\textsf{Preuve :} \\
Par le Lemme 6.2 du chapitre \ref{Ly12}, il existe un opérateur non nul \[ \varphi^+ : \ind_K^G \, V \to \ind_K^G \, V' \] de support $K \varphi_1 K$, ainsi qu'un opérateur non nul \[ \varphi^- : \ind_K^G \, V' \to \ind_K^G \, V  \] de support $K \varphi_1^r K$ pour $r \geq 1$ minimal tel que l'on ait $\varphi_1^{r+1} \in \widetilde{Z}_A \big( V_{U \cap K} \big)$. On dénote de la manière suivante les opérateurs qu'ils induisent après tensorisation par $\chi$:
\begin{equation} \label{entrelace} \ind_K^G \, V \otimes_\chi \Fpbar \xrightarrow{\overline{\varphi}^+} \ind_K^G \, V' \otimes_\chi \Fpbar \xrightarrow{\overline{\varphi}^-} \ind_K^G \, V \otimes_\chi \Fpbar .\end{equation} 
En effet, comme $\chi$ est un caract\` ere de $\Hh(A,A \cap K, V) = \Hh \big( A,A \cap K, V'_{U \cap K} \big)$, on peut voir cela apr\` es application de la transform\' ee de Satake; or, la structure de module \` a gauche et \` a droite sur $\Hh(A,A \cap K, V,V'_{U \cap K})$ \' etant la m\^ eme une fois identifi\' ees $\Hh(A,A \cap K, V)$ et $\Hh \big( A,A \cap K, V'_{U \cap K} \big)$, $\varphi^+$ et $\varphi^-$ induisent bien respectivement $\overline{\varphi}^+$ et $\overline{\varphi}^-$ comme voulu. \\
On veut montrer que ces morphismes dans (\ref{entrelace}) sont tous deux non nuls: on va le faire pour $\overline{\varphi}^-$, le cas de $\overline{\varphi}^+$ étant identique. Il suffit de voir que l'image de $[1,v']$ pour $v' \in V'_{U \cap K}$ non nul par $\varphi^-$ n'est pas dans\footnote{on utilise les notations de l'appendice \ref{appsylva}} $\Gamma(V,\chi)$. Notons $p_U : V \twoheadrightarrow V_{U \cap K}$ la projection canonique. Si on fixe un isomorphisme de $\Fpbar$-droites $V'_{U \cap K} \xrightarrow{\sim} V_{U \cap K}$ envoyant $v'$ sur $v$, et que l'on écrit $K \varphi_1^r K = \coprod_{i \in \I} K \varphi_1^r k_i$, le support $X_{\varphi^-}$ de \[ \varphi^-([1,v']) = \sum\nolimits_{i \in \I} [k_i^{-1} \varphi_1^{-r},p_U(k_i v)] \] vérifie $e_\Z(X_{\varphi^-}) =1$ et $\delta_\T(X_{\varphi^-}) = 2r$ (pour ce diam\` etre, l'argument est le m\^ eme que dans le lemme \ref{support1demi} en appendice). Par la proposition \ref{mainsupport}.(i), $\varphi^-([1,v'])$ n'appartient pas à $\Gamma(V,\chi)$: $\overline{\varphi}^-$ est bien non nul comme voulu. \\
Comme $\varphi^- \circ \varphi^+$ est un élément de $\Hh(G,K,V)$, il existe un $\lambda \in \Fpbar$ (la valeur de $\chi$ en $\varphi^- \circ \varphi^+$) tel que $\overline{\varphi}^- \circ \overline{\varphi}^+$ soit $\lambda \, \Id$. De même, $\overline{\varphi}^+ \circ \overline{\varphi}^-$ est $\lambda' \Id$ pour un certain $\lambda' \in \Fpbar$; montrons $\lambda' = \lambda$. On a 
\[ \lambda' \overline{\varphi}^+ = (\overline{\varphi}^+ \circ \overline{\varphi}^-) \circ \overline{\varphi}^+ = \overline{\varphi}^+ \circ (\overline{\varphi}^- \circ \overline{\varphi}^+) = \lambda \overline{\varphi}^+ .\]
Comme $\overline{\varphi}^+$ est non nul, on a bien $\lambda' = \lambda$. Montrons que $\lambda$ est inversible, et la preuve sera terminée. Par la proposition \ref{cosocles}.(ii).(a) et la discussion dans la preuve, ainsi que le théorème \ref{GL2preuve1}, $\ind_K^G \, V' \otimes_\chi \Fpbar$ est isomorphe à une série principale. De ce fait, $\overline{\varphi}^+$ est surjective, $\overline{\varphi}^-$ est injective et $\overline{\varphi}^- \circ \overline{\varphi}^+$ ne peut être nulle. On a ainsi $\lambda \neq 0$ et $\overline{\varphi}^+$ est l'isomorphisme cherché. \hfill$\Box$ \\

Parmi les $V'$ telles que $V_{U \cap K}$ et $V_{U \cap K}'$ sont conjuguées par $\varphi_1^{\Z}$, une seule au plus n'est pas régulière; en effet, dans ce cas-là, en particulier, $V'_{U \cap K}$ est un caractère de la forme $\chi \otimes \chi$ et les conjugués de $V_{U \cap K} = \chi_1 \otimes \chi_2$ par $\varphi_1^\Z$ étant les $\chi_1' \otimes \chi_2$ pour $\chi_1'$ parcourant les conjugués de $\chi_1$ par $\varpi^\Z$, un seul au plus est de la forme $\chi_2 \otimes \chi_2$. Alors, la proposition \ref{cosocles}.(ii).(a) et la proposition \ref{chgtpds2} (ou bien la proposition \ref{cosocles}.(ii).(b)) impliquent directement l'énoncé suivant.

\begin{cor}
Soient $V$ une $K$-représentation irréductible et $\chi$ un caractère de $\Hh \big( A, A \cap K, V_{U \cap K}\big)$. On note encore $\chi$ le caractère de $\Hh(G,K,V)$ que l'on obtient en le précomposant par ${}' \Sss_G$. Soit $V'$ une $K$-représentation irréductible telle que $V_{U \cap K}'$ et $V_{U \cap K}$ sont conjuguées\footnote{en particulier, grâce à (\ref{HeckeA}), on peut identifier $\Hh \big( A, A \cap K, V_{U \cap K}' \big)$ et $\Hh \big( A, A \cap K, V_{U \cap K} \big)$, et donc faire de $\chi$ un caractère de $\Hh(G,K,V')$}. Supposons $\widetilde{Z}_A \big( V_{U \cap K} \big) \neq A_\Lambda$ ou bien ($\widetilde{Z}_A \big( V_{U \cap K} \big) = A_\Lambda$ et $\chi(\tau^1_0) \neq \chi(\tau^0_1)$). Alors on a un isomorphisme \[ \ind_K^G \, V \otimes_\chi \Fpbar \xrightarrow{\sim} \ind_K^G \, V' \otimes_\chi \Fpbar \] de $G$-représentations. 
\end{cor}
\textsf{Preuve :} \\
Par hypoth\` ese, il existe $a,b \in \Z$ tels que $V'_{U \cap K}$ et $(\varphi_1^a \varpi^b).V_{U \cap K}$ sont des $(A \cap K)$-repr\' esentations isomorphes. Comme l'action par conjugaison de $\varpi^{d\Z}$ stabilise chaque caract\` ere de $\Oo_D^\times$, on peut supposer $a \geq 0$, ce que l'on va donc faire. Quitte \` a \' echanger les r\^ oles de $V$ et $V'$, on peut supposer que $V'$ n'est pas un caract\` ere de $K$. Soit $V''$ la $K$-repr\' esentation qui a la m\^ eme r\' egularit\' e que $V$ et qui est telle que $V''_{U \cap K}$ et $\varpi^b. V_{U \cap K}$ sont isomorphes. On a alors un op\' erateur non nul $\ind_K^G \, V \to \ind_K^G \, V''$ port\' e par la classe de $\varpi^b K$. Parce que ce dernier est port\' e par une simple classe, il est un isomorphisme d'inverse un op\' erateur port\' e par la classe $\varpi^{-b} K$ (en ajustant les scalaires de sorte \` a ce que la compos\' ee fasse bien l'identit\' e). Tout ceci passe au quotient pour induire 
\[ \ind_K^G \, V \otimes_\chi \Fpbar \xrightarrow{\sim} \ind_K^G \, V'' \otimes_\chi \Fpbar .\]
Il reste \` a \' etablir que $\ind_K^G \, V'' \otimes_\chi \Fpbar$ et $\ind_K^G \, V' \otimes_\chi \Fpbar$ sont isomorphes. Soit $V_1$ l'unique $K$-repr\' esentation r\' eguli\` ere telle que $(V_1)_{U \cap K}$ et $\varphi_1. V''_{U \cap K}$ sont des $(A \cap K)$-repr\' esentations isomorphismes. Deux cas se pr\' esentent \` a nous: ou bien $V''$ est un caract\` ere de $K$ et l'isomorphisme de $G$-repr\' esentations $\ind_K^G \, V'' \otimes_\chi \Fpbar \xrightarrow{\sim} \ind_K^G \, V_1 \otimes_\chi \Fpbar$ suit directement de la proposition \ref{chgtpds2}. Ou bien $V''$ n'est pas un caract\` ere et on peut utiliser la preuve de la proposition \ref{chgtpds2} (qui est encore valable en plus facile car ni $V''$, ni $V_1$ ne sont des caract\` eres). Dans les deux cas, apr\` es $b$ \' etapes, on obtient
\[ \ind_K^G \, V'' \otimes_\chi \Fpbar \xrightarrow{\sim} \ind_K^G \, V_1 \otimes_\chi \Fpbar \xrightarrow{\sim} \dots \xrightarrow{\sim} \ind_K^G \, V_b \otimes_\chi \Fpbar = \ind_K^G \, V' \otimes_\chi \Fpbar \] 
puisque $V'$ n'est pas un caract\` ere. Le r\' esultat est prouv\' e. \hfill$\Box$\\

On remarque que la proposition \ref{cosocles} implique que l'on connaît le $K$-socle des représentations de Steinberg, c'est-à-dire leur plus grande sous-$K$-représentation semisimple. Ce fait est établi en toute généralité dans le chapitre \ref{Ly11} mais on l'obtient ici à moindres frais.

\begin{cor} \label{socleSt}
Soit $\rho_0$ un caractère de $F^\times$. On factorise sa restriction à $\Oo_F^\times$ par la réduction modulo $(\varpi_F)$:
\[ \rho_0|_{\Oo_F^\times} : \Oo_F^\times \twoheadrightarrow k_F^\times \xrightarrow{\overline{\rho}_0} \Fpbar^\times ;\]
et on note 
\begin{equation} \label{defeta} \eta: k_D^\times \xrightarrow{N_{k_D/k_F}} k_F^\times \xrightarrow{\overline{\rho}_0} \Fpbar^\times .\end{equation}
Alors le $K$-socle de $\St_B \rho_0$ est irréductible et isomorphe à $(\eta \circ \overline{\det}) \otimes \Sym^{\vec{p-1}} \, \Fpbar^2$.
\end{cor}
\textsf{Preuve :} \\
Supposons le $K$-socle de $\St_B \rho_0$ non irréductible. Alors il existe deux $K$-représentations irréductibles $V_1$ et $V_2$ vérifiant $V_1 \oplus V_2 \subseteq \big( \St_B \rho_0 \big)|_K$. Mais alors, par le lemme \ref{coinvV}, $\big( \St_B \rho_0 \big)^{I(1)}$ contiendrait la somme directe $V_1^{I(1)} \oplus V_2^{I(1)}$ de deux droites; et ceci contredit la proposition \ref{StI(1)}.(i). L'hypothèse initiale est donc fausse et le $K$-socle de $\St_B \rho_0$ est irréductible. \\
Soit $V = (\eta \circ \overline{\det}) \otimes \Sym^{\vec{p-1}} \, \Fpbar^2$ avec $\eta$ comme en (\ref{defeta}). Par le lemme \ref{coinvV}, $\widetilde{Z}_A \big(V_{U \cap K} \big)$ est $A_\Lambda$. Soit alors $\chi$ le caractère de $\Hh \big( A, A \cap K, V_{U \cap K} \big)$ défini par $\chi(\tau^1_0) = \chi(\tau^0_1) = \rho_0(\Nrm \, \varpi)$. Par les \textit{(i)} et \textit{(ii)} de la proposition \ref{cosocles}, on obtient alors un isomorphisme de $G$-représentations \[ \ind_K^G \, V \otimes_\chi \Fpbar \xrightarrow{\sim} \Ind_B^G \big( \rho_0 \circ \det\nolimits_G \big) ,\] et donc une surjection \[ \ind_K^G \, V \otimes_\chi \Fpbar \twoheadrightarrow \St_B \rho_0 .\]
En particulier, on a un morphisme non nul de $K$-représentations $V \to \St_B \rho_0$; c'est une injection car $V$ est irréductible. Il en résulte que $V$ est le $K$-socle de $\St_B \rho_0$. \hfill$\Box$\\

Il reste uniquement à examiner le cas où $V$ est un caractère de $K$ (en particulier $V|_{A \cap K} = V_{U \cap K})$ vérifiant $\widetilde{Z}_A(V) = A_\Lambda$ et où $\chi$ vérifie $\chi(\tau^1_0) = \chi(\tau^0_1)$. C'est l'énoncé suivant.

\begin{prop} \label{cosocles2}
Soient $V$ un caractère de $K$ et $\chi$ un caractère de $\Hh \big( A, A \cap K, V \big)$. On note encore $\chi$ le caractère $\chi \circ {}' \Sss_G$ de $\Hh(G,K,V)$. Supposons $\widetilde{Z}_A(V) = A_\Lambda$ et $\chi(\tau^1_0) = \chi(\tau^0_1)$. 
\begin{itemize}
\item[(i)] Tout quotient irréductible de dimension finie de $\ind_K^G \, V \otimes_\chi \Fpbar$ est isomorphe à $\xi \circ \det_G$ où $\xi$ est un caractère de $F^\times$ uniquement déterminé par $V$ et $\chi$. Réciproquement, tout caractère de $G$ est quotient d'un $\ind_K^G \, V \otimes_\chi \Fpbar$ pour une paire $(V,\chi)$ comme précédemment.
\item[(ii)] Tout quotient irréductible admissible de $\ind_K^G \, V \otimes_\chi \Fpbar$ est de dimension finie.
\end{itemize}
\end{prop}
\textsf{Remarque :} 
On pourrait en fait prouver que $\ind_K^G \, V \otimes_\chi \Fpbar$ s'inscrit dans une suite exacte non scindée de $G$-représentations 
\[ 0 \to \St_B \xi \to \ind_K^G \, V \otimes_\chi \Fpbar \to \xi \circ \det\nolimits_G \to 0 .\]
Il faudrait alors retravailler la combinatoire de \cite{BarLiv95} mais cela permettrait d'éviter l'appel à la Proposition 14.1 du chapitre \ref{Ly12}. \\
\textsf{Preuve :} \\
Soit $\pi$ un quotient irréductible de dimension finie de $\ind_K^G \, V \otimes_\chi \Fpbar$. Par le lemme \ref{factordet} et parce que $\Fpbar$ est algébriquement clos, $\pi$ est un caractère de $G$: il existe un caractère $\xi$ de $F^\times$ vérifiant $\pi = \xi \circ \det_G$. A cause de la relation (\ref{detNrm}), $\xi \circ \det_G$ est une sous-représentation de $\Ind_B^G \, (\xi \circ \Nrm) \otimes (\xi \circ \Nrm)$; et on a donc un morphisme non nul de $G$-représentations \[ \ind_K^G \, V \otimes_\chi \Fpbar \to \Ind_B^G \, (\xi \circ \Nrm) \otimes (\xi \circ \Nrm) .\] Par le lemme \ref{parampos}, on a alors
\begin{equation} \label{xientier} V_{U \cap K} = (\xi \circ \Nrm)|_{\Oo_D^\times} \otimes (\xi \circ \Nrm)|_{\Oo_D^\times} .\end{equation} 
Comme $\Nrm : \Oo_D^\times \to \Oo_F^\times$ est surjective, cela détermine $\xi|_{\Oo_F^\times}$. De plus, (\ref{xientier}) implique que $\widetilde{Z}_A \big( V_{U \cap K} \big)$ est $A_\Lambda$ et le lemme \ref{parampos} donne ensuite \[ \chi(\tau^1_0) = \xi(\Nrm \, \varpi) = \xi \big( (-1)^d \varpi^d \big) .\] 
Comme on a $F^\times = \Oo_F^\times \varpi^{d \Z}$, $\xi$ est uniquement déterminé. \\
Réciproquement, soit $\pi$ un caractère de $G$: il s'écrit $\pi = \xi \circ \det_G$ pour un certain caractère $\xi$ de $F^\times$. Factorisons la restriction de $\xi|_{\Oo_F^\times}$ par la réduction modulo $(\varpi_F)$ 
\[ \xi|_{\Oo_F^\times} : \Oo_F^\times \twoheadrightarrow k_F^\times \xrightarrow{\overline{\xi}} \Fpbar^\times ,\] et notons \[ \eta: k_D^\times \xrightarrow{N_{k_D/k_F}} k_F^\times \xrightarrow{\overline{\xi}} \Fpbar^\times .\]
Soit $V = \eta \circ \overline{\det}$, qui est un caractère de $K$. En particulier, on a $\widetilde{Z}_A \big( V_{U \cap K} \big) = \widetilde{Z}_A \big( V|_{A \cap K} \big) = A_\Lambda$. Soit $\chi$ le caractère de $\Hh(A, A \cap K, V)$ défini par \[ \chi(\tau^1_0) = \chi(\tau^0_1) = \xi(\Nrm \, \varpi) .\] 
Il s'agit maintenant de voir que l'on a une surjection 
\begin{equation} \label{surjcaract} \ind_K^G \, V \otimes_\chi \Fpbar \twoheadrightarrow \xi \circ \det\nolimits_G .\end{equation} 
Or on a la surjection $G$-équivariante 
\[ \begin{array}{ccc} \ind_K^G \, V & \to & \xi \circ \det\nolimits_G \\ {} [g,1] & \mapsto & \xi \circ \det\nolimits_G(g) \end{array} ,\]
qui passe au quotient pour donner (\ref{surjcaract}) par le début de la preuve. Le \textit{(i)} est terminé. \\
Prouvons maintenant \textit{(ii)}. Soit $\pi$ une représentation irréductible admissible qui est quotient d'une telle $\ind_K^G \, V \otimes_\chi \Fpbar$. Ecrivons $V = \xi_K \circ \overline{\det}$ pour $\xi_K$ un caractère de $k_D^\times$. Parce que $\widetilde{Z}_A(V)$ est $A_\Lambda$, $\xi_K$ possède la factorisation \[ \xi_K : k_D^\times \xrightarrow{N_{k_D/k_F}} k_F^\times \xrightarrow{\xi'} \Fpbar^\times .\] Par inflation de $\xi'$ à $\Oo_F^\times$, on obtient un caractère $\Oo_F^\times \to \Fpbar^\times$ que l'on note encore $\xi'$, puis on prolonge $\xi'$ en $\xi_1: F^\times \to \Fpbar^\times$ par $\xi_1|_{\Oo_F^\times} = \xi'$ et $\xi_1((-1)^d \varpi_F) = \chi(\tau^1_0)$. Par le Lemme 4.5 du chapitre \ref{Ly12}, la représentation admissible $\pi' = ( \xi_1^{-1} \circ \det_G) \otimes \pi$ de $G$ est quotient de $\ind_K^G \, \id \otimes_{\chi'} \Fpbar$ où $\chi' : \Hh(G,K) \to \Fpbar$ est défini par $\chi'(\tau^1_0) = \chi'(\tau^0_1) = 1$. Par la Proposition 14.1 du chapitre \ref{Ly12} appliquée à $\pi'$ et à $B^-$, ou bien $\pi'$ contient la représentation triviale, ou bien on a $\Ord_{B^-} \pi' \neq 0$ (voir \cite{Eme10}, Definition 3.1.9, et \cite{Vig12}, Definition 1.8). Dans le premier cas de figure, $\pi$ contiendrait le caractère $\xi_1 \circ \det_G$ et lui serait égal par irréductibilité; alors $\pi$ serait de dimension finie comme voulu. \\
Montrons que le cas $\Ord_{B^-} \pi' \neq 0$ est exclu. En effet, dans ce cas-là, par \cite{HenVig11b}, point 2 du cheminement de la preuve de la Proposition 7.9 et Lemma 7.10, $\Ord_{B^-} \pi'$ contient une représentation irréductible admissible de $A$ que l'on note $\rho = \rho_1 \otimes \rho_2$. Parce que $\rho_i = \rho_i^{D(1)}$ est de dimension finie, $\rho_i$ est décrite par la section \ref{Dtimes}, pour $i \in \{1,2\}$. Par la propriété d'adjonction de $\Ord_{B^-}$ (voir \cite{Vig12}, Proposition 1.7 et Remark 1.9.(ii), et \cite{Eme10}, Theorem 4.4.6) on a alors \[ \Hom_G \big( \Ind_B^G \, \rho_1 \otimes \rho_2, \pi' \big) \simeq \Hom_A \big( \rho_1 \otimes \rho_2 , \Ord_{B^-} \pi' \big) .\] De ce fait, par le corollaire \ref{St2} et le théorème \ref{GL2preuve1}, $\pi'$ est une série principale ou une représentation de Steinberg. Le cas $\pi' = \Ind_B^G \, \rho$ série principale est exlu par le lemme \ref{parampos} puisque l'on a $\id_{U \cap K} = \id$, $\widetilde{Z}_A(\id) = A_\Lambda$ et $\chi'(\tau^1_0) = \chi'(\tau^0_1)$. Aussi, comme le $K$-socle d'une représentation de Steinberg est irréductible régulier par le corollaire \ref{socleSt}, le cas Steinberg est aussi exclu. Au final, on ne peut avoir $\Ord_{B^-} \pi' \neq 0$, et $\pi'$ est de dimension finie. \hfill$\Box$ \\

On est maintenant à même de donner un théorème de classification. 

\begin{theo}
Toute représentation irréductible admissible de $G$ est isomorphe à l'une des représentations suivantes:
\begin{itemize}
\item[(a)] $\xi \circ \det_G$ où $\xi$ est un caractère de $F^\times$;
\item[(b)] $\St_B \xi$ où $\xi$ est un caractère de $F^\times$;
\item[(c)] $\Ind_B^G \, \rho_1 \otimes \rho_2$ où $\rho_1, \rho_2$ sont deux représentations irréductibles de $D^\times$ satisfaisant (a) ou (b) du corollaire \ref{GL2preuve3};
\item[(d)] une représentation supersingulière.
\end{itemize}
De plus, il n'y a pas d'entrelacement entre deux représentations quelconques de la liste ci-dessus.
\end{theo}
\textsf{Preuve :} \\
Le fait que toute représentation irréductible admissible $\pi$ de $G$ soit de l'une des formes mentionnées suit de la discussion au début du paragraphe et des propositions \ref{cosocles} et \ref{cosocles2}. \\
Aussi, la classe d'isomorphisme de $\pi$ possède comme invariant l'ensemble $\Pi(\pi)$ des paires $(V,\chi)$ où $V$ est une $K$-représentation irréductible et $\chi : \Hh(G,K,V) \to \Fpbar$ est un caractère de $\Fpbar$-algèbres tels que l'on ait une surjection $\ind_K^G \, V \otimes_\chi \Fpbar \twoheadrightarrow \pi$. Notons $\Pi_{\ns}(\pi)$ le sous-ensemble \og non supersingulier \fg\ de $\Pi(\pi)$, c'est-à-dire 
\[ \Pi_{\ns}(\pi) = \{ (V,\chi) \in \Pi(\pi) \ | \ \chi \textrm{ se factorise par } {}' \Sss_G \} .\] 
Par les propositions \ref{cosocles} et \ref{cosocles2}, on connaît les quotients irréductibles admissibles des $\ind_K^G \, V \otimes_\chi \Fpbar$ pour tout $\chi$ se factorisant par ${}' \Sss_G$ et on remarque que les ensembles $\Pi_{\ns}(\pi)$ sont deux à deux distincts pour $\pi$ dans les cas \textit{(a)}, \textit{(b)} ou \textit{(c)}. Aussi, par définition, on a $\Pi_{\ns}(\pi) = \varnothing$ si et seulement si $\pi$ est dans le cas \textit{(d)}. Cela conclut quant à l'absence d'entrelacement. \hfill$\Box$ \\

On finit par énoncer la conséquence suivante. On rappelle qu'une représentation irréductible de $G$ est dite \textit{supercuspidale} si elle n'est pas sous-quotient d'une $\Ind_B^G \, \rho$ avec $\rho$ irréductible admissible.

\begin{cor}
Soit $\pi$ une représentation irréductible admissible de $G$. Alors $\pi$ est supersingulière si et seulement si elle est supercuspidale.
\end{cor}

\section{Appendice: un peu de sylviculture} \label{appsylva}

Soit $V$ une $K$-représentation irréductible. Pour étudier l'action de l'algèbre de Hecke $\Hh(G,K,V)$ sur $\ind_K^G \, V$, on va regarder un peu la combinatoire de $K \backslash G$. L'espace $K \backslash G$ se plonge dans une infinité dénombrable de copies\footnote{c'est pour cela que l'on a envie de parler de forêt...} d'un arbre de Bruhat-Tits (au sens de \cite{Ser77}, avec action de $G$ via $g.(Kx) = Kx g^{-1}$) pour $\PGL_2$, et c'est ce que l'on commence par établir. On définit le morphisme de groupes $v_\beta = \dfrac{1}{d} \, v_D \circ \det_G : G \to \Z$.

\begin{lem} \label{foret}
$\phantom{,}$
\begin{itemize}
\item[(i)] L'application \[ \begin{array}{cccc} \iota_{\T \Z} : & K \backslash G & \to & K \varpi^\Z \backslash G \times \Z \\ & g & \mapsto & (\beta(g),v_\beta(g)) \end{array} \] est une injection ensembliste.
\item[(ii)] L'espace $K \varpi^\Z \backslash G$ peut être naturellement muni de la structure d'arbre de Bruhat-Tits de $\PGL(2,D)$, que l'on notera $\T$.
\end{itemize}
\end{lem}
\textsf{Preuve :} \\
Montrons \textit{(i)}. Soient $g_1, g_2 \in G$ tels que $K g_1$ et $K g_2$ ont même image par $(\beta, v_\beta)$. Alors il existe $k \in K$ et $i \in \Z$ tels que l'on ait $g_2 = k \varpi^i g_1$. Parce que $\det_G$ envoie $K$ sur $\Oo_F^\times$, on obtient $v_\beta(g_2) = v_\beta(g_1)+i$, et donc $i=0$. Le \textit{(i)} est prouvé. \\
Pour \textit{(ii)}, on se reportera à \cite{Ser77}, qui ne suppose pas la commutativité du corps de base et rentre ainsi dans notre cadre. \hfill$\Box$\\

On remarque que l'injection définie au lemme \ref{foret}.(i) identifie $K \backslash G$ à un sous-ensemble discret de l'espace métrique qu'est l'immeuble élargi de $G$ (voir \cite{Lan00}, paragraphe 1.3). Cependant la métrique usuelle sur l'immeuble élargi ne nous intéresse pas vraiment et on préfère considérer des invariants qui distinguent bien $K \varpi^\Z \backslash G$ et $\Z$: c'est ce que l'on présente maintenant. \\ Soit $d_\T$ la distance usuelle sur l'espace $\T_0 = K \varpi^\Z \backslash G$ des sommets de $\T$: deux sommets $x$ et $y$ de $\T$ sont à distance $1$ si et seulement si il existe une arête de $\T$ qui les relie. On a, pour tout $n \geq 0$, 
\begin{equation} \label{dist1} d_\T(x,y)=n \Leftrightarrow x y^{-1} \in K \varphi_1^n \varpi^\Z K .\end{equation}
On se sert de $d_\T$ pour définir un diamètre $\delta_\T$ sur la \og forêt \fg\ $\T \times \Z$ de sommets $K \backslash G$ par $(\beta,v_\beta)^{-1}$. \\ Soit $X$ une partie finie non vide de $K \backslash G$; on note \[ \delta_\T(X) = \max \, \{ d_\T( \beta(x),\beta(y)) \ | \ x,y \in X , \ v_\beta(x) = v_\beta(y) \} .\]
Aussi, on définit l'étendue $e_\Z$ de $X$ par: 
\[ e_\Z(X) = \max \, \{ v_\beta(x) \ | \ x \in X \} - \min \, \{ v_\beta(x) \ | \ x \in X \} +1 .\]
Par convention, on prendra $\delta_\T(\varnothing) = 0$ et $e_\Z(\varnothing) = 0$. \\

Si on munit $\T \times \Z$ de l'action de $G$ définie par 
\begin{equation} \label{actionforet} g. \big( K \varpi^\Z h,a \big) = \big( K \varpi^\Z h g^{-1} , a - v_\beta(g) \big) ,\end{equation}
alors l'injection $\iota_{\T \Z}$ est $G$-équivariante. On remarque que l'action (\ref{actionforet}) de $G$ sur $\T$ est $d_\T$-isométrique. Cela implique facilement l'énoncé suivant.

\begin{lem}
Soit $X$ une partie finie de $K \backslash G$. Alors, pour tout $g \in G$, on a $e_\Z(gX) = e_\Z(X)$ et $\delta_\T(gX) = \delta_\T(X)$.
\end{lem}
\textsf{Preuve :} \\
Soit $g \in G$. Pour tous $x, y \in X$, on a $v_\beta(g.x) = v_\beta(x) - v_\beta(g)$ ; et donc on a \[ v_\beta(g.x) = v_\beta(g.y) \Leftrightarrow v_\beta(x) = v_\beta(y) .\] La preuve est terminée. \hfill$\Box$\\

Soit $V$ une représentation irréductible de $K$. Le support d'une fonction de $\ind_K^G \, V$ est une partie finie de $K \backslash G$, et l'objectif de ce paragraphe est de donner des conditions nécessaires sur $\delta_\T(X)$ et $e_\Z(X)$ pour qu'une fonction de support $X$ soit dans l'image de certains opérateurs de Hecke de $\Hh(G,K,V)$. \\
Soit $\chi : \Hh(G,K,V) \to \Fpbar$ un caractère de $\Fpbar$-algèbres qui se factorise par ${}' \Sss_G$. On définit le $\Fpbar$-espace vectoriel \[ \Gamma(V,\chi) = \ker \big( \ind_K^G \, V \to \ind_K^G \, V \otimes_\chi \Fpbar \big) .\]
On suppose à présent $\widetilde{Z}_A \big( V_{U \cap K} \big)$ de la forme $\left( \begin{array}{cc} \varpi^{d_0 \Z} & 0 \\ 0 & \varpi^{d_0 \Z} \end{array} \right)$ pour un diviseur $d_0 \geq 1$ de $d$. Cela nous assure en particulier que l'algèbre de Hecke $\Hh(G,K,V)$ est engendrée par les opérateurs suivants: l'opérateur $T$ de support $K \varphi_1^{d_0} K$ et induisant l'identité sur $V_{U \cap K}$, l'opérateur $Z$ de support $\varpi^{d_0} K$ et induisant l'identité sur $V$, et son inverse $Z^{-1}$ (voir \cite{HenVig11}, paragraphe 7.3). De ce fait, une famille génératrice de $\Gamma(V,\chi)$ est donnée par les $T^i Z^j [g,v] - \chi(T^i Z^j) [g,v]$ pour $i \geq 0$, $j \in \Z$, $g \in G$ et $v \in V$. 

\begin{lem} \label{support1}
Soient $V, \chi, d_0$ comme précédemment. Soit $x$ un élément de $\Gamma(V,\chi)$ engendré par les $Z^j [g,v] - \chi(Z^j)[g,v]$. Alors le support $X_x$ de $x$ est d'étendue $e_\Z(X_x)$ nulle ou strictement supérieure à $2 d_0$.
\end{lem}
\textsf{Preuve :} \\
Soit $x$ un élément de $\Gamma(V,\chi)$ engendré par les $Z^j[g,v] - \chi(Z^j) [g,v]$. On écrit 
\begin{equation} \label{xenZ} x = \sum\nolimits_{i \in \I} Z^{j_i} [g_i,v_i] - \chi(Z^{j_i})[g_i,v_i] ,\end{equation}
où $\I$ est un ensemble fini, les $g_i \in G$ et $v_i \in V$ pour tout $i \in \I$. On suppose l'écriture (\ref{xenZ}) minimale, au sens de $|\I|$ minimal. L'étendue $e_\Z(X_x)$ du support $X_x$ de $x$ est nulle si et seulement si $x$ est nul. Supposons $x$ non nul et montrons que $e_\Z(X_x)$ est strictement supérieure à $2 d_0$. \\
L'ensemble $\I$ est non vide puisque $x$ est non nul. On prend alors $i \in \I$ et le terme $Z^{j_i}[g_i,v_i] - \chi(Z^{j_i})[g_i,v_i]$ apparaît dans l'écriture (\ref{xenZ}). On allège les notations en $g =g_i, \, v = v_i, \, j= j_i$ pour la suite. On a l'alternative suivante: ou bien $K g^{-1}$ et $K \varpi^{jd_0} g^{-1}$ sont tous deux dans $X_x$, ou bien l'un au moins ne l'est pas. Dans le premier cas, on a \[ e_\Z(X_x) \geq 2 j d_0 +1 \geq 2 d_0+1 \] puisque $j$ est non nul par minimalité de $|\I|$; la preuve est terminée. \\
Montrons que le second cas aboutit à une contradiction. En effet, disons que $Kg^{-1}$ n'est pas dans le support de $x$, le cas de $K \varpi^{j d_0} g^{-1}$ étant similaire grâce à l'identité 
\begin{equation} \label{inverseZ} Z^j [g,v] - \chi(Z^j) [g,v] = - \chi(Z^j) \, \big( Z^{-j} [g \varpi^{-j d_0}, v] - \chi(Z^{-j}) [g \varpi^{-j d_0}, v] \big) .\end{equation}
Alors il existe $i_1, \dots, i_n \in \I$ (avec $n \geq 1$) et $g_{i_1}, \dots, g_{i_n} \in G$ vérifiant $K g_{i_k}^{-1} = K g^{-1}$ pour tout $k$, et tels que si on écrit $[g,v_k'] = [g_{i_k},v_{i_k}]$ (c'est-à-dire $v_k' = g^{-1} g_{i_k} v_{i_k}$) alors on a 
\begin{equation} \label{annulZ} \chi(Z^j) \, v + \sum_{k=1}^n \chi(Z^{j_{i_k}}) \, v_k' = 0 \end{equation}
(quitte à utiliser (\ref{inverseZ}) pour les $j_{i_k}$ au lieu de $j$). On veut alors raccourcir l'écriture de 
\begin{align*} y & = Z^j [g,v] - \chi(Z^j) [g,v] + \sum_{k=1}^n \big( Z^{j_{i_k}} [g_{i_k}, v_{i_k}] - \chi(Z^{j_{i_k}}) [g_{i_k}, v_{i_k}] \big) \\ & = Z^j [g,v] - \chi(Z^j) [g,v] + \sum_{k=1}^n \big( Z^{j_{i_k}} [g, v_k'] - \chi(Z^{j_{i_k}}) [g, v_k'] \big) .\end{align*}
A cause de (\ref{annulZ}), on a \[ y = Z^j [g,v] + \sum_{k=1}^n Z^{j_{i_k}} [g,v_k'] .\]
On réutilise à nouveau (\ref{annulZ}) pour développer $Z^j [g,v]$ et obtenir:
\begin{align*} y & = - \sum_{k=1}^n \chi(Z^{-j}) \chi(Z^{j_{i_k}}) \, Z^j [g,v_k'] + \sum_{k=1}^n Z^{j_{i_k}} [g,v_k'] \\ & = \sum_{k=1}^n \big( Z^{j_{i_k}-j} [g \varpi^{-j d_0}, v_k'] - \chi(Z^{j_{i_k}-j})[g \varpi^{-jd_0}, v_k'] \big) .\end{align*}
Cette expression de $y$ est plus courte que l'écriture initiale; et si on injecte cela dans (\ref{xenZ}), cela contredit la minimalité de $|\I|$. C'est ce qu'on voulait et la preuve est terminée. \hfill$\Box$\\

On note $p_U: V \to V_{U \cap K}$ la projection canonique.

\begin{lem} \label{support1demi}
Soient $V, \chi, d_0$ comme précédemment. Soient $g \in G$, $v \in V$ non nul et \[ y(g,v) = T[g,v] - \chi(T) [g,v] , \quad y'(g,v) = Z^{-1} T [g, \chi(Z) v] - \chi(T)[g,v] .\] 
\begin{itemize}
\item[(i)] Le support $X_y$ de $y(g,v)$ est d'étendue $e_\Z(X_y)$ égale à $d_0+1$ et de diamètre $\delta_\T(X_y)$ égal à $2 d_0$.
\item[(ii)] Le support $X_y'$ de $y'(g,v)$ est d'étendue $e_\Z(X_y')$ égale à $d_0+1$ et de diamètre $\delta_\T(X_y')$ égal à $2 d_0$.
\end{itemize}
Si on écrit $K \varphi_1^{d_0} K = \coprod_{i=1..r} K \varphi_1^{d_0} k_i$ pour des $k_i \in K$, alors on pose 
\[ z(g,v) = y(g,v) + \chi(T)^{-1} \sum_{i=1}^r y' \big( g k_i^{-1} \varphi_1^{-d_0} , p_U(k_i v) \big) .\]
\begin{itemize}
\item[(iii)] Le support $X_z$ de $z(g,v)$ est d'étendue $e_\Z(X_z)$ égale à $1$ et de diamètre $\delta_\T(X_z)$ égal à $4 d_0$.
\end{itemize}
\end{lem}
\textsf{Preuve :} \\
Si on écrit $K \varphi_1^{d_0} K = \coprod_{i=1,.., r} K \varphi_1^{d_0} k_i$ pour des $k_i \in K$, alors on a 
\[ y(g,v) = - \chi(T) [g,v] + \sum_{i=1,..,r} [g k_i^{-1} \varphi_1^{-d_0} , p_U(k_i v) ] .\]
Chacun des supports des fonctions $[g k_i^{-1} \varphi_1^{-d_0} , p_U(k_i v) ]$ est un sommet $x_i \in K \backslash G$ avec $v_\beta(x_i) = v_\beta(g)+d_0$. Voyons que suffisamment de $p_U(k_i v)$ sont non nuls. En effet, le noyau de $p_U$ est engendré linéairement par les $v' - u.v'$ pour $v' \in V$ et $u \in U \cap K$. En particulier, $\Ker \, p_U$ est stable par l'Iwahori $I$. Supposons que $k \notin I$ envoie un élément $v' \in \Ker \, p_U$ non nul dans $\Ker \, p_U$. Alors, par décomposition de Bruhat, $k$ est dans $I s I$; et le sous-groupe de $K$ engendré par $k$ et $I$ est $K$. De cette manière, $K .v' \subseteq \Ker \, p_U$ est une sous-$K$-représentation qui est de codimension au moins $1$ dans $V$, et cela contredit l'irréductibilité de $V$. On en déduit que $k.v , k'.v \in \Ker \, p_U$ implique $k' k^{-1} \in I$, et au plus $|K \cap \varphi_1^{-d_0} K \varphi_1^{d_0} \backslash I | = q^{d(d_0-1)}$ parmi les $|K \cap \varphi_1^{-d_0} K \varphi_1^{d_0} \backslash K| = q^{d(d_0-1)}(q^d+1)$ vecteurs $p_U(k_iv)$ peuvent s'annuler. \\
A cause de cela, et parce que $\chi(T)$ est non nul, $X_y$ est d'étendue $d_0+1$. Enfin, parmi les $q^{d(d_0-1)}(q^d+1)$ sommets à distance $d_0$ de $\beta(K g^{-1})$, au moins $q^{d d_0}$ d'entre eux portent une fonction non nulle: il en existe alors deux qui réalisent le diamètre $\delta_\T(X_y) = 2 d_0$. En effet, si on fixe $\beta(K(g')^{-1})$ \` a distance $d_0$ de $\beta(Kg^{-1})$, il y a $q^{d(d_0-1)}$ sommets \` a distance $d_0$ de $\beta(Kg^{-1})$ et \` a distance strictement inf\' erieure \` a $2 d_0$ de $\beta(K(g')^{-1})$; comme on dispose d'au moins $q^{d d_0}$ sommets, l'affirmation est v\' erifi\' ee. \\
Pour $y'(g,v)$, il s'agit de voir que l'on a 
\[ y'(g,v) = y(g,v) + Z^{-1} T [g, \chi(Z) v] - \chi(Z^{-1}) \, T[g, \chi(Z) v] ,\]
et on obtient \[ X_y' \smallsetminus \{ K g^{-1} \} = \varpi^{-d_0} \big( X_y \smallsetminus \{ K g^{-1} \} \big) .\]
Le \textit{(ii)} s'en déduit. \\
Notons que $z(g,v)$ est bien défini parce que $T$ l'est. On écrit 
\[ z(g,v) = -\chi(T) [g,v] + \sum_{1 \leq i, j \leq r} \chi(Z) \chi(T)^{-1} \big[ g k_i^{-1} \varphi_1^{-d_0} k_j^{-1} \varphi_2^{d_0} , p_U(k_j p_U(k_iv)) \big] \] 
où on remarque que chacun des supports des fonctions de la somme de droite est un sommet $x_{ij} \in K \smallsetminus G$ vérifiant $v_\beta(x_{ij}) = v_\beta(g)$. En particulier, $X_z$ est d'étendue au plus $1$: en prouvant que son diamètre est $4 d_0$, on prouvera alors que $z(g,v)$ est non nul si $v$ est non nul, de sorte que l'étendue de $X_z$ sera exactement égale à $1$. On voit que si $i$ est diff\' erent de $i'$, les sommets $\beta(x_{ij}) = K \varpi^{\Z} \varphi_2^{-d_0} k_j \varphi_1^{d_0} k_i g^{-1}$ et $\beta(x_{i'j'})$ sont \` a distance $4d_0$ l'un de l'autre: en effet, comme $K \varpi^{\Z} \backslash G$ est muni d'une structure d'arbre, il suffit de voir que le chemin 
\[ K \varpi^{\Z} \varphi_2^{-d_0} k_j \varphi_1^{d_0} k_i g^{-1} \leadsto K \varpi^{\Z} \varphi_2^{-d_0+1} k_j \varphi_1^{d_0} k_i g^{-1} \leadsto \dots \leadsto K \varpi^{\Z} \varphi_1^{d_0} k_i g^{-1} \]
\[ \leadsto K \varpi^{\Z} \varphi_1^{d_0-1} k_i g^{-1} \leadsto \dots \leadsto K \varpi^{\Z} g^{-1} \leadsto K \varpi^{\Z} \varphi_1 k_{i'} g^{-1} \leadsto \dots  \]
\[ \phantom{xxxxxxxxxxxxxx}\leadsto K \varpi^{\Z} \varphi_1^{d_0} k_{i'} g^{-1} \leadsto \dots \leadsto K \varpi^{\Z} \varphi_2^{-d_0} k_{j'} \varphi_1^{d_0} k_{i'} g^{-1} \]
ne poss\` ede pas d'aller-retour (c'est-\` a-dire de cycle de longueur $2$). Et c'est pr\' ecis\' ement ce qu'assure la condition $i \neq i'$. On notera aussi que ce chemin est de longueur strictement inf\' erieure \` a $4d_0$ s'il y a un aller-retour (au moins) ; en particulier on a d'ores et d\' ej\` a $\delta_\T(x_z) \leq 4d_0$. \\
Aussi, par l'argument du $(ii)$, au plus $q^{d(2d_0-1)}$ des vecteurs $p_U(k_i v)$ peuvent s'annuler: on prend deux sommets $K \varpi^{\Z} \varphi_1^{d_0} k_i g^{-1}$ et $K \varpi^{\Z} \varphi_1^{d_0} k_{i'} g^{-1}$ tels que $p_U(k_i v)$ et $p_U(k_{i'}v)$ ne s'annulent pas. Parmi les sommets \` a distance $2 d_0$ de $K \varpi^{\Z} \varphi_1^{d_0} k_i g^{-1}$ (resp. $K \varpi^{\Z} \varphi_1^{d_0} k_{i'} g^{-1}$) seuls $q^{d(2d_0-1)}$ peuvent porter un vecteur $p_U(k_j p_U(k_i v))$ (resp. $p_U(k_{j'} p_U(k_{i'} v))$) qui s'annule. Deux tels $\beta(x_{ij})$ et $\beta(x_{i'j'})$ portant des vecteurs qui ne s'annulent pas assurent l'\' etendue voulue puisqu'ils v\' erifient $v_\beta(x_{ij}) = v_\beta(g) = v_\beta(x_{i'j'})$. La preuve est termin\' ee. \hfill$\Box$

\begin{lem} \label{support2}
Soient $V, \chi, d_0$ comme précédemment. Soit $x$ un élément non nul de $\Gamma(V,\chi)$. Supposons que le support $X_x$ de $x$ vérifie $e_\Z(X_x) \leq d_0 + 1$. Il existe $a \in \Z$, des ensembles finis disjoints $\I_1, \I_2, \I_3$ et des $g_i \in G$ avec $v_\beta(K g_i^{-1}) = a$ pour $i \in \I_1$, $v_\beta(K g_i^{-1}) = a+ d_0$ pour $i \in \I_2$ et $v_\beta(K g_i^{-1}) \in ]a,a+d_0[ \cap \Z$ pour $i \in \I_3$, des $v_i \in V$ pour tout $i \in \I_1 \cup \I_2 \cup \I_3$ tels que l'on ait \[ x = \sum_{i \in \I_1} y(g_i,v_i) + \sum_{i \in \I_2} y'(g_i,v_i) + \sum_{i \in \I_3} z(g_i,v_i).\]
\end{lem}
\textsf{Preuve :} \\
Ecrivons $x$ comme une somme finie \[ x = \sum\nolimits_i T^{r_i} Z^{s_i} [g_i,v_i] - \chi(T^{r_i} Z^{s_i}) [g_i,v_i] \] pour des $r_i \geq 0$, $s_i \in \Z$, $g_i \in G$ et $v_i \in V$. En utilisant successivement la relation 
\[ T^i Z^j [g,v] - \chi(T^i Z^j) [g,v] = \big( T^i [g \varpi^{-j},v] - \chi(T^i)[g \varpi^{-j},v] \big) + \big( Z^j [g, \chi(T^i) v] - \chi(Z^j) [g, \chi(T^i) v] \big) ,\]
on peut réécrire $x$ sous la forme
\[ \sum_{i \in \I} \big( T^{r_i}[g_i,v_i] - \chi(T^{r_i}) [g_i,v_i]\big) + \sum_{j \in \J} \big( Z^{s_j} [g_j,v_j] - \chi(Z^{s_j})[g_j,v_j] \big) .\]
Quitte à développer chacun des $T^{r_i}[g_i,v_i] - \chi(T^{r_i})[g_i,v_i]$ récursivement grâce à la formule (pour $k \geq 1$)
\[ T^{k+1} [g,v] - \chi(T^{k+1})[g,v] = (T - \chi(T)) (T^k [g,v]) + T^k [g, \chi(T) v ]- \chi(T^k) [g, \chi(T) v] ,\]
on peut supposer $r_i = 1$ pour tout $i \in \I$: c'est ce qu'on fera par la suite. Notons $a$ le minimum des $v_\beta(Kg)$ pour $Kg \in X_x$. Quitte à translater les sommets à l'aide de relations du type $Z^s [g,v] - \chi(Z^s) [g,v]$, on peut supposer $v_\beta(K g_i^{-1}) \in [a,a+2d_0-1]$ et on le suppose par la suite. Aussi, en réutilisant à nouveau des relations $Z^s [g,v] - \chi(Z^s) [g,v]$ pour transformer les $y(g_i,v_i)$ en $y'(g_i,v_i)$ lorsque l'on a $a+d_0 \leq v_\beta(K g_i^{-1})< a+ 2 d_0$, on écrit $x$ sous la forme \[ x = \sum_{i \in \I_1} y(g_i,v_i) + \sum_{i \in \I_2} y'(g_i,v_i) + \sum_{j \in \J} \big( Z^{s_j} [g_j,v_j] - \chi(Z^{s_j}) [g_j,v_j] \big) ,\]
où on a $a \leq v_\beta(K g_i^{-1})< a+ d_0$ pour tout $i \in \I_1$ et $a+d_0 \leq v_\beta(K g_i^{-1})< a+ 2 d_0$ pour tout $i \in \I_2$. 
Parce que le support de $x$ est d'étendue inférieure ou égale à $d_0$, avec \[ \min \, \{ v_\beta(Kg) \ | \ Kg \in X_x \} = a ,\] \[ x' =x - \sum_{i \in \I_1} y(g_i,v_i) - \sum_{i \in \I_2} y'(g_i,v_i) = \sum_{j \in \J} \big( Z^{s_j} [g_j,v_j] - \chi(Z^{s_j}) [g_j,v_j] \big) \]
a un support ayant une étendue inférieure ou égale à $2 d_0$ (par le lemme \ref{support1demi}). Par le lemme \ref{support1}, $x'$ est nul. \\
Regardons les sommets $K g_i^{-1}$ dans $\I_1$ avec $v_\beta(K g_i^{-1}) = c \in]a,a+d_0[$ fix\' e (si cet intervalle n'intersecte pas $\Z$, il n' y a rien \` a faire). On consid\` ere le \og diamètre des centres \fg\ des \' el\' ements de $\I_1$:
\[ \delta_c = \max\limits_{v_\beta(K g_i^{-1}) = v_\beta(K g_j^{-1}) = c} \, d_\T(\beta(K g_i^{-1}), \beta(K g_j^{-1})) .\]
Deux $K g_i^{-1}$ et $K g_j^{-1}$ réalisant ce diamètre $\delta_c$ possèdent (voir la preuve du lemme \ref{support1demi} à propos de $\Ker \, p_U$) des éléments à distance $\delta_c + 2d_0$ dans leur support qui sont extrémaux et sont donc dans $X_x$, \` a moins qu'ils ne s'annulent avec des sommets de $\I_2$. Dans le cas o\` u ils ne s'annulent pas avec des sommets de $\I_2$, cela contredit $e_\Z(X_x) \leq d_0+1$: ils sont donc contraints \` a s'annuler avec des sommets de $\I_2$. De ce fait, pour chaque sommet $K g_i^{-1}$ de $\I_1$ avec $v_\beta(K g_i^{-1}) = c$, on utilise l'identit\' e 
\[ y(g_i,v_i) = z(g_i,v_i) - \sum_{j=1}^r y' \big( g k_j^{-1} \varphi_1^{-d_0} , p_U (k_j v_i ) \big) :\]
on d\' eplace les tels indices $i$ dans un ensemble $\I_3$ (qui indice les $z(g_i,v_i)$) et on rajoute ceux correspondants aux $y'\big( g k_j^{-1} \varphi_1^{-d_0} , p_U (k_j v_i ) \big)$ dans $\I_2$ (en renum\' erotant convenablement). En effectuant cela pour tous les $c \in ]a,a+d_0[ \cap \Z$, il ne reste plus dans $\I_1$ que des sommets v\' erifiant $v_\beta(K g_i^{-1}) = a$. \\
Pour les sommets $K g_i^{-1}$ pour $i \in \I_2$, de deux choses l'une: ou bien ils font partie du support $X_x$, et comme $X_x$ est d'\' etendue $d_0+1$ au plus, on a alors $v_\beta(K g_i^{-1}) = a+d_0$. Ou bien ils s'annulent avec un sommet indic\' e dans $\I_1$, mais alors ils v\' erifient aussi $v_\beta(Kg_i^{-1}) = a+d_0$ (puisque ceux de $\I_1$ v\' erifient $v_\beta(Kg_i^{-1}) = a$). Dans les deux cas, on a bien l'\' ecriture voulue. \hfill$\Box$

\begin{prop} \label{mainsupport}
Soient $V$ une représentation irréductible de $K$ et $d_0 \geq 1$ avec $\widetilde{Z}_A \big( V_{U \cap K} \big)$ de la forme $\left( \begin{array}{cc} \varpi^{d_0 \Z} & 0 \\ 0 & \varpi^{d_0 \Z} \end{array} \right)$ et $\chi$ un caractère de $\Hh(G,K,V)$ se factorisant par ${}' \Sss_G$. Soit $x$ un élément non nul de $\Gamma(V,\chi)$.
\begin{itemize}
\item[(i)] Le support $X_x$ de $x$ vérifie $e_\Z(X_x) \geq d_0+1$ ou $\delta_\T(X_x) \geq 2 d_0$.
\item[(ii)] Si $X_x$ satisfait $e_\Z(X_x) = d_0 +1$ et $\delta_\T(X_x) = 2 d_0$, alors il existe $g,g' \in G$ et $v,v' \in V$ (avec éventuellement l'un d'eux nul) tels que l'on ait $x = y(g,v) + y'(g',v')$. Si de plus $v$ et $v'$ sont non nuls, alors on a $v_\beta(K (g')^{-1}) = v_\beta(K g^{-1})+d_0$. 
\end{itemize}
\end{prop}
\textsf{Preuve :} \\
Parce que $x$ est non nul, on a $e_\Z(X_x) \geq 1$. Supposons $e_\Z(X_x) \leq d_0$ et montrons $\delta_\T(X_x) \geq 2 d_0$. Par le lemme \ref{support2}, on prend $a \in \Z$, des ensembles finis $\I_1 , \I_2, \I_3$, des $v_i \in V$ tous non nuls et des $g_i \in G$ avec $v_\beta(K g_i^{-1})= a$ pour tout $i \in \I_1$, $v_\beta(K g_i^{-1})= a + d_0$ pour tout $i \in \I_2$ et $v_\beta(K g_i^{-1}) \in ]a,a+d_0[$ pour tout $i \in \I_3$ pour écrire 
\begin{equation} \label{decompx} x = \sum\nolimits_{i \in \I_1} y(g_i,v_i) + \sum\nolimits_{i \in \I_2} y'(g_i,v_i) + \sum\nolimits_{i \in \I_3} z(g_i,v_i) .\end{equation}
Dans une telle \' ecriture, on suppose le cardinal de $\I_3$ minimal. Si $\I_3$ est non vide, on va voir $\delta_\T(X_x) \geq  4d_0$. En effet, dans ce cas-l\` a, il existe un $c \in ]a,a+d_0[$ tel que $\I_3$ poss\` ede un indice $i$ correspondant au sommet $K g_i^{-1}$ avec $v_\beta(K g_i^{-1}) = c$. On prend une nouvelle fois le diam\` etre des centres sur $\I_3$:
\[ \delta_3 = \max\limits_{v_\beta(K g_i^{-1}) = v_\beta(K g_j^{-1}) = c} \, d_\T(\beta(K g_i^{-1}), \beta(K g_j^{-1})) .\]
Deux $K g_i^{-1}$ et $K g_j^{-1}$ réalisant ce diamètre $\delta_3$ possèdent (voir la preuve du lemme \ref{support1demi} à propos de $\Ker \, p_U$, \` a appliquer deux fois) des éléments à distance $\delta_3 + 4d_0$ dans le support. Comme ces derniers sont extr\' emaux, ils ne peuvent pas \^ etre annul\' es par d'autres de $\I_3$ et on a $\delta_\T(X_x) \geq 4d_0$ comme annonc\' e. \\
Supposons donc \` a pr\' esent $\I_3 = \varnothing$ dans l'\' ecriture (\ref{decompx}). On considère le diamètre des centres \[ \delta_{12} = \max\limits_{i,j \in \I_1 \cup \I_2} d_\T(\beta(K g_i^{-1}),\beta(K g_j^{-1})) .\]
Deux $K g_i^{-1}$ et $K g_j^{-1}$ réalisant ce diamètre possèdent des éléments dans le support de $y(g_i,v_i)$ ou $y'(g_i,v_i)$ (selon $i \in \I_1$ ou $i \in \I_2$) et $y(g_j,v_j)$ ou $y'(g_j,v_j)$ qui subsistent dans $X_x$ par maximalité. Si $i$ et $j$ sont tous deux dans le même ensemble d'indices $\I_1$ ou $\I_2$, alors on a $\delta_\T(X_x) = \delta_{12} + 2 d_0 \geq 2 d_0$. Cela fournit déjà \textit{(i)} et le cas d'égalité implique alors $\delta_{12} = 0$. Ceci nous assure $K g_i^{-1} = K g_j^{-1}$ grâce à $v_\beta(K g_i^{-1}) = v_\beta(K g_j^{-1})$ et au lemme \ref{foret}. La preuve est alors terminée dans ce cas-là, puisque pour la même raison $\I_2$ est vide ou un singleton. \\
Supposons maintenant $i \in \I_1$ et $j \in \I_2$ (le cas inverse étant identique). On prend $k \in \I_1$ tel que $K g_k^{-1}$ est à distance maximale \[ \delta_1(i) = \max\limits_{k \in \I_1} d_\T(\beta(K g_i^{-1}),\beta(K g_k^{-1})) \] de $K g_i^{-1}$. Ou bien, deux\footnote{il faut défausser celui situé sur le chemin entre $K g_i^{-1}$ et $K g_k^{-1}$} au moins des sommets du support de $T[g_k,v_k]$ sont dans $X_x$, ou bien tous sauf un au plus sont annulés par des sommets indicés par $\I_2$. Dans le premier cas, on a ce que l'on voulait: $\delta_\T(X_x)$ est supérieur ou égal à $\delta_1(i) + 2 d_0$, et le cas d'égalité donne $K g_i^{-1} = K g_k^{-1}$. Dans le second cas, les sommets de $T[g_k,v_k]$ sont annulés par des $[g_a,v_a]$ pour $a \in \I_2$ de support des sommets à distance $2 d_0$ les uns des autres. Parce que $K g_k^{-1}$ est extrémal (au sens qu'il réalise $\delta_1(i)$), les sommets des $Z^{-1}T[g_a,v_a]$ ne peuvent pas tous être annulés par des $[g_c,v_c]$ avec $c \in \I_1$, et on a alors $\delta_\T(X_x) \geq 4 d_0$. En particulier, il n'y a pas d'égalité $\delta_\T(X_x) = 2d_0$ dans ce cas-là. \\
On fait de même pour $j \in \I_2$ et cela termine la preuve. \hfill$\Box$

\chapter{Représentations de Steinberg modulo $p$ pour un groupe réductif sur un corps local} \label{Ly11}

\section{Introduction}

Au début des années 1950, Steinberg introduit dans \cite{Ste51} de nouvelles représentations (à coefficients complexes) pour le groupe général linéaire sur un corps fini. Quelques années plus tard, Curtis donne une formule agréable pour leur caractère (\cite{Cur66}). C'est dans l'esprit de cette dernière que sont aujourd'hui définies les représentations de Steinberg généralisées. \\
Dans \cite{Grk09}, Grosse-Klönne établit leur admissibilité et leur irréductibilité lorsque $G$ est un groupe classique déployé sur $F$ de caractéristique nulle et $R$ est un corps de caractéristique $p >0$. Ensuite, dans \cite{Her11b}, Herzig utilise les résultats préliminaires de \cite{Grk09} et sa machinerie propre pour étendre ces propriétés à tout groupe réductif déployé $G$, toujours sur $F$ de caractéristique nulle (mais avec $R$ algébriquement clos). \\
On notera de fait que le résultat principal de \cite{Her11b} pour $G$ le groupe général linéaire déployé met en emphase l'importance des représentations de Steinberg généralisées puisqu'avec les représentations dites supersingulières, elles représentent les \og briques fondamentales \fg\ pour construire toutes les représentations lisses admissibles irréductibles modulo $p$ de $G$. \\

L'objet de cette note est d'étendre les résultats de \cite{Grk09} et de \cite{Her11b} sur les Steinberg généralisées pour $G$ un groupe réductif quelconque. On développe ainsi l'analogue des paragraphes 2 et 3 de \cite{Grk09} et du paragraphe 7 de \cite{Her11b} dans le cas non déployé, mais le lecteur aura conscience que les arguments ne sont pas nouveaux.

\section{Contexte général} \label{contextgene}

Soient $p$ un nombre premier et $\Fpbar$ une clôture algébrique fixée de $\Fp$; tout corps fini de caractéristique $p$ sera vu comme un sous-corps de $\Fpbar$. Toute représentation considérée ici sera lisse. \\

Soit $F$ un corps local non archimédien localement compact à corps résiduel fini $k_F$ de caractéristique $p$. Soient $F^\sep$ une clôture séparable de $F$ et $F^\un \subseteq F^\sep$ la sous-extension maximale non ramifiée de $F$. On note $\I = \Gal(F^\sep/F^\un)$ le sous-groupe d'inertie du groupe de Galois absolu de $F$ et $\sigma \in \Gal(F^\un/F)$ le générateur topologique correspondant à un Frobenius arithmétique. \\

Soit $G$ un groupe réductif connexe sur $F$. Au paragraphe 7 de \cite{Kot97}, Kottwitz définit un morphisme fonctoriel et surjectif \[ \kappa_G: G(F^\un) \to X^* \big( Z(\widehat{G})^{\I} \big) ,\] où $\widehat{G}$ désigne le dual de Langlands connexe de $G$ et $Z(\widehat{G})$ son centre. On note $\B$ l'immeuble de Bruhat-Tits du groupe adjoint $G_{\mathrm{ad}}(F^\un)$. Un sous-groupe parahorique de $G$ est un groupe de la forme\footnote{on a noté $\Fix \, \Ff$ le fixateur point par point des simplexes de dimension $0$ composant $\Ff$} \[ \ker \, \kappa_G \cap G(F) \cap \Fix \, \Ff \] pour une facette $\sigma$-invariante $\Ff$ de $\B$ (voir 5.2.6 de \cite{BruTit84}, paragraphe 3.3 de \cite{HenVig11} ou paragraphe 8 de \cite{Hai09}). Si $\Ff$ est une chambre, on parle de \textit{sous-groupe d'Iwahori}. \\
Si $H$ est un sous-groupe parahorique de $G$ associé à une facette $\sigma$-invariante $\Ff$ de $\B$, on lui associe un groupe $H(F) \leq \widetilde{H} \leq G(F)$ défini par: \[ \widetilde{H} := \{ g \in G(F) \cap \Fix \, \Ff \ | \ \kappa_G(g) \textrm{ est de torsion} \} .\]

Soient $K$ un sous-groupe parahorique maximal spécial de $G(F)$ et $K(1)$ le pro-$p$-radical de $K$ (voir paragraphe 3.6 de \cite{HenVig11}). Le quotient $K/K(1)$ est le groupe des $k_F$-points d'un groupe réductif $\overline{G}$. \\
Soient $T$ un tore maximal parmi les tores $F$-déployés de $G$ et $A$ le centralisateur de $T$ dans $G$. Soient $B$ un parabolique minimal de $G$ de composante de Levi $A$ et $U$ son radical unipotent. Lorsque $Q$ est un parabolique contenant $B$, on notera $Q^-$ son parabolique opposé au sens du Théorème 4.15 de \cite{BorTit65}. \\
Pour tout sous-groupe $H$ de $G$ défini sur $F$, on note \[ \overline{H} := (H(F) \cap K)/(H(F) \cap K(1)) \] le sous-groupe de $\overline{G}(k_F)$ correspondant. \\
On confondra par abus $G$ et ses sous-groupes paraboliques (resp. composante de Levi, radical unipotent de sous-groupes paraboliques) avec leurs $F$-points. De même pour $\overline{G}$ et ses sous-groupes paraboliques (resp. composante de Levi, radical unipotent de sous-groupes paraboliques) avec leurs $k_F$-points. 

\section{Définitions et résultats}

Soit $R$ un anneau commutatif unitaire. Soit $Q$ un sous-groupe parabolique standard (c'est-à-dire qui contient\footnote{cette hypothèse n'est en fait utile que lorsque l'on veut utiliser la comparaison parabolique-compacte de \cite{HenVig11b}} $B$) de $G$. On définit la $G$-représentation à coefficients dans $R$ suivante : \[ \St_Q R := \frac{\Ind^{G}_{Q} \, \id}{\sum_{Q' \gneq Q} \Ind^{G}_{Q'} \, \id} .\] 
Ici, commme dans toute la suite, on a noté $\Ind$ le foncteur d'induction lisse et on fait agir $G$ sur $\Ind_Q^G \, \id$ par translation à droite.

\begin{theo} \label{Stirred}
Soit $R$ un corps de caractéristique $p$. La représentation de Steinberg généralisée $\St_Q R$ est irréductible et admissible.
\end{theo}
\textsf{Remarque :} 
Lorsque $F$ est de caractéristique $0$, l'admissibilité suit automatiquement de \cite{Vig07}. Par contre, lorsque $F$ est de caractéristique $p$, à ma connaissance on ne sait pas se passer de l'argument de ce papier. \\

On va présenter une preuve de cet énoncé dans les paragraphes qui suivent. Commençons par énoncer un corollaire (on dira un mot de la preuve dans le paragraphe \ref{Stfiltre}).

\begin{cor} \label{Stgeneral}
Les constituants de Jordan-Hölder de $\Ind^{G}_{Q} \, \id$ sont les $\St_{Q'} R$ pour les sous-groupes paraboliques\footnote{comme $Q'$ contient $Q$ standard, il l'est automatiquement aussi} $Q'$ contenant $Q$. Ils sont deux à deux non isomorphes et de multiplicité $1$.
\end{cor}

On donnera aussi la filtration par les cosocles de $\Ind^{G}_{Q} \, \id$. On note tout de suite le corollaire immédiat suivant.

\begin{cor} 
Soit $R$ un corps de caractéristique $p$. La représentation de Steinberg $\St_B R$ est irréductible et admissible.
\end{cor}

\section{Sous-groupe d'Iwahori et sous-groupes radiciels}

Soient $\Phi$ le système de racines de $G$ associé à $T$ et $\Phired$ le système réduit associé: \[ \Phired = \{ a \in \Phi \ | \ a/2 \notin \Phi \} .\] Le groupe parabolique minimal $B$ nous fournit un sous-ensemble de racines positives $\Phired^+ \subseteq \Phired$ et un système de racines simples $\Delta$. On note $\Phired^- := \Phired \smallsetminus \Phired^+$ et, pour $\alpha \in \Phired$, on appelle $s_\alpha$ la réflexion correspondante. On note $W$ le groupe de Weyl fini (déterminé par $T$) et $l : W \to \N$ la longueur (déterminée par $\Delta$). Soit $w_0$ l'élément le plus long de $W$. \\
Pour $w \in W$, on notera encore $w$ un relèvement (fixé une fois pour toutes) de $w$ dans le normalisateur $N_G(T) \cap K$ de $T$ dans $K$ (ce qui est possible car $K$ est spécial) ou bien l'image de ce relèvement dans $\overline{G}$. \\

On fait remarquer que le paragraphe 1 de \cite{Grk09} ne concerne que des groupes de réflexions cristallographiques avec système de racines réduit. Il est donc valable lorsque l'on travaille avec $\Phired$ et le lecteur ne devra pas être surpris quand on y fera référence. \\

Soit $I$ le sous-groupe d'Iwahori de $G$ suivant: si $x_0$ est le point spécial de l'immeuble de $G_{\mathrm{ad}}(F^{\un})$ fixe par $K$, et si $\Cc$ est la chambre de sommet $x_0$ et fixe par $B$, alors $I$ est le parahorique fixant $\Cc$ (voir paragraphe \ref{contextgene}). On dispose alors de la décomposition d'Iwasawa (voir \cite{BruTit72}, Proposition 7.3.1)\footnote{ici, comme dans tout ce qui suit, l'appel \` a \cite{BruTit72} n\' ecessite de faire attention que cela est bien loisible: c'est l'objet du Th\' eor\` eme 5.1.20 de \cite{BruTit84}, comme expliqu\' e dans son introduction. Par la suite, on gardera cet \' enonc\' e en t\^ ete sans le rappeler \` a chaque fois, mais on se permet d'insister que son importance est cruciale.}: \[ G = \coprod\nolimits_{w \in W} B w \It .\] Aussi, on a des injections naturelles 
\[ A \cap \Kt / A \cap K \hookrightarrow \It/I \hookrightarrow \Kt/K .\]
La composée est un isomorphisme par \cite{HenVig11}, Lemma 6.2.(iii); la première flèche est donc un isomorphisme \[ A \cap \Kt / A \cap K \xrightarrow{\sim} \It/I .\] On a donc finalement 
\begin{equation} \label{Iwasawa} G = \coprod\nolimits_{w \in W} B w (A \cap \Kt) I = \coprod\nolimits_{w \in W} B w I. \end{equation}
Pour tout sous-groupe $H$ de $G$, on pose $H^0 := H \cap I$. \\
Pour $\alpha \in \Phi$, on note $U_\alpha$ le sous-groupe radiciel associé. Comme on a l'inclusion $U_{2\alpha} \subsetneq U_\alpha$ si $\{\alpha,2\alpha\} \subseteq \Phi$, il convient de remarquer aussi $U^0_{2\alpha} \subseteq U^0_\alpha$. Par la Proposition 6.1.6 de \cite{BruTit72}, on a donc \[ \prod\nolimits_{\alpha \in \Phired^+} U_\alpha = U , \quad \prod\nolimits_{\alpha \in \Phired^+} U_\alpha^0 = U^0 ,\] quel que soit l'ordre choisi sur $\Phired^+$. \\

Soit $J$ un sous-ensemble de $\Delta$. Les $s_\alpha$ pour $\alpha \in J$ engendrent un sous-groupe $W_J$ de $W$. On note aussi \[ W^J := \{ w \in W \ | \ \forall \alpha \in J, \ l(w s_\alpha)> l(w) \} .\] On a, par \cite{Hum92}, Lemma 1.6 et Corollary 1.7, 
\begin{equation} \label{WJaltern} W^J = \{ w \in W \ | \ w(J) \subseteq \Phired^+ \} .\end{equation}
Aussi, gr\^ ace \` a \cite{Hum92}, Proposition 1.10 et paragraphe 5.12, on a le fait important suivant: l'ensemble $W^J$ est un système de représentants de $W/W_J$ contenant l'élément le plus court de chaque classe. \\
On note le sous-ensemble de $W^J$ constitué des éléments primitifs \[ W^J_{\pr} := W^J \smallsetminus \bigcup_{\alpha \in \Delta \smallsetminus J} W^{J \cup \{\alpha\}} = \{ w \in W^J \ | \ w(\Delta \smallsetminus J) \subseteq \Phired^- \} ,\] de sorte que l'on a $W = \coprod W^J_{\pr}$, lorsque $J$ parcourt les sous-ensembles de $\Delta$. \\
On définit aussi le sous-ensemble suivant de $\Phired$ : \[ W_J. J := \{ w \alpha \ | \ w \in W_J, \, \alpha \in J\} .\]

\section{Détermination de $(\St_Q R)^I$}

Soit $R$ un anneau commutatif unitaire. \\ 
Soit $J$ un sous-ensemble propre de $\Delta$. Pour $w \in W^{J \cup \{ \alpha\}}$ et $\alpha \in \Delta \smallsetminus J$, on définit \[ \partial(w) := \sum_{w' \in W^J , \, w' W_J \subseteq w W_{J \cup \{ \alpha \}}} w' \in R \big[W^J \big] ,\] où $R[W^J]$ désigne le $R$-module libre de base les éléments de $W^J$. En prolongeant par $R$-linéarité, on a la suite exacte 
\begin{equation} \label{MJdef} \bigoplus_{\alpha \in \Delta \smallsetminus J} R \big[W^{J \cup \{\alpha\}} \big] \xrightarrow{\partial} R \big[W^J \big] \xrightarrow{\nabla} \M_J(R) \to 0 , \end{equation} 
qui définit l'application linéaire $\nabla$ et le $R$-module $\M_J(R)$. Ce dernier module est un objet essentiel pour la compréhension du module des $I$-invariants de $\St_J R$. Et Grosse-Klönne a démontré (\cite{Grk09}, Proposition 1.3.(a)) que $\M_J(R)$ est libre de rang $|W^J_{\pr}|$. \\
Par la Proposition 2.4 de \cite{BusHen06} appliqué\footnote{dans cette référence toutes les représentations sont à coefficients complexes mais cela n'influe pas sur la preuve de la proposition en question} à $H = \{ 1 \}$ et au groupe profini $I$, le foncteur des fonctions localement constantes $C^\infty(I, \cdot)$ est exact. En appliquant $C^\infty(I, \cdot)$ à (\ref{MJdef}), on obtient la suite exacte \[ C^\infty \Big( I, \bigoplus_{\alpha \in \Delta \smallsetminus J} R \big[ W^{J \cup \{\alpha\}} \big] \Big) \to C^\infty \big( I,R \big[ W^J \big] \big) \to C^\infty(I,\M_J(R)) \to 0 .\] Par abus, on note encore ces flèches $\partial$ et $\nabla$ respectivement. \\
Notons, pour $J \subseteq \Delta$ et $w \in W$ : \begin{equation} \label{PJdef} P_J := B W_J B , \quad {}^w P_J := w P_J w^{-1} .\end{equation} On remarque que ${}^w P_J$ ne dépend que de la classe de $w$ dans $W/W_J$. Par définition de $P_J$ son radical unipotent est égal à \[ N_J = \prod_{\alpha \in \Phired^+ \smallsetminus ( \Phired^+ \cap W_J. J)} U_\alpha .\] On fixe aussi un sous-groupe de Levi $M_J$ de $P_J$ contenant $A$ et $W_J$. Enonçons tout de suite une inclusion entre les ${}^w P_J$ qui nous sera utile dans très peu de temps.

\begin{lem} \label{inclusPJw}
Soit $\alpha \in \Delta \smallsetminus J$. Soient $w$ et $w'$ des éléments de $W$ vérifiant $w' W_J \subseteq w W_{J \cup \{\alpha\}}$. On a les inclusions\footnote{on omet les parenthèses pour alléger les notations, mais ${}^w P_J^0$ doit se lire $( {}^w P_J )^0$} \[ {}^{w'} P_J \subseteq {}^w P_{J \cup \{\alpha\}} ; \quad {}^{w'} P_J^0 \subseteq {}^w P_{J \cup \{\alpha\}}^0 .\]
\end{lem}
\textsf{Preuve :} \\
Par hypoth\` ese, il existe un \' el\' ement $\sigma \in W_{J \cup \{\alpha\}}$ v\' erifiant $w' = w \sigma$. La premi\` ere inclusion vient imm\' ediatement: 
\[ {}^{w'} P_J = w \sigma P_J \sigma^{-1} w^{-1} \subseteq {}^w P_{J \cup \{\alpha\}} .\]
La seconde suit par intersection avec $I$. \hfill$\Box$\\
 
Pour $\alpha \in \Delta \smallsetminus J$, en notant $C^\infty({}^w P^0_{J \cup\{\alpha\}} \backslash I, R)$ le sous-espace de $C^\infty(I,R)$ constitué des fonctions ${}^w P^0_{J \cup\{\alpha\}}$-invariantes à gauche, on a une injection \[ \bigoplus_{w \in W^{J \cup \{\alpha \}}} C^\infty \big( {}^w P^0_{J \cup \{\alpha\}} \backslash I, R \big) \hookrightarrow C^\infty \big( I , R \big[W^{J \cup \{\alpha\}} \big] \big) \] donnée par la flèche \[ (f_{\alpha,w})_w \mapsto \sum\nolimits_w f_{\alpha, w} \, w .\] On a de même une injection \[ \bigoplus\nolimits_{w \in W^J} C^\infty \big( {}^w P^0_J \backslash I, R \big) \hookrightarrow C^\infty \big(I,R[W^J] \big) .\] On verra dorénavant les injections précédentes comme des inclusions. On a alors le diagramme commutatif suivant: 
\begin{equation} \label{diagramE2} \xymatrix{ C^\infty \big( I, \bigoplus\nolimits_{\alpha \in \Delta \smallsetminus J} R \big[ W^{J \cup \{\alpha\}} \big] \big) \ar[r]^{\partial} & C^\infty \big( I, R \big[ W^J \big] \big) \ar[r]^{\nabla} & C^\infty \big( I, \M_J(R) \big) \ar[r] & 0 \\ \bigoplus\limits_{\alpha \in \Delta \smallsetminus J, \, w \in W^{J \cup \{\alpha\}}} C^\infty \big( {}^w P^0_{J \cup \{\alpha\}} \backslash I , R \big) \ar[r]^{\phantom{xxxxx} \partial} \ar@{^(->}[u] \ar[r] & \bigoplus\limits_{w \in W^J} C^\infty \big( {}^w P^0_J \backslash I, R \big) \ar[r]^{\nabla} \ar@{^(->}[u] & C^\infty \big( I, \M_J(R) \big) \ar@{=}[u] & } .\end{equation}
Vérifions que l'image de $\bigoplus\nolimits_{\alpha, w} C^\infty \big( {}^w P^0_{J \cup \{\alpha\}} \backslash I , R \big)$ par $\partial$ est bien incluse dans $\bigoplus\nolimits_{W^J} C^\infty \big( {}^w P^0_J \backslash I, R \big)$. Fixons pour cela $\alpha \in \Delta \smallsetminus J$. On a 
\begin{equation} \label{def2delta} \partial \Big( \sum_{w \in W^{J \cup \{\alpha\}}} f_{\alpha, w} \, w \Big) = \sum_{w' \in W^J} \Big( \sum_{w \in W^{J \cup \{\alpha\}}, \, w' W_J \subseteq w W_{J \cup \{\alpha\}}} f_{\alpha, w} \Big) w';\end{equation} 
il s'agit donc de voir que $\sum\limits_{w' W_J \subseteq w W_{J \cup \{\alpha\}}} f_{\alpha,w}$ est ${}^{w'} P^0_J$-invariante à gauche. En fait, chacun des termes de cette somme l'est. En effet, chaque $f_{\alpha,w}$ est ${}^w P_{J \cup \{\alpha\}}^0$-invariant à gauche; par le lemme \ref{inclusPJw}, il est aussi ${}^{w'}P_J^0$-invariant et $\partial$ envoie bien $\bigoplus\nolimits_{\alpha, w} C^\infty \big( {}^w P^0_{J \cup \{\alpha\}} \backslash I , R \big)$ dans $\bigoplus\nolimits_{W^J} C^\infty \big( {}^w P^0_J \backslash I, R \big)$. 

\begin{prop} \label{suiteE2}
La ligne du bas du diagramme (\ref{diagramE2}) est exacte.
\end{prop}

L'introduction de quelques objets est nécessaire avant d'aborder la preuve de cette proposition. Remarquons que, par \cite{Hum92}, Corollary 1.5, le sous-système $\Phi_J \subseteq \Phired$ engendré par $J$ vérifie \[ \Phi_J^- = \Phired^- \cap W_J . J .\] 
On dit que $D \subsetneq \Phired$ est $J$-\textit{quasi-parabolique} s'il est l'intersection de certains $w(\Phired^- \smallsetminus \Phi_J^-)$. \\
Enonçons maintenant un lemme/définition particulièrement utile.

\begin{lem} \label{U0Ddef}
Soit $D \subsetneq \Phired$ une partie $J$-quasi-parabolique. Le produit $\prod\nolimits_D U^0_\alpha$ est indépendant de l'ordre choisi sur $D$: il forme un sous-groupe de $G$ que l'on notera $U^0_D$.
\end{lem}
\textsf{Preuve :} \\
Comme $D$ est $J$-quasi-parabolique, on voit qu'il est suffisant de prouver que $\prod_{w(\Phired^- \smallsetminus \Phi_J^-)} U^0_\alpha$ est indépendant de l'ordre sur $w(\Phired^- \smallsetminus \Phi_J^-)$ pour un $w \in W$ tel que $D$ est contenu dans $w(\Phired^- \smallsetminus \Phi_J^-)$. Quitte à conjuguer par $w$, on suppose à présent $w=1$. \\ 
Il reste donc à voir la condition de commutateurs sur la partie $\Phired^- \smallsetminus \Phi_J^- \subseteq \Phired^-$ pour pouvoir appliquer la Proposition\footnote{on l'applique \` a $Y_a = U_a \cap I \supseteq Y_{2a} = U_{2a} \cap I$ et $X_a = U_a \cap I$} 6.1.6 de \cite{BruTit72}. Il s'agit de voir dans un premier temps que si $a$ et $b$ sont des éléments de $\Phired^- \smallsetminus \Phi_J^-$ alors on a 
\begin{equation} \label{relcommnonentier} (U_a,U_b) \subseteq \langle U_{na + mb} \ | \ na+mb \in \Phired^- \smallsetminus \Phi_J^- , \ m,n \in \N^* \rangle .\end{equation} 
Pour voir que ceci est automatique, écrivons \[ a = a_J' + a_J , \quad a_J' = \sum_{\alpha \in \Delta \smallsetminus J} n_\alpha^{(a)} \alpha , \ a_J = \sum_{\alpha \in J} n_\alpha^{(a)} \alpha \] avec des $n_\alpha^{(a)} \leq 0$; et de même $b = b_J' + b_J$. On a alors, pour $n, m > 0$ : \[ na + mb = (na_J' + m b_J') + (n a_J + m b_J) .\] Parce que l'on a $a_J' \neq 0$ et $b_J' \neq 0$, on a $n a_J' + m b_J' \neq 0$ et donc $na + mb$ n'est pas un élément de $\Phi_J^-$. En prenant les intersections avec $I$, parce que $I \cap \prod_{\alpha \in \Phired^-} U_\alpha$ est égal à $\prod_{\alpha \in \Phired^-} U_\alpha^0$ par la Proposition\footnote{cette dernière s'occupe de $U \cap \widetilde{I}$, mais en prenant l'intersection avec le morphisme de Kottwitz, on voit l'égalité $U \cap \widetilde{I} = U \cap I$} 5.2.32 de \cite{BruTit72}, (\ref{relcommnonentier}) devient 
\[ (U^0_a,U^0_b) \subseteq \langle U^0_{na + mb} \ | \ na+mb \in \Phired^- \smallsetminus \Phi_J^- , \ m,n \in \N^* \rangle .\] 
Le lemme en découle. $\hfill\Box$ \\

Pour $w \in W$, on pose \[ D_w = w(\Phired^- \smallsetminus \Phi_J^-) :\] le groupe $U^0_{D_w}$ est l'intersection de $I$ avec le radical unipotent de ${}^w P_J^-$ (défini en (\ref{PJdef})) (voir \cite{Dem11b}, Proposition 1.12). On remarque que $D_w$ ne dépend que de la classe de $w$ dans $W/W_J$, et que, pour tout $\alpha \in \Delta \smallsetminus J$, on a $w(\Phired^- \smallsetminus \Phi_{J \cup \{\alpha\}}^-) \subseteq w(\Phired^- \smallsetminus \Phi_J^-)$. \\
L'introduction des ensembles $J$-quasi-paraboliques s'explique alors par le fait que l'on va s'intéresser à l'intersection de tels $U^0_{D_w}$ pour différents $w \in W$ et à la décomposition d'Iwahori qui suit. 

\begin{lem} \label{homeoIwahori}
Soient $J \subseteq \Delta$ et $w \in W$. Le produit ${}^w P^0_J \times U^0_{D_w} \to I$ est un homéomorphisme.
\end{lem}
\textsf{Preuve :} \\
Par la Proposition 5.2.32 de \cite{BruTit72}, on a la d\' ecomposition d'Iwahori 
\[ \It = \Big( \prod\nolimits_{\beta \in w \Phired^+} U^0_\beta \Big) \big( A \cap \Kt \big) \Big( \prod\nolimits_{\beta \in w \Phired^-} U^0_\beta \Big) ,\] 
où les produits sur $w \Phired^+$ et $w \Phired^-$ sont indépendants de l'ordre par la Proposition 6.1.6 de \cite{BruTit72}. En l'intersectant avec le noyau du morphisme de Kottwitz $\kappa_G$, parce que l'on a $A^0 = A \cap K$, cela permet d'\' ecrire 
\[ I = \Big( \prod\nolimits_{\beta \in w \Phired^+} U^0_\beta \Big) A^0 \Big( \prod\nolimits_{\beta \in w \Phired^-} U^0_\beta \Big) =: I^+ A^0 I^- .\] 
La d\' ecomposition 
\begin{equation} \label{Iwa4} {}^w P^0_J = \Big( \prod\nolimits_{\beta \in w \Phired^+} U^0_\beta \Big) A^0 \Big( \prod\nolimits_{\beta \in w \Phi_J^-} U^0_\beta \Big) \end{equation}
suit par intersection avec ${}^w P_J$: $I^+$ et $A^0$ sont inclus dans ${}^w P_J$ et donc ${}^w P^0_J$ est égal à $I^+ A^0 (I^- \cap {}^w P_J)$. Il en résulte la bijection
\begin{equation} \label{U0alpha} U^0_{D_w} = {}^w U^0_{\Phired^- \smallsetminus \Phi_J^-} \xrightarrow{\sim} {}^w P^0_J \backslash I .\end{equation}
On en déduit que le produit ${}^w P^0_J \times U^0_{D_w} \to I$ est un homéomorphisme: c'est une bijection par (\ref{U0alpha}), le produit est continu et une bijection continue d'un espace compact vers un espace séparé est un homéomorphisme. \hfill$\Box$\\

Pour un sous-ensemble $D$ de $\Phired$, on note 
\begin{equation} \label{defWJD} W^J(D) := \{ w \in W^J  \ | \ D \subseteq w(\Phired^- \smallsetminus \Phi_J^-) \} .\end{equation} 

\begin{lem} \label{gbiendef}
Soient $D,D' \subseteq \Phired$ deux ensembles $J$-quasi-paraboliques. Soient $\alpha \in \Delta \smallsetminus J$ et $w \in W^{J \cup \{\alpha\}}(D)$. On a l'égalité d'ensembles
\begin{equation} \label{gbiendef2} {}^w P^0_{J \cup \{\alpha\}} U^0_D \cap {}^w P^0_{J \cup \{\alpha\}} U^0_{D'} = {}^w P^0_{J \cup \{\alpha\}} \big( U^0_D \cap U^0_{D'} \big) .\end{equation}
\end{lem}
\textsf{Preuve :} \\
L'inclusion $\supseteq$ est évidente; prouvons $\subseteq$. En appliquant le lemme \ref{U0Ddef}, le produit $\prod\limits_{\beta \in D} U^0_\beta$ est indépendant de l'ordre. Choisissons un ordre sur $D$ tel que tout \' el\' ement de $D \smallsetminus (D \cap D')$ soit plac\' e avant tout \' el\' ement de $D \cap D'$. On forme le produit $\prod\limits_{\beta \in D \smallsetminus (D \cap D')} U^0_\beta$ suivant cet ordre et on le note $U^0_{D \smallsetminus D'}$, de sorte que l'on a $U^0_D = U^0_{D \smallsetminus D'} U^0_{D \cap D'}$ ($D \cap D'$ est $J$-quasi-parabolique et la notation est celle du lemme \ref{U0Ddef}). Il est équivalent de voir $\subseteq$ dans (\ref{gbiendef2}) et \[ {}^w P^0_{J \cup \{\alpha\}} U^0_{D \smallsetminus D'} \cap U^0_{D'} \subseteq {}^w P^0_{J \cup \{\alpha\}} \big( U^0_D \cap U^0_{D'} \big) .\] On va même montrer que ${}^w P^0_{J \cup \{\alpha\}} U^0_{D \smallsetminus D'} \cap U^0_{D'}$ est inclus dans ${}^w P^0_{J \cup \{\alpha\}}$. Notons 
\begin{equation} \label{3.9new} \Phi' = \Phired \smallsetminus w(\Phired^- \smallsetminus \Phi_{J \cup \{\alpha\}}^-) = \{ \beta \in \Phired \ | \ U_\beta \subseteq {}^w P^0_{J \cup \{\alpha\}} \} .\end{equation}
Puisque $w$ est dans $W^{J \cup \{\alpha\}}(D)$, l'intersection de $D$ et de $\Phi'$ est vide; ou autrement dit, la partie $(J \cup \{\alpha\})$-quasi-parabolique \[ D_{\alpha,w} = w(\Phired^- \smallsetminus \Phi_{J \cup \{\alpha\}}^-) \] contient $D$. Par le lemme \ref{U0Ddef} pour $D_{\alpha,w}$, le produit $\prod\limits_{\beta \in D_{\alpha,w}} U^0_\beta$ est ind\' ependant de l'ordre: on choisit un ordre sur $D_{\alpha,w}$ tel que sa restriction \` a $D$ co\" incide avec l'ordre pr\' ec\' edemment choisi sur $D$ et tel que tout \' el\' ement de $D$ pr\' ec\` ede tout \' el\' ement de $D_{\alpha,w} \smallsetminus D$. On forme le produit $\prod\limits_{\beta \in D_{\alpha,w} \smallsetminus D} U^0_\beta$ suivant cet ordre et on le note $U^0_{D_{\alpha,w} \smallsetminus D}$, de sorte que l'on a $U^0_{D_{\alpha,w}} = U^0_D U^0_{D_{\alpha,w} \smallsetminus D}$.
Par (\ref{U0alpha}), on a alors la décomposition en produit direct d'ensembles
\begin{equation} \label{Iwa1} I = {}^w P^0_{J \cup \{\alpha\}} U^0_{D \smallsetminus D'} U^0_{D \cap D'} U^0_{D_{\alpha,w} \smallsetminus D} .\end{equation}
Par le (\ref{Iwa4}) dans la preuve du lemme \ref{homeoIwahori}, on a la décomposition d'Iwahori 
\begin{equation} \label{Iwa2} {}^w P^0_{J \cup \{\alpha\}} = \Big( \prod_{\beta \in w \Phired^+} U^0_\beta \Big) A^0 \Big( \prod_{\beta \in w \Phi_{J \cup \{\alpha\}}^-} U^0_\beta \Big) .\end{equation}
Parce que $w \Phired^+$ vérifie la condition des commutateurs, on peut appliquer la Proposition\footnote{ici encore, on prend $Y_a = U_a \cap I \supseteq Y_{2a} = U_{2a} \cap I$ et $X_a = U_a \cap I$} 6.1.6 de \cite{BruTit72} et dire que $\prod\limits_{\beta \in w \Phired^+} U^0_\beta$ est un produit indépendant de l'ordre sur $w \Phired^+$, que l'on notera $U^0_{w \Phired^+}$. On choisit un ordre sur $w \Phired^+$ tel que tout élément de $w \Phired^+ \cap D'$ précède tout élément de $w \Phired^+ \smallsetminus (D' \cap w \Phired^+)$. On forme $U^0_{w \Phired^+ \cap D'} = \prod\limits_{\beta \in w \Phired^+ \cap D'} U^0_\beta$, $U^0_{w \Phired^+ \smallsetminus D'} = \prod\limits_{\beta \in w \Phired^+ \smallsetminus (D' \cap w \Phired^+)} U^0_\beta$, et on a $U^0_{w \Phired^+} = U^0_{w \Phired^+ \cap D'} U^0_{w \Phired^+ \smallsetminus D'}$. De la même manière, on choisit un ordre sur $w \Phi_{J \cup \{\alpha\}}^-$ tel que tout élément de $w \Phi_{J \cup \{\alpha\}}^-$ qui appartient à $D'$ précède tout élément qui n'y appartient pas: on obtient $U^0_{w \Phi_{J \cup \{\alpha\}}^-} = U^0_{w \Phi_{J \cup \{\alpha\}}^- \cap D'} U^0_{w \Phi_{J \cup \{\alpha\}}^- \smallsetminus D'}$. L'identité (\ref{Iwa2}) devient 
\[ {}^w P^0_{J \cup \{\alpha\}} = U^0_{w \Phired^+ \cap D'} U^0_{w \Phired^+ \smallsetminus D'} A^0 U^0_{w \Phi_{J \cup \{\alpha\}}^- \cap D'} U^0_{w \Phi_{J \cup \{\alpha\}}^- \smallsetminus D'} ;\]
et (\ref{Iwa1}) devient le produit direct 
\begin{equation} \label{Iwa3} I = U^0_{w \Phired^+ \cap D'} U^0_{w \Phired^+ \smallsetminus D'} A^0 U^0_{w \Phi_{J \cup \{\alpha\}}^- \cap D'} U^0_{w \Phi_{J \cup \{\alpha\}}^- \smallsetminus D'} U^0_{D \smallsetminus D'} U^0_{D \cap D'} U^0_{D_{\alpha,w} \smallsetminus D} .\end{equation}
Grâce à (\ref{Iwa3}), un élément de $U^0_{D'} \cap {}^w P^0_{J \cup \{\alpha\}} U^0_{D \smallsetminus D'}$ est dans le produit \[ U^0_{w \Phired^+ \cap D'} U^0_{w \Phi_{J \cup \{\alpha\}}^- \cap D'} \subseteq {}^w P^0_{J \cup \{\alpha\}} \] et le lemme est prouvé. \hfill$\Box$ \\ 

\textsf{Preuve de la proposition \ref{suiteE2}:} \\
On commence par numéroter tous les sous-ensembles $J$-quasi-paraboliques de $\Phired$: $D_0, D_1, D_2, \dots$ par ordre croissant de taille, c'est-à-dire avec 
\[ n<m \Rightarrow |D_n| \leq |D_m| .\] 
Soit $f \in C^\infty \big( I, R[W^J] \big)$, image de $(f_w)_{w \in W^J} \in \bigoplus\limits_{w \in W^J} C^\infty \big( {}^w P^0_J \backslash I , R \big)$ par la flèche verticale du diagramme (\ref{diagramE2}), tel que l'on a $\nabla f = 0$, c'est-à-dire $\nabla (f(x)) = 0$ pour tout $x \in I$. On cherche $g \in C^\infty \big( I, \bigoplus\limits_{\alpha \in \Delta \smallsetminus J} R \big[W^{J \cup \{\alpha\}} \big] \big)$, image de $(g_{\alpha,w})_{\alpha,w} \in \bigoplus\limits_{\alpha \in \Delta \smallsetminus J, \, w \in W^{J \cup \{\alpha\}}} C^\infty \big( {}^w P^0_{J \cup \{\alpha\}} \backslash I , R \big)$, vérifiant $f = \partial g$. \\
On va montrer par récurrence l'existence d'une telle fonction $g$ satisfaisant $f = \partial g$ sur $\bigcup_{n \geq 0} U^0_{D_n}$. Ceci implique $f = \partial g$ car $f$ et $\partial g$ proviennent de la ligne du bas du diagramme (\ref{diagramE2}) et que l'on a ${}^w P_J^0 \big( \bigcup_{n \geq 0} U^0_{D_n} \big) = I$ pour tout $w \in W^J$ par le lemme \ref{homeoIwahori}. La démonstration se fera en deux étapes:
\begin{itemize}
\item il existe $g$ telle que $f$ est égale à $\partial g$ sur $U^0_{D_0}$;
\item si $f$ est nulle sur $\bigcup_{n<m} U^0_{D_n}$, alors il existe $g$ vérifiant $f = \partial g$ sur $\bigcup_{n \leq m} U^0_{D_n}$ ($m \geq 1$).
\end{itemize}
La situation de l'étape d'initiation est la suivante: on a \[ D_0 = \varnothing = (\Phired^- \smallsetminus \Phi_J^-) \cap w_0(\Phired^- \smallsetminus \Phi_J^-) , \quad U^0_{D_0} = \{1\} .\] La fonction $f$ vérifie $(\nabla f)(1) = \nabla(f(1)) = 0$. Comme la suite (\ref{MJdef}) est exacte, on choisit, une famille d'éléments $(g^{(1)}_{\alpha,w})_{\alpha,w} \in \bigoplus_{\alpha \in \Delta \smallsetminus J} R \big[W^{J \cup \{ \alpha\}} \big]$ telle que l'on a $\partial ((g^{(1)}_{\alpha,w})_{\alpha,w}) = f(1)$. Soit, pour tout $\alpha \in \Delta \smallsetminus J$ et tout $w \in W^{J \cup \{\alpha\}}$, $g_{\alpha,w}$ une fonction de $C^\infty \big( {}^w P^0_{J \cup \{\alpha\}} \backslash I , R \big)$ valant $g^{(1)}_{\alpha,w}$ sur ${}^w P^0_{J \cup \{\alpha\}}$. Alors l'image $g$ de  
\[ (g_{\alpha,w})_{\alpha,w} \in \bigoplus\limits_{\alpha \in \Delta \smallsetminus J, w \in W^{J \cup \{\alpha\}}} C^\infty \big( {}^w P^0_{J \cup \{\alpha\}} \backslash I, R \big) \]
dans $C^\infty \big( I, \bigoplus\nolimits_{\alpha \in \Delta \smallsetminus J} R \big[ W^{J \cup \{\alpha\}} \big] \big)$ vérifie $\partial g = f$ sur $U^0_{D_0} = \{1 \}$. L'étape d'initiation est terminée. \\ 
Montrons maintenant la propriété de propagation de la récurrence. Soit, pour tout $w \in W^J$, $f_w \in C^\infty({}^w P^0_J \backslash I, R)$ nulle sur $\bigcup_{n<m} U^0_{D_n}$. \\
Si $w$ est un élément de $W^J \smallsetminus W^J(D_m)$, $f_w$ est nulle sur $U^0_{D_m}$ puisque l'on a ${}^w P^0_J U^0_{D_m} = {}^w P^0_J U^0_{D_n}$ pour $n<m$ vérifiant 
\[ D_m \cap w(\Phired^- \smallsetminus \Phi_J^-) = D_n .\]
En effet, il suffit de remarquer qu'un tel $n<m$ existe par définition de $W^J(D_m)$ (voir (\ref{defWJD})) et que l'on a alors $U^0_{D_m} = (U^0_{D_m} \cap {}^w P_J^0) U^0_{D_n}$. \\
Nous allons trouver la fonction $g$ comme image de 
\[ (g_{\alpha,w})_{\alpha,w} \in \bigoplus_{\alpha \in \Delta \smallsetminus J, \, w \in W^{J \cup \{\alpha\}}(D_m)} C^\infty( {}^w P^0_{J \cup \{\alpha\}} \backslash I, R) \]
telle que $(f - \partial g)_{w'}$ est nulle sur $\bigcup_{n \leq m} U^0_{D_n}$ pour tout $w' \in W^J(D_m)$. Comme $(\partial g)_{w'}$ est nulle\footnote{en (\ref{def2delta}), la somme sur les $w$ vérifiant $w' W_J \subseteq w W_{J \cup \{\alpha\}}$ ne fait intervenir que des $w \notin W^{J \cup \{\alpha\}}(D_m)$ car $w' W_J \subseteq w W_{J \cup \{\alpha\}}$ implique ${}^{w'} P_J \subseteq {}^w P_{J \cup \{\alpha\}}$ par le lemme \ref{inclusPJw}, et donc $w(\Phired^- \smallsetminus \Phi^-_{J \cup \{\alpha\}}) \subseteq w'(\Phired^- \smallsetminus \Phi_J^-)$ par (\ref{3.9new}) (voir définition (\ref{defWJD}) aussi)} pour $w' \in W^J \smallsetminus W^J(D_m)$, $(f- \partial g)_{w'}$ sera nulle sur $\bigcup_{n \leq m} U^0_{D_n}$ pour tout $w' \in W^J$. \\
La fonction $f$ est localement constante sur le compact $\bigcup_{n \leq m} U^0_{D_n}$ et nulle sur $\bigcup_{n<m} U^0_{D_n}$. Soient $(C_i)_{0 \leq i \leq r}$ les ouverts disjoints de $\bigcup_{n \leq m}U^0_{D_n}$ vérifiant $\bigcup_i C_i = \bigcup_{n \leq m} U^0_{D_n}$ et tels que $f_w$ est constant sur $C_i$, égal à $f_w^{(i)}$, et $f_w^{(0)} = 0$ pour tout $w \in W^J(D_m)$ et $(f_w^{(i)})_{w \in W^J(D_m)} \neq (f_w^{(i')})_{w \in W^J(D_m)}$ si $i \neq i'$. En particulier, on notera que $C_0$ contient $\bigcup_{n<m} U^0_{D_n}$ et que $U^0_{D_m}$ est l'union disjointe de $C_0 \cap U^0_{D_m}$ et des $C_i$ pour $1 \leq i \leq r$. \\
La suite extraite de (\ref{MJdef})
\begin{equation} \label{MJdefnew} \bigoplus_{\alpha \in \Delta \smallsetminus J} R \big[ W^{J \cup \{\alpha\}}(D_m)\big] \xrightarrow{\partial} R \big[ W^J(D_m) \big] \xrightarrow{\nabla} \M_J(R) \end{equation}
est encore exacte (voir \cite{Grk09}, Proposition 1.3.(b)). Par ailleurs, $(f_w^{(i)})_{w \in W^J(D_m)}$ appartient au noyau de $\nabla$ dans (\ref{MJdefnew}) car on a $\nabla f = 0$ et $f_w = 0$ pour $w \in W^J \smallsetminus W^J(D_m)$. Il existe alors, pour tout $0 \leq i \leq r$,
\[ g^{(i)} = \big( g^{(i)}_{\alpha,w} \big)_{\alpha \in \Delta \smallsetminus J,w \in W^{J \cup \{\alpha\}}(D_m)} \in \bigoplus_{\alpha \in \Delta \smallsetminus J} R \big[ W^{J \cup \{\alpha\}}(D_m) \big] \]
tel que l'on a $g^{(0)} = 0$ et
\begin{equation} \label{deltagi} \partial (g^{(i)}) = (f_w^{(i)})_{w \in W^J(D_m)} \in R \big[ W^J(D_m) \big] .\end{equation}
Posons $D_{\alpha,w} = w(\Phired^- \smallsetminus \Phi_{J \cup \{\alpha\}}^-)$ pour $\alpha \in \Delta \smallsetminus J$ et $w \in W$. Pour tout $w \in W$, le produit ${}^w P^0_{J \cup \{\alpha\}} \times U^0_{D_{\alpha,w}} \to I$ est un homéomorphisme par le lemme \ref{homeoIwahori}, donc on a 
\[ C^\infty \big( {}^w P^0_{J \cup \{\alpha\}} \backslash I, R \big) \simeq C^\infty(U^0_{D_{\alpha,w}},R) .\]
Soient $\alpha \in \Delta \smallsetminus J$ et $w \in W^{J \cup \{\alpha\}}(D_m)$. Notons $U'$ la projection\footnote{remarquons que $U'$ dépend de $\alpha$ et de $w$} de $\bigcup_{n<m} U^0_{D_n}$ sur $U^0_{D_{\alpha,w}}$: c'est le sous-espace de $U^0_{D_{\alpha,w}}$ vérifiant 
\[ {}^w P^0_{J \cup \{\alpha\}} U' = {}^w P^0_{J \cup \{\alpha\}} \big( \bigcup\nolimits_{n<m} U^0_{D_n} \big) .\]
Supposons qu'il existe, pour $\alpha \in \Delta \smallsetminus J$ et $w \in W^{J \cup \{\alpha\}}(D_m)$, une fonction $g_{\alpha,w} \in C^\infty(U^0_{D_{\alpha,w}},R)$ nulle sur $U' \cup (C_0 \cap U^0_{D_m})$ et constante, égale à $g_{\alpha,w}^{(i)}$, sur $C_i$ pour $1 \leq i \leq r$. Soit alors $g \in C^\infty \big( I, \bigoplus\limits_{\alpha \in \Delta \smallsetminus J} R \big[W^{J \cup \{\alpha\}} \big] \big)$ l'image de 
\[ (g_{\alpha,w})_{\alpha,w} \in \bigoplus_{\alpha \in \Delta \smallsetminus J, \, w \in W^{J \cup \{\alpha\}}(D_m)} C^\infty( {}^w P^0_{J \cup \{\alpha\}} \backslash I, R) .\]
Soit $w' \in W^{J}(D_m)$. Grâce à (\ref{deltagi}), on a alors 
\[ (\partial g)_{w'} = \sum_{\alpha \in \Delta \smallsetminus J} \sum_{w \in W^{J \cup \{\alpha\}}(D_m), \, w' W_J \subseteq w W_{J \cup \{\alpha\}} } g_{\alpha,w}^{(i)} = f^{(i)}_{w'} = f_{w'} \] 
sur chaque $C_i$ pour $1 \leq i \leq r$; la même égalité est vraie sur $C_0 \cap U^0_{D_m}$, de sorte que l'on a $(\partial g)_{w'} = f_{w'}$ sur tout $U^0_{D_m}$.
Pour $x \in \bigcup_{n \leq m} U^0_{D_n} \smallsetminus U^0_{D_m}$, pour tout $\alpha \in \Delta \smallsetminus J$ et tout $w \in W^{J\cup\{\alpha\}}(D_m)$, il existe $p_{\alpha,w} \in {}^w P^0_{J \cup \{\alpha\}}$ et $u_{\alpha,w} \in U'$ vérifiant $x = p_{\alpha,w} u_{\alpha,w}$. On a alors 
\[ (\partial g)_{w'}(x) = \sum_{\alpha \in \Delta \smallsetminus J} \sum_{w \in W^{J \cup \{\alpha\}}(D_m), \, w' W_J \subseteq w W_{J \cup \{\alpha\}} } g_{\alpha,w}(u_{\alpha,w}) = 0 .\]
Il nous reste à vérifier l'existence des fonctions $g_{\alpha,w}$. Soient $\alpha \in \Delta \smallsetminus J$ et $w \in W^{J \cup \{\alpha\}}(D_m)$. Par définition de $U'$, on a 
\[ {}^w P^0_{J \cup \{\alpha\}} U' \cap {}^w P^0_{J \cup \{\alpha\}} U^0_{D_m} = {}^w P^0_{J \cup \{\alpha\}} \big( \bigcup\nolimits_{n<m} U^0_{D_n} \big) \cap {}^w P^0_{J \cup \{\alpha\}} U^0_{D_m} ,\]
et par le lemme \ref{gbiendef}, on a 
\[ {}^w P^0_{J \cup \{\alpha\}} \big( \bigcup\nolimits_{n<m} U^0_{D_n} \big) \cap {}^w P^0_{J \cup \{\alpha\}} U^0_{D_m} = {}^w P^0_{J \cup \{\alpha\}} \big( \big( \bigcup\nolimits_{n<m} U^0_{D_n} \big) \cap U^0_{D_m} \big) .\]
On sait alors que $U' \cap U^0_{D_m}$ est inclus dans $C_0 \cap U^0_{D_m}$: comme $g^{(0)}_{\alpha,w}$ est nul, les conditions imposées sur les valeurs de $g_{\alpha,w}$ sont compatibles. 
L'espace $U' \cup U^0_{D_m}$ est une union disjointe de $U' \cup (C_0 \cap U^0_{D_m})$ et des $C_i$ pour $1 \leq i \leq r$. Les $C_i$, pour $1 \leq i \leq r$, sont ouverts compacts dans $U^0_{D_m}$ et $U'$ est compact. Le complémentaire de $U'$ dans $U' \cup U^0_{D_m}$ contient $C_i$, donc $C_i$ est ouvert dans $U' \cup U^0_{D_m}$ pour tout $1 \leq i \leq r$. Comme $\bigcup_{1 \leq i \leq r} C_i$ est compact, son complémentaire $U' \cup (C_0 \cap U^0_{D_m})$ dans $U' \cup U^0_{D_m}$ est ouvert. \\
Il existe donc une fonction $g'_{\alpha,w} \in C^\infty(U' \cup U^0_{D_m},R)$ nulle sur $U' \cup (C_0 \cap U^0_{D_m})$ et constante, égale à $g^{(i)}_{\alpha,w}$, sur $C_i$ pour $1 \leq i \leq r$. Comme $U' \cup U^0_{D_m}$ est fermé dans $U^0_{D_{\alpha,w}}$, le morphisme de restriction \[ C^\infty(U^0_{D_{\alpha,w}},R) \to C^\infty(U' \cup U^0_{D_m},R) \] est surjectif: il existe une fonction $g_{\alpha,w}$ sur $U^0_{D_{\alpha,w}}$ prolongeant $g'_{\alpha,w}$ comme voulue.
La récurrence est donc terminée et la proposition prouvée. $\hfill\Box$ \\

Pour $J \subseteq \Delta$, on définit la $G$-représentation de Steinberg généralisée $\St_J R$ par la suite exacte \begin{equation} \label{Stdef} \bigoplus_{\alpha \in \Delta \smallsetminus J} C^\infty(P_{J \cup \{\alpha\}} \backslash G,R) \xrightarrow{\partial} C^\infty(P_J \backslash G,R) \to \St_J R \to 0 . \end{equation} On observe que cela coïncide avec la définition de $\St_Q R$ lorsque $Q = P_J$. 

\begin{cor} \label{Stlibre}
Le $R$-module $\St_J R$ est libre. Et il existe une injection $I$-équivariante \[ \iota_R : \St_J R \hookrightarrow C^\infty(I,\M_J(R)) \] dont la formation commute aux changements de base.
\end{cor}
\textsf{Preuve :} \\
La flèche $i \mapsto P_J w^{-1} i$ induit une bijection ensembliste \[ {}^w P^0_J \backslash I \xrightarrow{\sim} P_J \backslash P_J w^{-1} I .\] De plus, toute inclusion $w' W_J \subseteq w W_{J \cup \{\alpha\}}$ donne alors lieu à un diagramme commutatif comme suit (par le lemme \ref{inclusPJw}): \[\xymatrix{ {}^w P^0_J \backslash I \ar[r] \ar[d]^\sim & {}^{w'} P^0_{J \cup \{\alpha\}} \backslash I \ar[d]^\sim \\ P_J \backslash P_J w^{-1} I \ar[r] & P_{J \cup \{\alpha\}} \backslash P_{J \cup \{\alpha\}} (w')^{-1} I } .\]
Les $w^{-1}$ pour $w \in W^J$ (resp. $w' \in W^{J \cup \{\alpha\}}$) forment un système de représentants de $W_J \backslash W$ (resp. $W_{J \cup \{\alpha\}} \backslash W$). On a les décompositions d'Iwasawa suivantes (conséquences directes de (\ref{Iwasawa}) et du Corollaire 4.2.2 de \cite{BruTit72}): \[ P_J \backslash G = \coprod_{w \in W^J} P_J \backslash P_J w^{-1} I , \quad P_{J \cup \{\alpha\}} \backslash G = \coprod_{w' \in W^{J \cup \{\alpha\}}} P_{J \cup \{\alpha\}} \backslash P_{J \cup \{\alpha\}} (w')^{-1} I .\] On en déduit les sommes directes \[ C^\infty( P_J \backslash G, R) = \bigoplus_{w \in W^J} C^\infty(P_J \backslash P_J w^{-1} I , R) ,\] \[ C^\infty(P_{J \cup \{\alpha\}} \backslash G,R) = \bigoplus_{w' \in W^{J \cup \{\alpha\}}} C^\infty(P_{J \cup \{\alpha\}} \backslash P_{J \cup \{\alpha\}} (w')^{-1} I ,R ) .\] Il en découle le diagramme commutatif suivant: 
\begin{equation} \label{commStI} \xymatrix{ \bigoplus\limits_{\alpha \in \Delta \smallsetminus J} C^\infty( P_{J \cup \{\alpha\}}\backslash G,R) \ar[r] \ar[d]^\sim & C^\infty(P_J \backslash G,R) \ar[r] \ar[d]^\sim & \St_J R \ar[r] \ar@{.>}[d] & 0 \\ \bigoplus\limits_{\alpha \in \Delta \smallsetminus J, \, w' \in W^{J \cup \{\alpha\}}} C^\infty( {}^{w'}P^0_{J \cup \{\alpha\}} \backslash I,R) \ar[r] & \bigoplus\limits_{w \in W^J} C^\infty( {}^w P^0_J \backslash I,R) \ar[r] & C^\infty(I,\M_J(R)) & } .\end{equation}
La première ligne est exacte par définition de $\St_J R$ (voir (\ref{Stdef})) et la seconde l'est par la proposition \ref{suiteE2}. On en déduit que l'application \[ \iota_R : \St_J R \hookrightarrow C^\infty(I,\M_J(R)) \] est injective; elle est aussi $I$-équivariante car toutes les flèches pleines du diagramme le sont. \\ Etudions d'abord le cas $R = \Z$: parce que l'on a \[ \M_J(R) = \M_J(\Z) \otimes_{Z} R ,\] on a la propriété de changement de base voulue. Grâce à la Proposition 1.3.(a) de \cite{Grk09} et au lemme \ref{Ccompact}, le $\Z$-module $C^\infty(I,\M_J(\Z))$ est libre. Or un sous-module d'un module libre sur un anneau principal est libre. Ainsi $\St_J \Z$ est libre, et $\St_J R$ est libre sur $R$ par changement de base\footnote{on a bien $\St_J \Z \otimes_{\Z} R = \St_J R$ comme la suite exacte \[ 0 \to \sum_{\alpha \in \Delta \smallsetminus J} C^\infty(P_{J \cup \{\alpha\}} \backslash G,\Z) \to C^\infty(P_J \backslash G, \Z) \to \St_J \Z \to 0 \] reste exacte après tensorisation par $R$, les deux $\Z$-modules à gauche étant libres par l'appendice \ref{Ccompactlibre}}. $\hfill\Box$ \\

Avant de déterminer à proprement parler $(\St_J R)^I$, on commence par déterminer $C^\infty(P_J \backslash G,R)^I$.

\begin{lem} \label{rgtotal}
Le $R$-module $C^\infty(P_J \backslash G,R)^I$ est libre de rang $|W^J|$.
\end{lem}
\textsf{Preuve :} \\
La décomposition d'Iwasawa (voir (\ref{Iwasawa})) fournit \[ G = \coprod\nolimits_{W^J} P_J w^{-1} I .\] Grâce à cette décomposition, on définit, pour $w \in W^J$, la fonction $f_w$ de $C^\infty(P_J \backslash G,R)^I$ par $f_w(w^{-1}) = 1$ et $f_w(w') = 0$ si $w' \in (W^J)^{-1} \smallsetminus \{w^{-1}\}$. Le $R$-module $C^\infty(P_J \backslash G, R)^I$ est alors libre et engendré par $(f_w)_{w \in W^J}$. $\hfill\Box$

\begin{prop} \label{IinvVQ}
L'espace des $I$-invariants $\big( \St_J R \big)^I$ est un $R$-module libre de rang $|W^J_{\pr}|$.
\end{prop}
\textsf{Preuve :} \\
Appliquons le foncteur des $I$-invariants au carr\' e commutatif du diagramme de la preuve du corollaire \ref{Stlibre}. On obtient les applications de $R$-modules : 
\begin{equation} \label{Iinvline} R[W^J] = \bigoplus_{w \in W^J} C^\infty \big( {}^w P^0_J \backslash I, R \big)^I \to \big( \St_J R \big)^I \to C^\infty(I,\M_J(R))^I = \M_J(R) .\end{equation}
Parce que $\iota_R$ est injective et que le foncteur des $I$-invariants est exact \` a gauche, la fl\` eche de droite de (\ref{Iinvline}) est injective. Enfin la compos\' ee (\ref{Iinvline}) est surjective par d\' efinition de $\M_J(R)$. Alors on a un isomorphisme 
\[ \big( \St_J R \big)^I \xrightarrow{\sim} \M_J(R) \] 
de $R$-modules, ce qui fait de $(\St_J R)^I$ un $R$-module libre de rang $|W^J_\pr|$ (par \cite{Grk09}, Proposition 1.2.(a)). \hfill$\Box$

\begin{cor} \label{Stbase}
Les images $\overline{f}_w$ dans $\St_J R$ des fonctions $f_w$ (pour $w \in W^J_{\pr}$) définies dans la preuve du lemme \ref{rgtotal} forment une base explicite de $\big( \St_J R \big)^I$.
\end{cor}
\textsf{Preuve :} \\
En prenant les $I$-invariants du diagramme commutatif dans la preuve du corollaire \ref{Stlibre}, on obtient le diagramme commutatif suivant, o\` u les carr\' es du haut sont commutatifs.
\[\hspace{-0.3cm} \xymatrix{ \bigoplus_{\alpha \in \Delta \smallsetminus J} C^\infty \big( P_{J \cup \{\alpha\}} \backslash G , R \big)^I \ar[r] \ar[d]^\sim & C^\infty \big( P_J \backslash G, R \big)^I \ar[r] \ar[d]^\sim & \big( \St_J R \big)^I \ar[d]^\sim& \\
\bigoplus_{\alpha \in \Delta \smallsetminus J, \, w' \in W^{J \cup \{\alpha\}}} C^\infty \big( {}^{w'} P^0_{J \cup \{\alpha\}} \backslash I, R \big)^I \ar[r]^{\phantom{xxxxxxx} \partial} & \bigoplus_{w \in W^J} C^\infty \big( {}^w P^0_J \backslash I, R \big)^I \ar[r] & C^\infty \big( I, \M_J(R) \big)^I & \\ 
\bigoplus_{\alpha \in \Delta \smallsetminus J, \, w' \in W^{J \cup \{\alpha\}}} R[W^{J \cup \{\alpha\}}] \ar[u]_\sim \ar[r]^{\phantom{xxxxxxxxx} \partial} & R[W^J] \ar[r] \ar[u]_\sim & \M_J(R) \ar[r] \ar[u]_\sim & 0 }\]
On remarque que l'on a utilis\' e la proposition \ref{IinvVQ} pour affirmer l'isomorphisme de $R$-modules libres que constitue la fl\` eche verticale en haut \` a droite. Les carr\' es du bas du diagramme sont commutatifs en vertu de ce que la description des fonctions du lemme \ref{rgtotal} est compatible avec la discussion d'avant la proposition \ref{suiteE2}. Le r\' esultat suit. \hfill$\Box$\\

Une autre conséquence importante est l'assertion de l'admissibilité dans le théorème \ref{Stirred}.

\begin{cor} \label{Stadm}
Supposons que $R$ soit un corps de caractéristique $p$. Alors $\St_J R$ est admissible.
\end{cor}
\textsf{Preuve :} \\
En regardant le diagramme commutatif (\ref{commStI}), on voit que $\St_J R$ est l'image de $\bigoplus_{w \in W^J} C^\infty \big( {}^w P^0_J \backslash I, R \big)$ dans $C^\infty(I,\M_J(R))$. Comme chaque ${}^w P^0_J$ contient $A \cap I$, cette image est en fait incluse dans $C^\infty \big(A \cap I \backslash I,\M_J(R) \big) \subseteq C^\infty(I,\M_J(R))$. Si $I(1)$ désigne le pro-$p$-Sylow de $I$, on a $(A \cap I) I(1) = I$; il suit les \' egalit\' es $C^\infty \big( A \cap I \backslash I, \M_J(R) \big)^I = C^\infty \big( A \cap I \backslash I, \M_J(R) \big)^{I(1)}$, et donc $\big( \St_J R \big)^I = \big( \St_J R \big)^{I(1)}$. \\
Traitons d'abord le cas $R = \Fpbar$: le $\Fpbar$-espace vectoriel $(\St_J \Fpbar)^{I(1)} = (\St_J \Fpbar)^I$ est de dimension finie par la proposition \ref{IinvVQ}. Comme $I(1)$ contient $K(1)$, il est ouvert et on peut appliquer \cite{Pas04}, Theorem 6.3.2, et $\St_J \Fpbar$ est admissible. \\
Dans le cas $R = \Fp$, on remarque que l'on a $(\St_J \Fpbar)^H = (\St_J \Fp)^H \otimes_{\Fp} \Fpbar$ pour tout sous-groupe ouvert $H$ de $G$. L'admissibilité de $\St_J \Fpbar$ implique alors celle de $\St_J \Fp$. Enfin, pour $R$ un corps de caractéristique $p$, $\Fp$ est naturellement un sous-corps de $R$. On a alors $(\St_J R)^H = (\St_J \Fp)^H \otimes_{\Fp} R$ pour tout $H \leq G$ ouvert. Le résultat en découle. $\hfill\Box$ \\
\textsf{Remarque :} 
Marie-France Vignéras nous fait remarquer que l'on peut aussi prouver ce fait en se servant du corollaire \ref{Stlibre}. Alors, pour tout sous-groupe compact ouvert $H$ de $I$, l'espace d'invariants $(\St_J R)^H$ s'injecte dans le $R$-espace vectoriel de dimension finie $C(I/H,\M_J(R))$. Par ailleurs, un tel argument reste tout à fait valable pour $R$ un anneau principal (sans hypothèse sur la caractéristique).

\section{Comparaison avec le cas fini}

On retourne à $R$ un anneau commutatif unitaire. On cherche à comparer les représentations de Steinberg généralisées avec leur analogue dans le cas fini. \\

On remarque d'abord que comme on a choisi $K$ comme étant le parahorique associé à un sommet spécial, le groupe de Weyl associé à l'adhérence schématique $\widetilde{T}$ de $\overline{T}$ est encore $W$ (voir \cite{BruTit84}, Propositions 4.4.5 et 4.6.4; \cite{HaiRap08}, Proposition 12). Le système de racines $\overline{\Phi}$ associé à la paire $(\overline{B},\widetilde{T})$ est un sous-système de $\Phi$ et il n'est pas nécessairement réduit; on considère alors son réduit $\Phibarred$. Parce que $\Phibarred$ et $\Phired$ sont tous deux des sous-systèmes réduits de $\Phi$ qui contiennent une $\R$-base de $X^*(T) \otimes_{\Z} \R$, on a une correspondance bijective \[\begin{array}{ccc} \Phired & \xrightarrow{\sim} & \Phibarred \\ a & \mapsto & a \textrm{ ou } 2a \end{array} .\] Cette bijection exhibe un sous-ensemble de racines simples $\overline{\Delta}$ de $\overline{\Phi}$ en correspondance avec $\Delta$ de la même manière: \[ \overline{\Delta} := \bigcup\nolimits_{a \in \Delta} \{a,2a\} \cap \Phibarred .\] Pour $J$ un sous-ensemble de $\Delta$, on associe alors de cette manière $\overline{J} \subseteq \overline{\Delta}$. \\
De même que pour le paragraphe précédent, le lecteur ne s'étonnera pas des invocations du paragraphe 1 de \cite{Grk09} pour des résultats concernant $\Phibarred$.

\begin{lem} \label{parabcompatibles}
Soit $J$ un sous-ensemble de $\Delta$. La réduction $\overline{P}_J$ de $P_J$ est $P_{\overline{J}}$, le parabolique de $\overline{G}$ associé à $\overline{J}$.
\end{lem}
\textsf{Remarque :} 
On peut aussi se reporter au Corollaire 4.6.4 de \cite{BruTit84}. \\
\textsf{Preuve :} \\
La réduction de $P_J$ est \[ \overline{P}_J = P_J \cap K / P_J \cap K(1) \] et est engendrée par les images dans $\overline{G}$ de $A \cap K$ et des $U_a \cap K$ pour $a \in \Phired^+ \cup W_J . J$. En effet, par d\' ecomposition de Bruhat, il suffit de voir que $B \cap K$ et $W_J \subseteq K$ sont engendr\' es par $A \cap K$ et les $U_a \cap K$ pour $a \in \Phired^+ \cup W_J.J$. Pour $B \cap K$, cela suit de $\cite{HenVig11}$, Theorem 6.5, qui affirme $B \cap K = (A \cap K)(U \cap K)$; et donc $B \cap K$ est engendr\' e par $A \cap K$ et les $U_a \cap K$ pour $a \in \Phired^+$. Soient $\alpha \in J$ et $u \in U_{-\alpha} \cap K$ d'image non nulle dans $\overline{U}_{-\alpha}$. Par \cite{Car85}, Corollary 2.6.2, il existe alors $\overline{b_1}, \overline{b}_2 \in \overline{B}$ tels que l'on ait $\overline{u} = \overline{b}_1 s \overline{b}_2$ o\` u $s$ est la r\' eflexion dans $W_J \subseteq K$ correspondant \` a $\alpha$. En relevant $\overline{b}_1$ et $\overline{b}_2$ respectivement en $b_1$ et $b_2$ dans $B \cap K$, on voit que $s$ est bien dans le groupe engendr\' e par $A \cap K$ et les $U_a \cap K$ pour $a \in \Phired^+ \cup W_J.J$. Les tels $s$ formant un syst\` eme de g\' en\' erateurs de $W_J$, on a bien l'affirmation voulue sur $P_J \cap K$. \\ 
Le groupe $P_{\overline{J}}$ est quant à lui engendré par $\overline{A} = Z_{\overline{G}}(\overline{T})$ et les $U_{\overline{a}}$ pour $\overline{a} \in \overline{\Phi}^+ \cup W_{\overline{J}}. \overline{J}$. Ce sont bien là les mêmes groupes radiciels par la remarque suivant le Lemma 6.12 de \cite{HenVig11}. $\hfill\Box$ \\

Pour $\overline{J} \subseteq \overline{\Delta}$ (correspondant à $J \subseteq \Delta$), on garde la même définition de la représentation de Steinberg généralisée $\overline{\St}_{\overline{J}} R$; on rappelle la suite exacte de $\overline{G}$-représentations qui la définit : \[ \bigoplus_{\alpha \in \Delta \smallsetminus J} C(\overline{P}_{J \cup \{\alpha\}} \backslash \overline{G},R) \to C(\overline{P}_J \backslash \overline{G},R) \to \overline{\St}_{\overline{J}} R \to 0 .\]
Commençons par exhiber des éléments particuliers de $\big( \overline{\St}_{\overline{J}} R \big)^{\overline{B}}$: pour $w$ un élément de $W^{\overline{J}}_{\pr}$, notons $\overline{g}_w \in \overline{\St}_{\overline{J}} R$ l'image de la fonction caractéristique $g_w$ de \[ \overline{P}_J w^{-1} \overline{B} = \overline{P}_J w^{-1} \big( \overline{U} \cap w \overline{U}^- w^{-1} \big) \] dans $\Ind_{\overline{P}_J}^{\overline{G}} \, \id$. \\
Soit $w \in W$. On va commencer par prouver l'affirmation \[\overline{P}_J w^{-1} \overline{B} = \overline{P}_J w^{-1} (\overline{U} \cap w \overline{U}^- w^{-1}) .\] On a d'abord \begin{equation} \label{BwU} \overline{B} w^{-1} \overline{B} = \overline{B} w^{-1} \overline{A} \overline{U} = \overline{B} w^{-1} \overline{U} .\end{equation} Ensuite, par la Proposition 6.1.6 de \cite{BruTit72}, on a \begin{equation}\label{sUeqn} \overline{U} = ( \overline{U} \cap w \overline{U} w^{-1})(\overline{U} \cap w \overline{U}^- w^{-1}) . \end{equation} Il s'ensuit \[ \overline{B} w^{-1} \overline{U} = \overline{B} w^{-1} \big( \overline{U} \cap w \overline{U}^- w^{-1} \big) .\] La même égalité subsiste alors avec $\overline{P}_J \geq \overline{B}$ à gauche plutôt que $\overline{B}$, et l'égalité voulue est prouvée. \\
Si on fait subir à $(\overline{G}, \Phibarred,\overline{J})$ le raisonnement du paragraphe précédent, on trouve alors que $\big( \overline{\St}_{\overline{J}} R \big)^{\overline{B}}$ est un $R$-module libre, de base $(\overline{g}_w)_{w \in W^{\overline{J}}_{\pr}}$. Cependant, savoir ce fait n'est pas nécessaire tout de suite, et la preuve de la proposition \ref{mainGrk1} nous le redonnera à moindres frais. \\

On a une décomposition d'Iwasawa $G = P_J K$ (voir \cite{HaiRos10}, Corollary 9.1.2). L'application $P_J \backslash G \to \overline{P}_J \backslash \overline{G}$ est continue (car $P_J K(1)$ est fermé dans $G$) et surjective. On a donc l'injection naturelle ($k$ désigne ici un représentant dans $K$): \[ \begin{array}{ccc} C(\overline{P}_J \backslash \overline{G},R) & \to & C^\infty( P_J \backslash G,R)\\ f & \mapsto & \big( P_J k \mapsto f(\overline{P}_J \overline{k}) \big) \end{array} .\] 
Dès lors, on a le diagramme suivant de $K$-représentations. \[\xymatrix{ \bigoplus_{ \alpha \in \Delta \smallsetminus J} C(\overline{P}_{J \cup \{\alpha\}} \backslash \overline{G},R) \ar[r]^{\phantom{xxxxx} \varphi_1} \ar@{^(->}[d] & C(\overline{P}_J \backslash \overline{G},R) \ar[r] \ar@{^(->}[d] & \overline{\St}_{\overline{J}} R \ar[r] \ar@{.>}[d]^{\iota} & 0 \\ \bigoplus_{\alpha \in \Delta \smallsetminus J} C^\infty(P_{J \cup \{\alpha\}} \backslash G,R) \ar[r]^{\phantom{xxxxx} \varphi_2} & C^\infty(P_J \backslash G,R) \ar[r] & \St_J R \ar[r] & 0 }\]
Comme les deux flèches verticales de gauche sont des injections et parce que l'on a l'égalité \[ \varphi_2 \bigg( \bigoplus_{\alpha \in \Delta \smallsetminus J} C^\infty(P_{J \cup \{\alpha\}} \backslash G,R) \bigg) \bigcap C(\overline{P}_J \backslash \overline{G},R) = \varphi_1 \bigg( \bigoplus_{ \alpha \in \Delta \smallsetminus J} C(\overline{P}_{J \cup \{\alpha\}} \backslash \overline{G},R) \bigg) ,\] la flèche de droite est aussi injective. Dans l'\' egalit\' e pr\' ec\' edante, l'inclusion $\supseteq$ est \' evidente. Pour obtenir l'inclusion inverse $\subseteq$, soit $f$ une fonction de 
\[ \varphi_2 \bigg( \bigoplus_{\alpha \in \Delta \smallsetminus J} C^\infty(P_{J \cup \{\alpha\}} \backslash G,R) \bigg) \bigcap C(\overline{P}_J \backslash \overline{G},R) ,\]
que l'on voit comme une fonction sur $K$. Elle s'écrit $f = \sum_{\alpha \in \Delta \smallsetminus J} \varphi_2(f_\alpha)$ avec $f_\alpha \in C^\infty(P_{J \cup \{\alpha\}} \backslash G,R)$ pour tout $\alpha \in \Delta \smallsetminus J$. Fixons, pour tout $\alpha \in \Delta \smallsetminus J$, $\K_\alpha$ un ensemble fini de représentants de $\overline{P}_{J \cup \{\alpha\}} \backslash \overline{G}$ dans $P_{J \cup \{\alpha\}} \backslash G$. On définit alors, pour tout $\alpha \in \Delta \smallsetminus J$, $f'_\alpha \in C(\overline{P}_{J \cup \{\alpha\}} \backslash \overline{G},R)$ par $f'_\alpha(\overline{k}) = f_\alpha(k)$ pour tout $k \in \K_\alpha$ d'image $\overline{k} \in \overline{P}_{J \cup \{\alpha\}} \backslash \overline{G}$. Parce que le diagramme ci-dessus est commutatif et que $f$ est dans $C(\overline{P}_J \backslash G, R)$, on obtient 
\[ \sum\nolimits_\alpha \varphi_1(f'_\alpha) = \sum\nolimits_\alpha \varphi_2(f_\alpha) = f \]
et l'inclusion $\subseteq$ voulue. 

\begin{prop} \label{mainGrk1} 
$\phantom{,}$
L'injection $K$-équivariante \[ \iota: \overline{\St}_{\overline{J}} R \hookrightarrow \St_J R \] induit l'isomorphisme de $R$-modules \[ \big(\overline{\St}_{\overline{J}} R \big)^{\overline{B}} \xrightarrow{\sim} \big( \St_J R \big)^I .\] 
\end{prop}
\textsf{Preuve :} \\
L'injection $\iota : \overline{\St}_{\overline{J}} R \hookrightarrow \St_J R$ étant $K$-équivariante, on en déduit l'injection \[ \big(\overline{\St}_{\overline{J}} R \big)^{\overline{B}} \hookrightarrow \big( \St_J R \big)^I ,\] que l'on notera encore $\iota$. Par la proposition \ref{IinvVQ}, le $R$-module $(\St_J R)^I$ est libre de rang $|W^J_{\pr}|$. \\ Commençons par considérer le cas $R = \Z$. 
En tant que sous-module de $(\St_J \Z)^I$, $(\overline{\St}_{\overline{J}} \Z)^{\overline{B}}$ est un $\Z$-module libre de rang majoré par $|W^{\overline{J}}_{\pr}| = |W^J_{\pr}|$ (on a $|W_J| = |W_{\overline{J}}|$ pour tout $J \subseteq \Delta$). Pour prouver l'isomorphisme voulu, il suffit d'exhiber une base de $(\overline{\St}_{\overline{J}} \Z)^{\overline{B}}$ qui s'envoie sur la base $(\overline{f}_w)_{w \in W^J_{\pr}}$ de $(\St_J \Z)^I$, où $\overline{f}_w$ est l'image de $f_w \in C^\infty(P_J \backslash G, \Z)$ (voir corollaire \ref{Stbase}) dans $\St_J \Z$. Il reste à remarquer que, par le lemme \ref{parabcompatibles}, la famille $(\overline{g}_w)_{w \in W^{\overline{J}}_{\pr}}$ s'envoie par $\iota$ sur $(\overline{f}_w)_{w \in W^J_{\pr}}$. \\
Revenons à $R$ général: on a la situation suivante
\[ \eta : (\overline{\St}_{\overline{J}} \Z)^{\overline{B}} \otimes_\Z R \hookrightarrow (\overline{\St}_{\overline{J}} R)^{\overline{B}} \mathop{\hookrightarrow}\limits^{\iota} \big( \St_J R \big)^I .\]
Parce que la famille $(\overline{g}_w \otimes 1)_{w \in W^{\overline{J}}_{\pr}}$ d'éléments de $(\overline{\St}_{\overline{J}} \Z)^{\overline{B}} \otimes_\Z R$ s'envoie par $\eta$ sur $(\overline{f}_w)_{w \in W^J_{\pr}}$, $\eta$ est un isomorphisme et\footnote{cela garantit au passage que la formation de $\big( \overline{\St}_{\overline{J}} \cdot \big)^{\overline{B}}$ commute au changement de base} $\iota$ induit l'isomorphisme $(\overline{\St}_{\overline{J}} R)^{\overline{B}} \xrightarrow{\sim} \big( \St_J R \big)^I$. $\hfill\Box$

\begin{cor}
La famille $(\overline{g}_w)_{w \in W^{\overline{J}}_{\pr}}$ forme une base du $R$-module libre $\big(\overline{\St}_{\overline{J}} R \big)^{\overline{B}}$.
\end{cor}

\section{Représentations de Steinberg généralisées dans le cas fini}

Dans toute cette section, $G$ désignera un groupe réductif fini: cela permettra d'alléger les notations et d'éviter de surligner un nombre déraisonnable de symboles. Aussi, quand on fera référence à l'\og axiomatique des systèmes de Tits \fg\ ou \og des BN-paires \fg, on pensera au cadre des \og strongly split BN-pairs of characteristic $p$ \fg\ du chapitre I.2 de \cite{CabEng04}, ce qui est loisible grâce au paragraphe 5.2 de \cite{HenVig11}. \\

Soit $S =\{ s_\alpha \ | \ \alpha \in \Delta \}$. Comme pour les autres éléments de $W$, on ne fera pas de distinction entre un élément $s \in S$ et son relèvement fixé dans $N_G(T)$. \\
Supposons à partir de maintenant que $R$ est un anneau commutatif unitaire de caractéristique $p$. On cherche à comprendre un peu mieux la structure de $(\St_J R)^B$ en tant que $\Hh_R(G,B)$-module. En effet, si $\pi$ est une $G$-représentation, la réciprocité de Frobenius confère à son espace de $B$-invariants \[ \pi^B \simeq \Hom_B(\id,\pi) \simeq \Hom_G \big( \Ind_B^G \, \id, \pi \big) \] une structure de module à droite sur \[ \Hh_R(G,B) := \End_{G}\big( \Ind_B^G \, \id \big) .\] A tout $w \in W$, on peut associer l'opérateur de Hecke $T_w$ défini sur $\pi^B$ par \[ v \mapsto v T_w = \sum_{\gamma \in B \backslash B w B} \gamma^{-1} v = \sum_{u \in (B \cap w^{-1} B w) \backslash B} u^{-1} w^{-1} v .\]
Pour $w \in W$, on note\footnote{par \cite{Car85}, paragraphe 1.18, Corollary 2.5.17 et discussion suivant Proposition 2.6.3, l'ordre n'a pas d'importance} \[ U_w := U \cap w U^- w^{-1} = \prod_{\alpha \in \Phired^+ , \, w^{-1}(\alpha) \in \Phired^-} U_\alpha .\] 
Par la première étape de la preuve de la proposition \ref{mainGrk1} (notamment l'équation (\ref{sUeqn})), pour $s \in S$, on peut voir $U_s$ comme un ensemble de représentants de $(B \cap s^{-1} B s) \backslash B = (B \cap sBs^{-1}) \backslash B$ (car $s^2$ est un élément de $T$). Pour $\pi = \Ind_{P_J}^{G} \, \id$, $w \in W^{J}$ et $s \in S$, la formule d'action se réécrit 
\begin{equation} \label{Taction} g_w T_s = \sum\limits_{u \in U_s} u^{-1} s^{-1} g_w = \sum_{u \in U_s} \id_{ P_J w^{-1} U_w su}. \end{equation} C'est cette action qu'on va investiguer à travers les deux résultats techniques suivants. \\
Pour $w \in W$, on note $w^J$ le représentant de $w W_J$ dans $W^J$ (on rappelle que $W^J$ est un système de représentants de $W/W_J$).

\begin{lem} \label{Usu}
Soient $w \in W^J$ et $s \in S$.
\begin{itemize}
\item[(a)] Si $(sw)^J = w$, alors on a \[ P_J w^{-1} U_w su = P_J w^{-1} U_w \quad \textrm{si } u \in U_s .\]
\item[(b)] Si $l((sw)^J) > l(w)$, alors on a \[ P_J w^{-1} U_w s U_s = P_J w^{-1} s U_{sw} .\]
\item[(c)] Si $l((sw)^J) < l(w)$, alors on a $s = s_\beta$ avec $\beta \in \Delta$ et $w^{-1}(\beta) \in \Phired^-$.  Posons \[ U' = \prod_{\alpha \in \Phired^+ \smallsetminus \{\beta\} , \, w^{-1}(\alpha) \in \Phired^-} U_\alpha ;\] c'est un sous-groupe de $U_w$. On a \[ P_J w^{-1} U' u s U_s = P_J w^{-1} U_w \quad \textrm{si } u \in U_s \smallsetminus \{1\} ,\] \[ P_J w^{-1} U' su = P_J w^{-1} s U_{sw} \quad \textrm{si } u \in U_s .\]
\end{itemize}
De plus, toutes ces relations sont des égalités entre produits directs d'ensembles.
\end{lem}
\textsf{Remarque :} 
Ce sont des raffinements dans l'axiomatique des systèmes de Tits que l'on peut déjà trouver dans \cite{Grk09} pour le cas déployé (Lemma 3.1): on se permet de reproduire sa preuve ici en rajoutant quelques commentaires, notamment pour le fait que les produits d'ensembles sont directs. \\
\textsf{Preuve :} \\
Commençons par le fait que les produits sont directs: il s'agit tout d'abord de voir que le produit $P_J w^{-1} U_w$ est direct. Supposons pour cela \[q_1 w^{-1} u_1 = q_2 w^{-1} u_2 \quad \textrm{avec } q_1, q_2 \in P_J, \ u_1, u_2 \in U_w .\] On a alors \[ u_1 u_2^{-1} \in w P_J w^{-1} \cap w U^- w^{-1} \cap U .\] Regardons à quoi ressemblent les éléments de $w^{-1} U w \cap U^- \cap P_J$. Ils sont inclus dans \[ U^- \cap w^{-1} U w = \prod_{\alpha \in\Phired^- , \, w(\alpha) \in \Phired^+} U_\alpha .\] De plus, par définition de $w \in W^J$ (voir (\ref{WJaltern})), pour tout $\alpha$ négativement engendré par $J$, on a $w(\alpha) \in \Phired^-$. Dès lors, comme $P_J$ ne contient que les sous-groupes radiciels associés à la partie quasi-close (au sens du 3.8 de \cite{BorTit65}) $\Phired^+ \cup ( \Phired^- \cap W_J. J )$, l'intersection $w^{-1} U w \cap U^- \cap P_J$ est réduite à l'élément neutre et le produit $P_J w^{-1} U_w$ est direct. \\
Une fois que toutes les égalités seront prouvées, le fait que les autres produits sont directs se ramène à chaque fois au cas précédent. Par exemple, regardons le terme de gauche de la première égalité de (c). L'ensemble $P_J w^{-1} U' us U_s$ est de cardinal majoré par $|P_J| \cdot |U'| \cdot |U_s| = |P_J| \cdot |U_w|$. Or l'égalité avec $P_J w^{-1} U_w$, qui est un produit direct, nous dit que le cardinal de $P_J w^{-1} U' us U_s$ est exactement $|P_J| \cdot |U_w|$. C'est donc que le produit est aussi direct et cela prouve bien le fait voulu. \\
Montrons (a). On a d'abord \begin{equation}\label{UsuA} P_J w^{-1} U_w s = P_J w^{-1} B s \subseteq P_J w^{-1} B \cup P_J w^{-1} s B \end{equation} par les équations (\ref{BwU}) et (\ref{sUeqn}) et l'axiomatique des BN-paires. Parce que l'on a $(sw)^J = w$, ces deux dernières doubles classes sont égales et leur union vaut simplement \[ P_J w^{-1} B = P_J w^{-1} U_w .\] L'inclusion (\ref{UsuA}) devient une égalité puisqu'en réappliquant $s$ on obtient: \[ P_J w^{-1} U_w = P_J w^{-1} U_w s^2 \subseteq P_J w^{-1} U_w s \subseteq P_J w^{-1} U_w .\] Pour finir, $P_J w^{-1} B$ est bien entendu invariant par translation à droite par $U_s$. \\
Attaquons nous à (b). On a dans un premier temps \[ P_J w^{-1} U_w s U_s = P_J w^{-1} B s U_s = P_J w^{-1} B s B .\] Par définition de $P_J$, on obtient \[ P_J w^{-1} B s B = \bigcup\nolimits_{v \in W_J} B v B w^{-1} B s B .\] L'axiomatique des systèmes de Tits nous dit que $B w^{-1} B s B$ est exactement $B w^{-1} s B$ car on a \[ l(w^{-1} s) = l(sw) > l(w) = l(w^{-1}) \] par \cite{Grk09}, Lemma 1.3.(a). On termine alors: \[ P_J w^{-1} B s B = \bigcup\nolimits_{v \in W_J} B v B w^{-1} s B = P_J w^{-1} s B = P_J w^{-1} s U_{sw} .\] 
Enfin pour (c), on remarque d'abord que l'on a $(s(sw)^J)^J = w^J =w$ et $l((sw)^J)<l(w)$. Par \cite{Grk09}, Lemma 1.4.(c), on a $l(sw)<l(w)$ et donc $w^{-1}(\beta) \in \Phired^-$. On observe aussi \[ s^{-1} U_{sw} s = s^{-1} U s \cap w U^- w^{-1} = \prod_{s(\alpha) \in \Phired^+, \, w^{-1}(\alpha) \in \Phired^-} U_\alpha .\] Or on a (voir Proposition 1.4 de \cite{Hum92}) \[ s(\Phired^+) = (\Phired^+ \smallsetminus \{\beta\}) \cup \{-\beta\} , \quad w^{-1}(-\beta) \in \Phired^+ .\] Il en résulte que la condition sur les indices du produit se réécrit \[ \alpha \in \Phired^+ \smallsetminus \{\beta\}, \ w^{-1}(\alpha) \in \Phired^- ;\] et on en déduit $s^{-1} U_{sw} s = U'$. Cela implique directement $U' s = s U_{sw}$ et comme on a \[ u \in U_s \subseteq B , \quad P_J w^{-1} s U_{sw} = P_J w^{-1} s B ,\] la dernière égalité en découle. Pour la première égalité, écrivons l'inclusion (le détail est identique au (b)) \[ P_J w^{-1} U_w s \subseteq P_J w^{-1} U_w \cup P_J w^{-1} s U_{sw} ;\] cette union est disjointe car on a $sw W_J \neq w W_J$. De l'égalité $P_J w^{-1} U' s = P_J w^{-1} s U_{sw}$ que l'on vient de prouver, et parce que le produit $P_J w^{-1} U_w$ est direct, on déduit \[ P_J w^{-1} (U_w \smallsetminus U') s \subseteq P_J w^{-1} U_w .\] Lorsque $u$ est un élément de $U_s \smallsetminus \{1\} = U_\beta \smallsetminus \{1\} \subseteq U_w$, on a l'inclusion $U' u \subseteq U_w \smallsetminus U'$; il s'ensuit \[ P_J w^{-1} U' u s U_s \subseteq P_J w^{-1} (U_w \smallsetminus U') s U_s \subseteq P_J w^{-1} U_w U_s = P_J w^{-1} U_w .\] Voyons l'inclusion inverse: par \cite{Car85}, Corollary 2.6.2, il existe $u_1 \in U_s$ et $b \in U_s T \subseteq B$ tel que l'on ait la décomposition $sus = bsu_1$. On a alors \[ P_J w^{-1} s B sus U_s = P_J w^{-1} s B s U_s = P_J w^{-1} s B s B .\] Il s'ensuit que \[ P_J w^{-1} U' = P_J w^{-1} s U_{sw} s = P_J w^{-1} s B s \] est inclus dans \[ P_J w^{-1} s B sus U_s = P_J w^{-1} s U_{sw} sus U_s = P_J w^{-1} U' us U_s .\] Au final, on a bien l'égalité voulue (gr\^ ace \` a $U_w = U' U_s$) $\hfill\Box$\\

Le cas fini peut se voir de manière similaire au cas $p$-adique traité dans le lemme \ref{rgtotal}: le $R$-module $C(P_J \backslash G, R)$ est libre et une base est donnée par les fonctions $g_w$ pour $w$ parcourant $W^J$. On tâche d'investiguer sa structure en tant que $\Hh_R(G,B)$-module à droite lorsque $R$ est de caractéristique $p$.

\begin{lem} \label{lemT}
Soient $w \in W^J$ et $s \in S$.
\begin{itemize}
\item[(a)] Si $(sw)^J = w$, alors on a $g_w T_s = 0$.
\item[(b)] Si $l((sw)^J) > l(w)$, alors on a $g_w T_s = g_{(sw)^J}$.
\item[(c)] Si $l((sw)^J) < l(w)$, alors on a $g_w T_s = - g_w$.
\end{itemize}
\end{lem}
\textsf{Preuve :} \\
Le (a) est conséquence immédiate de (\ref{Taction}), du lemme \ref{Usu}.(a) et de l'égalité $|U_s| = 0$ dans $R$ de caractéristique $p$ (voir \cite{Car85}, page 74). Le (b) suit de (\ref{Taction}) et de la décomposition en produit direct du lemme \ref{Usu}.(b) aussi. \\
Intéressons nous au (c): on utilise la décomposition en produit direct $U_w = U' U_s$. On a alors \[ g_w T_s = \sum_{u \in U_s} \sum_{u' \in U_s} \id_{P_J w^{-1} U' u' s u} = \sum_{u \in U_s} \id_{P_J w^{-1} U' su} + \sum_{u \in U_s} \sum_{u' \neq 1} \id_{P_J w^{-1} U' u' su } .\] Par le lemme \ref{Usu}.(c) et $|U_s| = 0$ dans $R$, le premier terme vaut $0$ et le second \[ \sum_{u' \in U_s \smallsetminus \{1\}} \id_{P_J w^{-1} U_w} = -g_w .\] Le lemme en découle. $\hfill\Box$\\

Les actions précédemment étudiées dans $\big( \Ind_B^G \, \id \big)^B$ sont compatibles à celles du quotient $\big( \St_J R \big)^B$.

\begin{prop} \label{mainGrk2}
Supposons $R$ de caractéristique $p$. Tout sous-$\Hh_R(G,B)$-module non nul de $\big( \St_J R \big)^B$ contient l'élément $\overline{g}_{z^J}$ de $\St_J R$, où $z^J$ désigne l'élément de longueur maximale\footnote{un élément de longueur maximale est aussi maximal pour $<_J$ par \cite{Grk09}, Lemma 1.4.(d); l'existence d'un unique élément $<_J$-maximal est ensuite donnée par \cite{Grk09}, Lemma 1.4.(e)} de $W^J$.
\end{prop}
\textsf{Preuve :} \\
Soit $E$ un sous-$\Hh_R(G,B)$-module non nul de $(\St_J R)^B$. Par la proposition \ref{mainGrk1}, $E$ contient un élément non nul $h = \sum\limits_{w \in W^J_{\pr}} \alpha_w(h) \overline{g}_w$ avec $\alpha_w(h) \in R$. On commence par rappeler la définition de l'ordre $<_J$ introduit dans \cite{Grk09}: on note $w <_J w'$ s'il existe $s_1, \dots, s_r$ dans $S$ tels que $w^{(i)} = (s_i \dots s_1 w)^J$ vérifie $l(w^{(i)}) > l(w^{(i-1)})$ pour tout $1 \leq i \leq r$ et $w^{(r)} = w'$. On fixe une énumération $w_1, w_2, \dots$ des éléments de $W^J_{\pr}$ vérifiant $w_j <_J w_i \Rightarrow i <j$: en particulier, on a $w_1 = z^J$. On veut montrer qu'il existe $h \in E$ non nul tel que \[ t(h) := \min \, \{ i \geq 1 \ | \ \forall j>i, \ \alpha_{w_j}(h) = 0 \}\] soit égal à $1$, c'est-à-dire $\overline{g}_{z^J} \in E$. Supposons le contraire et donc on a \[ t := \min \, \{ t(h) \ | \ h \in E \smallsetminus \{0\}\} \geq 2 .\] Par \cite{Grk09}, Lemma 1.5, il existe $w' \in W^J_{\pr}$ et $s \in S$ tel que $w_t <_J w'$, $l((sw_t)^J) < l(w_t)$ et $l(w') \leq l((sw')^J)$. Par définition de $<_J$, il existe $s_1, \dots, s_r$ dans $S$ tels que $w^{(i)} = (s_i \dots s_1 w_t)^J$ vérifie $l(w^{(i)}) > l(w^{(i-1)})$ pour tout $1 \leq i \leq r$ et $w^{(r)} = w'$. \\ Soit $h \in E \smallsetminus \{0\}$ avec $t(h)=t$. Commençons par remarquer que, grâce au lemme \ref{lemT}, on a $\alpha_{w^{(r)}}(h T_s) = 0$. On peut donc considérer $h \in E \smallsetminus \{0\}$ avec $t(h)=t$ et tel que $k(h) \geq 0$ est minimal, égal à $k$, où $k(h)$ est l'entier minimal de $[0,r]$ défini par $\alpha_{w^{(k(h))}}(h)=0$. Si on a $k(h)=0$, alors $\alpha_{w_t(h)}$ est nul, ce qui est contradictoire. On suppose ainsi $k>0$, et on va le faire diminuer: considérons $h' = h T_{s_k}$, et observons 
\[ \alpha_{w^{(k-1)}}(h') = - \alpha_{w^{(k)}}(h)=0 , \quad \alpha_{w^{(k)}}(h') = \alpha_{w^{(k-1)}}(h) \neq 0 .\]
Ceci nous assure $h' \neq 0$ et $k(h')<k$. Cela contredit la minimalité de $k$ et donc l'hypothèse initiale: on vient donc de montrer l'existence de $h \in E \smallsetminus \{0\}$ avec $t(h) = 1$. C'est le résultat voulu. $\hfill\Box$ 

\section{Paramètres de Hecke-Satake et preuve du théorème \ref{Stirred}} \label{Stfinal}

Soit $R$ un corps de caractéristique $p$. On commence par \' etudier le $K$-socle d'une repr\' esentation de Steinberg g\' en\' eralis\' ee. 

\begin{prop}
Soient $R$ un corps de caractéristique $p$ et $J \subseteq \Delta$. Le $K$-socle de la Steinberg généralisée $\St_J R$ est irréductible.
\end{prop}
\textsf{Preuve :} \\
Soit $V$ une sous-$K$-représentation irréductible de $\St_J R$. L'injection de la proposition \ref{mainGrk1} nous permet de voir $V^{U \cap K}$ comme un sous-espace de 
\[ \big(\St_J R \big)^{I(1)} = \big(\St_J R \big)^{I} \simeq \big(\overline{\St}_{\overline{J}} R \big)^{\overline{B}} \] 
(voir d\' ebut de preuve du corollaire \ref{Stadm}); de cette manière, $V^{U \cap K} = V^{B \cap K}$ est une droite stable par l'action de $\Hh_R(\overline{G},\overline{B})$. Or, par la proposition \ref{mainGrk2}, tout sous-$\Hh_R(\overline{G},\overline{B})$-module non nul de $\big( \overline{\St}_{\overline{J}} R \big)^{\overline{B}}$ contient l'image $\overline{g}_J$ de l'élément $g_J := \id_{\overline{P}_J z^{J} \overline{B}} \in \Ind_{\overline{P}_J}^{\overline{G}} \, \id$ dans $\overline{\St}_{\overline{J}} R$. Ainsi $V$ est généré par $\overline{g}_J$ en tant que $K$-représentation, et le $K$-socle de $\St_J R$ est irréductible. \hfill$\Box$ \\

Soit $V$ une $K$-représentation irréductible. On va noter $\Hh_R(G,K,V)$ l'algèbre de Hecke-Satake $\End_G(\ind_K^G \, V)$. Les éléments de $\Hh_R(G,K,V)$ sont des opérateurs à support fini parmi les doubles classes $K \backslash G / K$, donc on va fixer un système de représentants dominants dans $A$ (au sens qu'ils contractent $B$) que l'on notera $\Sigma_+$. \\ La transformée de Satake (voir \cite{HenVig11}, paragraphe 7.3) est un isomorphisme de $R$-algèbres 
\begin{equation} \label{Satake} \Sss : \Hh_R(G,K,V) \xrightarrow{\sim} \Hh_R^+ \big(A,A \cap K, V_{U \cap K} \big) \end{equation} 
où $\Hh_R^+ \big(A,A \cap K, V_{U \cap K} \big)$ désigne la sous-algèbre de engendrée par les opérateurs portés par la classe $z(A \cap K)$ pour $z \in \Sigma_+$ quand ils existent. \\

Soit $J$ un sous-ensemble de $\Delta$. On va particulièrement s'intéresser au cas où $V$ est $V_J$, l'unique $K$-représentation irréductible $M_J$-coréguliere\footnote{cela veut dire que le stabilisateur dans $K$ de la représentation associée à $V_J^{U^- \cap K}$ est inclus dans $(P_J^- \cap K) K(1)$} telle que l'action de $M_J$ sur la droite $(V_J)_{N_J \cap K}$ est triviale (par la Proposition 3.11 de \cite{HenVig11b}). Dans ce cas-là, en particulier, $V_{U \cap K}$ est la représentation triviale $\id$ de $A \cap K$ et toute classe $z(A \cap K)$ pour $z \in \Sigma_+$ porte un unique opérateur de Hecke $\tau_z$ de $\Hh_R^+ \big(A,A \cap K, V_{U \cap K} \big)$, envoyant $z$ sur $\Id_{V_{U \cap K}}$ (voir \cite{HenVig11}, paragraphe 7.3); et ils en constituent une base en tant que $R$-espace vectoriel. De plus, l'algèbre $\Hh_R^+ \big(A,A \cap K, \id \big)$ est alors commutative et de type fini sur $R$ (voir paragraphe 7.2 de \cite{HenVig11}). On en déduit alors que $\Hh_R(G,K,V_J)$ est aussi commutative et de type fini par l'isomorphisme de Satake (\ref{Satake}). \\

Soit $\pi$ une $G$-représentation admissible à coefficients dans $R$. Alors, comme $K(1)$ est un pro-$p$-groupe ouvert, le sous-espace des $K(1)$-invariants $\pi^{K(1)}$ est non nul et de dimension finie sur $R$. Il possède donc une sous-$K$-représentation irréductible $V$. En particulier, $\Hom_K(V,\pi)$ est non trivial et de dimension finie. De plus, c'est un module à droite sur l'algèbre $\Hh_R(G,K,V)$. \\ On suppose à présent que $V$ est un $V_J$ pour un certain sous-ensemble $J$ de $\Delta$, et que $R$ est algébriquement clos. Alors $\Hh_R(G,K,V_J)$ étant commutative, elle possède un sous-espace propre associé à un caractère $\chi: \Hh_R(G,K,V_J) \to R$ dans son action sur $\Hom_K(V_J,\pi)$. Dit autrement, on a un morphisme non nul de $G$-représentations \[ \ind_K^G \, V_J \otimes_{\Hh,\chi} R \to \pi .\] On tâche de déterminer les caractéristiques d'un tel caractère $\chi$ lorsque $\pi$ est $\St_J R$. Pour être précis, notons $\chi^{(A)}$ la composée \[ \Hh_R^+ \big(A,A \cap K, (V_J)_{U \cap K} \big) \mathop{\xrightarrow{\Sss^{-1}}}\limits_{\sim} \Hh_R(G,K,V_J) \xrightarrow{\chi} R ;\] c'est $\chi^{(A)}$ que l'on va expliciter.

\begin{prop} \label{paramHecke}
Soit $R$ un corps algébriquement clos de caractéristique $p$. Soit $J$ un sous-ensemble de $\Delta$.
\begin{itemize}
\item[(i)] Le $K$-socle de $\St_J R$ s'identifie \` a $V_J$.
\item[(ii)] Il existe un morphisme non nul de $G$-représentations \[ \ind_K^G \, V_J \otimes_{\Hh,\chi} R \to \St_J R \] si et seulement si $\chi^{(A)}$ est le caractère envoyant $\tau_z$ sur $1$ pour tout $z \in \Sigma_+$.
\end{itemize}
\end{prop}
\textsf{Remarque :} 
En particulier, $\St_J R$ n'est pas supersingulier au sens de \cite{HenVig11b}. \\
\textsf{Preuve :} \\
On note $\Hh = \Hh_R(G,K,V_J)$ et $\Hh_M = \Hh_R(M_J, M_J \cap K, \id)$. Soient $\chi : \Hh \to R$ et $\chi_M : \Hh_M \to R$ les caractères d'algèbres associés \` a $\chi^{(A)}$, caractère envoyant tout $\tau_z$  sur $1$ (pour $z \in \Sigma_+$). Par \cite{HenVig11b}, Theorem 1.2, parce que $V_J$ est $M_J$-corégulière et que $(V_J)_{N_J \cap K}$ est la $(M_J \cap K)$-représentation triviale, on a un morphisme surjectif de $G$-représentations 
\begin{equation} \label{VJsurj} \ind_K^G \, V_J \otimes_\chi R \xrightarrow{\sim} \Ind_{P_J}^G \big( \ind_{M_J \cap K}^{M_J} \, \id \otimes_{\chi_M} R \big)\twoheadrightarrow \Ind_{P_J}^G \, \id .\end{equation}
On en déduit\footnote{jusqu'à présent, hormis dans l'introduction, on a défini $\St_J R$ à partir des $C^\infty(P_J \backslash G,R)$. Remarquons que la définition peut se faire de manière équivalente à partir d'induites paraboliques: l'application $R$-linéaire $C^\infty(P_J \backslash G,R) \to \Ind_{P_J}^G \, \id$ est un isomorphisme car la projection canonique $G \to P_J \backslash G$ possède une section continue} l'existence d'un morphisme surjectif (en particulier non nul) de $G$-représentations $\ind_K^G \, V_J \otimes_\chi R \to \St_J R$. Il s'ensuit que $V_J$ génère $\St_J R$ en tant que $G$-représentation. Parce que $\St_J R$ est de $K$-socle irréductible, $V_J$ est l'unique $K$-représentation irréductible contenue dans $\St_J R$. Cela prouve (i) et le sens indirect de (ii). \\
Prouvons maintenant la seconde implication de (ii). Comme $K \cap P_J \backslash K \to P_J \backslash G$ est un homéomorphisme, on a par réciprocité de Frobenius \[ \Hom_K(V_J, \Ind_{P_J}^G \, \id) = \Hom_{K \cap P_J} (V_J, \id)\] et donc, puisque l'on a $(V_J)_{N_J \cap K} = \id$, on déduit \[ \Hom_K(V_J, \Ind_{P_J}^G \, \id) = \Hom_{K \cap M_J} (\id,\id) .\] Ce dernier espace est donc un $R$-espace vectoriel de dimension $1$. Le morphisme surjectif $\Ind_{P_J}^G \, \id \twoheadrightarrow \St_J R$ de $G$-représentations induit le morphisme 
\[ \psi : \Hom_K \big( V_J, \Ind_{P_J}^G \, \id \big) \to \Hom_K \big( V_J, \St_J R \big) \]
de $R$-espaces vectoriels. On vient de voir que $\Hom_K \big( V_J, \Ind_{P_J}^G \, \id \big)$ est de dimension $1$, et $\Hom_K \big( V_J, \St_J R \big)$ est aussi de dimension $1$ par (i) et lemme de Schur. De ce fait, $\psi$ est un isomorphisme si et seulement si il est non nul. Or le fait que (\ref{VJsurj}) induise $\ind_K^G \, V_J \otimes_\chi R \twoheadrightarrow \St_J R$ implique la non nullité de $\psi$: $\psi$ est un isomorphisme. \\
Apr\` es application de la r\' eciprocit\' e de Frobenius, on a un isomorphisme 
\[ \Hom_G \big( \ind_K^G \, V_J, \Ind_{P_J}^G \, \id \big) \xrightarrow{\sim} \Hom_G \big( \ind_K^G \, V_J, \St_J R \big) \]
de $R$-espaces vectoriels. Ce dernier est $\Hh_R(G,K,V_J)$-\' equivariant car $\psi$ est juste induit par la projection d\' efinissant $\St_J R$. De ce fait, toute fl\` eche non nulle $\ind_K^G \, V_J \otimes_{\Hh, \chi} R \to \St_J R$ pour un certain $\chi$ se factorise \` a travers $\ind_K^G \, V_J \otimes_{\Hh, \chi} R \to \Ind_{P_J}^G \, \id$. Et par le diagramme (4) de \cite{HenVig11b}, l'isomorphisme de r\' eciprocit\' e de Frobenius 
\[ \Hom_G \big( \ind_K^G \, V_J, \Ind_{P_J}^G \, \id \big) \xrightarrow{\sim} \Hom_{M_J} \big( \ind_{M_J \cap K}^{M_J} \, \id, \id \big) \]
est $\Hh_R(G,K,V_J)$-\' equivariant, o\` u l'action au but se fait \` a travers la transform\' ee de Satake partielle $\Sss_G^{M_J} : \Hh_R(G,K,V_J) \hookrightarrow \Hh_R(M_J, M_J \cap K, \id)$. Il s'agit donc de d\' eterminer les valeurs propres de Hecke possibles pour l'action de $\Hh_R(M_J, M_J \cap K, \id)$ sur $\Hom_{M_J} \big( \ind_{M_J \cap K}^{M_J} \, \id, \id \big)$.
Cet espace est unidimensionnel puisqu'il est isomorphe à \[ \Hom_{M_J \cap K} (\id, \id) = \Hom_{A \cap K}(\id, \id) = \Hom_A( \ind_{A \cap K}^A \, \id, \id) ,\] et c'est à travers ce dernier que l'action de $\Hh_R(A, A \cap K, \id)$ se lit naturellement. On conclut donc qu'elle se fait à travers le caractère $\chi^{(A)}$ envoyant chaque $\tau_z$ sur $1$, d'où l'implication qu'il restait à prouver pour (ii).
\hfill$\Box$ \\

On exhibe de la preuve pr\' ec\' edente le fait plus précis suivant, qui va facilement impliquer le théorème \ref{Stirred}. 

\begin{cor}
Soient $R$ un corps algébriquement clos de caractéristique $p$ et $J$ un sous-ensemble de $\Delta$. Le $K$-socle $V_J$ de $\St_J R$ l'engendre en tant que $G$-repr\' esentation.
\end{cor}

On a alors l'irréductibilité de $\St_J R$ dans le cas $R$ algébriquement clos de caractéristique $p$ comme suit. Soit $\pi \subseteq \St_J R$ une sous-représentation non nulle. Alors $\pi$ contient une sous-$K$-représentation irréductible qui, par l'argument précédent, est donc $V_J$. Mais on sait que $V_J$ génère $\St_J R$; donc on a $\pi = \St_J R$ et l'irréductibilité voulue. Le théorème \ref{Stirred} en découle par le corollaire \ref{Stadm} et le lemme \ref{Ralgclos} tout à fait général.

\section{Induites paraboliques de Steinberg généralisées} \label{Stfiltre}

Commençons par un mot sur la preuve du corollaire \ref{Stgeneral}, notamment sur le fait que deux $J, J' \subseteq \Delta$ distincts engendrent des $\St_J R$ et $\St_{J'} R$ non isomorphes. Cela suit immédiatement de ce que l'on vient de faire puisque l'on a alors $V_J \neq V_{J'}$. \\ 

Définissons la filtration suivante sur $\Ind_{P_J}^G \, \id$ : \[ \Fil^i = \begin{cases} \sum_{J' \supseteq J, \ |J' \smallsetminus J| \geq i} \Ind_{P_{J'}}^G \, \id & \textrm{pour } 0 \leq i \leq |\Delta \smallsetminus J| , \\ 0 & \textrm{pour } i > |\Delta \smallsetminus J| . \end{cases} \]
Et montrons que c'est la filtration par les cosocles de $\Ind_{P_J}^G \, \id$, c'est-à-dire la filtration descendante définie par $\Fil^0 = \Ind_{P_J}^G \, \id$ et $\Fil^i$ est telle que $\grr^{i-1} := \Fil^{i-1}/\Fil^i$ est le cosocle de $\Fil^{i-1}$ pour $i \geq 1$. On a le diagramme commutatif suivant, pour $i \geq 0$ et $J' \supseteq J$ avec $|J' \smallsetminus J| = i$ (en particulier $i \leq |\Delta \smallsetminus J|$): 
\[\xymatrix{ & \sum_{J'' \supsetneq J'} \Ind_{P_{J''}}^G \, \id \ar[r] \ar@{^(->}[d] & \Ind_{P_{J'}}^G \ar[r] \ar@{^(->}[d] & \St_{J'} R \ar[r] \ar@{.>}[d] & 0 \\ 0 \ar[r] & \Fil^{i+1} \ar[r] & \Fil^i \ar[r] & \grr^i \ar[r] & 0 } .\]
Si la flèche $\St_{J'} R \to \grr^i$ était triviale, cela voudrait dire que l'image de $\Ind_{P_{J'}}^G \, \id \hookrightarrow \Fil^i$ serait incluse dans $\Fil^{i+1}$, ce qui serait absurde par le lemme \ref{StFilabs} ult\' erieur. De ce fait, et parce que l'on a $\St_{J_1} R \nsim \St_{J_2} R$ pour $J_1 \neq J_2$, on a une injection \[ \bigoplus_{J' \supseteq J , \ |J' \smallsetminus J| = i} \St_{J'} R \hookrightarrow \grr^i .\] 
Pour voir que c'est bien le cosocle de $\Fil^i$, comme on connaît les composantes de Jordan-Hölder de $\Ind_{P_J}^G \, \id$ (ce sont les $\St_{J'} R$ pour $J' \supseteq J$, avec\footnote{ce sont les telles repr\' esentations de Steinberg comme on peut le voir par r\' ecurrence descendante sur $|J|$ \` a partir de la d\' efinition des $\St_{J'} R$. La décomposition $G =\coprod_{w \in W^J} P_J w^{-1} B$ fait que les restrictions respectives de $C^\infty(P_J \backslash G, R)$ et $\St_J R$ à $B$ possèdent des filtrations telles que $\bigoplus_i \gr^i$ sont respectivement égales à $\bigoplus_{w \in W^J} C^\infty(P_J \backslash P_J w^{-1} B, R)$ et $\bigoplus_{w \in W^J_\pr} C^\infty(P_J \backslash P_J w^{-1} B, R)$. Comme on a de plus $W^J = \coprod_{J' \supseteq J} W^{J'}_\pr$, les $\St_{J'} R$ apparaissent avec multiplicité $1$ comme voulu} multiplicit\' e $1$), il suffit de voir que toute flèche $\Ind_{P_{J'}}^G \, \id \to \St_{J''} R$ est triviale dès que l'on a $J'' \supsetneq J'$ ou $J'' \nsupseteq J'$. 
Le second cas est clair et on veut montrer que tout morphisme $f : \Ind_{P_{J'}}^G \, \id \to \St_{J''} R$ de $G$-représentations est trivial si $J'' \supsetneq J'$. Si ce n'était pas le cas, $f$ serait surjectif par irréductibilité de $\St_{J''} R$, de noyau contenant $\St_{J'} R$. Mais comme $V_{J'} \subseteq (\St_{J'} R)|_K$ génère $\Ind_{P_{J'}}^G \, \id$ par (\ref{VJsurj}), toute injection $\St_{J'} R \hookrightarrow \Ind_{P_{J'}}^G \, \id$ se doit d'être un isomorphisme. C'est absurde, donc $f$ est nul et $\St_{J'} R$ est bien le plus gros quotient semi-simple de $\Ind_{P_{J'}}^G \, \id$. Il en r\' esulte que $(\Fil^i)_i$ est bien la filtration par les cosocles comme annoncé.

\begin{lem} \label{StFilabs}
Soient $J \subseteq \Delta$ un ensemble, $i \in [0, |\Delta \smallsetminus J|]$ un entier et $J' \supseteq J$ un sous-ensemble de $\Delta$ avec $|J' \smallsetminus J| = i$. Alors l'image de $\Ind_{P_{J'}}^G \, \id \hookrightarrow \Fil^i$ n'est pas incluse dans $\Fil^{i+1}$.
\end{lem}
\textsf{Preuve :} \\
Le cas $i=0$ r\' esulte de ce que $\St_J R$ est non triviale. On suppose donc \` a pr\' esent $i \geq 1$. \\
Supposons par l'absurde que $\Ind_{P_{J'}}^G \, \id$ est incluse dans $\Fil^{i+j}$ pour un $j \geq 1$ maximal, c'est-\` a-dire avec $\Ind_{P_{J'}}^G \, \id \nsubseteq \Fil^{i+j+1}$ (un tel $j$ existe bien puisque la filtration devient nulle au bout d'un certain rang). On commence par \' etablir 
\begin{equation} \label{IndFil} \Ind_{P_{J'}}^G \, \id = \Big( \Ind_{P_{J'}}^G \, \id \cap \Fil^{i+j+1} \Big) + \sum_{J'' \supsetneq J'} \Ind_{P_{J''}}^G \, \id .\end{equation}
Soit $f$ un \' el\' ement de $\Ind_{P_{J'}}^G \, \id \Big/ \Big( \Big( \Ind_{P_{J'}}^G \, \id \cap \Fil^{i+j+1} \Big) + \sum_{J'' \supsetneq J'} \Ind_{P_{J''}}^G \, \id \Big)$, que l'on voit dans $\Ff = \Fil^{i+j} \Big/ \Big( \Big( \Ind_{P_{J'}}^G \, \id \cap \Fil^{i+j+1} \Big) + \sum_{J'' \supsetneq J'} \Ind_{P_{J''}}^G \, \id \Big)$. Alors $f$ s'\' ecrit $\sum_{J''} f_{J''}$ o\` u les $J''$ parcourent les $J'' \supseteq J$ avec $|J'' \smallsetminus J| \geq i+j$, $J '' \nsupseteq J'$ et $f_{J''}$ appartient \` a l'image de $\Ind_{P_{J''}}^G \, \id$ dans $\Ff$. On prend ensuite un $\alpha \in J' \smallsetminus J$ (cet ensemble est non vide car $i$ est non nul) et $s \in W_{J'}$ la r\' eflexion correspondante. Parce que l'on a $J' \nsubseteq J''$, et que l'on a quotient\' e par $\Ind_{P_{J'}}^G \, \id \cap \Fil^{i+j+1}$, l'\' ecriture $f = \sum f_{J''}$ est en fait invariante \` a gauche par $s$. En effectuant de m\^ eme pour toute racine de $J' \smallsetminus J$, on voit que chaque $f_{J''}$ est en fait nulle dans $\Ff$ (car $j \geq 1$), ce qui donne la nullit\' e de $f$. Et donc (\ref{IndFil}) comme annonc\' e: mais cela implique que l'on a une surjection naturelle $\Ind_{P_{J'}}^G \, \id \cap \Fil^{i+j+1} \twoheadrightarrow \St_{J'} R$. De ce fait, ou bien il existe un $J'' \supsetneq J'$ avec $\Ind_{P_{J''}}^G \, \id \nsubseteq \Fil^{i+j+1}$, ce qui est exclu par l'inclusion $\Ind_{P_{J'}}^G \, \id \subseteq \Fil^{i+j}$. Ou bien $\Ind_{P_{J'}}^G \, \id \cap \Fil^{i+j+1}$ est tout $\Ind_{P_{J'}}^G \, \id$, c'est-\` a-dire que $\Ind_{P_{J'}}^G \, \id$ est inclus dans $\Fil^{i+j+1}$, ce qui contredit la maximalit\' e de $j \geq 1$. Toutes les possibilit\' es am\` enent \` a des contradictions: le lemme est prouv\' e. \hfill$\Box$ 

\begin{cor} \label{JHsteinberg}
Soient $J' \subseteq J$ des sous-ensembles de $\Delta$. Alors la représentation\footnote{où on note $\St^{M_J}_{J'} R$ une représentation de Steinberg généralisée du groupe réductif $M_J$} $\Ind_{P_J}^G (\St^{M_J}_{J'} R)$ est de longueur finie, de constituants de Jordan-Hölder les $\St_{J''} R$ avec $J'' \subseteq \Delta$ vérifiant $J \cap J'' = J'$, chacun apparaissant avec multiplicité $1$. 
\end{cor}
\textsf{Preuve :} \\
Commençons tout d'abord par remarquer que, puisque l'on a $J' \subseteq J$, la condition $J \cap J'' = J'$ est équivalente à $J'' \supseteq J'$ et $J \smallsetminus J' \subseteq \Delta \smallsetminus J''$. C'est sous cette seconde forme que nous allons l'utiliser au cours du raisonnement qui suit. \\
Prouvons-le par récurrence descendante sur $J' \subseteq J$. L'étape d'initiation $J' = J$ est juste le corollaire \ref{Stgeneral}. Soit $J' \neq J$ et supposons le résultat vrai pour tout parabolique $J_0 \subseteq J$ contenant strictement $J'$. Par définition de la représentation de Steinberg généralisée, on a une suite exacte de $M_J$-représentations \begin{equation} \label{Stexact} 0 \to \Ker \to \Ind_{P_{J'}}^{M_J} \, \id \to \St_{J'}^{M_J} R \to 0 ,\end{equation} où $\Ker$ est par là-même définie et a constituants de Jordan-Hölder les $\St_{J_0} \id$ pour $J' \subsetneq J_0 \subseteq J$ par le corollaire \ref{Stgeneral}. Appliquons le foncteur exact $\Ind_{P_J}^G$ (voir \cite{Vig12}, Proposition 1.1) à (\ref{Stexact}) pour obtenir: \[ 0 \to \Ind_{P_J}^G(\Ker) \to \Ind_{P_J}^G \big( \Ind_{P_{J'}}^{M_J} \, \id \big) \to \Ind_{P_J}^G (\St_{J'}^{M_J} R) \to 0 .\] Comme on a \[ M_J/(M_J \cap P_{J'}) \xrightarrow{\sim} M_J B /(M_J \cap P_{J'})B = P_J / (M_J \cap P_{J'}) B = P_J / P_{J'} \] par décomposition de Levi, le terme central est $\Ind_{P_{J'}}^G \, \id$, de constituants de Jordan-Hölder les $\St_{J''} R$ avec $J'' \supseteq J'$ par le corollaire \ref{Stgeneral}. Par hypothèse de récurrence, les constituants de $\Ind_{P_J}^G(\Ker)$ sont les $\St_{J''} R$ avec $J'' \subseteq \Delta$ vérifiant $J'' \supseteq J_0$ et $J \smallsetminus J_0 \subseteq \Delta \smallsetminus J''$ pour un certain $J' \subsetneq J_0 \subseteq J$. \\
Mais alors, soit $J''$ tel que $\St_{J''} R$ est un constituant de Jordan-Hölder de $\Ind_{P_J}^G(\St_{J'}^{M_J} R)$. On a $J'' \supsetneq J'$, et aussi, $J \smallsetminus J' \subseteq \Delta \smallsetminus J''$. En effet, supposons cette seconde inclusion fausse, c'est-à-dire $J'' \cap (J \smallsetminus J') \neq \varnothing$: on considère $J_0 \subseteq \Delta$ avec \[ J_0 := J' \cup \big( J'' \cap ( J \smallsetminus J') \big) \supsetneq J' \] et on a $J_0 \subseteq J''$, $J \smallsetminus J_0 \subseteq \Delta \smallsetminus J''$, donc $\St_{J''} R$ apparaît déjà dans $\Ind_{P_J}^G (\Ker)$ par l'assertion de multiplicité $1$ dans le corollaire \ref{Stgeneral}. C'est absurde. Enfin, la décomposition disjointe (qu'il est plus facile de voir avec la condition équivalente $J \cap J'' = J'$ dans le terme de droite) \[ \{ J'' \supseteq J' \} = \coprod_{J_0 \supseteq J'} \{ J'' \supseteq J_0 \ | \ J \smallsetminus J_0 \subseteq \Delta \smallsetminus J'' \}\] nous permet de dire que ce sont les seuls constituants qui interviennent. La récurrence est terminée. \hfill$\Box$

\section{Appendice: de la liberté de $C^\infty(X,\Z)$} \label{Ccompactlibre}

Parce qu'on se sert constamment du fait suivant, on se permet de rappeler sa preuve probablement bien connue. \\
Pour $X$ un espace topologique, on note $\U$ la famille des recouvrements de $X$ par un nombre fini d'ouverts disjoints. On remarque que $\U$ un ensemble ordonné par $U \leq V$ si pour tout $A \in U$ il existe des $B_i \in V$ vérifiant $A = \bigcup B_i$. 

\begin{lem} \label{Ccompact}
Soit $X$ un espace topologique compact et totalement discontinu. Supposons que $\U$ possède un sous-ensemble $\V$ dénombrable et cofinal. Alors l'espace $C^\infty(X,\Z)$ des fonctions localement constantes sur $X$ à valeurs dans $\Z$ est un $\Z$-module libre.
\end{lem}
\textsf{Remarque :} 
Par changement de base $\Z \to R$, si $R$ est un anneau commutatif unitaire, $C^\infty(X,R)$ est alors aussi un $R$-module libre. \\
\textsf{Preuve :} \\
Que $C^\infty(X,\Z)$ est un $\Z$-module est évident; il s'agit de voir qu'il est libre. Soit $f$ un élément de $C^\infty(X,\Z)$. Pour tout $x \in X$, on fixe un ouvert $V_x$ contenant $x$ tel que $f$ est constant sur $V_x$. Le compact $X$ est alors recouvert par les $V_x$ pour $x$ parcourant $X$. Par compacité, on peut en extraire un recouvrement fini $X = \bigcup X_i$ par des ouverts $X_i$ sur lesquels $f$ est constant. Si $X_i$ et $X_j$ sont non disjoints, on peut les remplacer par $X_i \cap X_j$, $X_i \smallsetminus X_i \cap X_j$ et $X_j \smallsetminus X_i \cap X_j$, qui sont trois ouverts deux à deux disjoints d'union $X_i \cup X_j$. En répétant le procédé, on peut supposer que l'union $X = \bigcup X_i$ est une union finie disjointe. \\
On note que $\U$ est filtrant puisque si $U, V \in \U$, on peut, à partir du procédé précédent appliqué à $U \cup V$, obtenir un recouvrement $W$ vérifiant $U \leq W$ et $V \leq W$. Pour tout $U = \{X_i\}$ élément de $\U$, on note $C_{(U)}(X,\Z)$ le sous-$\Z$-module de $C^\infty(X,\Z)$ constitué des fonctions constantes sur chaque $X_i$. On vient de voir que tout élément de $C^\infty(X,\Z)$ vit dans un $C_{(U)}(X,\Z)$ pour un $U \in \U$ convenable. Et si on considère les flèches d'inclusion $C_{(U)}(X,\Z) \to C_{(V)}(X,\Z)$ pour $U \leq V$, on obtient un système inductif et on écrit $C^\infty(X,\Z)$ comme limite inductive de modules libres: \[ C^\infty(X,\Z) = \varinjlim_{U \in \U} C_{(U)}(X,\Z) \simeq \varinjlim_{U \in \U} \Z^{|U|} .\]
On se sert de l'hypothèse sur $\U$ et on peut numéroter un sous-ensemble cofinal $\V$ de $\U$ (quitte \` a enlever des \' el\' ements du $\V$ de l'\' enonc\' e) $\V = \{ V_n \ | \ n \in \N \}$ en respectant l'ordre: on demande $V_i \leq V_j \Rightarrow i \leq j$. On réécrit alors \[ C^\infty(X,\Z) = \varinjlim_{V \in \V} C_{(V)}(X,\Z) = \varinjlim_n \sum_{k \leq n} C_{(V_k)}(X,\Z) .\] On construit par récurrence sur $n \in \N$ une base de $C_n := \sum_{k \leq n} C_{(V_k)}(X,\Z)$. L'étape d'initiation $n=0$ consiste simplement à choisir une base $(b_0, \dots , b_{i_0})$ du $\Z$-module libre $C_{(V_0)}(X,\Z)$. Supposons que l'on a construit une base $(b_0, \dots, b_{i_n})$ de $C_n$. Parce que $C_n \cap C_{(V_{n+1})}(X,\Z)$ est un $C_{(V_{n+1}')}(X,\Z)$ pour un certain $V_{n+1}' \leq V_{n+1}$ ($V_{n+1}'$ non n\' ecessairement dans $\V$), $C_{(V_{n+1})}(X,\Z) / \big( C_n \cap C_{(V_{n+1})}(X,\Z) \big)$ est sans torsion. Alors $C_n$ est un facteur direct du $\Z$-module libre $C_{n+1}$. On peut alors compléter $(b_0, \dots, b_{i_n})$ en une base $(b_0, \dots, b_{i_{n+1}})$ de $C_{n+1}$. La récurrence est alors prouvée et le résultat suit. $\hfill\Box$

\section{Appendice: de la non nécessité de $R$ algébriquement clos} \label{parRnonalgclos}

Soient $R$ un corps et $R^\al$ une clôture algébrique fixée de $R$. Soit $G$ un groupe topologique. On se restreint à la catégorie des représentations lisses de $G$ à coefficients dans $R$ (respectivement $R^\al$).

\begin{lem} \label{Ralgclos}
Soit $\pi$ une représentation de $G$ définie sur $R$. Supposons que $\pi \otimes_R R^\al$ soit une représentation irréductible. Alors $\pi$ est irréductible. 
\end{lem}
\textsf{Preuve :} \\ 
Toute suite exacte courte \[ 0 \to \rho \to \pi \to \pi/\rho \to 0 \] induit la suite exacte courte \[ 0 \to \rho \otimes R^\al \to \pi \otimes R^\al \to \pi/\rho \otimes R^\al \to 0 \] dont on sait que l'un des termes extrémaux est trivial. $\hfill\Box$

\section{Appendice: sur l'ordre $<_J$}

On rappelle la définition de $<_J$: pour $w,w' \in W^J$, on note $w <_J w'$ s'il existe $s_1, \dots, s_r \in S$ tels que $w^{(i)} = (s_i s_{i-1} \cdots s_1 w)^J$ vérifie $l(w^{(i)}) > l(w^{(i-1)})$ pour tout $1 \leq i \leq r$ et $w^{(r)} = w'$. \\
On établit la caractérisation suivante de $<_J$, qui dit en particulier que $<_J$ est un raffinement de la restriction à $W^J$ de l'ordre fort dans un groupe de Coxeter fini.

\begin{prop}
Soient $w,w' \in W^J$. On a $w<_J w'$ si et seulement si il existe $s_1, \dots, s_r \in S$ tels que $w^{(i)} = s_i \cdots s_1 w$ soit un élément de $W^J$ de longueur $l(w)+i$ pour tout $1 \leq i \leq r$ et $w^{(r)} = w'$.
\end{prop}
\textsf{Remarque :}
Le cas $r=1$ est déjà présent dans \cite{Grk09}, Lemma 1.4.(b). \\
\textsf{Preuve :} \\
Le sens $(\Leftarrow)$ est immédiat puisque $w \in W^J$ implique $w^J =w$. Supposons donc $w <_J w'$ et prenons $s_1, \dots, s_r$ comme dans la définition de $<_J$. Prouvons d'abord $s_1 w \in W^J$ avec $l(s_1w) = l(w)+1$. On a \[ l(w) < l((s_1w)^J) \leq l(s_1 w) \leq l(w)+1 ,\] où la première inégalité suite de la définition de $<_J$ et la deuxième de celle de $W^J$. Mais alors, on a $l((s_1w)^J)=l(s_1w)=l(w)+1$. Cela dit en particulier que $s_1w$ est de longueur minimale dans $s_1 w W_J$ et on a $(s_1w)^J = s_1 w \in W^J$. Par une récurrence immédiate, $w^{(i)} = s_i \cdots s_1 w$ est un élément de $W^J$ de longueur $l(w)+i$ pour tout $1 \leq i \leq r$. Le résultat est prouvé. $\hfill\Box$

\section{Appendice: de l'irréductibilité de la Steinberg généralisée dans le cas fini}

Soit $R$ un corps algébriquement clos de caractéristique $p$. Le travail effectué nous permet de découvrir ou redécouvrir quelques résultats sur les Steinberg généralisées pour un groupe réductif fini $\overline{G}$. 

\begin{prop}
La plus grande sous-$\overline{G}$-représentation irréductible de $\overline{\St}_{\overline{J}} R$ est $V_J$.
\end{prop}
\textsf{Remarque :} 
En particulier, si $\overline{\St}_{\overline{J}} R$ est irréductible, alors on a $V_J = \overline{\St}_{\overline{J}} R$. \\
\textsf{Preuve :} \\
On a une inclusion $\overline{\St}_{\overline{J}} R \subseteq \St_J R$, et on sait que le $K$-socle de $\St_J R$ est irréductible, égal à $V_J$: $V_J$ est donc aussi le $K$-socle de $\overline{\St}_{\overline{J}} R$. Comme $K(1)$ agit trivialement sur $V_J \subseteq \overline{\St}_{\overline{J}} R$ et $\overline{\St}_{\overline{J}} R$, le résultat se traduit en termes de $\overline{G}$-représentations. $\hfill\Box$

\begin{prop}
Supposons $\Phibarred$ irréductible et $\overline{J} \notin \{\varnothing, \overline{\Delta} \}$. Alors $\overline{\St}_{\overline{J}} R$ n'est pas irréductible.
\end{prop}
\textsf{Remarque :} 
$\St_{\overline{\Delta}} R = \id$ est bien sûr irréductible; quant à la Steinberg $\St_\varnothing R$, en utilisant \cite{CabEng04}, Theorem 6.10, Theorem 6.12 et Definition 6.13, on voit qu'elle est aussi irréductible. \\
\textsf{Preuve :} \\
Par \cite{CabEng04}, Theorem 6.12, si $V$ est une $\overline{G}$-représentation irréductible alors son espace de $\overline{U}$-invariants est de dimension $1$. De ce fait, si $V$ est une représentation avec $\dim \, V^{\overline{U}} \geq 2$, alors $V$ n'est pas irréductible. Par les propositions \ref{IinvVQ} et \ref{mainGrk1}, et le d\' ebut de la preuve du corollaire \ref{Stadm}, $(\overline{\St}_{\overline{J}} R)^{\overline{U}} = (\overline{\St}_{\overline{J}} R)^{\overline{B}}$ est un $R$-espace vectoriel de dimension $|W^J_\pr|$. Il s'agit d'examiner la cardinalité de $W^J_\pr$ et le résultat suit par le lemme \ref{WJpr}. $\hfill\Box$

\begin{lem} \label{WJpr}
Supposons $\Phired$ irréductible. Alors on a $|W^J_\pr| \geq 1$, avec égalité si et seulement si $J$ est $\varnothing$ ou $\Delta$.
\end{lem}
\textsf{Preuve :} \\
Rappelons la définition suivante de $W^J_\pr$: \[ W^J_\pr = \{ w \in W \ | \ \forall \alpha \in J, \ l(w s_\alpha) > l(w) \ ; \ \forall \beta \in \Delta \smallsetminus J , \ l(w s_\beta) < l(w) \} .\]
Notons $w_{\Delta \smallsetminus J}$ l'élément le plus long de $W_{\Delta \smallsetminus J}$. C'est un élément de $W^J_\pr$, et de ce fait on a la minoration voulue. Il reste à déterminer le cas d'égalité. D'abord, on remarque $W^\varnothing_\pr = \{w_\Delta\}$ et $W^\Delta_\pr = \{1\}$, de sorte qu'on veut maintenant montrer que si $J$ n'est pas $\varnothing$ ou $\Delta$, alors $W^J_\pr$ contient un autre élément que $w_{\Delta \smallsetminus J}$. \\
Supposons $J \neq \varnothing, \Delta$. On cherche un élément $w \in W_J \smallsetminus \{1\}$ vérifiant $l(w w_{\Delta \smallsetminus J} s_\alpha) > l(w w_{\Delta \smallsetminus J})$ pour tout $\alpha \in J$ . Parce que $\Phired$ est irréductible, on peut choisir $\beta \in J$ tel que $(\Delta \smallsetminus J ) \cup \{\beta\}$ engendre un sous-système de $\Phired$ avec au plus autant de composantes irréductibles que celui engendré par $\Delta \smallsetminus J$. Montrons que $s_\beta$ est l'élément $w \in W^J \smallsetminus \{1\}$ cherché. \\ En effet, remarquons d'abord que l'on a \[ l(w_{\Delta \smallsetminus J}) = l (s_\beta w_{\Delta \smallsetminus J} s_\alpha) < l(s_\beta w_{\Delta \smallsetminus J}) \] pour tout $\alpha \in \Delta \smallsetminus J$. Ensuite, supposons qu'il existe un élément $\gamma \in J$ avec $l(s_\beta w_{\Delta \smallsetminus J} s_\gamma) < l(s_\beta w_{\Delta \smallsetminus J})$, alors cela veut dire que $s_\beta w_{\Delta \smallsetminus J}$ possède une écriture qui se termine par $s_\gamma$, disons $w' s_\gamma$ avec $l(w') = l(w_{\Delta \smallsetminus J})$ et $w'$ ne se terminant pas par $s_\gamma$. Maintenant on a $W_J w_{\Delta \smallsetminus J} W_J = W_J w' W_J$, ce qui force $w' = w_{\Delta \smallsetminus J}$. Cela implique $s_\gamma = w_{\Delta \smallsetminus J} s_\beta w_{\Delta \smallsetminus J} \in \W_{(\Delta \smallsetminus J) \cup \{\beta\}}$ et donc $\gamma = \beta$. Mais dans ce cas-là, c'est que $s_\beta w_{\Delta \smallsetminus J}$ est l'élément le plus long de $W_{(\Delta \smallsetminus J) \cup \{\beta\}}$. Mais on sait aussi que la longueur de $w_{(\Delta \smallsetminus J) \cup \{\beta\}}$ est égale à $|\Phi^+_{(\Delta \smallsetminus J) \cup \{\beta\}}|$ (\cite{Hum92}, I.4.8), et on a alors \[ |\Phi^+_{(\Delta \smallsetminus J) \cup \{\beta\}}| = |\Phi^+_{(\Delta \smallsetminus J)}| + 1 .\] Il s'ensuit \[ \Phi^+_{(\Delta \smallsetminus J) \cup \{\beta\}} = \Phi^+_{(\Delta \smallsetminus J)} \coprod \{\beta\} \] et cela contredit le fait que $\Phi_{(\Delta \smallsetminus J) \cup \{\beta\}}$ a moins (éventuellement le même nombre) de composantes irréductibles que $\Phi_{\Delta \smallsetminus J}$. C'est absurde, et le résultat est prouvé. $\hfill\Box$ \\
\textsf{Remarque :} 
Marie-France Vignéras nous fait remarquer que l'élément $z^J = w_\Delta w_J \in W^J$ de la proposition \ref{mainGrk2} convient aussi en tant qu'élément distinct de $w_{\Delta \smallsetminus J}$ dans $W^J_\pr$ (voir aussi \cite{Grk09}, Lemma 1.4.(e)). En effet, il est de longueur maximale $|\Phi^+| - |\Phi_J^+| \neq |\Phi^+_{\Delta \smallsetminus J}|$ et est donc différent de $w_{\Delta \smallsetminus J}$. Et sa longueur excède aussi celle de tout élément de $W^{J'}$ pour $J' \supsetneq J$ et il est donc primitif.

\section{Appendice: représentations de Steinberg généralisées pour le groupe dérivé}

Dans cette section uniquement, on distinguera le groupe r\' eductif $\underline{G}$ d\' efini sur $F$ de ses $F$-points $G = \underline{G}(F)$. De m\^ eme, $B$ et $P$ seront respectivement les $F$-points de $\underline{B}$ et $\underline{P}$. \\
On note $\underline{D}(G)$ le groupe dérivé de $G$, c'est-à-dire le faisceau fppf des commutateurs de $\underline{G}$. C'est un groupe semi-simple (voir \cite{Dem11}, Théorème 6.2.1.(iv)) et on notera $D(G)$ pour son groupe des $F$-points. Le but de ce paragraphe est de comparer les représentations de Steinberg généralisées pour $G$ et pour $D(G)$. \\
Soient $R$ un corps de caractéristique $p$ et $J \subseteq \Delta$. Pour les distinguer, on notera $\St^{(G)}_J R$ et $\St^{(D(G))}_J R$ la représentation de Steinberg généralisée respectivement pour $G$ et $D(G)$, par rapport à $J$.

\begin{prop} \label{StD(G)}
La restriction de $\St_J^{(G)} R$ à $D(G)$ est isomorphe à $\St^{(D(G))}_J R$.
\end{prop}

A partir de là, on peut déduire l'irréductibilité de $\St^{(D(G))}_J R$ en utilisant toute la machinerie de ce papier pour $\underline{D}(G)$, y compris l'appendice \ref{parRnonalgclos}. \\

\textsf{Preuve de la proposition \ref{StD(G)} :} \\
Le groupe $\underline{B} \cap \underline{D}(G)$ est un parabolique minimal de $\underline{D}(G)$ (voir \cite{Dem11}, Proposition 6.2.8.(ii)), et on appelle standard un parabolique de $\underline{D}(G)$ contenant $\underline{B} \cap \underline{D}(G)$. La flèche $\underline{P} \mapsto \underline{P} \cap \underline{D}(G)$ est une bijection entre l'ensemble des paraboliques standards de $\underline{G}$ et celui des paraboliques standards de $\underline{D}(G)$. Soient $\underline{P}$ un parabolique standard de $\underline{G}$ et $J \subseteq \Delta$ l'ensemble v\' erifiant $\underline{P} = \underline{P}_J$. On veut voir que l'injection $P \cap D(G) \backslash D(G) \hookrightarrow P \backslash G$ est une bijection. Pour cela, utilisons la d\' ecomposition de Bruhat pour $D(G)$:
\[ P \cap D(G) \backslash D(G) = \coprod_{w \in W^J} P \cap D(G) \backslash (P \cap D(G)) w^{-1} (B \cap D(G)) ,\] 
o\` u on a relev\' e chaque $w \in W^J$ en un \' el\' ement de $D(G)$. En notant que $U_w = U \cap w U^- w^{-1}$ est inclus dans $D(G)$, elle se r\' e\' ecrit encore
\begin{equation} \label{BruDG} P \cap D(G) \backslash D(G) = \coprod_{w \in W^J} P \cap D(G) \backslash (P \cap D(G)) w^{-1} U_w .\end{equation}
De m\^eme, pour $G$ on \' ecrit (en gardant les m\^ emes rel\` evements pour $W^J$) 
\begin{equation} \label{BruG} P \backslash G = \coprod_{w \in W^J} P \backslash P w^{-1} B = \coprod_{w \in W^J} P \backslash P w^{-1} U_w \end{equation}
La comparaison de (\ref{BruG}) et de (\ref{BruDG}) nous donne que $P \cap D(G) \backslash D(G) \hookrightarrow P \backslash G$ est surjective, et donc bijective. \\
De ce fait, la restriction à $D(G)$ de l'induite $\Ind_P^G \, \id$ est $\Ind_{P \cap D(G)}^{D(G)} \, \id$. Par définition (\ref{Stdef}), la restriction à $D(G)$ de $\St_J^{(G)} R$ s'identifie à $\St_J^{(D(G))} R$. \hfill$\Box$ \\

\chapter{Des représentations modulo $p$ de $\mathrm{GL}(m,D)$, $D$ algèbre à division sur un corps local} \label{Ly12}

\section{Introduction}

L'histoire de l'étude des représentations lisses modulo $p$ de $\GL(m,D)$ commence avec Barthel-Livné en 1994-95 avec le cas déployé $D=F$ et $m=2$ (voir \cite{BarLiv94}, \cite{BarLiv95}). Récemment, Herzig a généralisé cette méthode au cas $D=F$, $m \geq 2$ quelconque (\cite{Her11b}). On continue ici la classification en abordant le cas non nécessairement déployé. C'est la première fois qu'un exemple de groupe non quasi-déployé est traité. En effet, Abe a étendu le travail de Herzig au cas de tout groupe réductif connexe déployé (\cite{Abe11}), et Abdellatif au cas d'un groupe quasi-déployé de rang $1$ (\cite{Abd11}). \\
Soient $G = \GL(m,D)$ et $P$ le parabolique standard de type $(m_1, \dots, m_r)$ avec $\sum_i m_i = m$. On commence par établir un critère d'irréductibilité. \\

Avant de l'énoncer, rappelons brièvement la définition d'une représentation de Steinberg généralisée. Il existe un morphisme de groupes $G \mathop{\to}\limits^{\det_G} F^\times$ qui coïncide avec la norme réduite $\Nrm$ pour $m=1$ (voir paragraphe \ref{pardeterminants}). La représentation de Steinberg relativement à $P$ est alors \[ \St_P \id := \frac{\Ind_P^G \, \id}{\sum_{P' \gneq P} \Ind_{P'}^G \, \id} ,\] où $P'$ parcourt l'ensemble des paraboliques contenant $P$; on définit aussi, pour $\rho : D^\times \mathop{\to}\limits^{\Nrm} F^\times \mathop{\to}\limits^{\rho^0} \Fpbar^\times$ un caractère lisse se factorisant par $\Nrm$, \[ \St_P \rho^0 := (\rho^0 \circ \det\nolimits_G) \otimes \St_P \id .\]
Les résultats que l'on obtient nécessitent l'hypothèse suivante:
\begin{itemize}
\item[(a)] $D$ est de degré $d^2$ sur $F$ avec $d$ un nombre premier, ou égal à $1$;
\item[(b)] ou $m$ est inférieur ou égal à $3$.
\end{itemize}
On la supposera donc satisfaite dans le reste de l'introduction. On espère que l'hypothèse est superflue mais on ne sait pas l'enlever pour le moment. \\
Aussi, les résultats sont conditionnels à une conjecture (voir paragraphe \ref{LKglm}) sur des valeurs de la transformée de Satake que l'on sait être vraie sous l'hypothèse (b). 

\begin{theo} \label{ABmain1}
Soit, pour tout $1 \leq i \leq r$, $\sigma_i$ une représentation irréductible admissible de $\GL(m_i,D)$ qui soit dans l'un des cas suivants:
\begin{itemize}
\item[(a)] $\sigma_i$ est supercuspidale\footnote{une représentation est dite \textit{supercuspidale} si elle n'est pas sous-quotient d'une induite parabolique propre d'une représentation irréductible admissible} avec $m_i > 1$;
\item[(b)] $\sigma_i$ est une représentation de $D^\times$ (i.e. $m_i = 1$) avec $\dim \, \sigma_i > 1$;
\item[(c)] $\sigma_i \simeq \St_{Q_i} \rho_i^0$ pour un parabolique standard $Q_i \leq \GL(m_i,D)$ et un caractère $\rho_i : D^\times \xrightarrow{\Nrm} F^\times \xrightarrow{\rho_i^0} \Fpbar^\times$.
\end{itemize}
Supposons $\rho_i \neq \rho_{i+1}$ s'il y a deux blocs adjacents dans le cas (c). Alors $\Ind_{P}^G(\sigma_1 \otimes \dots \otimes \sigma_r)$ est irréductible admissible.
\end{theo}

Dans un second temps, on prouve le théorème de classification suivant.

\begin{theo} \label{ABmain2}
Soit $\pi$ une représentation irréductible admissible de $G$. Alors il existe un parabolique standard $P$ de Levi $\prod\limits_{i=1}^m \GL(m_i,D)$ avec $\sum_i m_i = m$ et des représentations irréductibles admissibles $\sigma_i$ de $\GL(m_i,D)$ tels que l'on ait $\pi \simeq \Ind_P^G (\sigma_1 \otimes \cdots \otimes \sigma_r)$. De plus, pour tout $i$, on est dans l'un des cas (a), (b) ou (c) du théorème \ref{ABmain1}.
Aussi, s'il y a deux blocs adjacents dans le cas (c), alors on a $\rho_i \neq \rho_{i+1}$. \\ Enfin, si $\Ind_P^G (\sigma_1 \otimes \cdots \otimes \sigma_r)$ et $\Ind_{P'}^G (\sigma_1' \otimes \cdots \otimes \sigma_{r'}')$ sont deux telles représentations irréductibles admissibles isomorphes, alors on a $r=r'$, $M_i =M_i'$ et $\sigma_i = \sigma_i'$ pour tout $i$.
\end{theo}

En fait, comme dans \cite{Her11b}, on va effectuer le travail avec des blocs élémentaires qui seront des représentations supersingulières et ce ne sera qu'à la fin qu'on saura qu'une représentation irréductible admissible est supersingulière si et seulement si elle est supercuspidale. \\ 

On remarque enfin que tout ceci est compatible avec les calculs explicites de l'auteur sur le cas $\GL(2,D)$ (voir chapitre \ref{Ly11b}).

\section{Notations et généralités} \label{notations}

Soient $p$ un nombre premier et $\Fpbar$ une clôture algébrique fixée de $\Fp$; tout corps fini de caractéristique $p$ sera vu comme un sous-corps de $\Fpbar$. \\

L'objet de cette section est de donner les notations et de présenter des faits généraux bien connus sur la théorie des représentations ou bien sur $\GL(m,D)$.

\subsection{Représentations et algèbres de Hecke} \label{generalreps}

On pourra se reporter au paragraphe 2 de \cite{BarLiv94}. Dans tout ce qui suit, toute représentation considérée sera lisse (et on oubliera souvent de le mentionner), à coefficients dans un anneau commutatif $R$; en pratique, dans la suite $R$ sera $\Fpbar$ ou $\C$. Soient $G$ un groupe topologique et $H \leq G$ un sous-groupe fermé. On notera $\ind_H^G$ le foncteur d'induction compacte lisse et $\Ind_H^G$ celui d'induction lisse. L'action sur une induite se fera par translation à droite, à savoir $g. f : x \mapsto f(xg)$. On suppose $H$ ouvert. Soit $V$ une représentation de $H$. Pour tous $g \in G$ et $v \in V$, on définit $[g,v]$ comme étant la fonction de $\ind_H^G \, V$ à support dans $H g^{-1}$ et prenant pour valeur $v$ en $g^{-1}$. En particulier, si $h$ est un élément de $H$, on a $[gh,v] = [g,hv]$ et $g[1,v]=[g,v]$. \\ 
Pour toutes représentations $V_1$ et $V_2$ de $H$, on définit l'espace d'entrelacements \[ \Hh_R(G,H,V_1,V_2) := \Hom_G \big( \ind_H^G \, V_1, \ind_H^G \, V_2 \big) .\] On se permettra aussi de noter les algèbres de Hecke \[ \Hh_R(G,H,V) := \Hh_R(G,H,V,V), \quad \Hh_R(G,H) := \Hh_R(G,H,\id) ,\] où $\id$ désigne la représentation triviale de $H$. \\ 
Lorsque $H$ est ouvert dans $G$, par réciprocité de Frobenius, on a \[ \Hh_R(G,H,V_1,V_2) \simeq \Hom_H(V_1, \ind_H^G \, V_2) .\] 
Supposons que, pour tout $g \in G$, la double classe $HgH$ est union finie de classes à gauche (ou à droite), et que $V_1$ et $V_2$ sont finiment engendrées en tant que $H$-représentations. Alors on peut (voir \cite{HenVig11}, section 2.2) exhiber l'isomorphisme 
\[\begin{array}{ccc} \Hom_H(V_1, \ind_H^G \, V_2) & \simeq & \left\{ \begin{tabular}{r|l} \multirow{2}{*}{$G \xrightarrow{f} \Hom_{\Fpbar}(V_1,V_2)$} & $f(h_2 gh_1) = h_2 f(g) h_1$ \\ & $H \backslash \Supp \, f / H$ fini \end{tabular} \right\} \\ \big( v \mapsto (g \mapsto f(g)(v) )\big) & \mapsfrom & f \\ f & \mapsto & \big( g \mapsto (v_1 \mapsto f(v_1)(g)) \big) . \end{array} \] 
On va souvent assimiler les opérateurs de $\Hh_R(G,H,V_1,V_2)$ à des fonctions sur $G$ à travers cet isomorphisme. \\
Pour tous $V_1, V_2, V_3$ et tous $f_1 \in \Hh_R(G,H,V_1,V_2)$, $f_2 \in \Hh_R(G,H,V_2,V_3)$, on a la convolée $f_2 * f_1 \in \Hh_R(G,H,V_1,V_3)$ définie par 
\begin{equation} \label{GLmdefconvol} f_2 * f_1  : g \mapsto \sum_{x \in G/H} f_2(x) f_1(x^{-1} g) .\end{equation}
Aussi, quand on voit les opérateurs de $\Hh_R(G,H,V)$ comme fonctions sur $G$, la structure multiplicative est donnée par la convolution, qui est bien compatible à la multiplication par composition de la première définition (voir \cite{BarLiv94}, Proposition 5).

\subsection{Des notations pour $\GL(m,D)$}

Soit $F$ un corps local non archimédien localement compact à corps résiduel fini $k_F$ de caractéristique $p$; on notera $\Oo_F$ son anneau de valuation, $\m_F$ son idéal maximal et $q = p^f$ le cardinal de $k_F$.\\
Soit $D$ une algèbre à division de centre $F$: elle est de degré $d^2$ sur $F$ pour un certain entier $d \geq 1$, et d'invariant de Brauer $a_D/d$ avec $a_D$ un entier de $[1,d[$ premier à $d$. Alors $D$ possède un sous-corps commutatif maximal $E$ tel que $E/F$ est une extension non ramifiée de degré $d$, unique \` a conjugaison dans $D$ pr\` es (voir \cite{Ser67}, Appendix); on note $\Oo_E$ son anneau d'entiers. On peut munir $D$ d'une valuation discrète  $v_D$; on fixe une uniformisante $\varpi_D$ de $D$ (et on normalise à $v_D(\varpi_D) = 1$) de sorte que $\varpi := \varpi_D^d$ est une uniformisante de $F$. On note $\Oo_D$ l'anneau de valuation de $D$, $\m_D = (\varpi_D)$ son idéal maximal et $k_D = \Oo_D/\m_D$ son corps résiduel. On note $x \mapsto \overline{x}$ l'application de réduction de $\Oo_D$ à $k_D$. L'action de l'uniformisante sur le corps résiduel est $\varpi_D \overline{x} \varpi_D^{-1} = \overline{x}^{q^{a_D}}$ pour tout $x \in \Oo_D$. On remarque que $k_D$ est aussi le corps résiduel de $E$. Ainsi, on note $[ \phantom{x}]$ l'application de Teichmüller $k_D^\times \to \Oo_E^\times$, que l'on prolonge en $k_D \to \Oo_E$ en envoyant $0$ sur $0$. \\ Pour plus de précisions on pourra aller voir le chapitre 17 de \cite{Pie82} ou le chapitre 14 de \cite{Rei75}. \\

Soit $m \geq 1$ un entier. Soient $B$ le parabolique minimal de $G:=\GL(m,D)$ composé des triangulaires supérieures, $A$ sa composante de Levi diagonale et $U$ son radical unipotent, composé des matrices unipotentes supérieures. On définit aussi $B^- = A U^-$ le parabolique minimal opposé. De manière générale, si $P = MN$ désigne un parabolique de radical unipotent $N$, on note $P^- = M N^-$ le parabolique opposé. On dira que $P$ est un parabolique \textit{standard} s'il contient $B$ et \textit{antistandard} s'il contient $B^-$. \\ 
Soit $T$ le centre $Z(A)$ de $A$: c'est un tore $F$-déployé maximal (parmi les tores $F$-déployés) de $G$. On associe à $T$ le groupe de Weyl étendu $W_e$ et le groupe de Weyl fini $W$. La paire $(B,T)$ fournit un système de racines positives, ainsi qu'un sous-ensemble $\Delta$ de racines simples. On note $\{s_\alpha \ | \ \alpha \in \Delta\} = S_\Delta$, qui constitue un ensemble générateur de réflexions dans $W$. On dénote $w_0$ l'élément le plus long de $W$ pour la longueur définie par $S_\Delta$. \\ 
Pour un parabolique standard $P$ de composante de Levi $M$, le tore $T$ définit un sous-ensemble de racines simples $\Delta_M \subseteq \Delta$. Lorsque l'on veut éviter d'alourdir les notations (et d'introduire le Levi standard quand ce n'est pas nécessaire), on se permettra de noter $\Delta_P$ au lieu de $\Delta_M$. On note $W_M$ le sous-groupe de $W$ engendré par $\Delta_M$. \\
On note respectivement $X^*(T)$ et $X_*(T)$ les groupes des caractères et cocaractères algébriques de $T$ à valeurs dans $F^\times$. Ils sont en dualité par \[ \langle \phantom{x}, \phantom{x}\rangle : X_*(T) \times X^*(T) \to \Z .\] On dit qu'un cocaractère est \textit{antidominant} s'il appartient à l'ensemble \[X_*(T)^- := \{ x \in X_*(T) \ | \ \langle x, \alpha \rangle \leq 0 \textrm{ pour tout $\alpha \in \Delta$}\} .\] Pour $J \subseteq \Delta$, on dira que $x \in X_*(T)$ est \textit{antidominant hors de} $J$ si on a $\langle x,\alpha \rangle \leq 0$ pour tout $x \notin J$. Les mêmes notions existent avec des inégalités strictes et on parlera de \textit{stricte antidominance} par exemple; on notera $X_*(T)^{--}$ le monoïde correspondant. Aussi, il y a la notion opposée de \textit{dominance} définie par $\langle x , \alpha \rangle \geq 0$ pour tout $\alpha \in \Delta$. \\
Pour $\lambda \in X_*(T)^-$, on définit $P_\lambda$ comme étant le parabolique antistandard dont on va fixer le Levi $M_\lambda$ défini par le système de racines \[ \Delta_\lambda := \{ \alpha \in \Delta \ | \ \langle \lambda, \alpha \rangle = 0 \} ;\] $P_{-\lambda}$ désignera son parabolique opposé. \\

Soient $K := \GL(m,\Oo_D)$ le compact maximal de $G$, $K(1)$ son pro-$p$-radical égal à $1 + \varpi_D M_m(\Oo_D)$ et $\overline{G} := K/K(1) \simeq \GL(m,k_D)$. Pour $H$ un sous-groupe de $G$, $\overline{H}$ désignera le sous-groupe de $\overline{G}$ correspondant: \[ \overline{H} := (H \cap K)/(H \cap K(1)) ;\] par exemple $\overline{B}$ est le Borel de $\overline{G}$ composé des matrices triangulaires supérieures, et que l'on appellera standard. \\ Aussi, on note $K(h)= 1+\varpi_D^h M_m(\Oo_D)$ pour $h \geq 1$; et $\overline{\phantom{x}}^{(h)}$ désigne la réduction modulo $K(h)$.

\subsection{Des sous-groupes à un paramètre dans $A$} \label{oneparamsbgp}

On note $\Lambda := A/A \cap K$; c'est un groupe abélien libre de type fini de rang $m$. L'élévation à la puissance $d$ induit un isomorphisme $\Lambda \xrightarrow{\sim} T/T \cap K$, d'où on déduit l'isomorphisme de groupes 
\begin{equation} \label{versSclassic} \begin{array}{ccccc} \Lambda & \xrightarrow{\sim} & T/T \cap K & \xrightarrow{\sim} & X_*(T) \\ & & x(\varpi) & \mapsfrom & x \end{array} .\end{equation} 
Notons $A_\Lambda$ le sous-groupe de $A$ à coefficients diagonaux dans $\varpi_D^\Z$: c'est un système de représentants de $\Lambda$ dans $A$. La composition de la flèche précédente avec la projection de $A_\Lambda$ dans $\Lambda$ induit alors un isomorphisme de groupes \[ \iota_\Lambda : A_\Lambda \xrightarrow{\sim} \Lambda \xrightarrow{\sim} T/T \cap K \xrightarrow{\sim} X_*(T) .\] On remarque que toutes ces flèches sont $W$-équivariantes. \\ Pour un cocaractère $x \in X_*(T)$, on définit alors une application multiplicative $\widetilde{x} : F^\times \to A$ par \[ \widetilde{x}(\varpi) = \iota_\Lambda^{-1}(x) , \quad \widetilde{x}|_{\Oo_F^\times} = x|_{\Oo_F^\times} .\]
En particulier, on remarque que cette définition fait de $x \mapsto \widetilde{x}$ un morphisme de groupes de \[ \widetilde{\cdot} : X_*(T) \to \Hom_{\gr}(F^\times,A) .\] Enfin, on note $A_\Lambda^-$ le sous-monoïde de $A_\Lambda$ constitué des éléments antidominants, c'est-à-dire dont l'action par conjugaison contracte $U^-$. \\

Soit $P$ un parabolique standard de Levi $M$. On définit: \[ \widetilde{Z}_M = \iota_\Lambda^{-1} (X_*(Z(M)) \leq A_\Lambda ,\] où $X_*(Z(M))$ désigne le sous-groupe de $X_*(T)$ des cocaractères à valeurs dans le centre $Z(M)$ de $M$. \\ Pour des cocaractères $\lambda$ et $\mu$, on note $\lambda \leq \mu$ si $\mu-\lambda$ est somme de coracines positives. Lorsque $M$ est un Levi standard (c'est-à-dire contenant $A$), la paire $(B \cap M,T)$ définit la paire radicielle $(\Phi_M,\Phi_M^\vee)$ avec l'ordre $\leq_M$: on notera $\lambda \leq_M \mu$ lorsqu'on peut prendre les coracines dans $\Phi_M^\vee$.

\subsection{Représentations irréductibles de $\GL(m,k_D)$} \label{paragpoids}

On va suivre l'exposition de \cite{HenVig11b}. D'ailleurs, on se reportera à \cite{Cur70} et à \cite{HenVig11b} pour les preuves. \\

On commence par remarquer que $\overline{A}$ est le groupe des $k_F$-points du tore diagonal non $k_F$-déployé $(\Res_{k_D/k_F} \Gm)^m$: en d'autres termes, $\overline{A}$ est simplement constitué des matrices diagonales de $\GL(m,k_D)$. Soit $\chi: \overline{A} \to \Fpbar^\times$ un caractère de $\overline{A}$. 
Pour $\alpha \in \Delta$, on définit le conjugué de $\chi$ suivant: $\chi^\alpha = \chi (s_\alpha \cdot s_\alpha^{-1})$, o\` u on confond $s_\alpha \in W$ avec l'image dans $\overline{K}$ d'un rel\` evement dans $K$. On remarquera que comme $k_D^\times$ est ab\' elien, $\chi^\alpha$ ne d\' epend pas du rel\` evement $s_\alpha$ de $\alpha$. On note alors 
\begin{equation} \label{defDeltachi} \Delta_\chi = \{ \alpha \in \Delta \ | \ \chi^\alpha = \chi \} \end{equation} 
le sous-ensemble de $\Delta$ caractérisant la régularité du caractère $\chi$. On v\' erifie que notre d\' efinition de $\Delta_\chi$ est compatible avec celle de \cite{HenVig11}.

\begin{lem}
Soit $\chi : \overline{A} \to \Fpbar^\times$ un caract\` ere. On a \[ \Delta_\chi = \big\{ \alpha \in \Delta \ \big| \ \chi \textrm{ est trivial sur } \langle \overline{U}_\alpha, \overline{U}_{-\alpha} \rangle \cap \overline{A} \big\} .\]
\end{lem}
\textsf{Preuve :} \\
Notons $\chi = \chi_1 \otimes \dots \otimes \chi_m$ o\` u chaque $\chi_j$ est un caract\` ere $k_D^\times \to \Fpbar^\times$. Soient $\alpha \in \Delta$ et $i \in [1,m-1]$ l'entier tel que $\alpha$ est la racine entre $i$ et $i+1$. Alors $\chi^\alpha$ est le caract\` ere $\chi_1 \otimes \dots \otimes \chi_{i-1} \otimes \chi_{i+1} \otimes \chi_i \otimes \chi_{i+2} \otimes \dots \otimes \chi_m$, de sorte que $\alpha$ est dans $\Delta_\chi$ si et seulement si on a $\chi_i = \chi_{i+1}$. \\
Mais le groupe engendr\' e par $\overline{U}_\alpha$ et $\overline{U}_{-\alpha}$ est le sous-groupe de $\overline{K}$ des matrices diagonales par blocs de type $1 \times \dots \times 1 \times 2 \times 1 \times \dots \times 1$ valant $1$ sur les blocs de taille $1$ et \' etant une matrice de $\SL(2,k_D)$ sur le bloc de taille $2$. De ce fait, $\langle \overline{U}_\alpha, \overline{U}_{-\alpha} \rangle \cap \overline{A}$ est constitu\' e des matrices $\Diag(1, \dots, 1, x, x^{-1}, 1, \dots, 1)$ pour $x \in k_D^\times$, les coefficients \' eventuellement non \' egaux \` a $1$ \' etant en $i$ et $i+1$. Alors $\chi$ est trivial sur $\langle \overline{U}_\alpha, \overline{U}_{-\alpha} \rangle \cap \overline{A}$ si et seulement si $\chi_i$ est \' egal \` a $\chi_{i+1}$. Cela donne la compatibilit\' e voulue. \hfill$\Box$ \\

Soit $V$ une $\overline{G}$-représentation irréductible. Soit $\overline{P} = \overline{M} \, \overline{N}$ un parabolique standard de $\overline{G}$. On sait que $V^{\overline{N}}$ est une $\overline{M}$-représentation irréductible. En particulier, pour $\overline{P} = \overline{B}$, $V^{\overline{U}}$ et un certain sous-ensemble de $\Delta_{V^{\overline{U}}}$ suffisent à déterminer $V$ à isomorphisme près: c'est l'énoncé suivant.

\begin{prop} \label{irredGbar}
$\phantom{,}$
\begin{itemize}
\item[(i)] Soit $V$ une $\overline{G}$-représentation irréductible. La droite $V^{\overline{U}}$ est stabilisée par un parabolique standard $\overline{P}_V$ associé à un sous-ensemble $\Delta_V$ de $\Delta_{V^{\overline{U}}}$.
\item[(ii)] On a une bijection ensembliste : \[ \begin{array}{cccc} \Par_{\overline{B}} : & V & \mapsto & \big( V^{\overline{U}}, \Delta_V \big) \\ & \left\{ \begin{array}{c} \textrm{classes d'isomorphismes de} \\ \overline{G} \textrm{-représentations irréductibles }\end{array} \right\} & \xrightarrow{\sim} & \left\{ \begin{array}{c|l} (\chi,\Delta_V) & \chi: \overline{A} \to \Fpbar^\times \textrm{ caractère} , \\ & \Delta_V \subseteq \Delta_\chi \end{array} \right\} \end{array} .\]
\end{itemize}
\end{prop}
\textsf{Preuve :}
Voir \cite{Cur70}, Theorem 5.7, reformulé en \cite{HenVig11b}, Theorem 3.7. \hfill$\Box$\\

On remarque, par la proposition \ref{irredGbar}.(i), que $V^{\overline{N}}$ est de dimension $1$ si et seulement si $\overline{P}$ est inclus dans $\overline{P}_V$; et en particulier, $V$ est un caractère de $\overline{G}$ si et seulement si on a $\Delta_V = \Delta$. \\
Enfin, il est intéressant de savoir que l'on a accès aux paramètres, au sens de la proposition \ref{irredGbar}.(ii) de la représentation duale $V^* = \Hom_{\Fpbar}(\Fpbar,V)$. 

\begin{prop} \label{paramV*}
Soit $V$ une représentation irréductible de $\overline{G}$, de paramètres $(\chi,\Delta_V)$. Alors $V^*$ est irréductible, et on a \[ \Par_{\overline{B}^-}(V) = \big( w_0(\chi), w_0(\Delta_V) \big) , \quad \Par_{\overline{B}}(V^*) = \big( w_0(\chi^{-1}), -w_0(\Delta_V) \big) .\]
\end{prop}
\textsf{Preuve :}
Voir \cite{HenVig11b}, Lemma 3.11 et Lemma 3.12. \hfill$\Box$\\

Un autre point de vue intéressant à ce sujet est celui de la Proposition 2.2 de \cite{Her11b} et de la Proposition II.3.10 de \cite{Jan87}; cependant on ne l'utilisera pas ici et on ne le décrira donc pas. \\

Comme toute représentation irréductible $V$ de $K$ est égale à son espace (non nul) d'invariants par le pro-$p$-radical $K(1)$, l'action de $K$ sur $V$ se factorise en fait par $\overline{G} = K/K(1)$. Réciproquement toute représentation irréductible de $\overline{G}$ définit une représentation irréductible de $K$ par inflation. De ce fait, il n'y a pas lieu dans la suite de faire de distinction entre une $\overline{G}$-représentation irréductible et une $K$-représentation irréductible.

\subsection{Les \og déterminants \fg\ dans $\GL(m,D)$} \label{pardeterminants}

Lorsque $D$ n'est pas commutatif (c'est-à-dire $d>1$), il n'y a pas de déterminant canonique $\GL(m,D) \to D^\times$. Cependant, on va tout de même présenter un morphisme de groupes $\det_G : \GL(m,D) \to F^\times$ qu'on conviendra d'appeler déterminant. \\
Parce que $E$ déploie $D$, on a le diagramme commutatif suivant entre algèbres de matrices, qui nous permet de définir $\det_0$:
\[ \xymatrix{ M_m(D) \otimes_F E \ar[r]^{\phantom{xx} \sim} & M_{md}(E) \ar[d]^{\det} \\ M_m(D) \ar@{^{(}->}[u] \ar[r]^{\det_0} & E } .\]
L'isomorphisme de déploiement \[ \eta_D : M_m(D) \otimes_F E \xrightarrow{\sim} M_{md}(E) \] n'est pas $\Gal(E/F)$-équivariant mais on sait par application du théorème de Skolem-Noether que, pour tout élément $\sigma \in \Gal(E/F)$, il existe une matrice inversible $c_\sigma \in \GL(md,E)$ telle que pour tout $x \otimes 1 \in M_m(D) \otimes_F E$, on a : 
\[ \sigma \circ \eta_D (x \otimes 1) = c_\sigma \big( \eta_D(x \otimes 1) \big) c_\sigma^{-1} .\] 
En particulier, $\det_0(x)$ est fixe par $\Gal(E/F)$ pour tout $x \in M_m(D)$, et $\det_0$ est donc à valeurs dans $F$. On définit alors la restriction de $\det_0$ à $\GL(m,D)$, $\det_G : \GL(m,D) \to F^\times$, qui est un morphisme de groupes. Pour $m=1$, $\det_G$ devient le morphisme de groupes $\Nrm: D^\times \to F^\times$ appelé norme réduite. \\
On remarquera que si $A \in \GL(m,D)$ est une matrice triangulaire supérieure par blocs de blocs diagonaux $(A_1, \dots, A_r) \in \prod\limits_{i=1}^r \GL(m_i,D)$ avec $\sum_i m_i =m$, alors on a \[ \det\nolimits_G(A) = \prod_{i=1}^r \, \det\nolimits_{\GL(m_i,D)} (A_i) ;\] en particulier, pour des blocs de dimension $1$, on obtient 
\begin{equation} \label{GLmdetNrm} \det\nolimits_G(A) = \Nrm(A_1) \cdots \Nrm(A_m) .\end{equation}
Pour plus de détails, on pourra se reporter au \S 12 $\mathrm{n^o}$3 de \cite{Bou58II8}, au \S 4 de \cite{Bou58II8}, ainsi qu'au paragraphe 16 de \cite{Pie82}. \\

Un autre morphisme dérivé du déterminant va aussi entrer en jeu lorsque l'on s'intéressera aux représentations de $\GL(m,\Oo_D)$. Il s'agit de la composée \[ \overline{\det} : \GL(m,\Oo_D) \to \GL(m,k_D) \xrightarrow{\det} k_D^\times ,\] où la première application est la projection canonique $K \to K/K(1)$, qui correspond concrètement à la réduction modulo $(\varpi_D)$ coefficient par coefficient. \\ On fait la remarque que si on note $N_{k_D/k_F}$ la norme de l'extension galoisienne $k_D/k_F$, alors on peut vérifier la compatibilité que résume le diagramme commutatif suivant:
\[\xymatrix{ \GL(m,\Oo_D) \ar[d]_{\det_G} \ar[r]^{\phantom{xxx} \overline{\det}} &  k_D^\times \ar[d]^{N_{k_D/k_F}} \\ \Oo_F^\times \ar[r]^{\overline{\cdot}} & k_F^\times } .\] Elle sera prouvée et exploitée dans la preuve de la proposition \ref{conversemult}.

\section{Transformée de Satake, isomorphisme de comparaison} \label{rappels}

L'objet de cette section est de rappeler des faits établis dans d'autres références récentes sur la transformée de Satake et la comparaison entre induite parabolique et induite compacte. Dans un dernier paragraphe, on établira quelques résultats élémentaires sur la transformée de Satake de certains opérateurs de Hecke.

\subsection{Transformation de Satake} \label{parIntroSatake}

On se contente ici de rappeler les résultats utiles de \cite{HenVig11} et de \cite{HenVig11b}. Soient $P$ un sous-groupe parabolique\footnote{on remarquera que certains résultats de \cite{HenVig11} nécessitent que $P$ vérifie une décomposition d'Iwasawa $G = K P$, mais c'est ici automatique parce que $K$ est un parahorique maximal spécial (voir \cite{BruTit72}, Proposition 4.4.6)} de $G$, $M$ un sous-groupe de Levi de $P$ et $N$ son radical unipotent. \\
Soit $V$ une représentation irréductible de $\overline{G}$, inflatée à $K$. Lorsque l'on prend $P = B$ minimal, l'espace $V^{U \cap K}$ est de dimension $1$ et $A \cap K$ agit dessus à travers un caractère $\chi$. On définit alors son stabilisateur \[ Z_A(\chi) = \{ a \in A \ | \ \forall x \in A \cap K , \, \chi(a^{-1} x a) = \chi(x) \} \] et $Z_A^-(\chi)$ (resp. $Z_A^+(\chi)$) le sous-monoïde de $Z_A(\chi)$ des éléments antidominants (resp. dominants), c'est-à-dire dont l'action par conjugaison contracte $U^-$ (resp. contracte $U$). Pour éviter d'introduire $\chi$, on notera parfois $Z_A(V) := Z_A(\chi)$ et $Z_A^\pm(V) := Z_A^\pm(\chi)$. \\
On fait la remarque que la définition de $Z_A(\chi)$ est tout à fait explicite et se lit de la manière suivante: le caractère $\chi$ s'écrit $\chi_1 \otimes \dots \otimes \chi_m$ où chaque $\chi_i$ est un caractère $\Oo_D^\times \to k_D^\times \to \Fpbar^\times$, de stabilisateur dans $\varpi_D^{\Z}$ un $\varpi_D^{d_i \Z}$ pour un certain entier naturel $d_i \mid d$. L'intersection $A_\Lambda \cap Z_A(\chi)$ est diagonalement $\varpi_D^{d_1 \Z} \times \cdots \times \varpi_D^{d_m \Z}$, et comme l'on a $A = (A \cap K) \rtimes A_\Lambda$ et que chaque $\chi_i$ se factorise par $k_D^\times$, qui est abélien, on a 
\begin{equation} \label{ZAchi} Z_A(\chi) = (A \cap K) \rtimes (A_\Lambda \cap Z_A(\chi)) .\end{equation}
On définit grâce à cela $\widetilde{Z}_M(\chi) = \widetilde{Z}_M \cap Z_A(\chi)$ et $\widetilde{Z}_M^\pm(\chi) = \widetilde{Z}_M \cap Z_A^\pm(\chi)$. On remarquera en particulier que $\widetilde{Z}_M(\chi) \cap K$ est réduit à l'élément neutre. On note $\widetilde{Z}_M^{--}(\chi)$ le sous-monoïde de $\widetilde{Z}_M^-(\chi)$ constitué des éléments $t$ vérifiant $v_D(\beta(t))<0$ pour tout $\beta \in \Delta \smallsetminus \Delta_M$. Enfin, comme précédemment, on se permettra de noter $\widetilde{Z}_M^\ast(V) := \widetilde{Z}_M^\ast(\chi)$ pour $\ast \in \{ \varnothing, -, --, +, ++ \}$. \\

L'énoncé suivant définit la transformée de Satake (par rapport au parabolique standard minimal $B$).

\begin{prop} \label{HVsatake}
Le morphisme d'algèbres \begin{equation} \label{defSatake} \begin{array}{cccc} \Sss_G : & \Hh_{\Fpbar}(G,K,V) & \to & \Hh_{\Fpbar} \big( A, A \cap K, V^{U \cap K} \big) \\ & f & \mapsto & \Big( x \mapsto \sum_{U/U \cap K} f(xu) \Big) \end{array} \end{equation} est injectif. Son image est l'espace $\Hh_{\Fpbar}^- \big( A, A \cap K, V^{U \cap K} \big)$ des fonctions de support $Z^{-}_A(V)$.
\end{prop}
\textsf{Preuve :}
Voir \cite{HenVig11}, Theorem 7.1. \hfill$\Box$\\

Puisque l'on a $A = (A \cap K) \rtimes A_\Lambda$, on peut appliquer la deuxième proposition du paragraphe 2.10 de \cite{HenVig11} pour voir que $\Hh_{\Fpbar} \big( A, A \cap K , V^{U \cap K} \big)$ est commutative; il en est alors de même de $\Hh_{\Fpbar}(G,K,V)$. \\

On peut aussi définir une transformée de Satake partielle, lorsque le parabolique standard $P$ n'est plus supposé minimal: on note cette dernière \begin{equation} \label{defSatpartiel} \begin{array}{cccc} \Sss_G^M : & \Hh_{\Fpbar}(G,K,V) & \to & \Hh_{\Fpbar} \big( M, M \cap K, V^{N \cap K} \big) \\ & f & \mapsto & \Big( m \mapsto \sum_{N/N \cap K} f(mn) \Big) .\end{array} \end{equation} 
Le morphisme de $\Fpbar$-algèbres $\Sss_G^M$ est aussi injectif. \\
On aura besoin de la variante suivante de la transformée de Satake: notons $p_N$ la projection canonique $V \twoheadrightarrow V_{N \cap K}$ et définissons \[\begin{array}{cccc} ' \Sss_G^M : & \Hh_{\Fpbar}(G,K,V) & \to & \Hh_{\Fpbar} \big( M, M \cap K , V_{N \cap K} \big) \\ & f & \mapsto & \big( m \mapsto p_N \sum_{n \in N \cap K \backslash N} f(nm) \big) \end{array} .\] 
Comme avant, on va noter $' \Sss_G$ pour $' \Sss_G^A$. \\ 

Par la Proposition 2.3 de \cite{HenVig11b}, si on note \[ \begin{array}{cccc} \iota_G : & \Hh_{\Fpbar}(G,K,V^*) & \to & \Hh_{\Fpbar}(G,K,V)^{\op} \\ & \Phi & \mapsto & \big( g \mapsto \Phi(g^{-1})^t \big) \end{array} ,\] on a le diagramme commutatif suivant 
\[\xymatrix{ \Hh_{\Fpbar}(G,K,V^*) \ar[d]_{\iota_G}^\sim \ar@{^{(}->}[rr]^{\Sss_G^M} & & \Hh_{\Fpbar}\big( M, M \cap K, (V^*)^{N \cap K} \big) \ar[d]_{\iota_M}^\sim \\ \Hh_{\Fpbar}(G,K,V)^{\op} \ar@{^{(}->}[rr]^{('\Sss_G^M)^{\op}} & & \Hh_{\Fpbar} \big( M, M \cap K, V_{N \cap K} \big)^{\op} }\]
où $V^* = \Hom_{\Fpbar}(V,\Fpbar)$ désigne la $\overline{G}$-représentation duale. Parce que les algèbres du diagramme ci-dessus sont ici commutatives, on obtient l'énoncé suivant. 

\begin{prop} \label{SatakeDual}
Soient $V$ une $K$-représentation irréductible et $P = MN$ un parabolique standard. Alors on a le diagramme commutatif de $\Fpbar$-algèbres commutatives suivant:
\[\xymatrix{ \Hh_{\Fpbar}(G,K,V^*) \ar[d]_{\iota_G}^\sim \ar@{^{(}->}[rr]^{\Sss_G^M} & & \Hh_{\Fpbar}\big( M, M \cap K, (V^*)^{N \cap K} \big) \ar[d]_{\iota_M}^\sim \\ \Hh_{\Fpbar}(G,K,V) \ar@{^{(}->}[rr]^{' \Sss_G^M} & & \Hh_{\Fpbar} \big( M, M \cap K, V_{N \cap K} \big) } .\]
\end{prop}

Les isomorphismes de cette proposition justifient que l'on peut jongler invariablement entre $\Sss_G^M$ et $' \Sss_G^M$ à condition de prendre en compte $V^*$ au lieu de $V$ lorsque que c'est approprié. Cela explique aussi l'apparition de la représentation duale $V^*$ dans ce travail alors qu'il n'en était pas mention dans \cite{Her11b}. \\

Ajoutons un mot sur les opérateurs de Hecke de $\Hh_{\Fpbar}(G,K,V)$. Soit $t$ un élément de $A_\Lambda^-$. On sait que le sous-espace de $\Hh_{\Fpbar}(G,K,V)$ formé des fonctions à support dans $K t K$ est de dimension $1$ si on a $t \in \widetilde{Z}_A^-(V)$, et il est de dimension $0$ le cas contraire. De plus, dans le cas non réduit à $0$, il existe un unique opérateur $T_t$ de $\Hh_{\Fpbar}(G,K,V)$ de support $K t K$ tel que $T_t (t)$ induit l'identité sur $V^{U \cap K}$: c'est à ce dernier que l'on fera référence quand on parlera de l'opérateur de $\Hh_{\Fpbar}(G,K,V)$ porté par la double classe $K t K$. Aussi, lorsque l'on a $t = \widetilde{\lambda}(\varpi) \in \widetilde{Z}_A^-(V)$ pour un certain $\lambda \in X_*(T)^-$, on notera parfois $T_\lambda$ l'opérateur $T_t$. 

\begin{prop} \label{Prop2.11} 
L'application $'\Sss_G^M$ est injective et est une localisation. De plus, on a $'\Sss_G = {}' \Sss_M \circ {}' \Sss_G^M$.
\end{prop}
\textsf{Preuve :}
Voir \cite{HenVig11b}, Proposition 4.1, équation (9) et Proposition 4.5. \hfill$\Box$\\

Le morphisme injectif $'\Sss_G^M$ est en fait une localisation (dans un cadre commutatif) par rapport à une partie multiplicative que l'on peut engendrer par un unique élément. Explicitons ce dernier: fixons un élément $s_M \in A_\Lambda \cap Z(M)^+$ qui est strictement dominant hors de $\Delta_M$ (comme il est choisi central dans $M$, il est même tout à fait indépendant de $V$). Comme $s_M$ est central dans $M$, on peut considérer l'unique opérateur ${}' T_M$ de $\Hh_{\Fpbar} \big( M, M \cap K, V_{N \cap K} \big)$ de support $(M \cap K) s_M$ et valant $\Id_{V_{N \cap K}}$ en $s_M$. On remarquera que ${}' T_M$ est défini, comme dans la discussion précédente, par sa valeur en $V_{U \cap K}$; mais comme $s_M$ est central dans $M$, ${}' T_M(s_M)$ induit bien l'identité sur $V_{N \cap K}$. Alors ${}' T_M$ est dans l'image de ${}' \Sss_G^M$ et on a \begin{equation} \label{localH} \Hh_{\Fpbar} \big( M, M \cap K, V_{N \cap K} \big) = {}' \Sss_G^M \big( \Hh_{\Fpbar}(G,K,V) \big) \big[ {}' T_M^{-1} \big] .\end{equation}

Notons $T_M \in \Hh_{\Fpbar} \big( M, M \cap K, V^{N \cap K} \big)$ l'opérateur porté par la double classe $(M \cap K) s_M^{-1}$ et valant $\Id_{V^{N \cap K}}$ en $s_M^{-1}$. Cette définition est compatible avec celle de ${}' T_M$ au sens de l'énoncé suivant. 

\begin{lem}
On a $\iota_M \big( {}' T_M^{(V)} \big) = T_M^{(V^*)}$, où les indices supérieurs précisent la $K$-représentation pour laquelle l'opérateur de Hecke est considéré.
\end{lem}
\textsf{Preuve :} \\
Par définition, $\iota_M$ est la fonction \[ \iota_M : T \mapsto \big( g \mapsto T(g^{-1})^t \big) ,\] où l'indice $\cdot^t$ désigne la transposée. En particulier, $\iota_M \big( {}' T_M^{(V)} \big)$ a donc pour support $(M \cap K) s_M^{-1}$ et a valeur $\big( \Id_{V_{N \cap K}} \big)^t = \Id_{(V^*)^{N \cap K}}$ en $s_M^{-1}$. \hfill$\Box$ \\

En l'appliquant (\ref{localH}) à $V^*$ (puisque l'on a $(V^*)^* \simeq V$), on a donc l'énoncé similaire pour $\Sss_G^M$: 
\begin{equation} \label{eq5'} \Hh_{\Fpbar} \big( M, M \cap K, V^{N \cap K} \big) = \Sss_G^M \big( \Hh_{\Fpbar}(G,K,V) \big) \big[ T_M^{-1} \big] .\end{equation}

\subsection{Représentations régulières, induction parabolique et induction compacte} \label{rappelcomparaison}

Là encore, on continue de rappeler les résultats utiles de \cite{HenVig11} et de \cite{HenVig11b}. \\
Soient $V$ une $K$-représentation irréductible et $P = MN$ un parabolique standard de $G$. On dira que $V$ est $P$\textit{-régulière} si le stabilisateur $\overline{P}_V$ de la droite $V^{U \cap K}$ est contenu dans $\overline{P}$. On peut exhiber le lien suivant entre $K$-représentations irréductibles et $(M \cap K)$-représentations irréductibles.

\begin{lem} \label{Lem2.5}
Les représentations irréductibles $P$-régulières $V$ de $K$ sont en bijection avec les représentations irréductibles de $M \cap K$ via l'application $V \mapsto V^{N \cap K}$. 
\end{lem}
\textsf{Preuve :}
Voir \cite{HenVig11b}, Proposition 3.10. \hfill$\Box$\\

L'énoncé de comparaison suivant est un argument-clef dans toute la machinerie que constitue l'article \cite{Her11b} ou bien le travail qui va suivre: c'est le théorème principal de \cite{HenVig11b}. On prend les notations suivantes pour l'énoncé ci-dessous: $\Hh = \Hh_{\Fpbar}(G,K,V)$ et $\Hh_M = \Hh_{\Fpbar}(M,M \cap K, V_{N \cap K})$. 

\begin{prop} \label{Thm3.1}
Soit $V$ une représentation irréductible $P^-$-régulière de $K$. Soit $\chi : \Hh_M \to \Fpbar$ un morphisme d'algèbres, que l'on voit aussi grâce à la transformée de Satake ${}' \Sss_G^M$ comme un morphisme $\Hh \to \Fpbar$. L'application \[ \ind_K^G \, V \otimes_{\Hh_,\chi} \Fpbar \to \Ind_P^G \big( \ind_{M \cap K}^M \, V_{N \cap K} \otimes_{\Hh_M, \chi} \Fpbar \big) \] est un isomorphisme de $G$-représentations.
\end{prop}
\textsf{Preuve :}
Voir \cite{HenVig11b}, Theorem 1.2 et Corollary 1.3. \hfill$\Box$

\subsection{Calcul de quelques transformées de Satake}

On commence par énoncer un petit lemme sur la réduction. Soient $M$ un Levi standard de blocs $\GL(m_1, D) \times \cdots \GL(m_r,D)$ avec $\sum_i m_i=m$ et $t \in A_\Lambda$ un élément de coefficients diagonaux $(\varpi_D^{a_1}, \dots, \varpi_D^{a_1}, \dots, \varpi_D^{a_r},\dots, \varpi_D^{a_r})$ avec $m_i$ occurences du coefficient $\varpi_D^{a_i}$ pour des entiers relatifs $a_i$, $1 \leq i \leq r$. On définit une action de $A_\Lambda$ par conjugaison sur $\overline{M}$ en faisant agir $t$ par $x \mapsto \overline{\varpi_D^{a_i}[x] \varpi_D^{-a_i}}$ sur les coefficients du $i$-ème bloc diagonal de $\overline{M}$. \\
Soit $t = \widetilde{\lambda}(\varpi) \in A_\Lambda^-$ pour un certain $\lambda \in X_*(T)^-$. On note respectivement $P_t, M_t, N_t, \Delta_t$ pour $P_\lambda, M_\lambda, N_\lambda, \Delta_\lambda$; on fera attention que $P_t$ est antistandard.

\begin{lem} \label{Prop3.8}
Soit $t \in A_\Lambda^-$. Alors $K \cap t^{-1} K t$ se réduit en $\overline{P}_t$. De plus, on a \[ \big\{ (\overline{g}, \overline{tgt^{-1}}) \ \big| \ g \in K \cap t^{-1} K t \big\} = \big\{ (g_+,g_-) \in \overline{P}_t \times \overline{P}_t^- \ \big| \ \pi(g_-) = t \pi(g_+) t^{-1} \in \overline{M}_t \big\} ,\] où $\pi$ désigne la projection canonique de $\overline{P}_t^-$ (resp. $\overline{P}_t$) sur $\overline{M}_t$.
\end{lem}
\textsf{Preuve :} \\
C'est le pendant de la Proposition 3.8 de \cite{Her11}, qui est clair pour $\GL(m,D)$. Voir aussi la Proposition 6.13 de \cite{HenVig11}. \hfill$\Box$\\
 
On va énoncer un lemme qui va nous permettre grosso modo, dans le calcul de transformées de Satake d'opérateurs de Hecke, de nous ramener au cas où $V$ est un caractère de $K$. Aussi, on note $T_t^M$ l'opérateur de $\Hh_{\Fpbar} \big( M,M \cap K, V^{N \cap K} \big)$ de support $(M \cap K)t(M \cap K)$ et induisant l'identité sur $V^{U \cap K}$; et on écrira $\tau_t$ au lieu de $T_t^A$.

\begin{lem} \label{deGàMv}
Soient $P = MN$ un parabolique standard, $V$ une $K$-représentation irréductible $P$-régulière et $t \in \widetilde{Z}_A^-(V)$. Alors on a $\Sss_G^M(T_t) = T_t^M$.
\end{lem}
\textsf{Remarque :} 
Cela veut dire en particulier que si $V$ est le plus régulier possible (c'est-à-dire $P_V = B$), le calcul des transformées de Satake est immédiat. \\
\textsf{Preuve :} \\
Commençons par remarquer que si $g = k_1 tk_2 \in KtK$ satisfait $T_t(g)|_{V^{N \cap K}} \neq 0$, alors la réduction $\overline{k}_2 \in \overline{K}$ de $k_2$ est dans $\overline{P}_t \overline{P}$. En effet, $T_t(g)|_{V^{N \cap K}} \neq 0$ implique que l'image de $k_2 V^{N \cap K}$ dans $V_{N_t \cap K}$ est non nulle, et par \cite{HenVig11b}, Corollary 3.17, cela donne $\overline{k}_2 \in \overline{P}_t \overline{P}_V \overline{P}$. Parce que $V$ est $P$-régulière, on a $\overline{P}_V \subseteq \overline{P}$ et l'affirmation voulue. Notons $I_P$ le parahorique de $G$ correspondant à $P$, c'est-à-dire l'image réciproque de $\overline{P}$ par la projection canonique $K \twoheadrightarrow \overline{G}$. Par le lemme \ref{Prop3.8}, on peut retraduire $\overline{k}_2 \in \overline{P}_t \overline{P}$ en: $k_2 \in (K \cap t^{-1}Kt) I_P$ et donc $g \in KtI_P$. \\ 
En vue d'utiliser la définition (\ref{defSatpartiel}), prenons $m \in M$ et $n \in N$ et appliquons l'observation précédente à $g =mn \in P \cap KtK$. Si la restriction de $T_t(mn)$ à $V^{N \cap K}$ est non nulle, alors on a $mn \in KtI_P \cap P$. Par décomposition d'Iwahori\footnote{le pro-$p$-radical $I_P(1)$ de $I_P$ admet une décomposition d'Iwahori \[ I_P(1) = (I_P(1) \cap N^-)(I_P(1) \cap M)(I_P(1) \cap N) \] par \cite{Sec04}, Proposition 2.14; de plus l'application \[ I_P \cap M / I_P(1) \cap M \to I_P / I_P(1)\] est bijective et on peut écrire \[ I_P = (I_P(1) \cap N^-)(I_P(1) \cap M)(I_P(1) \cap N)(I_P \cap M).\] Parce que $I_P \cap M$ normalise $I_P(1) \cap N$, cela devient \[ I_P = (I_P(1) \cap N^-)(I_P \cap M)(I_P(1) \cap N) .\] On vérifie alors que c'est bien une décomposition.} \begin{equation} \label{decompIwahori} I_P = (I_P \cap N^-)(I_P \cap M)(I_P \cap N) ,\end{equation} on écrit \[ KtI_P = Kt(I_P \cap N^-)(I_P \cap P) .\] Parce que $t$ contracte $U^- \cap K \geq I_P \cap N^-$, on a \[ Kt(I_P \cap N^-)(P \cap K) = Kt(P \cap K) .\] Mais alors, $mn$ appartient à \[ Kt(P \cap K) \cap P = (P \cap K)t(P \cap K) .\] Comme $t^{-1}$ contracte $U \cap K \geq N \cap K$, on a finalement \[ mn \in (M \cap K)t(M \cap K)(N \cap K) .\] Il s'ensuit $m \in (M \cap K)t(M \cap K)$ et $n \in N \cap K$. La transformée $\Sss_G^M(T_t)$ est alors de support inclus dans $(M \cap K)t(M \cap K)$ et $\Sss_G^M(T_t)(t) = T_t(t)|_{V^{N \cap K}}$ induit l'identité sur $V^{U \cap K} \subseteq V^{N \cap K}$: il en résulte $\Sss_G^M(T_t) = T_t^M$. \hfill$\Box$\\

Enfin, le résultat suivant nous rappelle qu'il est facile de calculer la transformée de Satake d'un opérateur de Hecke supporté par une simple classe, c'est-à-dire lorsque l'on a $KtK = tK$.

\begin{lem} \label{Satakehomothetie}
Soient $V$ une représentation irréductible de $K$ et $t \in \widetilde{Z}_A(V)$ vérifiant $\Delta_t = \Delta$. Alors on a $\Sss_G(T_t) = \tau_t$.
\end{lem}
\textsf{Preuve :} \\
La condition $\Delta_t = \Delta$ nous assure que $t$ est une homothétie: $t$ normalise $K$ et on a $KtK = tK$. Mais alors, si $t'$ est un élément de $A_\Lambda^-$, $tK \cap t' U \neq \varnothing$ implique $t=t'$ et on a $tK \cap tU = t(U \cap K)$. Le résultat découle alors directement de la définition (\ref{defSatake}). \hfill$\Box$

\section{Paramètres de Hecke-Satake}

Dans ce paragraphe, on tâche de comprendre un peu les caractéristiques des caractères de $\Fpbar$-algèbres $\chi : \Hh_{\Fpbar}(G,K,V) \to \Fpbar$. \\
On note $\chi^{(A)}$ l'application $\chi \circ (\Sss_G)^{-1}$, qui est donc définie sur $\Sss_G \big( \Hh_{\Fpbar}(G,K,V) \big)$.

\begin{lem} \label{Lem4.3} 
Un morphisme de $\Fpbar$-algèbres $\chi : \Hh_{\Fpbar}(G,K,V) \to \Fpbar$ se factorise à travers $\Sss_G^M$ si et seulement si $\chi^{(A)}$ contient $\widetilde{Z}^-_M(V)$ dans son support.
\end{lem}
\textsf{Preuve :} \\ 
A cause de (\ref{eq5'}), $\chi$ se factorise à travers $\Hh_{\Fpbar} \big( G,K, V^{N \cap K} \big)$ si et seulement si on peut prolonger $\chi$ à $\big((\Sss_G^M)^{-1}(T_M) \big)^{-1}$, autrement dit si et seulement si on a $\chi \big((\Sss_G^M)^{-1}(T_M) \big) \neq 0$. Parce que $s_M^{-1}$ est central dans $M$, $\Sss_M(T_M)$ est l'opérateur porté par la classe $s_M^{-1} (A \cap K)$, et on le notera $\tau_{s_M^{-1}}$. En particulier, $\chi \big((\Sss_G^M)^{-1}(T_M) \big) \neq 0$ est équivalent à $\chi^{(A)}(\tau_{s_M^{-1}}) \neq 0$. Reste à voir que ceci équivaut à $\chi^{(A)}(\tau_t) \neq 0$ pour tout opérateur $\tau_t$ porté par un $t(A \cap K)$ pour $t \in \widetilde{Z}_M^-(V)$. Le sens $(\Leftarrow)$ est évident puisque l'on a $s_M^{-1} \in \widetilde{Z}_M^-(V)$. Réciproquement, supposons l'existence de $t \in \widetilde{Z}_M^-(V)$ vérifiant $\chi(\tau_t)=0$. Parce que $s_M^{-1}$ est strictement antidominant hors de $\Delta_M$, il existe $n \geq 1$ tel que $s_M^{-n} t^{-1}$ soit un élément de $\widetilde{Z}_M^-(V)$. On a alors 
\[ \chi \big((\Sss_G^M)^{-1}(T_M) \big) ^n = \chi^{(A)}(\tau_{s_M^{-n}}) = \chi^{(A)}(\tau_{s_M^{-n} t^{-1}}) \chi^{(A)}(\tau_t) = 0\]
 et on a prouvé ce que l'on voulait. Le morphisme $\chi$ se prolonge donc à $\Hh_{\Fpbar} \big( M,M \cap K, V^{N \cap K} \big)$ si et seulement si $\chi^{(A)}$ ne s'annule sur aucune des classes portées par des éléments de $\widetilde{Z}_M^-(V)$. \hfill$\Box$

\begin{lem} \label{Lem4.3demi}
Soit $\chi : \Hh_{\Fpbar}(G,K,V) \to \Fpbar$ un caractère de $\Fpbar$-algèbres. Il existe un unique Levi standard $M$ tel que la restriction à $A_\Lambda$ du support de $\chi^{(A)}$ est égal à $\widetilde{Z}_M^-(V)$.
\end{lem}
\textsf{Preuve :} \\
On prend un Levi standard minimal $M$ tel que l'on ait $\widetilde{Z}_M^-(V) \subseteq \Supp \, \chi^{(A)}$; d'abord, un tel Levi existe car $\Supp \, \chi^{(A)}$ contient au moins les éléments de $\widetilde{Z}_G^-(V)$. Ensuite, cela aura deux conséquences:
\begin{itemize}
\item[(i)] $M$ est unique pour la condition de minimalité et $\widetilde{Z}_M^-(V) \subseteq \Supp \, \chi^{(A)}$;
\item[(ii)] $M$ vérifie $\widetilde{Z}_M^-(V) = (\Supp \, \chi^{(A)}) \cap A_\Lambda$.
\end{itemize}
La conjonction de ces deux faits prouvera le lemme voulu. \\
Commençons par (i): supposons l'existence d'un autre Levi standard $M'$ minimal vérifiant $\widetilde{Z}_{M'}^-(V) \subseteq \Supp \, \chi^{(A)}$. Par définition, l'élément $s_M s_{M'}$ vérifie $\beta(s_M s_{M'}) >0$ pour toute racine $\beta \notin \Delta_{M \cap M'}$. Alors il existe $n \geq 1$ assez grand tel que $(s_M s_{M'})^{-n}$ s'écrive $s_{M \cap M'}^{-1} t$ pour un $t \in \widetilde{Z}_A^-(V)$. On a alors \[ \chi^{(A)}\big( (s_M s_{M'})^{-n} \big) = \chi^{(A)}(s_{M \cap M'}^{-1}) \, \chi^{(A)}(t) \neq 0 ,\] et donc $\chi^{(A)}(s_{M \cap M'}^-1) \neq 0$. Par la preuve du lemme \ref{Lem4.3}, on obtient l'inclusion $\widetilde{Z}_{M \cap M'}^-(V) \subseteq \Supp \, \chi^{(A)}$. Sauf à avoir $M=M'$, cela contredit la minimalité de $M$: (i) est prouvé. \\ 
Pour (ii), si $\widetilde{Z}_M^-(V) \subseteq (\Supp \, \chi^{(A)}) \cap A_\Lambda$ était une inclusion stricte, il existerait $t \in A_\Lambda$ avec $v_D(\beta(t)) < 0$ pour un $\beta \notin \Delta_M$. Fixons deux tels $\beta$, $t$. Mais alors, $(t s_M^{-1})^n$ se factoriserait par un $s^{-1}_{M'}$ pour un Levi $M' \lneq M$ et un $n \geq 1$ assez grand, et cela contredirait la minimalité de $M$. Le résultat suit. \hfill$\Box$

\begin{prop} \label{Prop4.1} 
On a une bijection ensembliste \[ \begin{array}{ccr} \Hom_{\Fpbar\textrm{-alg}}\big( \Hh_{\Fpbar}(G,K,V),\Fpbar \big) & \leftrightarrow & \left\{ \begin{tabular}{r|l} \multirow{2}{*}{$(M,\chi_M)$} & $M$ Levi standard, \\ & $\chi_M : \widetilde{Z}_M(V) \to \Fpbar^\times$ caractère \end{tabular} \right\} \end{array} .\]  A un couple $(M,\chi_M)$, on associe le morphisme 
\begin{equation} \label{Prop3.6} \begin{array}{ccc} \Hh_{\Fpbar}(G,K,V) & \to & \Fpbar \\ t & \mapsto & \sum_{z \in \widetilde{Z}_M(V)} \Sss_G(t)(z) \chi_M(z)^{-1} \end{array} .\end{equation}
\end{prop}
\textsf{Preuve :} \\ 
On commence par compléter chaque $M$ en un parabolique standard $P = M N$ avec un radical unipotent $N$. \\ Si $\chi$ est un élément de $\Hom_{\Fpbar\textrm{-alg}}(\Hh_{\Fpbar}(G,K,V),\Fpbar)$, par le lemme \ref{Lem4.3demi}, la restriction à $A_\Lambda$ du support de $\chi^{(A)}$ est de la forme $\widetilde{Z}_M^-(V)$ pour un certain Levi standard $M$. Par le lemme \ref{Lem4.3}, $\chi$ se factorise alors à travers un $\chi^{(M)} : \Hh_{\Fpbar}(M,M \cap K, V^{N \cap K}) \to \Fpbar$. \\
Soit $H_M$ la sous-algèbre de $\Hh_{\Fpbar}(M,M \cap K,V^{N \cap K})$ des fonctions $t$ de support vérifiant $\Supp \, t \cap A_\Lambda \subseteq \widetilde{Z}_M(V)$. L'id\' eal $I_M$ des fonctions $t$ avec $t|_{\widetilde{Z}_M(V)} = 0$ constitue un suppl\' ement de $H_M$ dans $\Hh_{\Fpbar} \big( M, M \cap K , V^{N \cap K }\big)$ et on veut donc param\' etriser 
\[ \Hom_{\Fpbar\textrm{-alg}} \big( \Hh_{\Fpbar} \big( M, M \cap K,V^{N \cap K} \big) ,\Fpbar \big) = \Hom_{\Fpbar\textrm{-alg}}(H_M, \Fpbar) .\]
La restriction \` a $A_\Lambda$ fait de $H_M$ une $\Fpbar$-alg\` ebre isomorphe \` a
\[ \{ t : \widetilde{Z}_M(V) \to \Fpbar \textrm{ à support fini} \} .\]
Gr\^ ace au crochet de dualit\' e $(f,g) = \sum\limits_{z \in \widetilde{Z}_M(V)} f(z) g(z)$, on voit que le dual lin\' eaire $\Hom_{\Fpbar}(H_M,\Fpbar)$ est isomorphe \` a l'espace vectoriel 
\[ H_M' = \{ t : \widetilde{Z}_M(V) \to \Fpbar \} .\] En testant \` a travers $( \cdot, f)$ o\` u $f$ est \` a support un singleton, on voit que $t \in H_M'$ est dans $\Hom_{\Fpbar\textrm{-alg}}(H_M, \Fpbar)$ si et seulement si $t$ est \` a valeurs dans $\Fpbar^\times$ et d\' efinit un morphisme de groupes $t : \widetilde{Z}_M(V) \to \Fpbar^\times$. On note $\chi_M$ l'inverse de celui qui correspond de la sorte \` a $\chi^{(M)}$. \\
Et, dans l'autre sens, $\chi_M$ détermine uniquement l'élément 
\[ \chi_M' : t \mapsto \sum_{z \in \widetilde{Z}_M(V)} t(z) \chi_M(z)^{-1} \] 
dans $\Hom_{\Fpbar\textrm{-alg}}(H_M,\Fpbar)$. Si $t$ est port\' e par un singleton $z_0 \in \widetilde{Z}_M(V)$ alors $\chi_M'(t) = \chi_M^{-1}(z_0)$ est \' egal \` a $\chi^{(M)}(z_0)$, ce qui permet de d\' efinir $\chi^{(M)} : H_M \to \Fpbar$ de mani\` ere compatible \` a la construction pr\' ec\' edente; on le prolonge \` a $\Hh_{\Fpbar} \big( M, M \cap K, V^{N \cap K} \big)$ en envoyant $I_M$ sur $0$. \\
Pour d\' efinir un morphisme d'alg\` ebres $\Hh_{\Fpbar}(G,K,V) \to \Fpbar$, il reste \` a pr\' ecomposer par $\Sss_G^M$ pour donner 
\begin{equation} \label{presqueProp4.1} t \mapsto \sum_{z \in \widetilde{Z}_M(V)} \Sss_G^M(t)(z) \chi_M(z)^{-1} .\end{equation}
Parce que l'on ne voit que la restriction de $\Sss_G^M(t)$ \` a $\widetilde{Z}_M(V) (M \cap K)$ dans la formule (\ref{presqueProp4.1}), le lemme \ref{Satakehomothetie} et la proposition \ref{Prop2.11} nous permettent d'y remplacer $\Sss_G^M(t)$ par $\Sss_G(t)$. \hfill$\Box$ \\

Pour tout caractère de $\Fpbar$-algèbres $\chi: \Hh_{\Fpbar}(G,K,V) \to \Fpbar$, on dira que $\chi$ a pour \textit{paramètres de Hecke-Satake} la paire $(M,\chi_M)$ qui lui associée par la proposition \ref{Prop4.1}. On fait la remarque que par rapport à la paramétrisation utilisée dans \cite{Her11b}, on ne tient pas compte ici d'une partie intégrale pour $\chi_M$: notre paramétrisation \og ne voit pas \fg\ ce qui est sur $A \cap K$. \\ 
On dira qu'une $K$-représentation irréductible $V$ est un ($K$-)\textit{poids} de $\pi$ si on a $\Hom_K(V,\pi) \neq 0$. Si $\pi$ est une représentation admissible, l'espace \[ \Hom_G \big( \ind_K^G \, V, \pi \big) \simeq \Hom_K(V,\pi) \] est de dimension finie sur $\Fpbar$ et en particulier l'action de l'algèbre commutative $\Hh_{\Fpbar}(G,K,V)$ dessus possède des vecteurs propres. Chaque vecteur propre $f \in \Hom_K(V,\pi)$ définit un caractère $\chi : \Hh_{\Fpbar}(G,K,V) \to \Fpbar$ de $\Fpbar$-algèbres et on dira qu'un tel $\chi$ est un \textit{caractère propre} sur $\Hom_K(V,\pi)$.

\begin{lem} \label{Lem4.6}
Soit $\pi$ une représentation admissible de $G$ de caractère central $\chi_Z$. Soient $V$ un poids de $\pi$ et $\chi: \Hh_{\Fpbar}(G,K,V) \to \Fpbar$ un caractère propre de $\Hom_K(V,\pi)$. Si $\chi$ est de paramètres de Hecke-Satake $(M,\chi_M)$, alors la restriction de $\chi_M$ \`a $Z(G) \cap A_\Lambda$ est $\chi_Z$, o\`u $Z(G)$ d\' esigne le centre de $G$.
\end{lem}
\textsf{Preuve :} \\
Comme les fonctions de support $K t K$ pour $t \in \widetilde{Z}_A^-(V)$ forment une famille génératrice de $\Hh_{\Fpbar}(G,K,V)$, il s'agit de voir que l'on a l'égalité $\chi_Z(t) = \chi_M(t)$ sous la condition $t \in Z(G) \cap A_\Lambda \subseteq \widetilde{Z}_A^-(V)$. Lorsque $t$ est un élément central de $G$, la double classe $K t K$ est en fait une classe simple $t K$ et l'image par $\Sss_G^M$ de l'opérateur $T_t$ de support porté par cette classe est celui porté par $t (M \cap K)$ (voir \cite{HenVig11}, paragraphe 7.3). Dans ce cas-là, par (\ref{Prop3.6}), on a $\chi(T_t) = \chi_M(t^{-1})$. Reste à établir $\chi(T_t) = \chi_Z(t^{-1})$ pour conclure. \\ 
Soit $f \in \Hom_K(V, \pi)$ un élément propre pour $\chi$, et regardons $f$ comme un élément de $\Hom_G \big( \ind_K^G \, V, \pi \big)$ par réciprocité de Frobenius. On a alors d'une part \[ (t.f) \circ T_t = \chi_Z (t) \, (f \circ T_t) = \chi_Z (t) \chi(T_t) f .\] Et d'autre part, pour $[g,v] \in \ind_K^G \, V$, parce que $t$ est central dans $G$, la formule d'action (9) de \cite{BarLiv94} donne \[ t. f \circ T_t([g,v]) = t f \big( [ g t^{-1}, v ] \big) = f([g,v]) .\] De ce fait, on a $\chi_Z (t) \chi(T_t) = 1$, d'où le lemme. \hfill$\Box$

\begin{lem} \label{paramtorsion}
Soient $\pi$ une représentation admissible de $G$ et $\xi$ un caractère de $F^\times$. Soient $V$ un poids de $\pi$ et $f \in \Hom_K(V,\pi)$ un vecteur propre de caractère propre $\chi : \Hh_{\Fpbar}(G,K,V) \to \Fpbar$ de paramètre de Hecke-Satake $(M,\chi_M)$. Alors $f$ induit un morphisme $K$-équivariant \[ f_\xi : (\xi \circ \det\nolimits_G) \otimes V \to (\xi \circ \det\nolimits_G) \otimes \pi ,\] qui est vecteur propre pour $\Hh_{\Fpbar} \big( G,K,(\xi \circ \det_G) \otimes V \big)$ de caractère propre $\chi_\xi : \Hh_{\Fpbar} \big( G,K,(\xi \circ \det_G) \otimes V \big) \to \Fpbar$ paramétrisé par $ \big( M, \chi_M (\xi \circ \det_G)|_{\widetilde{Z}_M(V)} \big)$.
\end{lem}
\textsf{Preuve :} \\
Le morphisme de $G$-représentations 
\[ \begin{array}{ccc} (\xi \circ \det\nolimits_G) \otimes \ind_K^G \, V  & \to & \ind_K^G \big( (\xi \circ \det\nolimits_G) \otimes V \big) \\ a \otimes f & \mapsto & \big( g \mapsto a \, \xi \circ \det\nolimits_G(g) f(g) \big) \end{array} \]
est un isomorphisme et $f$ induit bien $f_\xi$ comme désiré. Aussi, en identifiant $\End_{\Fpbar}(V)$ et $\End_{\Fpbar}((\xi \circ \det_G) \otimes V)$, on obtient un isomorphisme d'algèbres \[ \begin{array}{cccc} \eta_\xi: & \Hh_{\Fpbar}(G,K,V) & \to & \Hh_{\Fpbar} \big( G,K, (\xi \circ \det_G) \otimes V \big) \\ & \phi & \mapsto & \phi_\xi \end{array} \] et on a le diagramme commutatif suivant
\begin{equation} \label{commtorsion} \xymatrix{ \Hh_{\Fpbar} \big( G,K, (\xi \circ \det_G) \otimes V \big) \ar[rr]^{\phantom{xxxxxxxxxx} \chi_\xi} & & \Fpbar \\ \Hh_{\Fpbar}(G,K,V) \ar[u]^{\eta_\xi}_{\sim} \ar[urr]_{\chi} & & } .\end{equation}
De plus, on a \[ \Sss_G(\phi_\xi)(a) = \xi \circ \det\nolimits_G(a) \, \Sss_G(\phi)(a) \] pour tout $\phi \in \Hh_{\Fpbar}(G,K,V)$ et $a \in A$. Le résultat suit alors de la combinaison de (\ref{commtorsion}) et de (\ref{Prop3.6}). \hfill$\Box$ \\

Dans l'étude de l'irréductibilité d'induites paraboliques, on veut comprendre les caractères propres sur $\Hom_K(V,\pi)$ pour $\pi = \Ind_P^G \, \sigma$ où $P$ est un parabolique standard et $\sigma$ une représentation admissible de son Levi standard $M$. \\
Par la proposition \ref{SatakeDual}, $\iota_G$ induit une bijection \[ \iota_G^* : \Hom_{\Fpbar\textrm{-alg}} \big( \Hh_{\Fpbar}(G,K,V), \Fpbar \big) \xrightarrow{\sim} \Hom_{\Fpbar\textrm{-alg}} \big( \Hh_{\Fpbar}(G,K,V^*), \Fpbar \big) .\] Si $\chi : \Hh_{\Fpbar}(G,K,V) \to \Fpbar$ est un caractère de $\Fpbar$-algèbres, on dira que la paire $(M, \chi_M)$ associée à $\iota_G^*(\chi)$ par la proposition \ref{Prop4.1} est le \textit{paramètre HS-dual} de $\chi$; on insiste que $\chi_M$ est alors un caractère $\widetilde{Z}_M(V^*) \to \Fpbar^\times$.

\begin{lem} \label{paramHSinduite}
Soient $P = MN$ un parabolique standard de $G$ et $\sigma$ une représentation admissible de $M$. Soit $f \in \Hom_K(V,\Ind_P^G \, \sigma)$ un vecteur propre de caractère propre $\chi : \Hh_{\Fpbar}(G,K,V) \to \Fpbar$ de paramètre HS-dual $(L,\chi_L)$. Alors $f$ induit $f_{(N)} \in \Hom_{M \cap K}(V_{N \cap K}, \sigma)$ de caractère propre \[ \chi_{(N)} : \Hh_{\Fpbar}(G,K,V_{N \cap K}) \to \Fpbar \] de paramètre HS-dual $(L,\chi_L)$; et $f \mapsto f_{(N)}$ est une bijection entre les ensembles de vecteurs propres de $\Hom_K(V,\Ind_P^G \, \sigma)$ et de $\Hom_{M \cap K}(V_{N \cap K}, \sigma)$.
\end{lem}
\textsf{Remarque :}
En particulier, on voit l'inclusion $L \subseteq M$. \\
\textsf{Preuve :} \\
Pa réciprocité de Frobenius, on voit plutôt $f$ comme un élément propre dans $\Hom_G(\ind_K^G \, V, \Ind_P^G \, \sigma)$. Par le diagramme (4) de \cite{HenVig11b}, l'isomorphisme de réciprocité de Frobenius \[ \begin{array}{ccc} \Hom_G(\ind_K^G \, V,\Ind_P^G \, \sigma) & \xrightarrow{\sim} & \Hom_M(\ind_{M \cap K}^M \, V_{N \cap K}, \sigma) \\ f & \mapsto & f_{(N)} \end{array} \] est $\Hh_{\Fpbar}(G,K,V)$-équivariant (où l'action de $\Hh_{\Fpbar}(G,K,V)$ au but se fait à travers ${}' \Sss_G^M$) et on a le diagramme commutatif suivant: 
\begin{equation} \label{diagHSinduit} \xymatrix{ \Hh_{\Fpbar} \big( M, M \cap K, (V^*)^{N \cap K} \big) \ar[r]^{\iota_M} & \Hh_{\Fpbar}(M,M \cap K, V_{N \cap K}) \ar[dr]^{\chi_{(N)}} &  \\ \Hh_{\Fpbar}(G,K,V^*) \ar[r]^{\iota_G} \ar[u]_{\Sss_G^M} & \Hh_{\Fpbar}(G,K,V) \ar[r]^{\chi} \ar[u]_{{}' \Sss_G^M} & \Fpbar } .\end{equation} Le résultat s'en déduit. \hfill$\Box$

\section{Quelques représentations irréductibles}

\subsection{Représentations de Steinberg généralisées} \label{steinberg} 

Soient $Q$ un sous-groupe parabolique standard de $G$ et $\xi$ un caractère lisse de $F^\times$. On définit la $G$-représentation à coefficients dans $\Fpbar$ suivante : \[ \St_Q \xi := \frac{\Ind^{G}_{Q} (\xi \circ \det_G) }{\sum_{Q' \gneq Q} \Ind^{G}_{Q'} (\xi \circ \det_G)} ,\] où on a encore noté $\det_G$ sa restriction à un sous-groupe parabolique standard de $G$. En particulier, lorsque $\xi$ est le caractère trivial de $F^\times$, on obtient \[ \St_Q \id := \frac{\Ind^{G}_{Q} \, \id}{\sum_{Q' \gneq Q} \Ind^{G}_{Q'} \, \id} .\]

\begin{lem}
Soient $\xi$ un caractère lisse de $F^\times$ et $Q$ un parabolique standard de $G$. Alors on a l'isomorphisme de $G$-représentations \[ (\xi \circ \det\nolimits_G) \otimes \Ind_Q^G \, \id \xrightarrow{\sim} \Ind_Q^G(\xi \circ \det\nolimits_G) .\]
\end{lem}
\textsf{Preuve :} \\
Le morphisme de $G$-représentations \[ a \otimes f \mapsto \big( g \mapsto a \, \xi \circ \det\nolimits_G(g) f(g) \big) \] ci-dessus admet naturellement \[  f \mapsto 1 \otimes \big(g \mapsto \xi^{-1} \circ \det\nolimits_G(g) f(g) \big) \] pour inverse. \hfill$\Box$ \\

De ce fait, on obtient l'isomorphisme de $G$-représentations \[ \St_Q \xi \simeq (\xi \circ \det\nolimits_G) \otimes \St_Q \id .\] On peut donc utiliser les résultats du chapitre \ref{Ly11} pour en déduire immédiatement les énoncés suivants, et on se référera à l'article mentionné pour les preuves.

\begin{prop} \label{GLmStirred}
Soient $\xi$ un caractère lisse de $F^\times$ et $Q$ un parabolique standard de $G$. La représentation de Steinberg généralisée $\St_Q \xi$ est irréductible et admissible.
\end{prop}

\begin{cor} \label{GLmStgeneral}
Les constituants de Jordan-Hölder de $\Ind^{G}_{Q} (\xi \circ \det_G)$ sont les $\St_{Q'} \xi$ pour $Q' \geq Q$ parabolique\footnote{comme $Q'$ contient $Q$ standard, il l'est automatiquement aussi}. Ils sont deux à deux non isomorphes et de multiplicité $1$.
\end{cor}

On aura de plus besoin de faire l'induction parabolique des représentations de Steinberg généralisées et de comprendre les représentations irréductibles en sous-quotient.

\begin{cor} \label{GLmJHsteinberg}
Soient $P$ un parabolique standard de $G$ et $Q$ un parabolique standard du Levi standard $M$ de $P$. Alors la représentation $\Ind_P^G (\St_Q \xi)$ est de longueur finie, de constituants de Jordan-Hölder les $\St_R \xi$ avec $R$ parabolique standard vérifiant $\Delta_M \cap \Delta_R = \Delta_Q$, chacun apparaissant avec multiplicité $1$. 
\end{cor}
\textsf{Remarque :}
Puisque l'on a $\Delta_Q \subseteq \Delta_M$, la condition $\Delta_M \cap \Delta_R = \Delta_Q$ est équivalente à $\Delta_R \supseteq \Delta_Q$ et $\Delta_M \smallsetminus \Delta_Q \subseteq \Delta \smallsetminus \Delta_R$. Cette seconde forme est parfois plus agréable à manipuler. \\

Au cours de la preuve de la proposition \ref{GLmStirred}, on a en fait dû examiner les $K$-poids de $\St_Q \xi$. Comme cela sera utile par la suite, on rappelle l'énoncé. 

\begin{prop} \label{Stsocle}
Le $K$-socle de $\St_Q \xi$ est irréductible et il est isomorphe à l'unique $K$-représentation irréductible $V_Q^\xi$ de paramètre \[ \Par_{\overline{B}^-}(V_Q^\xi) = \big( (\xi \circ \det\nolimits_G)|_{A \cap K} , - \Delta_Q \big) .\]
\end{prop}

Plus tard, on voudra connaître les paramètres de Hecke-Satake des représentations de Steinberg généralisées. Grâce au lemme \ref{paramtorsion} et \` a la Proposition 8.2 du chapitre \ref{Ly11}, on a le résultat suivant.

\begin{prop} \label{paramSt}
Le paramètre HS-dual de $\St_Q \xi$ est donné par la paire \[ (A, \xi^{-1} \circ \det\nolimits_G : \widetilde{Z}_A \to \Fpbar^\times) .\]
\end{prop}

\subsection{Représentations supersingulières} \label{prgssing}

Maintenant que l'on a rappelé les propriétés voulues pour les représentations de Steinberg généralisées, il nous reste à exhiber un autre type de \og briques de base \fg\ pour les induites paraboliques considérées. Ce sont les représentations supersingulières; il faut garder en tête que l'on peut y penser comme des représentations supercuspidales aussi bien que l'équivalence entre les deux notions n'apparaîtra qu'à la toute fin du papier. \\

Suivant \cite{HenVig11b}, on dira qu'une représentation irréductible $\pi$ de $G$ est $K$\textit{-supersingulière} (et on oubliera toujours le $K$, celui-ci ayant été fixé une fois pour toutes) si on a \[ \Hh_{\Fpbar} \big( M , M \cap K, V_{N \cap K} \big) \otimes_{\Hh, {}' \Sss^M_G} \Hom_G \big( \ind_K^G \, V, \pi \big)= 0 \] pour toute représentation irréductible $V$ de $K$ et tout parabolique standard $P = MN$ propre de $G$. On fait la remarque que par convention on n'appellera pas supersingulière une représentation de $D^\times$ irréductible; on a bien conscience que cette convention n'est pas des plus r\' epandues mais elle est commode ici. \\

Aussi, une représentation admissible est $K$-supersingulière si et seulement si pour tout représentation irréductible $V$ de $K$ et tout caractère $\chi : \Hh_{\Fpbar}(G,K,V) \to \Fpbar$ propre sur $\Hom_K(V,\pi)$, le paramètre HS-dual de $\chi$ est un $(G,\chi_G)$. C'est le sujet de la discussion après la Definition 7.2 dans \cite{HenVig11b}.

\section{Transformée de Satake et conjugaison}

Deux représentations irréductibles $V_1$ et $V_2$ de $\overline{G}$ sont dites \textit{conjuguées} si les représentations $(\chi_1,V_1^{\overline{U}})$ et $(\chi_2,V_2^{\overline{U}})$ de $\overline{B}$ sont des caractères conjugués au sens suivant: il existe $a \in A_\Lambda$ tel que l'on ait, pour tout $x \in A \cap K$, $\chi_2(x) = \chi_1(a^{-1} x a)$. En particulier, si les représentations $V_1^{\overline{U}}$ et $V_2^{\overline{U}}$ sont données par un même caractère, alors les $K$-représentations irréductibles $V_1$ et $V_2$ sont conjuguées. \\

Notons $X_*(T)^-_V$ le sous-monoïde de $X_*(T)^-$ constitué des éléments $\lambda$ vérifiant $\widetilde{\lambda}(\varpi) \in Z_A(V)$. On a besoin de voir que ce monoïde ne dépend que de la classe de conjugaison de la $K$-représentation irréductible $V$ considérée. On remarquera que, dans le cas déployé, cela veut juste dire que l'image de la transformée de Satake ne dépend que de la classe d'isomorphisme de la droite $V^{U \cap K}$ en tant que $(A \cap K)$-représentation, ce qui va de soi.

\begin{lem} \label{droiteconj}
Soient $V_1$ et $V_2$ deux $K$-représentations irréductibles conjuguées. Alors on a \[X_*(T)^-_{V_1} = X_*(T)^-_{V_2} .\]
\end{lem}
\textsf{Preuve :} \\
Il s'agit de voir le lien entre $Z_A(V_1)$ et $Z_A(V_2)$ lorsque $(\chi_1,V_1^{\overline{U}})$ et $(\chi_2,V_2^{\overline{U}})$ vérifient $\chi_2 = \chi_1(a^{-1} \cdot a)$ pour un certain $a \in A_\Lambda$. Et par définition de $Z_A(V)$, on a l'égalité \[Z_A(V_2) = a^{-1} Z_A(V_1) a .\] Parce que $A_\Lambda$ est commutatif, le résultat suit alors de la discussion avant la description (\ref{ZAchi}). \hfill$\Box$\\

Soient $V_1$ et $V_2$ deux représentations irréductibles de $K$. On note respectivement $(\chi_1, \Delta_1)$ et $(\chi_2,\Delta_2)$ les paramètres $\Par_{\overline{B}}(V_1)$ et $\Par_{\overline{B}}(V_2)$. On cherche à comprendre quelques opérateurs d'entrelacement dans $\Hh_{\Fpbar}(G,K,V_1,V_2)$. \\
Soit $t = \widetilde{\lambda}(\varpi) \in A_\Lambda^-$ pour un certain $\lambda \in X_*(T)^-$.  

\begin{lem} \label{Heckepdsconj}
L'espace des fonctions de $\Hh_{\Fpbar}(G,K,V_1,V_2)$ de support inclus dans $KtK$ est de dimension majorée par $1$, avec égalité si et seulement si on a \begin{equation} \label{condpdsconj} \chi_2 = \chi_1 (t^{-1} \cdot t) \textrm{ sur } A \cap K  \quad \textrm{et} \quad \Delta_1 \cap \Delta_t = \Delta_2 \cap \Delta_t .\end{equation}
\end{lem}
\textsf{Remarque :} 
Un \' enonc\' e similaire est pr\' esent au paragraphe 7.7 de \cite{HenVig11} mais la preuve y est omise. \\
\textsf{Preuve :} \\
Soit $f \in \Hh_{\Fpbar}(G,K,V_1,V_2)$ un opérateur de support inclus dans la double classe $K t K$ pour un certain élément $t =\widetilde{\lambda}(\varpi)$ avec $\lambda \in X_*(T)^-$. Alors $f$ vérifie $k_2 f(t) = f(t) k_1$ pour tous $k_1, k_2 \in K$ vérifiant $k_2 t = t k_1$. Par le lemme \ref{Prop3.8}, on a alors le diagramme commutatif 
\begin{equation} \label{commpdsconj} \xymatrix{ (V_1)_{N_t \cap K} \ar[r]^{f(t)} \ar[d]^{k_1} & V_2^{N_t^- \cap K} \ar[d]^{t k_1 t^{-1}} \\ (V_1)_{N_t \cap K} \ar[r]^{f(t)} & V_2^{N_t^- \cap K} }\end{equation}
et par \cite{Hum06}, Theorem 5.1, $V_2^{N_t^- \cap K}$ et $(V_1)_{N_t \cap K}$ sont irréductibles en tant que $(M_t \cap K)$-représentations, de sorte que $f$ est dans un espace de dimension au plus $1$ par lemme de Schur. S'il existe $f \in \Hh_{\Fpbar}(G,K,V_1,V_2)$ non nul et de support $KtK$, le diagramme (\ref{commpdsconj}) et le lemme de Schur nous assurent que $V_1^{N_t^- \cap K}$ et $t^{-1}.V_2^{N_t^- \cap K}$ sont deux représentations irréductibles de $M_t \cap K$ isomorphes. La proposition \ref{irredGbar}.(ii) nous assure alors (\ref{condpdsconj}). \\
Réciproquement, supposons la condition (\ref{condpdsconj}) satisfaite. Alors, parce que l'on a (voir section 5.7 de \cite{HenVig11}) \[ \Par_{\overline{B}} \big( V_1^{N_t^- \cap K} \big) = (\chi_1 , \Delta_1 \cap \Delta_t ) \quad \textrm{et} \quad \Par_{\overline{B}} \big( V_2^{N_t^- \cap K} \big) = (\chi_2 , \Delta_2 \cap \Delta_t ) ,\] $V_1^{N_t^- \cap K}$ et $t^{-1}.V_2^{N_t^- \cap K}$ sont des $(M_t \cap K)$-représentations isomorphes, et on peut choisir pour $f(t)$ un tel isomorphisme (unique à un scalaire près par lemme de Schur). Cela termine la preuve. \hfill$\Box$\\

Soit $t = \widetilde{\lambda}(\varpi) \in A_\Lambda^-$. On suppose que $t, V_1, V_2$ vérifient (\ref{condpdsconj}). D'après la preuve du lemme \ref{Heckepdsconj}, il n'y a en fait qu'une seule manière (à une constante inversible près) de définir un opérateur d'entrelacement de $\Hh_{\Fpbar}(G,K,V_1,V_2)$ porté par une double classe $K t K$: sa valeur en $t$ sera choisie (et on ne le précisera plus par la suite) comme la composée 
\begin{equation} \label{defHecke} V_1 \twoheadrightarrow (V_1)_{N_t \cap K} \xrightarrow{\sim} t^{-1}.V_2^{N_t^- \cap K} \xrightarrow{t} V_2^{N_t^- \cap K} \hookrightarrow V_2 ,\end{equation}
où l'isomorphisme de $(M_t \cap K)$-représentations au milieu est choisi et fixé une fois pour toutes, et o\` u, comme dans la preuve du lemme \ref{Heckepdsconj}, la fl\` eche not\' ee $t$ est $v \mapsto tv$. Plus tard, on notera cet opérateur $\varphi_t^?$, ou encore $\varphi_\lambda^?$ si on a posé $t = \widetilde{\lambda}(\varpi)$ pour un certain $\lambda \in X_*(T)^-$ (avec $?$ indice variable selon les choix de $V_1$ et $V_2$). \\
Disons maintenant un mot de la transformée de Satake sur l'espace d'entrelacements $\Hh_{\Fpbar}(G,K,V_1,V_2)$. On pourra aussi se reporter aux 2.8 et 7.9 de \cite{HenVig11}.

\begin{prop} \label{Prop6.3}
Soient $V_1$, $V_2$ et $V_3$ trois représentations irréductibles de $K$. La transformée de Satake \[\begin{array}{cccc} \Sss_G : & \Hh_{\Fpbar}(G,K,V_1,V_2) & \to & \Hh_{\Fpbar} \big( A, A \cap K, V_1^{U \cap K}, V_2^{U \cap K} \big) \\  & f & \mapsto & \Big( x \mapsto \sum_{U/U \cap K} f(xu) \Big) \end{array}\] est injective et compatible aux compositions au sens que le diagramme suivant commute. 
\[ \xymatrix{ \Hh_{\Fpbar}(G,K,V_2,V_3) \times \Hh_{\Fpbar}(G,K,V_1,V_2) \ar[r] \ar[d]^{\Sss_G \times \Sss_G} & \Hh_{\Fpbar}(G,K,V_1,V_3) \ar[d]^{\Sss_G} \\ \Hh_{\Fpbar} \big( A, A \cap K, V_2^{U \cap K}, V_3^{U \cap K} \big) \times \Hh_{\Fpbar} \big( A, A \cap K, V_1^{U \cap K}, V_2^{U \cap K} \big) \ar[r] & \Hh_{\Fpbar} \big( A, A \cap K, V_1^{U \cap K}, V_3^{U \cap K} \big) }\] 
\end{prop}
\textsf{Preuve :} \\
Pour prouver l'injectivité, prenons $f$ un élément non nul de $\Hh(G,K,V_1,V_2)$ et montrons que $\Sss_G(f)$ est aussi non nul. Supposons d'abord que $f$ est la fonction caractéristique d'une double classe, disons $K \widetilde{\lambda}(\varpi) K$ pour un certain $\lambda \in X_*(T)^-$. Par la Proposition 4.4.4.(ii) de \cite{BruTit72}, on sait que si l'intersection $K \widetilde{\lambda}(\varpi) K \cap \widetilde{\mu}(\varpi) U$ est non vide, alors on a $\mu \geq_{\R} \lambda$, c'est-à-dire que $\mu - \lambda$ est combinaison linéaire à coefficients réels positifs de cocaractères, et $K \widetilde{\lambda}(\varpi) K \cap \widetilde{\lambda}(\varpi) U = \widetilde{\lambda}(\varpi) (U \cap K)$. La définition de la transformée de Satake nous donne alors directement $\Sss_G(f)(\widetilde{\lambda}(\varpi))=1$, et donc $\Sss_G(f) \neq 0$. Revenons maintenant à un $f$ quelconque non nul: écrivons $\coprod_i K \widetilde{\mu}_i(\varpi) K$ son support. 
Si on prend $\mu_j$ minimal pour $\geq_\R$ dans cette écriture, le calcul précédent nous montre que les transformées de Satake des fonctions caractéristiques des autres doubles classes $K \widetilde{\mu}_i(\varpi) K$ s'annule en $\widetilde{\mu}_j(\varpi)$ et donc seul intervient dans $\Sss_G(f)(\widetilde{\mu}_j(\varpi))$ la valeur de $f$ sur $\widetilde{\mu}_j(U \cap K)$, qui est non nulle. On en déduit l'injectivité voulue. \\
Pour voir la compatibilité aux compositions, prenons $f_1 \in \Hh_{\Fpbar}(G,K,V_1,V_2)$ et $f_2 \in \Hh_{\Fpbar}(G,K,V_2,V_3)$. Pour $a \in A$, on a alors: 
\begin{align*} \Sss_G(f_2 * f_1) (a) & = \sum_{n \in U/U \cap K} f_2 * f_1(an) \\ & = \sum_{x \in G/K} \sum_{n \in U/U \cap K} f_2(x) f_1(x^{-1} an) .\end{align*}
De plus, à cause de la bijection $B/B\cap K \xrightarrow{\sim} G/K$, on a finalement: 
\[ \Sss_G(f_2 *f_1) (a) = \sum_{x \in B/B \cap K} \sum_{n \in U/U \cap K} f_2(x) f_1(x^{-1}an) .\]
Aussi, d'un autre côté, on a: 
\begin{align*} \Sss_G(f_2) * \Sss_G(f_1)(a) & = \sum_{y \in A/A \cap K} \Sss_G(f_2) (y) \Sss_G(f_1)(y^{-1}a) \\ & = \sum_{y \in A / A \cap K} \sum_{u \in U/U \cap K} f_2(yu) \sum_{v \in U/U \cap K} f_1(y^{-1}av) \\ & = \sum_{y \in A/A \cap K} \sum_{u \in U/U \cap K} f_2(yu) \sum_{v \in U/U \cap K} f_1(u^{-1} y^{-1} a (a^{-1} y u y^{-1} a)v) .\end{align*}
En translatant $v$ à gauche par $a^{-1} y u y^{-1} a$, on obtient \[ \Sss_G(f_2) * \Sss_G(f_1)(a) = \sum_{y \in A / A \cap K} \sum_{u \in U/ U \cap K} f_2(yu) \sum_{v \in U/U \cap K} f_1(u^{-1} y^{-1}a v) .\] Enfin, parce que l'on a l'égalité ensembliste \[ \{ yu \ | \ y \in A / A \cap K, \, u \in U/U \cap K \} = B/B \cap K ,\] on réécrit \[ \Sss_G(f_2) * \Sss_G(f_1) (a) = \sum_{b \in B / B \cap K} f_2(b) \sum_{v \in U/ U \cap K} f_1(b^{-1}av) ;\] ce qui donne $\Sss_G(f_2 *f_1) = \Sss_G(f_2) * \Sss_G(f_1)$. \hfill$\Box$

\section{Formule de Lusztig-Kato revisitée} \label{LKglm}

Dans le cas déployé $D=F$, l'article \cite{Her11b} utilise la formule de Lusztig-Kato pour expliciter des calculs sur la transformée de Satake. On va comparer des transformées de Satake pour $\GL(m,D)$ et pour son sous-groupe $\GL(m,E)$, qui est déployé sur le sous-corps non ramifié maximal $E$, pour pallier l'absence de formule de Lusztig-Kato pour $\GL(m,D)$: c'est l'objet de ce paragraphe. \\

Par le paragraphe 7.1 de \cite{HenVig11} ou le Theorem 1.0.1 de \cite{HaiRos10}, on a l'isomorphisme d'algèbres communément appelé isomorphisme de Satake classique \[ \Sss: \Hh_{\C}(G,K) \xrightarrow{\sim} \C[A/A \cap K]^W .\] Par le paragraphe 7.11 de \cite{HenVig11}, on a $\Sss = \delta_G^{1/2} \Sss_G^{\C}$ où $\delta_G^{1/2}$ est le caractère module $A \to \C^\times$ de $B^-$ et\footnote{on fera attention que l'on a fait une somme sur l'unipotent opposé pour passer de dominant à antidominant; cela explique aussi pourquoi on a choisi le caractère module de $B^-$ et pas de $B$} 
\begin{equation} \label{SatakC} \begin{array}{cccc} \Sss_G^{\C} : & \Hh_{\C}(G,K) & \to & \Hh_{\C}(A, A \cap K) \\ & f & \mapsto & \big( x \mapsto \sum\limits_{u \in U^-/ U^- \cap K} f(xu) \big) \end{array} .\end{equation} 
Par (\ref{versSclassic}), on a un isomorphisme d'algèbres \[ \C[A/A \cap K] \xrightarrow{\sim} \C[X_*(T)] ,\] qui est de plus compatible avec l'action de $W$. On définit alors la composée \[\overline{\Sss} : \Hh_{\C}(G,K) \mathop{\xrightarrow{\sim}}\limits_{\Sss} \C[A/A \cap K]^W \xrightarrow{\sim} \C[X_*(T)]^W .\]
Pour $g \in G$, on note $\id_{KgK}$ l'opérateur de $\Hh_{\C}(G,K)$ porté par la double classe $K g K$ et valant $1$ en $g$. \\
On note $\GE = \GL(m,E) \leq \GL(m,D) = G$, $\KE$ le compact maximal $\GL(m,\Oo_E)$, $T_{(E)}$ le tore $E$-déployé diagonal de $G_{(E)}$ et $\Hh_{\C}(G_{(E)},K_{(E)})$ l'algèbre de Hecke-Satake comme définie au paragraphe \ref{generalreps}. On remarquera que $W$, le groupe de Weyl fini défini à l'aide de $T$ pour $G$, coïncide avec le groupe de Weyl défini par $T_{(E)}$ pour $G_{(E)}$; on les confond donc par la suite. La transformée de Satake classique pour $\GE$ est \[ \Sss_{(E)} : \Hh_{\C}(\GE,\KE) \xrightarrow{\sim} \C[X_*(T_{(E)})]^W ;\] parce que le groupe des cocaractères $X_*(T_{(E)})$ s'identifie canoniquement à $X_*(T)$, on note encore $\Sss_{(E)}$ l'isomorphisme de $\C$-algèbres \[ \Sss_{(E)} : \Hh_{\C}(\GE,\KE) \xrightarrow{\sim} \C[X_*(T)]^W .\] L'ingrédient important de cette section réside dans l'énoncé qui va suivre : c'est une méthode de corps proches (voir \cite{Kaz86}) que l'on applique (certes pour des représentations modulo $p$), qui est ici très explicite puisqu'en niveau zéro.

\begin{prop} \label{beforeLK}
Supposons $m \leq 3$. Soit $\lambda$ un élément de $X_*(T)^+$. On a \[ \Sss_{(E)} \big(\id_{\KE \lambda(\varpi) \KE} \big) = \overline{\Sss} \big(\id_{K \widetilde{\lambda}(\varpi) K} \big) .\]
\end{prop}
\textsf{Preuve :} \\
De la même manière qu'en (\ref{SatakC}), on écrit $\Sss_{(E)} = \delta_{(E)}^{1/2} \Sss_{(E)}^{\C}$ pour $G_{(E)}$. Par les équations (21) et (24) de \cite{Car79}, on a 
\[ \delta_G \big( \widetilde{\lambda}(\varpi) \big) = q^{- v_D (\det_G \, \widetilde{\lambda}(\varpi))} = \big( q^d \big)^{- \frac{v_D (\det \, \lambda(\varpi))}{d} } = \delta_{(E)} \big( \lambda(\varpi) \big) ,\]
o\` u $\det : G_{(E)} \to E^\times \subseteq D^\times$ est le d\' eterminant usuel. Il s'agit maintenant de comparer $\Sss_G^{\C}(T_\lambda)$ et $\Sss_{(E)}^\C(T_\lambda^{(E)})$ où on a noté $T_\lambda^{(E)}$ l'opérateur de $\Hh_\C(\GE,\KE)$ porté par $\KE \lambda(\varpi) \KE$. Pour $f = T_\lambda^{(E)}$, $\mu$ dans $X_*(T)^+$ et $x = \mu(\varpi)$, la formule (\ref{SatakC}) pour $\GE$ se lit\footnote{$\BE^-$ et $\UE^-$ sont les notations pour $B^- \cap \GE$ et $U^- \cap \GE$ respectivement}: \[ \Sss_{(E)}^\C \big( T_\lambda^{(E)} \big) (\mu(\varpi)) = \sum_{\UE^-/\UE^- \cap \KE} T_\lambda^{(E)}(\mu(\varpi)u) .\] La décomposition d'Iwasawa nous permet d'écrire 
\begin{equation} \label{egalIw} \KE \lambda(\varpi) \KE = \coprod_{i \in \I} \lambda_i(\varpi) u_i \KE ,\end{equation} 
où $\I$ est un ensemble fini et $\lambda_i \in X_*(T)$, $u_i \in \UE^-$ pour tout $i \in \I$. L'intersection $\lambda_i(\varpi) u_i \KE \cap \mu(\varpi) \UE^-$ est non vide si et seulement si on a $\lambda_i = \mu$; en effet, supposons $\lambda_i(\varpi) u_i k = \mu(\varpi) u$ pour un certain $k \in \KE$ et $u \in \UE^-$. On a alors 
\[ k = u_i^{-1} \lambda_i(\varpi)^{-1} \mu(\varpi) u \in \BE^- \cap \KE \quad \textrm{et} \quad \mu(\varpi)^{-1} \lambda_i(\varpi) \in \UE^- (\BE^- \cap \KE) \UE^- .\] 
D'où l'égalité $\lambda_i = \mu$; dans ce cas (voir paragraphe 6.9 de \cite{HenVig11}), on a \[ \lambda_i(\varpi) u_i \KE \cap \mu(\varpi) \UE^- = \mu(\varpi) u_i (\UE^- \cap \KE) .\] Le coefficient\footnote{$\tau_\mu^{(E)}$ désigne l'opérateur de $\Hh_\C(A \cap \GE, A \cap \KE)$ porté par la double classe $(A \cap \KE) \mu(\varpi) (A \cap \KE) = \mu(\varpi) (A \cap \KE)$} de $\tau_\mu^{(E)}$ dans l'expression de $\Sss_{(E)}^\C \big( T_\lambda^{(E)} \big)$ est alors tout simplement
\begin{equation} \label{coefunram} \sum_{\lambda_i = \mu} T_\lambda^{(E)} (\mu(\varpi) u_i) = \big| \{ i \in \I \ | \ \lambda_i = \mu \}\big| .\end{equation}
On cherche ainsi \` a \' evaluer le cardinal de 
\begin{equation} \label{CartanInterIwasawa}\big( \KE \lambda(\varpi) \KE \cap \mu(\varpi) \UE^- \big) \big/ \UE^- \cap \KE .\end{equation}
Lorsque ce dernier est égal à celui de
\begin{equation} \label{CartanInterIwasawaD} \big( K \widetilde{\lambda}(\varpi) K \cap \widetilde{\mu}(\varpi) U^- \big) \big/ U^- \cap K ,\end{equation}
le coefficient de $\tau_\mu$ dans $\Sss_G^{\C}(T_\lambda)$ est égal à celui dans (\ref{coefunram}) comme voulu. L'égalité des cardinaux de (\ref{CartanInterIwasawa}) et (\ref{CartanInterIwasawaD}) est le fruit du lemme \ref{egalcardinaux} ultérieur si $m$ est inférieur ou égal à $3$: le r\' esultat est prouv\' e. \hfill$\Box$

\begin{lem} \label{egalcardinaux}
Supposons $m \leq 3$. Soient $\lambda, \mu$ des éléments de $X_*(T)^-$. On a l'égalité de cardinaux 
\[ \big| \big( K \widetilde{\lambda}(\varpi) K \cap \widetilde{\mu}(\varpi) U \big) \big/ U \cap K \big| = \big| \big( \KE \lambda(\varpi) \KE \cap \mu(\varpi) \UE \big) \big/ \UE \cap \KE \big| .\]
\end{lem}
\textsf{Preuve du cas $m \leq 2$ :} \\
Le cas $m=1$ étant trivial, supposons $m=2$. Notons
\begin{equation} \label{notmum=2} \lambda(\varpi) = \left( \begin{array}{cc} \varpi^\alpha & 0 \\ 0 & \varpi^\beta \end{array} \right) , \quad \mu(\varpi) u = \left( \begin{array}{cc} \varpi^x & a \\ 0 & \varpi^y \end{array} \right) \in \mu(\varpi) \UE .\end{equation}
Notons $a \mapsto \widetilde{a}$ la bijection $E \xrightarrow{\sim} D$ envoyant $\sum\limits_{i \gg - \infty} \varpi^i [a_i]$ sur $\sum\limits_{i \gg - \infty} \varpi_D [a_i]$. On voit alors que 
\[ \left( \begin{array}{cc} \varpi^x & a \\ 0 & \varpi^y \end{array} \right) \mapsto \left( \begin{array}{cc} \varpi_D^x & \widetilde{a} \\ 0 & \varpi_D^y \end{array} \right) =: \widetilde{\mu}(\varpi) \widetilde{u} \] définit une bijection 
\[ \Theta : \mu(\varpi) \UE \xrightarrow{\sim} \widetilde{\mu}(\varpi) U ;\] $\Theta$ respecte les classes modulo $\UE \cap \KE$ et $U \cap K$ respectivement et passe donc au quotient pour induire 
\[ \mu(\varpi) \UE / \UE \cap \KE \xrightarrow{\sim} \widetilde{\mu}(\varpi) U / U \cap K \]
comme le montre le petit calcul
\[ \left( \begin{array}{cc} \varpi_D^x & \widetilde{a} \\ 0 & \varpi_D^y \end{array} \right) \left( \begin{array}{cc} 1 & \widetilde{b} \\ 0 & 1 \end{array} \right) = \left( \begin{array}{cc} \varpi_D^x & \widetilde{a} + \varpi_D^x \widetilde{b} \\ 0 & \varpi_D^y \end{array} \right) .\]
Il reste à prouver que $\Theta$ envoie bijectivement $\mu(\varpi) \UE \cap \KE \lambda(\varpi) \KE$ sur $\widetilde{\mu}(\varpi) U \cap K \widetilde{\lambda}(\varpi) K$ pour conclure. Or la théorie des diviseurs élémentaires nous donne (avec les notations de (\ref{notmum=2}), voir aussi appendice \ref{appdivelem})
\[ \mu(\varpi) u \in \KE \lambda(\varpi) \KE \Leftrightarrow \begin{cases} \min(x,y, d^{-1} v_D(a)) = \alpha \\ x+y = \alpha + \beta \end{cases} \Leftrightarrow \widetilde{\mu}(\varpi) \widetilde{u} \in K \widetilde{\lambda}(\varpi) K .\]
La preuve est terminée. \hfill$\Box$\\
\textsf{Preuve du cas $m=3$ :} \\
Nous allons construire une bijection
\begin{equation} \label{interclassesbij} \KE \lambda(\varpi) \KE \cap \mu(\varpi) \UE \xrightarrow{\sim} K \widetilde{\lambda}(\varpi) K \cap \widetilde{\mu}(\varpi) U \end{equation}
qui envoie un système de représentants de (\ref{CartanInterIwasawa}) sur un système de représentants de (\ref{CartanInterIwasawaD}), établissant alors l'égalité des cardinaux voulue. 
Ecrivons pour cela
\[ \mu(\varpi) u' = \left( \begin{array}{ccc} \varpi^x & a' & b' \\ 0 & \varpi^y & c' \\ 0 & 0 & \varpi^z \end{array} \right) \in \mu(\varpi) \UE , \quad \widetilde{\mu}(\varpi) u = \left( \begin{array}{ccc} \varpi_D^x & a & b \\ 0 & \varpi_D^y & c \\ 0 & 0 & \varpi_D^z \end{array} \right) \in \widetilde{\mu}(\varpi) U ,\]
où on écrit
\[ a = \varpi_D^{\overline{a}} a_0 + \varpi_D^{\overline{a}+1} a_1 + \dots , \quad b = \varpi_D^{\overline{b}} b_0 + \varpi_D^{\overline{b}+1} b_1 + \dots, \quad c= \varpi_D^{\overline{c}} c_0 + \varpi_D^{\overline{c}+1} c_1 + \dots \]
avec les $a_i, b_i, c_i \in \Oo_E \subseteq \Oo_D$ des représentants de Teichmüller et $a_0,b_0,c_0 \in \Oo_E^\times$. De même on écrit 
\[ a' = \varpi^{\overline{a}'} a_0' + \varpi^{\overline{a}'+1} a_1' + \dots , \quad b' = \varpi^{\overline{b}'} b_0' + \varpi^{\overline{b}'+1} b_1' + \dots, \quad c' = \varpi^{\overline{c}'} c_0' + \varpi^{\overline{c}'+1} c_1' + \dots .\]
De plus, on note 
\[ \lambda(\varpi) = \left( \begin{array}{ccc} \varpi^\alpha & 0 & 0 \\ 0 & \varpi^\beta & 0 \\ 0 & 0 & \varpi^\gamma \end{array} \right) , \quad \widetilde{\lambda}(\varpi) = \left( \begin{array}{ccc} \varpi_D^\alpha & 0 & 0 \\ 0 & \varpi_D^\beta & 0 \\ 0 & 0 & \varpi_D^\gamma \end{array} \right).\]
La théorie des diviseurs élémentaires nous dit 
\begin{equation}\hspace{-0.5cm} \label{divelemE} \mu(\varpi) u' \in \KE
  \lambda(\varpi) \KE \Leftrightarrow \left\{ \begin{array}{l} \min(x,y,z,
      \overline{a}', \overline{b}', \overline{c}') = \alpha \\ x+y+z =
      \alpha + \beta + \gamma \\ \min(x+y,y+z,x+z,z+\overline{a}',x+
      \overline{c}', d^{-1} v_D(a'c'-\varpi^y b')) = \alpha +
      \beta \end{array} \right. ,\end{equation}
ainsi que (voir appendice \ref{appdivelem})
\begin{equation}\hspace{-0.5cm} \label{divelemD} \widetilde{\mu}(\varpi) u
  \in K \widetilde{\lambda}(\varpi) K \Leftrightarrow
  \left\{ \begin{array}{l} \min(x,y,z, \overline{a}, \overline{b},
      \overline{c}) = \alpha \\ x+y+z = \alpha + \beta + \gamma \\
      \min(x+y,y+z,x+z,z+\overline{a},x+ \overline{c}, d^{-1} v_D \Big(
      \det_G \left( \begin{array}{cc} a & b \\ \varpi_D^y & c \end{array}
      \right) \Big) ) = \alpha + \beta \end{array} \right. .\end{equation}
Notons $\widetilde{\phantom{x}} : E \to D$ l'application qui envoie $a'$ sur $a$ si $a$ est tel que $\overline{a} = \overline{a}'$ et $a_i = a_i'$ pour tout $i \geq 0$. Les deux premières conditions de (\ref{divelemE}) et (\ref{divelemD}) sont identiques si l'on prend $a = \widetilde{a'}$, $c = \widetilde{c'}$ et $\overline{b} = \overline{b}'$. On compte ensuite le nombre de $b$ avec $\overline{b} = \overline{b}'$ qui sont tels que les troisièmes conditions de (\ref{divelemE}) et (\ref{divelemD}) sont identiques. L'application $\widetilde{\phantom{x}} : E \to D$ étant bijective, si on arrive à faire correspondre à chaque $b'$ un et un seul $b$ comme précédemment, cela nous fournira la bijection (\ref{interclassesbij}) voulue. Commençons par remarquer que si $c'=0$ (et donc $c= \widetilde{c'} = 0$), alors la dernière condition de (\ref{divelemE}) et (\ref{divelemD}) est respectivement 
\[ \min(x+y,y+z,x+z,z+\overline{a}, y + \overline{b}') = \alpha +  \beta \quad \textrm{et} \quad \min(x+y,y+z,x+z,z+\overline{a}, y + \overline{b}) = \alpha + \beta .\]
On prend alors $b' \mapsto \widetilde{b'}$ pour assurer la bijection voulue. On suppose à présent $c \neq 0$ (et donc $c' \neq 0$). Parce que l'on a 
\[ \left( \begin{array}{cc} a & b \\ \varpi_D^y  & c \end{array} \right) \left( \begin{array}{cc} 1 & 0 \\ -c^{-1} \varpi_D^y  & 1 \end{array} \right) = \left( \begin{array}{cc} a-bc^{-1} \varpi_D^y & b \\ 0 & c \end{array} \right) ,\]
il s'ensuit
\begin{equation} \label{calculdet2} \det\nolimits_G \Big( \left( \begin{array}{cc} a & b \\ \varpi_D^y  & c \end{array} \right) \big) = \Nrm (ac-bc^{-1} \varpi_D^y c) .\end{equation}
Lorsque l'on a $\overline{a}+\overline{c} \neq y + \overline{b}$, les troisièmes conditions de (\ref{divelemE}) et (\ref{divelemD}) sont équivalentes sous $a = \widetilde{a'}$, $c = \widetilde{c'}$ et $\overline{b} = \overline{b}'$. On prend alors à nouveau $b' \mapsto \widetilde{b'}$ dans ce cas-là. \\
On suppose à présent $c \neq 0$ et on va étudier le cas $\overline{a}+\overline{c} = y + \overline{b}$. On écrit 
\[bc^{-1} \varpi_D^y c = \varpi_D^{\overline{b}+y} x_0 + \varpi_D^{\overline{b}+y+1} x_1 + \varpi_D^{\overline{b}+y+2} x_2 + \dots \]
avec les $x_i$ des Teichmüller dans $\Oo_E$. Le point important est que chaque $x_i$ est la somme d'un terme $\varpi_D^y b_i \varpi_D^y x_0'$ avec $x_0' \in \Oo_E^\times$ et de termes ne faisant intervenir, parmi les coefficients de $b$, que des $b_j$ pour $j >i$. Cela nous permet d'affirmer que 
\begin{equation} \label{bijpourD} \begin{array}{cccc} \Phi_D: & \{ b \in D \ | \ v_D(b) \geq \overline{a}+ \overline{c} -y \} & \to & \{ x \in D \ | \ v_D(x) \geq \overline{a}+ \overline{c} \} \\ & b & \mapsto & ac - bc^{-1} \varpi_D^y c \end{array} \end{equation}
est une bijection. Et de la même manière on a
\begin{equation} \label{bijpourE} \begin{array}{cccc} \Phi_E : & \{ b \in E \ | \ d^{-1} v_D(b) \geq \overline{a}+ \overline{c} -y \} & \xrightarrow{\sim} & \{ x \in E \ | \ d^{-1} v_D(x) \geq \overline{a}+ \overline{c} \} \\ & b & \mapsto & ac - \varpi^y b \end{array} .\end{equation}
Grâce à (\ref{bijpourD}) et (\ref{bijpourE}), on a la bijection
\[ \begin{array}{cccc} \Phi: & \{ b' \in E \ | \ d^{-1} v_D(b') = \overline{a}+\overline{c}-y \} & \to & \{ b \in D \ | \ v_D(b) = \overline{a}+\overline{c}-y \} \\ & b' & \mapsto & \Phi_D^{-1}\circ \widetilde{\cdot} \circ \Phi_E (b') \end{array} ,\]
qui fait en sorte que les troisièmes conditions de (\ref{divelemE}) et (\ref{divelemD}) soient équivalentes pour $b'$ et $b = \Phi(b')$. De plus, cette bijection $\Phi$ vérifie 
\begin{equation} \label{compatibmodulo} b'_1 = b'_2 \mod \varpi^{\overline{a}+\overline{c}-y+k} \Leftrightarrow \Phi(b'_1) = \Phi(b'_2) \mod \varpi_D^{\overline{a}+\overline{c}-y+k} \end{equation}
pour tout $k \geq 1$ (cela vient du fait que (\ref{bijpourD}) (resp. (\ref{bijpourE})), vue comme application des coefficients du développement $\varpi_D$-adique (resp. $\varpi$-adique) de $b$, est donnée par une matrice (infinie) triangulaire supérieure: c'est le point important souligné précédemment). \\
La bijection $E \xrightarrow{\sim} D$ définie par ce $b' \mapsto \Phi(b')$ pour $\overline{b}' = \overline{a} + \overline{c}-y$ et par $b' \mapsto \widetilde{b'}$ autrement donne alors une bijection (en prenant $a' \mapsto \widetilde{a'}$ et $c' \mapsto \widetilde{c'}$)
\[ \Theta : \mu(\varpi) \UE \xrightarrow{\sim} \widetilde{\mu}(\varpi) U \]
qui se restreint en la bijection (\ref{interclassesbij}) annoncée. \\
Pour pouvoir comparer les cardinaux de (\ref{CartanInterIwasawa}) et de (\ref{CartanInterIwasawaD}), il reste à choisir un ensemble de représentants de (\ref{CartanInterIwasawa}) et de voir que son image par la bijection (\ref{interclassesbij}) constitue un ensemble de représentants de (\ref{CartanInterIwasawaD}).
Notons, pour tout $\zeta \in \Z$:
\[ \Rr^{(\zeta)}_{(E)} = \{ 0\} \cup \big\{ \varpi^{\overline{a}} [x_{\overline{a}}] + \varpi^{\overline{a}+1} [x_{\overline{a}+1}] + \dots + \varpi^{\zeta-1} [x_{\zeta-1}] \ \big| \ \overline{a}<\zeta , x_{\overline{a}} \in k_D^\times, x_{\overline{a}+1} , \dots, x_{\zeta-1} \in k_D \big\} ,\]
\[ \Rr^{(\zeta)} = \{ 0\} \cup \big\{ \varpi_D^{\overline{a}} [x_{\overline{a}}] + \varpi_D^{\overline{a}+1} [x_{\overline{a}+1}] + \dots + \varpi_D^{\zeta-1} [x_{\zeta-1}] \ \big| \ \overline{a}<\zeta , x_{\overline{a}} \in k_D^\times, x_{\overline{a}+1} , \dots, x_{\zeta-1} \in k_D \big\}.\]
Regardons la classe de $\widetilde{\mu}(\varpi) u$ dans $\widetilde{\mu}(\varpi) U / U \cap K$: prenons $g,h,i \in \Oo_D$ et calculons
\[ \left( \begin{array}{ccc} \varpi_D^x & a & b \\ 0 & \varpi_D^y & c \\ 0 & 0 & \varpi_D^z \end{array} \right) \left( \begin{array}{ccc} 1 & g & h \\ 0 & 1 & i \\ 0 & 0 & 1 \end{array} \right) = \left( \begin{array}{ccc} \varpi_D^x & a+ \varpi_D^x g & b + \varpi_D^x h + ai \\ 0 & \varpi_D^y & c + \varpi_D^y i \\ 0 & 0 & \varpi_D^z \end{array} \right) .\]
Le calcul similaire est valable pour $\mu(\varpi)u'$ et on peut fixer l'ensemble suivant de représentants de $\mu(\varpi) \UE / \UE \cap \KE$:
\[ \Rr^\mu_{(E)} := \Bigg\{ \left( \begin{array}{ccc} \varpi^x & a & b \\ 0 & \varpi^y & c \\ 0 & 0 & \varpi^z \end{array} \right) \ \Bigg| \ a \in \Rr^{(x)}_{(E)} , c \in \Rr^{(y)}_{(E)}, b \in \Rr^{(\min(x,\frac{v_D(a)}{d}))}_{(E)} \Bigg\}.\]
Parce que $\Theta$ respecte les valuations des coefficients surdiagonaux (notamment (\ref{compatibmodulo})), elle induit une bijection de $\Rr^\mu_{(E)}$ sur l'ensemble
\[ \Rr^{\widetilde{\mu}} := \Bigg\{ \left( \begin{array}{ccc} \varpi_D^x & a & b \\ 0 & \varpi_D^y & c \\ 0 & 0 & \varpi_D^z \end{array} \right) \ \Bigg| \ a \in \Rr^{(x)} , c \in \Rr^{(y)}, b \in \Rr^{(\min(x,v_D(a)))} \Bigg\} \]
de représentants de $\widetilde{\mu}(\varpi) U / U \cap K$. Parce que de plus $\Theta$ se restreint en la bijection (\ref{interclassesbij}), elle induit aussi 
\[ \Rr^\mu_{(E)} \cap \KE \lambda(\varpi) \KE \xrightarrow{\sim} \Rr^{\widetilde{\mu}} \cap K \widetilde{\lambda}(\varpi) K .\]
Cela nous donne alors l'égalité de cardinaux voulue. \hfill$\Box$\\

On espère que la proposition \ref{beforeLK} soit vraie sans restriction sur $m$. C'est ce que l'on formule dans la conjecture suivante; cette dernière est équivalente à l'égalité des cardinaux des intersections des doubles classes de Cartan et des classes d'Iwasawa comme dans le lemme \ref{egalcardinaux}.

\begin{conj}[formule de comparaison de transformées de Satake] \label{conjSatake}
Soit $\lambda$ un élément de $X_*(T)^+$. On a \[ \Sss_{(E)} \big(\id_{\KE \lambda(\varpi) \KE} \big) = \overline{\Sss} \big(\id_{K \widetilde{\lambda}(\varpi) K} \big) .\]
\end{conj}

On rappelle les définitions des groupes de Weyl étendus pour $(G,T)$ et $(G_{(E)},T_{(E)})$ respectivement: \[ W_e = N_G(T) / Z_G(T) \cap K, \quad W_e^{(E)} = N_{G_{(E)}}(T_{(E)}) / Z_{G_{(E)}}(T_{(E)}) \cap K_{(E)} ,\] où $N_G( \cdot)$ et $Z_G( \cdot)$ désignent respectivement normalisateur et centralisateur. 

\begin{lem} \label{Weylextend}
Les groupes $W_e^{(E)}$ et $W_e$ sont canoniquement isomorphes.
\end{lem}
\textsf{Preuve :} \\
On confond le groupe de Weyl fini $W$ avec le sous-groupe de $G_{(E)} \leq G$ des matrices de permutation. En particulier, il s'identifie au groupe de Weyl fini associé à $(G_{(E)},T_{(E)})$. Tout élément de $W_e^{(E)}$ possède un représentant de la forme $w x(\varpi) \in N_{G_{(E)}} (T_{(E)})$ pour $w \in W$ et $w \in X_*(T)$. La flèche $[w x(\varpi)] \mapsto [w \widetilde{x}(\varpi)]$ est alors bien définie (voir paragraphe \ref{oneparamsbgp}) et induit un isomorphisme canonique entre $W_e^{(E)}$ et $W_e$. \hfill$\Box$\\

Grâce à ce lemme, on peut faire hériter à $W_e$ des résultats bien connus pour $W_e^{(E)}$, que l'on rappelle maintenant. Le groupe de Weyl étendu admet une décomposition $W_e = \Lambda \rtimes W$, où $\Lambda$ peut se relever dans $G$ par $A_\Lambda$. L'action de $W$ sur $\Lambda$ dans ce produit semi-direct est la conjugaison. Aussi, $W_e$ admet une décomposition faisant appel au groupe de Weyl affine $W_a$ (voir \cite{BusKut93}, paragraphe 5.4), $W_e = W_a \rtimes \langle \Omega \rangle$, où $\Omega$ est un élément d'ordre $m$ de $W_e$. Le groupe de Weyl $W$ est un groupe de Coxeter fini admettant $S_\Delta$ comme système générateur minimal. Il existe une réflexion $s_\sim$ qui complète $S_\Delta$ en un système générateur minimal $S_\Delta^\sim = S_\Delta \cup \{ s_\sim \}$ du groupe de Coxeter $W_a$. Cela permet de définir une longueur $l$ sur $W_a$, que l'on étant à $W_e$ en définissant $l(w \Omega^x) = l(w)$ pour tout $w \in W_a$ et tout $x \in \Z$. \\
Soient $\rho$ la demi-somme des racines positives de $\Phi_G$ et $\lambda$ un copoids dominant. On rappelle que $w_0$ est l'élément le plus long de $W$. Alors $w_\lambda := \widetilde{w_0 \lambda}(\varpi) w_0$ est l'unique élément de longueur maximale de $W \widetilde{\lambda}(\varpi) W$ et il est de longueur $l(w_\lambda) = l(w_0)+2 \langle \rho, \lambda \rangle$ (voir la preuve de \cite{NelRam03}, Lemma 2.7). \\

Pour $x \leq y$ deux éléments de $W_e$, on note $P_{x,y}(v) \in 1 + v \Z[v]$ le polynôme de Kazhdan-Lusztig introduit dans \cite{KazLus79}. Si $\mu$ est un copoids dominant, on note $V_\mu$ le module de plus haut poids $\mu$ du groupe dual connexe $G^\vee = \GL(m,\C)$; et $\ch \, V_\mu$ désignera son caractère formel (voir \cite{Jan87} I.2.11, II.5.7): 
\begin{equation} \label{caracformel} \ch \, V_\mu = \sum_{\nu \in X_*(T)} (\dim \, V_\mu(\nu)) \, \tau_\nu ,\end{equation} 
où $V_\mu(\nu)$ désigne le sous-espace propre de poids $\nu$, et $\tau_\nu$ l'élément de $\C[X_*(T)]$ porté par $\nu$. \\
On peut maintenant énoncer la formule de Lusztig-Kato.

\begin{cor} \label{formuleLK}
Supposons la formule de comparaison de transformées de Satake. Soit $\mu$ un copoids dominant.  On a l'identité \[ \ch \, V_\mu = \sum_{\lambda \leq \mu} q^{-d \langle \rho, \mu \rangle} P_{w_\lambda,w_\mu}(q^d) \, \overline{\Sss} \big(\id_{K \widetilde{\lambda}(\varpi) K} \big) ,\] où $\lambda$ parcourt l'ensemble des copoids dominants $\leq$-inférieurs à $\mu$.
\end{cor}
\textsf{Preuve : } \\
Le groupe dual connexe de $\GE$ est aussi $\GE^\vee = \GL(m,\C)$. La formule de Lusztig-Kato pour $\GE$ s'écrit alors (voir \cite{HaiKotPra09}, Theorem 7.8.1): \[ \ch \, V_\mu = \sum_{\lambda \leq \mu} q^{-d \langle \rho, \mu \rangle} P_{w_\lambda,w_\mu}(q^d) \, \Sss_{(E)} \big(\id_{K_{(E)} \lambda(\varpi) K_{(E)}} \big) .\] Le résultat suit maintenant de la formule de comparaison de transformées de Satake. \hfill$\Box$

\section{Prolongement d'un caractère de $K$ à certaines doubles classes de Cartan} \label{prgcharacter}

Parce que $K$ est un compact maximal associé à un sommet spécial de $G$, la décomposition de Cartan nous donne (voir \cite{HaiRos10}, Theorem 1.0.3): 
\begin{equation} \label{decompCartan} G = \coprod_{t \in A_\Lambda^-} K t K .\end{equation} 
Tout caractère lisse de $K$ est trivial sur le pro-$p$-groupe $K(1)$ et se factorise donc par $\overline{K}$. Et comme tout caractère de $\overline{K}$ se factorise par le déterminant $\det$, il s'ensuit une correspondance bijective entre caractères lisses de $K$ et caractères de $k_D^\times$ via précomposition par $\overline{\det}$. \\

Quand on a un caractère $\chi_K : K \to \Fpbar^\times$, dans le cas déployé $D=F$, on peut l'étendre en un caractère de $G$ en utilisant le déterminant (c'est ce qui est fait au début de la preuve de la Proposition 5.1 de \cite{Her11b}). Ce n'est plus le cas dans le cas pour $\GL(m,D)$ en général, mais on peut se contenter du résultat suivant: on oublie que $\chi_K$ est un caractère et on peut le prolonger en une application ensembliste sur le support de l'algèbre de Hecke $\Hh_{\Fpbar}(G,K,\chi_K)$ de manière suffisamment naturelle pour que subsiste une propriété de multiplicativité (en un sens convenable) en petite dimension.

\begin{lem} \label{prolongmtcaract}
Soit $\chi$ un caractère $k_D^\times \to \Fpbar^\times$, que l'on voit comme un caractère $\chi_K$ de $K$ en précomposant par $\overline{\det}$. Soit $\rho$ un caractère de $\det_G(\widetilde{Z}_A(\chi_K))$. Alors $\chi_K$ s'étend naturellement en une application ensembliste \[ \widetilde{\chi}_\rho : \coprod_{t \in \widetilde{Z}_A^-(\chi_K)} K t K \to \Fpbar^\times \] par la formule 
\begin{equation} \label{defprlgt} \widetilde{\chi}_\rho ( k_1 t k_2 ) = \big( \chi_K (k_1 k_2) \big) \, \big( \rho \circ \det\nolimits_G(t) \big) \end{equation} pour tous $k_1, k_2 \in K$ et $t \in \widetilde{Z}_A^-(\chi_K)$.
\end{lem}
\textsf{Remarque :} 
Le stabilisateur de $\chi$ dans $D^\times$ est caractérisé par le sous-groupe \[ \big\{ z \in \Z \ \big| \ \forall x \in k_D^\times, \, \chi \big( \overline{\varpi_D^z [x] \varpi_D^{-z}} \big) = \chi(x) \big\} = d_0 \Z \] de $\Z$ (pour $d_0 \geq 1$ un entier). Commençons par remarquer que, par la discussion avant la proposition \ref{HVsatake}, un élément $t \in A_\Lambda^-$ est dans $\widetilde{Z}_A^-(\chi_K)$ si et seulement si tous ses coefficients sont des puissances de $\varpi_D^{d_0}$. En particulier, $\det_G(\widetilde{Z}_A(\chi_K))$ est le groupe monogène $((-1)^{d-1} \varpi)^{d_0 \Z}$ et $\rho$ est un caractère uniquement déterminé par sa valeur en $((-1)^{d-1} \varpi)^{d_0}$. \\  
\textsf{Remarque :} 
La source de $\widetilde{\chi}_\rho$ n'est pas toujours un sous-groupe de $G$, par exemple dans $\GL(2,D)$: 
\[ \left( \begin{array}{cc} 1 & 0 \\ 0 & \varpi_D^2 \end{array} \right) \left( \begin{array}{cc} 1 & \varpi_D \\ 0 & 1 \end{array} \right) \left( \begin{array}{cc} \varpi_D^2 & 0 \\ 0 & 1 \end{array} \right) = \left( \begin{array}{cc} \varpi_D^2 & \varpi_D \\ 0 & \varpi_D^2 \end{array} \right) = \left( \begin{array}{cc} 1 & 0 \\ \varpi_D & 1 \end{array} \right) \left( \begin{array}{cc} \varpi_D & 0 \\ 0 & \varpi_D^3 \end{array} \right) \left( \begin{array}{cc} \varpi_D & 1 \\ -1 & 0 \end{array} \right) .\]
En effet, si $\chi$ est tel que l'on a $d_0=2$, alors le terme de droite n'appartient pas à $K \widetilde{Z}_A^-(\chi_K) K$ puisque l'union est disjointe dans la décomposition de Cartan (\ref{decompCartan}). \\
\textsf{Preuve 1 :} \\
Soit $t$ un élément de $\widetilde{Z}_A^-(\chi_K)$. On veut étendre $\chi_K = \chi \circ \overline{\det}$ sur $K t K$, et on appellera $\widetilde{\chi}_\rho$ son prolongement. Naturellement, on a envie de définir $\widetilde{\chi}_\rho(k_1 t k_2) = \chi_K(k_1 k_2) \, \rho \circ \det_G(t)$ pour tous $k_1, k_2 \in K$: il faut vérifier que cette définition a un sens. Si on a $k_1 t k_2 = k_3 t k_4$ pour d'autres éléments $k_3, k_4 \in K$, alors on a \[ k_1^{-1} k_3 = t k_2 k_4^{-1} t^{-1} \in K \cap t K t^{-1} .\] 
Il s'agit donc de voir que si $k_0$ est un élément de $K \cap t K t^{-1}$ (dans la situation présente $k_0 = k_1^{-1} k_3$), alors on a $\chi_K(k_0) = \chi_K(t^{-1} k_0 t)$. La réduction $\overline{k}_0$ de $k_0 \in K \cap t K t^{-1}$ appartient à un parabolique standard de Levi standard $\prod_i \GL(m_i,k_D)$ avec $\sum\limits_{i=1}^r m_i = m$ par le lemme \ref{Prop3.8}. La projection de $\overline{k}_0$ sur ce Levi est alors un élément $(g_i) \in \prod_i \GL(m_i,k_D)$ et on a \[ \chi_K(k_0) = \prod\nolimits_i \chi \circ \det (g_i) .\] 
Remarquons que l'application \[\begin{array}{cccc} \sigma: & x & \mapsto & \overline{\varpi_D^{-d_0} [x] \varpi_D^{d_0}} \\ & k_D & \to & k_D \end{array}\] est un morphisme d'anneaux, de sorte que $\det: \GL(m,k_D) \to k_D^\times$ vérifie $\det \circ \sigma = \sigma \circ \det$. Ensuite, $\overline{K \cap t^{-1} K t}$ est le parabolique opposé à $\overline{K \cap t K t^{-1}}$ et on a alors le diagramme suivant, où la flèche verticale est $x \mapsto t^{-1} x t$. 
\[\xymatrix{ K \cap t K t^{-1} \ar[r] \ar[dd] & \overline{K \cap t K t^{-1}} \ar[dr] & \\ & & \prod\nolimits_i \GL(m_i,k_D) \\ t^{-1} K t \cap K \ar[r] & \overline{t^{-1} K t \cap K} \ar[ur] & }\] 
Cela correspond à l'identité \[ \chi_K(t^{-1} k_0 t) = \prod\nolimits_i \chi \circ \det \big( \overline{\varpi_D^{-d_0 a_i} [g_i] \varpi_D^{d_0 a_i}} \big) = \prod\nolimits_i \chi \big( \overline{\varpi_D^{-d_0 a_i} [\det \, g_i] \varpi_D^{d_0 a_i}} \big) ,\] où $t$ a pour coefficients diagonaux par blocs $\varpi_D^{d_0 a_1}$, \dots, $\varpi_D^{d_0 a_r}$. Comme $\varpi_D^{d_0}$ est dans le stabilisateur de $\chi$, on a finalement $\chi_K(t^{-1} k_0 t) = \chi_K(k_0)$. Il s'ensuit que $\widetilde{\chi}_\rho$ est bien défini. \hfill$\Box$\\
\textsf{Preuve 2 :} \\
Suite \` a la suggestion de l'un des rapporteurs, on rajoute une autre preuve de ce r\' esultat, qui tire \` a profit le fait que $\coprod\limits_{t \in \widetilde{Z}_A^-(\chi_K)} K t K$ est le support de l' alg\` ebre de Hecke $\Hh_{\Fpbar}(G,K,\chi_K)$. \\
On va d\' efinir s\' epar\' ement $\widetilde{\chi}_\rho$ sur chaque double classe $KtK$ pour $t \in \widetilde{Z}_A(\chi_K)$. Soient donc $t$ un \' el\' ement de $\widetilde{Z}_A(\chi_K)$ et $f_t$ l'unique fonction de $\Hh_{\Fpbar}(G,K,\chi_K)$ de support $KtK$ et valant $\rho \circ \det_G(t)$ en $t$ (par la discussion avant la proposition \ref{Prop2.11}). On d\' efinit la restriction de $\widetilde{\chi}_\rho$ \` a $KtK$ comme \' etant identiquement \' egale \` a $f_t|_{KtK}$: la bonne d\' efinition de l'op\' erateur $f_t$ nous assure celle de $\widetilde{\chi}_\rho|_{KtK}$ ainsi que la formule (\ref{defprlgt}). La preuve est termin\' ee. \hfill$\Box$ \\

Le prolongement $\widetilde{\chi}_\rho$ que l'on vient de définir, bien que naturel, a le défaut de ne pas être multiplicatif (en un sens convenable). C'est le cas pour $m=2$, mais des contre-exemples existent dès que l'on a $m>2$. 
Aussi, un cas facile pour lequel le prolongement $\widetilde{\chi}_\rho$ est multiplicatif est le suivant.

\begin{lem} \label{multsiNrm}
Soit $\chi_K$ un caractère lisse de $K$ qui se factorise par $\det_G$. Alors $\widetilde{Z}_A^-(\chi_K)$ est $A_\Lambda^-$ et, pour tout caractère $\rho$ de $A_\Lambda$, le prolongement $\widetilde{\chi}_\rho$ défini au lemme \ref{prolongmtcaract} est un caractère de $G$.
\end{lem}
\textsf{Remarque :}
Une réciproque à cet énoncé, qui ne servira pas pour la suite, sera étudiée en détail dans l'appendice \ref{appendixmult}. \\
\textsf{Preuve :} \\
On écrit $\chi_K = \chi' \circ \det_G$ sur $K$ pour $\chi'$ un caractère de $\det_G(K) = \Oo_F^\times$. Il suit de la définition que $\widetilde{Z}_A^-(\chi_K)$ est tout $A_\Lambda^-$, et $\det_G(\widetilde{Z}_A(\chi_K))$ est simplement $((-1)^d \varpi)^{\Z}$. On peut alors prolonger $\chi'$ en un caractère $\chi_G$ de $F^\times$ par $\chi_G |_{((-1)^d \varpi)^{\Z}} = \rho$ et la formule (\ref{defprlgt}) devient \[ \widetilde{\chi}_\rho(k_1 t k_2) = \chi_G \circ \det\nolimits_G(k_1 t k_2) .\] Mais alors $\widetilde{\chi}_\rho$ coïncide sur tout $G$ avec le caractère $\chi_G \circ \det_G$ et le résultat suit. \hfill$\Box$ \\

Soient $\chi$ un caractère de $k_D^\times$ et $\rho$ un caractère de $\det_G(\widetilde{Z}_A(\chi \circ \overline{\det}))$. On dira que l'application $\widetilde{\chi}_\rho$ est \textit{pseudo-multiplicative} si elle vérifie la propriété suivante: si $a$, $b$ et $ab$ sont dans $\coprod\limits_{t \in \widetilde{Z}_A^-(\chi \circ \overline{\det})} KtK$, alors on a $\widetilde{\chi}_\rho(ab) = \widetilde{\chi}_\rho(a) \widetilde{\chi}_\rho(b)$. \\
On remarquera que si la source de $\widetilde{\chi}_\rho$ est un groupe et que $\widetilde{\chi}_\rho$ est un caractère, alors en particulier $\widetilde{\chi}_\rho$ est pseudo-multiplicatif. C'est ce qui se passe dans le cas $m=1$.

\begin{lem}
Supposons $m=1$. Pour tout caractère $\chi$ de $\Oo_D^\times$ et tout caractère $\rho$ de $\Nrm(\widetilde{Z}_A(\chi))$, le support de $\widetilde{\chi}_\rho$ est un sous-groupe de $D^\times$ et $\widetilde{\chi}_\rho$ est un caractère.
\end{lem}
\textsf{Preuve :} \\
Soit $d_0 \geq 1$ le diviseur de $d$ vérifiant \[ \big\{ z \in \Z \ \big| \ \forall x \in k_D^\times, \, \chi \big( \overline{\varpi_D^z [x] \varpi_D^{-z}} \big) = \chi(x) \big\} = d_0 \Z .\] Parce que la valuation $v_D: D^\times \to \Z$ est un morphisme de groupes, le support de $\widetilde{\chi}_\rho$ est le sous-groupe des éléments de $D^\times$ de valuation dans $d_0 \Z$. Pour tous $x_1, x_2, x_3, x_4 \in \Oo_D^\times$ et $y_1, y_2 \in \Z$, on a par définition de $\widetilde{\chi}_\rho$: 
\begin{equation} \label{chitilde1} \widetilde{\chi}_\rho (x_1 \varpi_D^{d_0 y_1} x_2) = \chi( \overline{x_1 x_2}) \, \rho \big( ((-1)^d \varpi)^{d_0} \big)^{y_1} ,\end{equation} \begin{equation} \label{chitilde2} \widetilde{\chi}_\rho (x_3 \varpi_D^{d_0 y_2} x_4) = \chi( \overline{x_3 x_4}) \, \rho \big( ((-1)^d \varpi)^{d_0} \big)^{y_2} .\end{equation} On écrit ensuite \[ x_1 \varpi_D^{d_0 y_1} x_2 x_3 \varpi_D^{d_0 y_2} x_4 = x_1 \varpi_D^{d_0(y_1+y_2)} (\varpi_D^{-d_0 y_2} x_2 x_3 \varpi_D^{d_0 y_2}) x_4 ,\] pour avoir l'expression suivante de $\widetilde{\chi}_\rho(x_1 \varpi_D^{d_0 y_1} x_2 x_3 \varpi_D^{d_0 y_2} x_4)$: \begin{equation} \label{chitilde3} \chi \big( \overline{x_1  (\varpi_D^{-d_0 y_2} x_2 x_3 \varpi_D^{d_0 y_2}) x_4} \big) \rho \big( ((-1)^d \varpi)^{d_0} \big)^{y_1 + y_2} .\end{equation} Puisque $\varpi_D^{d_0 \Z}$ stabilise $\chi$, l'expression (\ref{chitilde3}) se réécrit enfin \[ \chi(\overline{x_1}) \, \chi(\overline{x_2 x_3}) \, \chi(\overline{x_4}) \, \rho \big( ((-1)^d \varpi)^{d_0} \big)^{y_1+y_2} .\] 
C'est bien là le produit des expressions dans (\ref{chitilde1}) et (\ref{chitilde2}). \hfill$\Box$\\

Pour $m=2$, la source de $\widetilde{\chi}_\rho$ n'est plus nécessairement un groupe comme l'atteste la remarque après le lemme \ref{prolongmtcaract}. Cependant $\widetilde{\chi}_\rho$ a le bon goût d'être pseudo-multiplicatif et ce sera essentiel pour notre travail pour $\GL(3,D)$. \\
On note $H_K$ le sous-groupe normal de $K$ engendré par les transvections dans $K$. Le fait que $H_K$ soit inclus dans le noyau de $\overline{\det}$ explique pourquoi il entre en jeu dans la compréhension de $\widetilde{\chi}_\rho$. 

\begin{lem} \label{t0m=2}
Supposons $m=2$. Soient $t_1$, $t_2$ des éléments de $A_\Lambda^-$ et $k_3$ un élément de $K$. Il existe $t_3^0 \in A \cap K$ et un unique $t_3 \in A_\Lambda^-$ vérifiant \[ t_1 k_3 t_2 \in H_K t_3 t_3^0 H_K \] et la propriété suivante: pour tout caractère $\chi_K$ de $K$, si $t_1$, $t_2$ et $t_3$ sont des éléments de $\widetilde{Z}_A^-(\chi_K)$, alors on a $\chi_K(t_3^0) = \chi_K(k_3)$.
\end{lem}
\textsf{Preuve :} \\
L'unicité provient de l'écriture en union disjointe de la décomposition de Cartan. \\ 
Commençons par écrire l'identité suivante qui nous sera de nombreuses fois utiles au cours de la preuve 
\begin{equation} \label{generateswap} s := \left( \begin{array}{cc} 0 & -1 \\ 1 & 0 \end{array} \right) = \left( \begin{array}{cc} 1 & -1 \\ 0 & 1 \end{array} \right) \left( \begin{array}{cc} 1 & 0 \\ 1 & 1 \end{array} \right) \left( \begin{array}{cc} 1 & -1 \\ 0 & 1 \end{array} \right) ,\end{equation} 
et nous indique que $s$ est un élément de $H_K$. Remarquons ensuite que tout élément $\left( \begin{array}{cc} \alpha & \beta \\ \gamma & \delta \end{array} \right)$ de $G$ avec $\alpha \neq 0$ possède la décomposition unique suivante 
\begin{equation} \label{bigcell} \left( \begin{array}{cc} \alpha & \beta \\ \gamma & \delta \end{array} \right) = \left( \begin{array}{cc} 1 & 0 \\ u & 1 \end{array} \right) \left( \begin{array}{cc} \eta & 0 \\ 0 & \nu \end{array} \right) \left( \begin{array}{cc} 1 & v \\ 0 & 1 \end{array} \right) ,\end{equation} avec $\eta = \alpha$, $u = \gamma \alpha^{-1}$, $v = \alpha^{-1} \beta$ et $\nu = \delta - \gamma \alpha^{-1} \beta$. Après cette parenthèse, écrivons explicitement les éléments \[ k_3 = \left( \begin{array}{cc} \alpha & \beta \\ \gamma & \delta \end{array} \right), \quad t_1 = \left( \begin{array}{cc} \varpi_D^a & 0 \\ 0 & \varpi_D^b \end{array} \right), \quad t_2 = \left( \begin{array}{cc} \varpi_D^x & 0 \\ 0 & \varpi_D^y \end{array} \right) ,\] avec $\alpha, \beta, \gamma, \delta \in \Oo_D$ et $a, b, x, y \in \Z$. On va déterminer $t_3$ et $t_3^0$ pour \[ t_1 k_3 t_2 = \left( \begin{array}{cc} \varpi_D^a \alpha \varpi_D^x & \varpi_D^a \beta \varpi_D^y \\ \varpi_D^b \gamma \varpi_D^x & \varpi_D^b \delta \varpi_D^y \end{array} \right) .\] On rappelle que l'on a supposé $t_1, t_2 \in A_\Lambda^-$, c'est-à-dire $a \leq b$ et $x \leq y$. \\ 
On a deux cas distincts à envisager. D'abord, lorsque $\alpha$ est non nul, on profite de (\ref{bigcell}) pour écrire 
\begin{equation} \label{speccell} \left( \begin{array}{cc} 1 & 0 \\ - \varpi_D^b \gamma \alpha^{-1} \varpi_D^{-a} & 1 \end{array} \right) t_1 k_3 t_2 \left( \begin{array}{cc} 1 & - \varpi_D^{-x} \alpha^{-1} \beta \varpi_D^y \\ 0 & 1 \end{array} \right) = \left( \begin{array}{cc} \varpi_D^a \alpha \varpi_D^x & 0 \\ 0 & \varpi_D^b (\delta - \gamma \alpha^{-1} \beta) \varpi_D^y \end{array} \right).\end{equation} 
On remarque que les matrices de transvection dans (\ref{speccell}) sont dans $H_K$ lorsque l'on impose $b-a - v_D(\alpha) \geq 0$ et $y-x - v_D(\alpha) \geq 0$. Sous ces conditions, on note \[ t_3 = \left( \begin{array}{cc} \varpi_D^{a+x+v_D(\alpha)} & 0 \\ 0 & \varpi_D^{b+y-v_D(\alpha)} \end{array} \right) \in A_\Lambda^- ,\] \begin{equation} \label{t31} t_3^0 = \left( \begin{array}{cc} \varpi_D^{-x-v_D(\alpha)} \alpha \varpi_D^x & 0 \\ 0 & \varpi_D^{-y+v_D(\alpha)} (\delta - \gamma \alpha^{-1} \beta) \varpi_D^y \end{array} \right) \in A \cap K .\end{equation}
On traite maintenant les cas restants: 
\begin{itemize}
\item[(a)] $\alpha = 0$;
\item[(b)] $v_D(\alpha) +x-y>0$ et $\alpha \neq 0$;
\item[(c)] $v_D(\alpha) +a-b>0$ et $\alpha \neq 0$.
\end{itemize}
On remarque tout de suite que, comme $k_3$ est un élément de $K$ et que l'on a $\alpha = 0$ ou bien $v_D(\alpha)>0$, cela implique $\gamma \in \Oo_D^\times$. On commence par écrire \[ g_3 := t_1 k_3 t_2 \left( \begin{array}{cc} 1 & - \varpi_D^{-x} \gamma^{-1} \delta \varpi_D^y \\ 0 & 1 \end{array} \right) ,\] où la matrice de transvection est bien dans $H_K$ car on a $\gamma \in \Oo_D^\times$ et $y-x \geq 0$. On note \[ g_3 = \left( \begin{array}{cc} \varpi_D^a \alpha \varpi_D^x & \varpi_D^a (\beta - \alpha \gamma^{-1} \delta) \varpi_D^y \\ \varpi_D^b \gamma \varpi_D^x & 0 \end{array} \right) =: \left( \begin{array}{cc} \gamma_1 & \gamma_2 \\ \gamma_3 & 0 \end{array} \right) .\]
Les cas (a), (b) et (c) correspondent alors respectivement à \[ \gamma_1 = 0, \quad v_D(\gamma_1) > v_D(\gamma_2), \quad v_D(\gamma_1) > v_D(\gamma_3) .\] 
Dans le cas (b), on utilise l'identité \begin{equation} \label{lem1casb} \left( \begin{array}{cc} \gamma_1 & \gamma_2 \\ \gamma_3 & 0 \end{array} \right) \left( \begin{array}{cc} 1 & 0 \\ - \gamma_2^{-1} \gamma_1 & 1 \end{array} \right) = \left( \begin{array}{cc} 0 & \gamma_2 \\ \gamma_3 & 0 \end{array} \right) ,\end{equation}
et dans le cas (c) l'identité
\begin{equation} \label{lem1casc} \left( \begin{array}{cc} 1 & - \gamma_1 \gamma_3^{-1} \\ 0 & 1 \end{array} \right) \left( \begin{array}{cc} \gamma_1 & \gamma_2 \\ \gamma_3 & 0 \end{array} \right) = \left( \begin{array}{cc} 0 & \gamma_2 \\ \gamma_3 & 0 \end{array} \right) .\end{equation}
Enfin, dans tous les cas, on utilise $s \in H_K$ (voir (\ref{generateswap})) et l'identité \[ \left( \begin{array}{cc} 0 & \gamma_2 \\ \gamma_3 & 0 \end{array} \right) \left( \begin{array}{cc} 0 & -1 \\ 1 & 0 \end{array} \right) = \left( \begin{array}{cc} \gamma_2 & 0 \\ 0 & -\gamma_3 \end{array} \right) \] pour voir que $t_1 k_3 t_2$ est dans la même classe que \[ \left( \begin{array}{cc} \varpi_D^a (\beta - \alpha \gamma^{-1} \delta) \varpi_D^y & 0 \\ 0 & -\varpi_D^b \gamma \varpi_D^x \end{array} \right) \] dans $H_K \backslash G / H_K$. On pose alors \[ t_3 = \left( \begin{array}{cc} \varpi_D^{a+y} & 0 \\ 0 & \varpi_D^{b+x} \end{array} \right) \in A_\Lambda ,\] \begin{equation} \label{t32} t_3^0 = \left( \begin{array}{cc} \varpi_D^{-y} (\beta - \alpha \gamma^{-1} \delta) \varpi_D^y & 0 \\ 0 & -\varpi_D^{-x} \gamma \varpi_D^x \end{array} \right) \in A \cap K \end{equation}
Lorsque l'on a $b+x-a-y \geq 0$, la paire $(t_3^0,t_3)$ fait l'affaire; sinon, on prend $(- s t_3^0 s, - s t^3 s)$. \\
Vérifions maintenant que $t_3^0$ possède la propriété escomptée; il en est alors de même de $-s t_3^0 s$. Soit $\chi : k_D^\times \to \Fpbar^\times$ un caractère de normalisateur $\varpi_D^{d_0 \Z}$ avec $d_0 \geq 1$ divisant $d$, $\chi_K = \chi \circ \overline{\det}$ et supposons $t_1, t_2, t_3 \in \widetilde{Z}_A^-(\chi_K)$. \\ 
Dans le cas (\ref{t31}), on a \[ \chi \circ \overline{\det} (t_3^0) = \chi \big( \overline{\varpi_D^{-x-v_D(\alpha)} \alpha \varpi_D^x} \big) \, \chi \big( \overline{\varpi_D^{-y+v_D(\alpha)} (\delta - \gamma \alpha^{-1} \beta) \varpi_D^y} \big) .\] Parce que l'on a $t_1, t_2, t_3 \in \widetilde{Z}_A^-(\chi_K)$, $x, y$ et $v_D(\alpha)$ sont tous trois dans $d_0 \Z$. Mais alors, on a \[ \chi \circ \overline{\det} (t_3^0) = \chi \big( \overline{\alpha \varpi_D^{-v_D(\alpha)}} \big) \, \chi \big( \overline{\varpi_D^{v_D(\alpha)} (\delta - \gamma \alpha^{-1} \beta)} \big) ,\] et il s'ensuit \[ \chi \circ \overline{\det} (t_3^0) = \chi \big( \overline{\alpha(\delta-\gamma \alpha^{-1} \beta)}\big) = \chi \circ \overline{\det} (k_3) .\]
Dans le cas (\ref{t32}), on a de la même manière \[ \chi \circ \overline{\det} (t_3^0) = \chi \big( \overline{\varpi_D^{-y} (\alpha \gamma^{-1} \delta - \beta) \varpi_D^y} \big) \, \chi \big( \overline{\varpi_D^{-x} \gamma \varpi_D^x} \big) ,\] qui donne bien \[ \chi \circ \overline{\det} (t_3^0) = \chi \big( \overline{(\alpha \gamma^{-1} \delta - \beta)\gamma}\big) = \chi \circ \overline{\det} (k_3) .\]
Le lemme est donc bien prouvé. \hfill$\Box$

\begin{prop} \label{multm=2}
Supposons $m=2$. Soient $\chi$ un caractère de $k_D^\times$ et $\rho$ un caractère de $\det_G(\widetilde{Z}_A(\chi \circ \overline{\det}))$. L'application $\widetilde{\chi}_\rho$ est pseudo-multiplicative.
\end{prop}
\textsf{Remarque :}
Cette proposition est uniquement un corollaire du lemme \ref{t0m=2}, au sens où si on arrive à prouver le lemme sous des hypothèses autres que \og $m=2$ \fg , alors la conclusion de la proposition suit sous les mêmes hypothèses. \\
\textsf{Preuve :} \\
On note $\chi_K = \chi \circ \overline{\det}$. Comme $\widetilde{\chi}_\rho$ est bien définie sur chaque double classe $K t K$ avec $t \in \widetilde{Z}_A^-(\chi_K)$, il s'agit de montrer que si $t_1, t_2, t_3 \in \widetilde{Z}_A^-(\chi_K)$ et $k_1, k_2, k_3 \in K$ vérifient \begin{equation} \label{egaldoubleclasses} t_1 k_3 t_2 = k_1 t_3 k_2 \end{equation}
alors on a \begin{equation} \label{multfinale} \chi_K(k_3) \, \rho \circ \det\nolimits_G (t_1 t_2) = \chi_K (k_1 k_2) \, \rho \circ \det\nolimits_G (t_3) .\end{equation} 
Parce que $\det_G : G \to F^\times$ est un morphisme de groupes et qu'il envoie $K$ dans $\Oo_F^\times$, (\ref{egaldoubleclasses}) implique \[ \det\nolimits_G (t_1 t_2) \in \det\nolimits_G (t_3) \, \Oo_F^\times .\] Mais, parce que $t_1,t_2,t_3$ sont dans $A_\Lambda$ et à cause de la définition de $\det_G$, on a $\det_G(t_1 t_2) = \det_G(t_3)$. Reste à voir que l'on a $\chi_K (k_3) = \chi_K (k_1 k_2)$ sous (\ref{egaldoubleclasses}). \\
Cela découle du lemme \ref{t0m=2} de la manière suivante: soient $t_3' \in A_\Lambda^-$ et $t_3^0 \in A \cap K$ fournis par le lemme \ref{t0m=2}. Comme la décomposition de Cartan est une union de doubles classes disjointes, on a $t_3 = t_3'$. De plus, $H_K$ est normalisé par $K$ et il existe donc (toujours par le lemme \ref{t0m=2}) des $h_1, h_2 \in H_K$ et $t_1^0, t_2^0 \in A \cap K$ vérifiant $k_1 = h_1 t^0_1$ et $k_2 = t^0_2 h_2$. En particulier, on a \[ \overline{\det}(k_1) = \overline{\det}(t^0_1) \quad \textrm{et} \quad \overline{\det}(k_2) = \overline{\det}(t^0_2) .\] On a alors, pour certains $h_3, h_4 \in H_K$: \[ t_3 t_3^0 = h_3 t_1^0 t_3 t_2^0 h_4 = h_3 t_3 (t_3^{-1} t_1^0 t_3 t_2^0) h_4 .\] Parce que $\widetilde{\chi}_\rho$ est bien défini (voir lemme \ref{prolongmtcaract}), on a alors \[ \chi_K (t_3^0) = \chi_K (t_3^{-1} t_1^0 t_3) \, \chi_K (t_2^0) = \chi_K (k_1 k_2) .\] Et cette quantité est bien $\chi_K (k_3)$ par la propriété fondamentale de $t_3^0$ (voir lemme \ref{t0m=2}). \hfill$\Box$ \\

On établit maintenant les propriétés que l'on va vraiment utiliser du prolongement $\widetilde{\chi}_\rho$.

\begin{cor} \label{utilchirho}
Soient $\chi$ un caractère de $k_D^\times$ et $\rho$ un caractère de $\det_G \big( \widetilde{Z}_A( \chi \circ \overline{\det} )\big)$. Soient $t \in \widetilde{Z}_A^-( \chi \circ \overline{\det} )$ et $u \in U$.
\begin{itemize}
\item[(i)] On a $t^{-1} \in K \widetilde{Z}_A^-( \chi \circ \overline{\det} ) K$ et $\widetilde{\chi}_\rho(t^{-1}) = \widetilde{\chi}_\rho(t)^{-1}$.
\item[(ii)] Supposons $\widetilde{\chi}_\rho$ pseudo-multiplicative. Si $tu$ est dans $K \widetilde{Z}_A^-( \chi \circ \overline{\det} ) K$, alors on a $\widetilde{\chi}_\rho(tu) = \widetilde{\chi}_\rho(t)$.
\end{itemize}
\end{cor}
\textsf{Preuve :} \\
On note encore $w_0$ pour le représentant $\left( \begin{array}{ccc} & & 1 \\ & \cdot^{\cdot^{\cdot}} & \\ 1 & & \end{array} \right)$ de l'élément le plus long du groupe de Weyl fini $W$. C'est un élément de $K$, et on a \[ t^{-1} = w_0(w_0 t^{-1} w_0) w_0 \in K(w_0 t^{-1} w_0)K .\] Aussi, si on note \[ t = \Diag(\varpi_D^{a_1}, \varpi_D^{a_2}, \dots, \varpi_D^{a_m}) \in \widetilde{Z}_A^-(\chi \circ \overline{\det}) ,\] alors $w_0 t^{-1}w_0 = \Diag(\varpi_D^{-a_m}, \dots, \varpi_D^{-a_1})$ est aussi dans $\widetilde{Z}_A^-(\chi \circ \overline{\det})$. Et il s'ensuit \[ \widetilde{\chi}_\rho(t^{-1}) = \chi \circ \overline{\det}(w_0^2) \, \rho \circ \det\nolimits_G(w_0 t^{-1} w_0) = \rho \circ \det\nolimits_G(t^{-1}) = \widetilde{\chi}_\rho(t)^{-1} .\]
Prouvons maintenant \textit{(ii)}. Il existe $z \in A_\Lambda^+$ à coefficients diagonaux dans $\varpi_D^{d\Z}$ tel que $zuz^{-1}$ soit dans $U \cap K$. Ainsi $z^{-1}$ et $w_0 z w_0$ sont dans $\widetilde{Z}_A^-(\chi \circ \overline{\det})$, et on a \[ \widetilde{\chi}_\rho(tu) = \widetilde{\chi}_\rho \big( (t z^{-1})(z u z^{-1}) w_0 (w_0 z w_0) w_0 \big) .\]
Par pseudo-multiplicativité de $\widetilde{\chi}_\rho$ et le \textit{(i)}, on a : \[ \widetilde{\chi}_\rho(tu) = \widetilde{\chi}_\rho \big( (t z^{-1}) zuz^{-1} \big) \, \widetilde{\chi}_\rho \big( w_0 (w_0 z w_0) w_0 \big) = \widetilde{\chi}_\rho(tz^{-1}) \, \widetilde{\chi}_\rho(z^{-1})^{-1} = \widetilde{\chi}_\rho(t) .\] La preuve est terminée. \hfill$\Box$

\section{Inversion partielle de la transformée de Satake}  \label{invSatake}

Pour $\zeta \geq 1$ un entier, on note \begin{equation} \label{defR} \Rr_{(\zeta)} := \{ [x_0] + \varpi_D [x_1] + \cdots + \varpi_D^{\zeta-1} [x_{\zeta-1}] \ | \ x_0, \dots, x_{\zeta-1} \in k_D \} .\end{equation} C'est un ensemble de représentants de $\Oo_D/(\varpi_D^\zeta)$ et on a $\Rr_{(\zeta)} \subseteq \Rr_{(\zeta+1)}$. \\
A $\lambda \in X_*(T)^-$, on a associé un parabolique antistandard $P_\lambda$ et un sous-groupe de Levi $M_\lambda$. Si on écrit $t= \widetilde{\lambda}(\varpi)$, on va noter $W_t$ le sous-groupe de $W$ engendré par $\Delta_{M_\lambda}$. On définit aussi \[ W^t := \{ w \in W \ | \ \forall \alpha \in \Delta_{M_\lambda}, \ l(w s_\alpha) > l(w) \};\] c'est un système de représentants de $W/W_t$ contenant l'élément le plus court de chaque classe (voir chapitre \ref{Ly11}, Lemme 4.1). On confondra $W^t$ et un ensemble de relèvements de $W^t$ dans $K$ à coefficients dans $\{0,1\}$. 

\begin{lem} \label{reprIwasawa}
Soient $\zeta \geq 1$ et $t \in A_\Lambda^-$ de la forme $\Diag(1, \dots, 1, \varpi_D^\zeta, \dots, \varpi_D^\zeta)$ avec $r \in [1,m-1]$ occurrences de $1$. Soit $\Rr_t$ l'ensemble composé des matrices $w \left( \begin{array}{cc} I_r & 0 \\ M & I_{m-r} \end{array} \right)$ pour $M \in M_{m-r,r}(\Rr_{(\zeta)})$ et $w \in W^t$. On a alors \[ K = \Rr_t (K \cap t K t^{-1}) .\]
\end{lem}
\textsf{Remarque :} 
Attention, il peut y avoir dans $\Rr_t$ deux représentants de la même classe (comme le montre une simple comparaison des cardinaux). \\
\textsf{Preuve :} \\
On a tout d'abord (voir lemme \ref{Prop3.8}): \[ K \cap t K t^{-1} = \left( \begin{array}{cc} \GL(r,\Oo_D) & M_{r,m-r}(\Oo_D) \\ \varpi_D^\zeta M_{m-r,r}(\Oo_D) & \GL(m-r,\Oo_D) \end{array} \right) =: I^-_{P_t} .\] On remarque tout de suite que la réduction $\overline{I}^-_{P_t}$ se confond avec $\overline{P}^-_t$.
Ecrivons la décomposition de Bruhat pour la réduction modulo $(\varpi_D)$ de $K$ (ce qui est possible puisque c'est alors un $\GL(m)$ sur le corps fini $k_D$). On a \[ \overline{K} = \coprod_{w \in W} \overline{B} w \overline{B} = \coprod_{w \in W} w (w^{-1} \overline{B}w) \overline{B} .\] On peut alors remplacer dans cette écriture $w^{-1} \overline{B}w$ par un ensemble de représentants de \[ (w^{-1} \overline{B}w) / (\overline{B} \cap w^{-1} \overline{B}w) \simeq (w^{-1} \overline{U} w) / (\overline{U} \cap w^{-1} \overline{U} w) .\] 
Par \cite{Car85}, Proposition 2.5.12, on a \[ w^{-1} \overline{U}w = (w^{-1} \overline{U}w \cap \overline{U}^-) (w^{-1} \overline{U}w \cap \overline{U}) \] et donc \[ (w^{-1} \overline{B}w) /( w^{-1} \overline{B}w \cap \overline{B}) \simeq (w^{-1} \overline{U}w \cap \overline{U}^-) .\] Ce dernier groupe est visiblement inclus dans $\overline{P}_t$ et on obtient 
\begin{equation} \label{decompKbar} \overline{K} = \bigcup_{w \in W} w \overline{P}_t \overline{B} = \bigcup_{w \in W^t} w \overline{P}_t \overline{B} = \bigcup_{w \in W^t} w \overline{P}_t \overline{P}_t^- .\end{equation} 
Comme $K_{\zeta} = 1+ \varpi_D^\zeta M_m(\Oo_D)$ est normal dans $K \geq I_{P_t}$, on a \[ K/(K \cap tKt^{-1}) \xrightarrow{\sim} \overline{K}^{(\zeta)} \! / (\overline{P}_t^-)^{(\zeta)},\] où rappelle que l'on a noté $\overline{\phantom{x}}^{(\zeta)}$ la réduction modulo $K(\zeta)$. Dans le cas $\zeta=1$, le lemme est donc prouvé grâce à l'égalité (\ref{decompKbar}). Lorsque l'on a $\zeta>1$, on procède par récurrence comme suit: supposons alors l'énoncé du lemme vrai pour $\zeta-1 \geq 1$. On écrit $t_0 = \Diag(1, \dots, 1, \varpi_D, \dots, \varpi_D)$ de sorte que l'on ait $t=t_0^\zeta$. Par hypothèse de récurrence, on a alors 
\begin{equation} \label{reprfiltr} K = \Rr_{t_0^{\zeta-1}} (K \cap t_0^{\zeta-1} K t_0^{1-\zeta}) .\end{equation} 
On cherche ensuite à évaluer $(K \cap t_0^{\zeta-1}K t_0^{1-\zeta})/(K \cap tKt^{-1})$. Notons $\Nn$ le sous-ensemble des matrices de $M_{m-r,r}(\Rr_{(\zeta)})$ des matrices à coefficients dans $\varpi_D^{\zeta-1}[k_D] = \Rr_{(\zeta-1)} \backslash \Rr_{(\zeta)}$. Parce que l'on a \[ \left( \begin{array}{cc} I_r & 0 \\ E & I_{m-r} \end{array} \right) \left( \begin{array}{cc} A & B \\ \varpi_D^\zeta C & D \end{array} \right) = \left( \begin{array}{cc} A & B \\ EA+\varpi_D^\zeta C & EB+D \end{array} \right) \] pour $E \in \Nn, A \in \GL(r,\Oo_D), B \in M_{r,m-r}(\Oo_D) , C \in M_{m-r,r}(\Oo_D), D \in \GL(m-r,\Oo_D)$, on obtient l'injection ensembliste 
\[ \left( \begin{array}{cc} I_r & 0 \\ \Nn & I_{m-r} \end{array} \right) \hookrightarrow (K \cap t_0^{\zeta-1} K t_0^{1-\zeta})/(K \cap tKt^{-1}) .\] 
Parce que tous deux sont de cardinal $|k_D|^{r(m-r)}$, cette flèche est en fait une bijection. On revient ensuite à (\ref{reprfiltr}) et on obtient \[ K= \Rr_{t_0^{\zeta-1}} \left( \begin{array}{cc} I_r & 0 \\ \Nn & I_{m-r} \end{array} \right) (K \cap t K t^{-1}) = \Rr_t (K \cap tKt^{-1}) .\] La récurrence est terminée et le résultat prouvé. \hfill$\Box$\\

Lorsque $t$ est une homothétie valant $1$ ou $\varpi_D^\zeta$, on prend alors $\Rr_t = \{1\}$, de sorte que le résultat du lemme \ref{reprIwasawa} reste vrai pour $r \in [0,m]$. 

\begin{lem} \label{miracle}
Soient $t = \Diag(1, \dots, 1, \varpi_D^\zeta, \dots, \varpi_D^\zeta) \notin \{1,\varpi_D^\zeta\}$ et $t' \in A_\Lambda^-$ vérifiant $\det_G(t) = \det_G(t')$. Supposons que les coefficients diagonaux de $t'$ sont dans $\varpi_D^{\zeta \Z}$. Alors $KtK \cap t' U \neq \varnothing$ implique $t = t'$. Dans ce cas-là, on a $Kt K \cap t U = t(U \cap K)$.
\end{lem}
\textsf{Preuve :} \\
Soit $t' \in A_\Lambda^-$ de coefficients diagonaux $\varpi_D^{a_1}, \dots, \varpi_D^{a_m}$ dans $\varpi_D^{\zeta \Z}$. On va chercher à montrer que la condition $KtK \cap t' U \neq \varnothing$ implique $a_1 \geq 0$ et $a_m = \zeta$. Comme on a $\det_G(t) = \det_G(t')$ avec $t' \in A_\Lambda^-$, cela donnera tout de suite $t=t'$. \\
Commençons par montrer $a_m = \zeta$. On a, par le lemme \ref{reprIwasawa}, $Kt K = \Rr_t t K$. Soit $r \in \Rr_t$. Comme $A \cap K$ normalise $H_K$ et $U$, et que l'on a $(A \cap K) H_K = K$, si $rtK \cap t' U$ est non vide, alors il en est de même pour $rt H_K \cap t'(A \cap K) U$. Supposons donc l'existence d'un $h \in H_K$ tel que l'on ait $(t')^{-1} r t h \in (A \cap K) U$. Observons 
\begin{equation} \label{matrixgcd} (t')^{-1} rt = \Diag(\varpi_D^{a_1}, \varpi_D^{a_2}, \dots, \varpi_D^{a_m})^{-1} w \left( \begin{array}{cc} I_b & 0 \\ M & \varpi_D^{\zeta} I_{m-b} \end{array} \right) ,\end{equation} 
où les coefficients de $M$ sont dans $\Rr_{(\zeta)}$: ils sont donc soit nuls, soit de valuation dans $[0,\zeta-1]$. Aussi, dans cette écriture, on a $a_1 \leq a_2 \leq \dots \leq a_m$ dans $\zeta \Z$ et $w \in W^t$. Puisque le produit par une transvection à droite dans $K$ ne change pas le plus grand commun diviseur d'une ligne, la condition $(t')^{-1}rth \in (A \cap K) U$ implique que la dernière ligne de (\ref{matrixgcd}) possède un coefficient de valuation minimale dans $\Oo_D^\times$. On en déduit $a_m = \zeta$  (si on avait permut\' e la derni\` ere ligne de $\left( \begin{array}{cc} I_b & 0 \\ M & \varpi_D^{\zeta} I_{m-b} \end{array} \right)$ avec l'une de ses $b$ premi\` eres lignes, on aurait $a_m=0$, et donc $\det_G(t') \leq 0$, ce qui est absurde). \\
Montrons maintenant $a_1 \geq 0$. Si on avait $a_1 \leq - \zeta$, alors le pgcd des coefficients de la première ligne de $(t')^{-1} rt$, donc de $(t')^{-1} rth$, serait supérieur ou égal à $\zeta$. Cela contredirait $(t')^{-1} rth \in (A \cap K)U$, et on a donc $a_1>-\zeta$. Comme on a supposé que les coefficients non nuls de $t'$ sont dans $\varpi_D^{\zeta \Z}$ et aussi $a_1 \leq a_m = \zeta$, on a de fait $a_1 \in \{0,\zeta\}$. Comme de plus, on a $t' \in A_\Lambda^-$ et $\det_G(t) = \det_G(t')$, cela impose $t=t'$ comme voulu. Cela termine la preuve de la première partie du lemme. \\
Pour la seconde assertion, on peut faire appel à la Proposition 6.9 de \cite{HenVig11}, ou bien présenter la preuve ad hoc suivante. On a facilement \[ t(U \cap K) \subseteq KtK \cap tU, \quad tK \cap tU = t(U \cap K) .\] 
On veut montrer que si $r \in \Rr_t$ est tel que $rtK \cap tU \neq \varnothing$, alors on a $r=1$, ce que l'on va faire par récurrence sur $m$ comme suit. Pour $m=1$, on a $\Rr_t = \{1\}$ et c'est trivialement vrai. Supposons le résultat vrai pour $m-1 \geq 1$ et montrons le au rang $m$. Ecrivons $r = w \left( \begin{array}{cc} I_b & 0 \\ M & I_{m-b} \end{array} \right)$ pour un $w \in W^t$ et $M \in M_{m-b,b}(\Rr_{(\zeta)})$. Parce que la dernière ligne de $t^{-1}rt$ doit avoir un coefficient de valuation nulle, la dernière ligne de $r$ doit appartenir aux lignes de $\left( \begin{array}{cc} M & I_{m-b} \end{array} \right)$, disons $(m_{i 1}, m_{i 2}, \dots, m_{i b}, 0 , \dots, 0, 1, 0, \dots, 0)$, où le $1$ est à la colonne $(b+i)$. Aussi cela impose de plus $m_{i 1} = m_{i 2} = \dots = m_{i b}= 0$. Si $b=m-1$, cela prouve déjà $r=1$; sinon, on peut appliquer l'hypothèse de récurrence en considérant pour $r_{(m-1)}$ la matrice où on a enlevé la dernière ligne et la $(b+i)$-ème colonne et pour $t_{(m-1)}$ l'élément $\Diag(1, \dots, 1, \varpi_D^\zeta, \dots, \varpi_D^\zeta)$ avec une occurrence de $\varpi_D^\zeta$ de moins que pour $t$. De ce fait, on obtient $r_{(m-1)} = 1$, et donc $r$ est de la forme 
\[ \left( \begin{array}{ccccccc} 
\multicolumn{3}{c}{\multirow{2}{*}{$I_{b+i-1}$}} & c_1 & \multicolumn{3}{c}{\multirow{2}{*}{$0$}} \\ 
\multicolumn{3}{c}{} & \multirow{2}{*}{\vdots} & \multicolumn{3}{c}{} \\ 
\multicolumn{3}{c}{\multirow{2}{*}{$0$}} & & \multicolumn{3}{c}{\multirow{2}{*}{$I_{m-b-i}$}} \\ 
\multicolumn{3}{c}{} & c_{m-1} & \multicolumn{3}{c}{} \\ 
0 & \cdots & 0 & 1 & 0 & \cdots & 0 
\end{array} \right) .\]
Aussi, $r = w \left( \begin{array}{cc} I_b & 0 \\ M & I_{m-b}\end{array} \right)$ a une $(b+i)$-ème colonne avec un seul coefficient non nul, ce qui force $c_1 = \dots = c_m = 0$. Enfin, parce que $w$ est un élément de $W^t$, aucune de ses écritures réduites ne peut se terminer par un $s_\alpha$ pour $\alpha \in \Delta_t$: il s'ensuit $b+i = m$, et donc $r=1$. La récurrence est terminée. \hfill$\Box$\\

On veut maintenant énoncer un analogue pour $\GL(m,D)$ de la Proposition 5.1 de \cite{Her11b}. Pour tout $\lambda \in X_*(T)_V^-$, on rappelle que $T_\lambda$ désigne l'opérateur de $\Hh_{\Fpbar}(G,K, V)$ porté par la double classe $K \widetilde{\lambda}(\varpi) K$ et $\tau_\lambda$ celui de $\Hh_{\Fpbar}\big(A, A \cap K, V^{U \cap K}\big)$ porté par la classe $\widetilde{\lambda}(\varpi)(A \cap K)$. \\
L'énoncé voulu va traduire une différence importante entre le cas déployé et le cas qui ne l'est pas. En effet, pour une $K$-représentation irréductible $\chi_K = \chi \circ \overline{\det}$ qui se factorise par $\det_G$, on arrive, comme dans le cas déployé, à inverser complètement la transformée de Satake. Si $\chi_K$ ne se factorise pas par $\det_G$, nous ne savons plus inverser la transformée de Satake. Cependant le résultat de la proposition \ref{Prop5.1} (en ce sens, partiel) suivant nous permettra tout de même de mener à bout le travail. \\

Pour $V$ une représentation irréductible de $K$, la proposition \ref{irredGbar} nous assure que la droite $V^{U \cap K}$ possède un stabilisateur de la forme $(P_V \cap K) K(1)$ où $P_V$ est un parabolique standard de $G$ admettant $P_V = M_V N_V$ pour décomposition de Levi standard. Ecrivons $M_V = \prod_i \GL(m_i,D)$ avec $\sum_i m_i = m$: l'action de $A \cap K$ sur $V^{U \cap K}$ est alors donnée par un caractère $\chi_1^{\otimes m_1} \otimes \chi_2^{\otimes m_2} \otimes \cdots \otimes \chi_r^{\otimes m_r}$. Notons, pour tout $1 \leq i \leq r$, $d_i \geq 1$ le diviseur de $d$ vérifiant \[ \big\{ z \in \Z \ \big| \ \forall x \in k_D^\times, \, \chi_i \big( \overline{\varpi_D^z [x] \varpi_D^{-z}} \big) = \chi_i(x) \big\} = d_i \Z .\]
On dira que la droite $V^{U \cap K}$ est un caractère de $A \cap K$ qui se \textit{factorise par $\Nrm$} si $\chi_i$ se factorise par $\Nrm$ pour tout $i$; c'est équivalent à $d_i = 1$ pour tout $i$, et lorsque l'on a $M_V = G$, cela revient à demander que $V$ se factorise par $\det_G$. \\
On va en fait avoir besoin d'une condition plus faible, qui laisse de côté les blocs de taille $1$ dans $M_V$. On dira que la droite $V^{U \cap K}$ est un caractère de $A \cap K$ qui se \textit{factorise presque partout par $\Nrm$} si on a \[ \forall i = 1, \dots, r, \quad m_i > 1 \Rightarrow d_i = 1 .\]

\begin{prop} \label{Prop5.1}
Supposons la formule sur les valeurs de la tranformée de Satake (conjecture \ref{conjSatake}) satisfaite. Soit $V$ une représentation irréductible de $K$. 
\begin{itemize}
\item[(i)] Supposons que $V^{U \cap K}$ se factorise presque partout par $\Nrm$. Alors, pour tout élément $\mu$ de $X_*(T)_V^-$, on a\footnote{on rappelle que $\geq_{M_V}$ a \' et\' e d\' efini \` a la fin du paragraphe \ref{oneparamsbgp}} \[ \tau_\mu = \sum_{\lambda \in X_*(T)_V^-, \, \lambda \geq_{M_V} \mu} \Sss_G(T_\lambda) .\]
\item[(ii)] Soit $-\lambda$ un copoids fondamental minuscule associé à une racine de $\Delta_V$, située dans le bloc $\GL(m_i,D)$ pour un certain $1 \leq i \leq r$, choisi de mani\` ere \` a avoir
\[ \lambda(\varpi) = \Diag(1, \varpi^{a_2}, \dots, \varpi^{a_m}) .\] 
Supposons $d_i = d_{i+1} = \cdots = d_r$. Alors on a \[ \tau_{d_i \lambda} = \Sss_G(T_{d_i \lambda}) .\]
\item[(iii)] Supposons $m = 2$ et que $V$ n'est pas un caractère qui se factorise par $\Nrm$. Alors, pour tout élément $\mu$ de $X_*(T)_V^-$, on a \[ \tau_\mu = \Sss_G(T_\mu) .\]
\end{itemize}
\end{prop}
\textsf{Remarque :} 
L'hypothèses nécessaire pour $(ii)$ est en particulier satisfaite lorsque $V$ est un caractère de $K$. \\
\textsf{Début de preuve (commun à $(i), (ii), (iii)$) :} \\
Par le lemme \ref{deGàMv}, on est ramen\' e \` a calculer $\Sss_{M_V} \big( T_{d_i \lambda}^{M_V} \big)$. On \' ecrit $M_V = \prod_j G_j$ avec $G_j = \GL(m_j,D)$. Aussi $V^{N_V \cap K}$ est une repr\' esentation irr\' eductible de $M_V \cap K$ qui se d\' ecompose en $\bigotimes_j V_j$ o\` u chaque $V_j$ est une repr\' esentation irr\' eductible de dimension finie de $G_j \cap K$ (vu dans $M_V \cap K$) par la Proposition 3.4.2 de \cite{Bum97}. On a ainsi \[ \Hh_{\Fpbar} \big( M_V, M_V \cap K, V^{N_V \cap K} \big) = \bigotimes\nolimits_j \Hh_{\Fpbar} (G_j, G_j \cap K, V_j) ,\]
et on \' ecrit $T^{M_V}_{d_i \lambda} = t_1 \otimes \dots \otimes t_r$. Sur toutes les composantes sauf une, le calcul se ram\` ene au lemme \ref{Satakehomothetie}. Pour celle qui reste, on est ramen\' e \` a un $\GL(m_j,D)$: on suppose donc $M_V = G$ par la suite. \\
Parce que l'on a supposé $M_V=G$, il existe un caractère $\chi : k_D^\times \to \Fpbar^\times$ tel que l'on ait $V = \chi \circ \overline{\det}$. Soit $\rho$ un caractère de $\det_G(\widetilde{Z}_A(V))$; par le lemme \ref{prolongmtcaract}, on a une application \[\widetilde{\chi}_\rho : \coprod_{t \in \widetilde{Z}_A^-(V)} K t K \to \Fpbar^\times \] dont la restriction à $K$ coïncide avec $V$. Soient $\lambda, \mu \in X_*(T)_V^-$ et $t := \widetilde{\lambda}(\varpi)$, $t' := \widetilde{\mu}(\varpi)$. Pour tous $k_1, k_2 \in K$, on a \[ T_\lambda(k_1 t k_2) = \chi \circ \overline{\det}(k_1 k_2) \in \End_{\Fpbar}(V) = \Fpbar .\] De ce fait, on a 
\begin{equation} \label{Tlambdaagit} T_{\lambda}(k_1 t k_2) = \widetilde{\chi}_\rho(k_1 t k_2) \widetilde{\chi}_\rho(t)^{-1} .\end{equation} 
La décomposition d'Iwasawa (voir \cite{HenVig11}, Proposition 6.5) nous permet d'écrire 
\[ K t K = \coprod\nolimits_{i \in I} t_i u_i K ,\] 
où $I$ est un ensemble fini, $\lambda_i \in X_*(T)$, $t_i = \widetilde{\lambda}_i(\varpi)$ et $u_i \in U$ pour tout $i \in I$. De plus, l'intersection $t_i u_i K \cap t' U$ est non vide si et seulement si $t_i = t'$; en effet, écrivons $t_i u_i k = t' u$ pour un certain $k \in K$ et un $u \in U$. On a alors \[ k = u_i^{-1} t_i^{-1} t' u \in B \cap K \quad \textrm{et} \quad (t')^{-1} t_i \in U (B \cap K) U .\] D'où l'égalité $t' = t_i$; dans ce cas, $t_i u_i K \cap t' U$ vaut exactement $t' u_i (U \cap K)$. Il en résulte que le coefficient de $\tau_\mu$ intervenant dans $\Sss_G(T_\lambda)$ (voir (\ref{defSatake}) et (\ref{Tlambdaagit})) est 
\begin{equation} \label{coefSatake} \sum_{\lambda_i = \mu} T_\lambda (t' u_i) = \sum_{\lambda_i = \mu} \widetilde{\chi}_\rho(t' u_i) \widetilde{\chi}_\rho(t)^{-1} .\end{equation} 
\textsf{Preuve de $(i)$ :} \\
Dans le cas où l'hypothèse du $(i)$ est satisfaite, puisque l'on a $M_V = G$, le prolongement $\widetilde{\chi}_\rho$ se factorise par $\det_G$ et est donc multiplicatif (voir lemme \ref{multsiNrm}); en particulier, on a $\widetilde{\chi}_\rho(t' u_i) = \widetilde{\chi}_\rho(t')$. L'égalité (\ref{coefSatake}) décrivant le coefficient de $\tau_\mu$ dans $\Sss_G(T_\lambda)$ devient alors tout simplement: \[ \sum_{\lambda_i = \mu} T_\lambda (t' u_i) = \big| \{ i \in I \ | \ \lambda_i = \mu \}\big| \cdot \widetilde{\chi}_\rho(t') \widetilde{\chi}_\rho(t)^{-1} .\]
De plus, si $K t K \cap t' U$ est non vide, on choisit un élément $x$ de cet ensemble. Le calcul de  son déterminant donne alors $\det_G (t) = \det_G (t')$. En particulier, on a l'égalité $\widetilde{\chi}_\rho(t') = \widetilde{\chi}_\rho(t)$ et le coefficient de $\tau_\mu$ dans $\Sss_G(T_\lambda)$ se réduit à $| \{ i \in I \ | \ \lambda_i = \mu \} |$, qui est une quantité indépendante de $\chi$ ou de $\rho$. \\
On peut donc prendre $V = \id$, pour lequel le calcul est identique à la fin de la preuve de \cite{Her11b}; on le rappelle tout de même ici. A partir de maintenant et jusqu'à la fin de la preuve, $\mu$ et $\lambda$ seront des copoids dominants. On fait ce changement antidominant/dominant pour alléger l'écriture, et à la fin on appliquera l'identité prouvée (\ref{inversfinal}) à $w_0 \lambda$ et $w_0 \mu$, qui seront bien antidominants (comme voulu dans l'énoncé à prouver). La formule de Lusztig-Kato s'écrit (voir corollaire \ref{formuleLK}): \[ \ch \, V_\mu = \sum_{\lambda \leq \mu} q^{-d \langle \mu, \rho \rangle} P_{w_\lambda, w_\mu}(q^d) \, \overline{\Sss} \big( \id_{K \widetilde{\lambda}(\varpi) K}\big).\] Si on regarde dans cette égalité le coefficient devant $\tau_{\mu'}$ pour $\mu' \in X_*(T)$, on a par définition de $\overline{\Sss}$ \[ \dim \, V_\mu(\mu') = \sum_{\lambda \leq \mu} q^{-d \langle \mu, \rho \rangle} P_{w_\lambda, w_\mu}(q^d) \sum_{u \in U/ U \cap K} \delta_G^{1/2} \big( \widetilde{\mu}'(\varpi) \big) \id_{K \widetilde{\lambda}(\varpi) K} \big( \widetilde{\mu}'(\varpi) u \big) ,\] où $V_\mu(\mu')$ désigne le sous-espace de $V_\mu$ de poids $\mu'$ et $\delta_G$ est le caractère module relativement à $B$. On a par définition \[ \delta_G^{1/2} \big( \widetilde{\mu}'(\varpi) \big) = q^{d \langle w_0 \mu', \rho \rangle} ,\] de sorte que l'on peut réécrire: \[ q^{d \langle \mu - w_0 \mu' , \rho \rangle} \dim \, V_\mu(\mu') = \sum_{\lambda \leq \mu} P_{w_\lambda, w_\mu}(q^d) \sum_{u \in U / U \cap K} \id_{K \widetilde{\lambda}(\varpi) K} \big( \widetilde{\mu}'(\varpi) u \big) .\] Parce que l'on a (voir \cite{Jan87}, Proposition II.2.2) \[ \dim \, V_\mu(\mu') \neq 0 \Rightarrow w_0 \mu \leq \mu' \leq \mu ,\] $P_{w_\lambda,w_\mu} (q^d) \in 1+q^d \Z[q^d]$, et que $q^{ d \langle \mu - w_0 \mu', \rho \rangle}$ est divisible par $p$ si $\mu \neq w_0 \mu'$, en réduisant modulo $p$, on obtient \[ \dim \, V_\mu(w_0 \mu) = \sum_{\lambda \leq \mu} \sum_{u \in U / U \cap K} \id_{K \widetilde{\lambda}(\varpi) K} \big( \widetilde{w_0 \mu}(\varpi) u \big) ,\] \[ 0 = \sum_{\lambda \leq \mu} \sum_{u \in U / U \cap K} \id_{K \widetilde{\lambda}(\varpi) K} \big( \widetilde{\mu}'(\varpi) u \big) \quad \textrm{pour } \mu' \neq w_0 \mu .\] La dimension de $V_\mu(w_0 \mu)$ vaut $1$ par les II.1.19 et II.2.2 de \cite{Jan87}, et donc en renormalisant via $\mu \mapsto \widetilde{\mu}$, on a 
\begin{equation} \label{inversfinal} \tau_{w_0 \mu} = \sum_{\lambda \leq \mu} \Sss_G (T_{w_0 \lambda}) .\end{equation} 
Le résultat en découle. \\
\textsf{Preuve de $(ii)$ :} \\
Ici, on rappelle que $V$ ne se factorise plus nécessairement par $\det_G$. On cherche \` a d\' eterminer $KtK \cap t'U$. Lorsque $t'$ est tel que cette intersection est non vide, on en prend un \' el\' ement et le calcul de son d\' eterminant donne $\det_G(t) = \det_G(t')$. De plus, l'hypoth\` ese $d_i = d_{i+1} = \dots = d_r$ nous place exactement dans les conditions d'application du lemme \ref{miracle}. On en d\' eduit $t=t'$ et $KtK \cap tU = t(U \cap K)$. De ce fait, par (\ref{defSatake}) et (\ref{coefSatake}), on a le r\' esultat voulu. \\
\textsf{Preuve de $(iii)$ :} \\
Dans le cas $m=2$ et $V$ un caractère de $K$, pour tout caractère $\rho$ de $\det_G(\widetilde{Z}_A(V))$, le prolongement $\widetilde{\chi}_\rho$ défini est pseudo-multiplicatif par la proposition \ref{multm=2}. On utilise alors le corollaire \ref{utilchirho} pour dire que (\ref{coefSatake}) est égal à $| \{ i \in I \ | \ \lambda_i = \mu \} |$, et le coefficient de $\tau_\mu$ dans $\Sss_G(T_\lambda)$ est indépendant de $\chi$ ou de $\rho$ comme au \textit{(i)}. On va utiliser \textit{(i)} pour effectuer des calculs pour $V = \id$. Pour $y \leq z$ dans $\Z$, notons $t^y_z = \left( \begin{array}{cc} \varpi_D^y & 0 \\ 0 & \varpi_D^z \end{array} \right) \in A_\Lambda^-$, $\tau^y_z$ l'opérateur de $\Hh_{\Fpbar}(A,A \cap K)$ de support $t^y_z(A \cap K)$ et valant $1$ en $t^y_z$, et $T^y_z$ l'opérateur de $\Hh_{\Fpbar}(G,K)$ de support $K t^y_z K$ et valant $1$ en $t^y_z$. On note $\Delta = \{\alpha\}$ et $d_0 \geq 1$ l'entier divisant $d$ vérifiant \[ \{ z \in \Z \ | \ \forall x \in k_D^\times, \, \chi \big( \overline{\varpi_D^z [x] \varpi_D^{-z}} \big) = \chi(x) \} = d_0 \Z ,\] et $\mu : x \mapsto \left( \begin{array}{cc} x^a & 0 \\ 0 & x^b \end{array} \right)$ avec $a \leq b$ dans $d_0 \Z$. 
Parce que $V$ ne se factorise pas par $\Nrm$, on a $d_0>1$. Aussi, on peut supposer $a<b$, le cas $a=b$ étant facile puisque la classe double $K t^a_a K$ est égale à la classe simple $t^a_a K$. On écrit, pour $V = \id$: 
\begin{equation} \label{invers1} \tau_\mu = \tau^a_b = \Sss_G(T^a_b) + \Sss_G(T^{a+1}_{b-1}) + \dots + \Sss_G \big( T^{a+\lfloor \frac{b-a}{2} \rfloor}_{b- \lfloor \frac{b-a}{2} \rfloor}\big) ,\end{equation}
\begin{equation} \label{invers2} \tau_{\mu + \alpha^\vee} = \tau^{a+1}_{b+1} = \Sss_G(T^{a+1}_{b-1}) + \Sss_G (T^{a+2}_{b-2}) + \dots + \Sss_G \big( T^{a+\lfloor \frac{b-a}{2} \rfloor}_{b- \lfloor \frac{b-a}{2} \rfloor}\big) ;\end{equation}
en retranchant (\ref{invers2}) à (\ref{invers1}), on obtient 
\begin{equation} \label{invers3} \Sss_G(T_\mu) = \tau_\mu - \tau_{\mu + \alpha^\vee} .\end{equation}
Revenons à notre $V$ non trivial: parce que $d_0$ est différent de $1$, $\mu + \alpha^\vee$ n'appartient pas à $X_*(T)_V^-$, et l'identité (\ref{invers3}) se réduit à $\Sss_G(T_\mu) = \tau_\mu$. Le \textit{(iii)} est prouvé. $\phantom{,}$ \hfill$\Box$ \\

On remarque que la proposition \ref{Prop5.1}.(iii) peut se calculer à la main, sans faire appel à la machinerie précédente. Les calculs sont toutefois réminiscents de ceux du lemme \ref{multm=2}. \\

\textsf{Seconde preuve de la proposition \ref{Prop5.1}.(iii) :} \\ 
On suppose encore $M_V = G$, le cas contraire étant résolu par les lemmes \ref{deGàMv} et \ref{Satakehomothetie}. La représentation $V$ est alors un caractère $\chi \circ \overline{\det}$ et on note $d_0 \geq 1$ le diviseur de $d$ satisfaisant 
\[ \{ z \in \Z \ | \ \forall x \in k_D^\times, \, \chi \big( \overline{\varpi_D^z [x] \varpi_D^{-z}} \big) = \chi(x) \} = d_0 \Z .\]
Pour tout $a \leq b$ dans $d_0 \Z$, on note $t^a_b = \left( \begin{array}{cc} \varpi_D^a & 0 \\ 0 & \varpi_D^b\end{array} \right)$ et $T^a_b$ l'opérateur de $\Hh_{\Fpbar}(G,K,V)$ de support $K t^a_b K$ et valant l'identité en $t^a_b$. Soient $\alpha \leq \beta$ dans $d_0 \Z$. On veut calculer 
\[ \Sss_G(T^\alpha_\beta)(t^a_b) = \sum_{x \in D/\Oo_D} T^\alpha_\beta \Big( t^a_b \left( \begin{array}{cc} 1 & x \\ 0 & 1 \end{array} \right) \Big) \]
pour $a \leq b$ dans $d_0 \Z$. L'élément $t^a_b \left( \begin{array}{cc} 1 & x \\ 0 & 1 \end{array} \right)$ appartient à $K t^a_b K$ pour $x \in \Oo_D$. Supposons $x \in D \smallsetminus \Oo_D$ et écrivons $v_D(x) =-v$ avec $v>0$. On a alors 
\begin{equation} \label{sataked2} t^a_b \left( \begin{array}{cc} 1 & x \\ 0 & 1 \end{array} \right) = \left( \begin{array}{cc} \varpi_D^a & \varpi_D^a x \\ 0 & \varpi_D^b \end{array} \right) = \left( \begin{array}{cc} 1 & 0 \\ \varpi_D^b x^{-1} \varpi_D^{-a} & 1 \end{array} \right) \left( \begin{array}{cc} \varpi_D^{a-v} & 0 \\ 0 & \varpi_D^{b+v} \end{array} \right) \left( \begin{array}{cc} \varpi_D^v & \varpi_D^v x \\ - \varpi_D^{-v} x^{-1} & 0 \end{array} \right) ,\end{equation}
et $t^a_b \left( \begin{array}{cc} 1 & x \\ 0 & 1 \end{array} \right)$ appartient à $K t^\alpha_\beta K$ si et seulement si $(\alpha,\beta) = (a-v,b+v)$. En particulier, aucun terme pour $x \notin \Oo_D$ ne contribue à $\Sss_G(T^\alpha_\beta)(t^\alpha_\beta)$ et donc on a $\Sss_G(T^\alpha_\beta)(t^\alpha_\beta) = 1$. \\
De plus, si $(\alpha,\beta) = (a-v,b+v)$, alors $v$ est aussi dans $d_0 \Z$; dans ce cas-là, par (\ref{sataked2}), on a \[ T^\alpha_\beta \Big( t^a_b \left( \begin{array}{cc} 1 & x \\ 0 & 1 \end{array} \right) \Big) = \chi( \varpi_D^{-v} x^{-1} \varpi_D^v x) .\]
Ce terme ne dépend que du coefficient de plus petite valuation de $x$, et si on écrit \[ x = \varpi_D^{-v} x_0 + \varpi_D^{-v+1} x_1 + \dots + \varpi_D^{-1} x_{v-1} + \Oo_D \in D/\Oo_D ,\]
en regroupant les termes de valuation $v$, on obtient 
\begin{equation} \label{regroupT} \sum_{v_D(x) = -v} T^\alpha_\beta \Big( t^a_b \left( \begin{array}{cc} 1 & x \\ 0 & 1 \end{array} \right) \Big) = \sum_{x_0 \in k_D^\times} \sum_{x_1, \dots, x_{v-1} \in k_D} \chi( \varpi_D^{-v} x_0^{-1} \varpi_D^v x_0) .\end{equation}
Parce que $V$ ne se factorise pas par $\Nrm$, on a $d_0 \geq 2$, et donc $v \geq 2$. Ainsi (\ref{regroupT}) est une somme de $|k_D|^{v-1}$ termes identiques et vaut donc $0$ modulo $p$. \\
De ce fait, on a $\Sss_G(T^\alpha_\beta)(t^a_b) = 0$ pour $(a,b) \neq (\alpha,\beta)$. Ceci démontre bien que $\Sss_G(T^\alpha_\beta)$ est l'opérateur de $\Hh_{\Fpbar}(G,K,V)$ de support $t^\alpha_\beta(A \cap K)$ et valant l'identité en $t^\alpha_\beta$. \hfill$\Box$

\section{Changement de poids} \label{prgchgtpds}

Dans tout ce paragraphe, on supposera la formule de comparaison de tranformées de Satake (conjecture \ref{conjSatake}) satisfaite. \\

Soit $\alpha$ une racine simple. On note alors $-\lambda$ le copoids fondamental minuscule associé à $\alpha$ v\' erifiant que le premier coefficient diagonal de $\lambda(x)$ est \' egal \` a $1$ pour tout $x$ (ce qui est possible pour $\GL(m,D)$). Notamment, cela signifie que l'on a $\langle -\lambda, \alpha \rangle =1$ et $\langle -\lambda, \beta \rangle =0$ pour toute racine simple $\beta \neq \alpha$. On pose $t := \widetilde{\lambda}(\varpi) \in A_\Lambda^-$. \\ 
Par exemple, pour $\GL(4,D)$, on peut avoir la situation suivante: 
\[ \alpha^\vee : x \mapsto \left( \begin{array}{cccc} 1 & & & \\ & x & & \\ & & x^{-1} & \\ & & & 1 \end{array} \right) , \quad \lambda : x \mapsto \left( \begin{array}{cccc} 1 & & & \\ & 1 & & \\ & & x & \\ & & & x \end{array} \right) .\]

Soit $V$ une représentation irréductible de $K$. Par le lemme \ref{irredGbar}, $V$ est paramétré (par rapport à $\overline{B}$) par une paire $(\chi,\Delta_V)$ où $\chi$ est un caractère de $A \cap K$ et $\Delta_V$ est un sous-ensemble de $\Delta_\chi$. On suppose que $V$ vérifie l'hypothèse \ref{VtoV'} suivante pour toute la suite du paragraphe \ref{prgchgtpds}.

\begin{hypo} \label{VtoV'}
La racine $\alpha$ appartient à $\Delta_V$.
\end{hypo}

On remarque que $\Delta_V \smallsetminus \{ \alpha\}$ est alors inclus dans $\Delta_{\chi(t^{-1} \cdot t)}$ (voir (\ref{defDeltachi})). Par la proposition \ref{irredGbar}, on prend pour $V'$ la représentation irréductible de $K$ définie par \begin{equation} \label{defV'} \Par_{\overline{B}}(V') = \big( \chi( t^{-1} \cdot t), \Delta_V \smallsetminus \{ \alpha \}\big) .\end{equation} On notera $(\chi',\Delta_{V'})$ cette paire de paramètres. Comme on a $\Delta_\lambda = \Delta \smallsetminus \{ \alpha \}$, grâce à l'hypothèse \ref{VtoV'}, $V$ et $V'$ vérifient les conditions (\ref{condpdsconj}) et le lemme \ref{Heckepdsconj} nous assure que l'espace $\Hh_{\Fpbar}(G,K,V,V')$ possède un opérateur non nul de support $KtK$; on en fixe donc un en suivant (\ref{defHecke}) que l'on notera $\varphi_\lambda^+$. \\

Soit $P_{(\alpha)} = M_{(\alpha)} N_{(\alpha)}$ le Levi standard associé à $\{\alpha\} \subseteq \Delta$. Parce que $N_{(\alpha)}$ contient le radical unipotent $N_V$ du stabilisateur $P_V \cap K$ de $V^{U \cap K}$, $V^{N_{(\alpha)} \cap K}$ est un caractère de $M_{(\alpha)} \cap K$: on écrit \[ V^{N_{(\alpha)} \cap K} = \chi_1 \otimes \cdots \otimes \chi_{i-1} \otimes \chi_0 \circ \overline{\det} \otimes \chi_{i+2} \otimes \cdots \otimes \chi_m ,\] où les $\chi_j$ sont des caractères de $\Oo_D^\times$ pour $j \neq 0$ et $\chi_0$ est un caractère de $k_D^\times$. Soit $d_0 \geq 1$ le diviseur de $d$ qui vérifie \begin{equation} \label{defd0} \big\{ z \in \Z \ \big| \ \forall x \in k_D^\times, \, \chi_0 \big( \overline{\varpi_D^z [x] \varpi_D^{-z}} \big) = \chi_0(x) \big\} = d_0 \Z .\end{equation}
Lorsque l'on a $d_0 > 1$, on fait alors l'hypothèse suivante.

\begin{hypo} \label{V'toV}
L'espace $\Hh_{\Fpbar}(G,K,V',V)$ possède un opérateur non nul de support $K t^{d_0-1} K$ si $d_0$ est strictement supérieur à $1$.
\end{hypo}
Si, pour $j \geq i+2$, on note $d_j \geq 1$ l'entier tel que l'on ait \[ \big\{ z \in \Z \ \big| \ \forall x \in k_D^\times, \, \chi_j \big( \overline{\varpi_D^z [x] \varpi_D^{-z}} \big) = \chi_j(x) \big\} = d_j \Z ,\] alors combinatoirement l'hypothèse \ref{V'toV} se traduit simplement par $d_j \mid d_0$ pour tout $j \geq i+2$. En effet, dans ce cas-là, on a $\chi(t^{-d_0} \cdot t^{d_0}) = \chi$, et l'équivalence entre les deux formulations suit du lemme \ref{Heckepdsconj}. \\

Dans le cadre de l'hypothèse \ref{V'toV}, lorsque l'on a $d_0 > 1$, on note $\varphi_\lambda^\sim$ le tel opérateur tel que $\varphi_\lambda^\sim(t^{d_0-1}) \circ \varphi_\lambda^+(t)$ induise l'identité sur $V^{N_t^- \cap K}$, ce qui est loisible par la proposition \ref{Prop6.3} et (\ref{defHecke}). \\ Lorsque l'on a $d_0 = 1$, parce que $t^d$ est central dans $A$, on sait grâce au lemme \ref{Heckepdsconj} que l'espace d'entrelacements $\Hh_{\Fpbar}(G,K,V',V)$ possède un opérateur non nul de support\footnote{ce $2d-1$ a été choisi pour que l'exposant soit strictement supérieur à $0$ même dans le cas déployé $d=1$; dans ce cas-là, ce choix est d'ailleurs compatible avec celui dans \cite{Her11b}} $K t^{2d-1} K$; on fixe et on note $\varphi_\lambda^-$ celui qui fait de $\varphi_\lambda^-(t^{2d-1}) \circ \varphi_\lambda^+(t)$ l'identité sur $V^{N_t^- \cap K}$. Enfin, lorsque l'on a $d_0 > 1$ (en particulier $d>1$), on note de plus $\varphi_\lambda^=$ l'opérateur de $\Hh_{\Fpbar}(G,K,V',V)$ tel que $\varphi_\lambda^=(t^{d-1}) \circ \varphi_\lambda^+(t)$ induise l'identité sur $V^{N_t^- \cap K}$. \\
Pour $? \in \{ \sim, -, = \}$ (et lorsque l'hypothèse \ref{V'toV} est satisfaite si besoin), $\varphi_\lambda^?$ est un opérateur de Hecke de $\Hh_{\Fpbar}(G,K,V)$ et on cherche à évaluer sa transformée de Satake. Le cas de $\varphi_\lambda^=$ ne nous intéresse que pour $m=3$; il sera donc repoussé plus loin dans ce paragraphe.

\begin{lem} \label{Prop6.8}
Soient $V, \alpha, \lambda, t$ comme précédemment, et supposons l'hypothèse \ref{VtoV'} satisfaite. Soit $V'$ comme en (\ref{defV'}).
\begin{itemize}
\item[(i)] Supposons $d_0 > 1$ et l'hypothèse \ref{V'toV} satisfaite. L'opérateur $\varphi_\lambda^\sim * \varphi_\lambda^+$ de $\Hh_{\Fpbar}(G,K,V)$ est $T_{d_0 \lambda}$. De plus, on a $\Sss_G(\varphi_\lambda^\sim * \varphi_\lambda^+) = \tau_{d_0 \lambda}$.
\item[(ii)] Supposons $d_0 = 1$. L'opérateur $\varphi_\lambda^- * \varphi_\lambda^+$ de $\Hh_{\Fpbar}(G,K,V)$ est $T_{2d \lambda}$. De plus, on a $\Sss_G(\varphi_\lambda^- * \varphi_\lambda^+) = \tau_{2d \lambda} - \tau_{2d \lambda + \alpha^\vee}$.
\end{itemize}
\end{lem}
\textsf{Preuve :} \\
Parce que le début de preuve est identique pour $(i)$ et $(ii)$, posons \[ (\delta_0,?) = \begin{cases} (d_0,\sim) & \textrm{si } d_0 > 1 , \\ (2d,-) & \textrm{si } d_0=1 .\end{cases} \] On cherche d'abord à montrer que $\varphi_\lambda^? * \varphi_\lambda^+$ a pour support $K t^{\delta_0} K$. Soient $\lambda' \in X_*(T)^-$ et $t' = \widetilde{\lambda}'(\varpi)$. On va donc montrer que si $\varphi_\lambda^? * \varphi_\lambda^+(t')$ est non nul alors on a $t' = t^{\delta_0}$. On a la bijection d'ensembles (voir le lemme \ref{Prop3.8}) 
\begin{equation} \label{KtK} \begin{array}{ccccc} K \backslash K t K & \mathop{\leftarrow}\limits^{\sim} & (K \cap t^{-1} K t) \backslash K & \simeq & \overline{P}_{\lambda} \backslash \overline{G} \\ K t k & \mapsfrom & (K \cap t^{-1} K t) k & & \end{array} .\end{equation} 
Par décomposition de Bruhat, c'est aussi \[ \coprod_{w \in W_\lambda \backslash W} \overline{P}_{\lambda} \backslash \overline{P}_{\lambda} w \overline{U}^- .\] 
En utilisant la flèche (\ref{KtK}), il en résulte la décomposition \[ K t K = \coprod_{w \in W_\lambda \backslash W} K t w (U^- \cap K) .\] La formule (\ref{GLmdefconvol}) se réécrit \[ \varphi_\lambda^? * \varphi_\lambda^+(t') = \sum_{x \in K \backslash KtK} \varphi_\lambda^?(t' x^{-1}) \varphi_\lambda^+(x) .\] Tout terme de $\varphi_\lambda^? * \varphi_\lambda^+(t')$ s'écrit alors $\varphi_\lambda^?(t' u^{-1} w^{-1} t^{-1}) \varphi_\lambda^+(twu)$ pour un $u \in U^- \cap K$ et un $w \in W$. Comme $t' u^{-1} (t')^{-1}$ est un élément de $U^- \cap K$, ce terme est non nul si et seulement si $\varphi_\lambda^?(t' w^{-1} t^{-1} w) \varphi_\lambda^+(w^{-1} tw)$ est non nul. Cela implique $\varphi_\lambda^?(t' w^{-1} t^{-1}w) \neq 0$ et donc l'existence d'un $w' \in W$ avec $\lambda' -w^{-1} \lambda = (\delta_0-1) w' \lambda$. Mais dans ce cas-là, on a: \[ \varphi_\lambda^?(t^{\delta_0-1}) (w')^{-1} w^{-1} \varphi_\lambda^+ (t) \neq 0 .\] Cela implique \[ p_{N_\lambda} \big( (w')^{-1} w^{-1} (V')^{N_\lambda^-} \big) \neq 0 ,\] où $p_{N_\lambda} : V' \twoheadrightarrow (V')_{N_\lambda}$ la projection canonique. Par \cite{HenVig11b}, Corollary 3.19, on a $(w')^{-1} w^{-1} \in W_\lambda W_{\Delta_{V'}} W_\lambda \subseteq W_\lambda$ et donc l'identité $\lambda' = \delta_0 w' \lambda$. Parce que $\lambda$ et $\lambda'$ sont tous deux antidominants, cela force $w' \lambda = \lambda$, et ainsi $\lambda' = \delta_0 \lambda$. On a bien l'affirmation voulue. Prenons donc $t' = t^{\delta_0}$ et regardons de nouveau quels termes contribuent à $\varphi_\lambda^? * \varphi_\lambda^+(t')$. Par l'analyse précédente, on a $(w')^{-1} w^{-1} \in W_\lambda$ et $w' \in W_\lambda$, donc $w$ appartient à $W_\lambda$ si le terme $\varphi_\lambda^?(t' u^{-1} w^{-1} t^{-1}) \varphi_\lambda^+(twu)$ a une contribution non nulle. Mais alors $Kt$ est la seule classe de $K \backslash KtK$ à intervenir. Il s'ensuit que $\varphi_\lambda^? * \varphi_\lambda^+$ est exactement l'opérateur $T_{\delta_0 \lambda}$ de $\Hh_{\Fpbar}(G,K,V)$. \\ 
Pour déterminer explicitement $\Sss_G(\varphi_\lambda^? * \varphi_\lambda^+)$, on commence par utiliser le lemme \ref{deGàMv} pour écrire \[ \Sss_G(\varphi_\lambda^? * \varphi_\lambda^+) = \Sss_{M_V} \big( T^{M_V}_{\delta_0 \lambda} \big) .\] Décomposons $M_V$ en $M_1 \times \cdots \times M_s$ où chaque $M_i$ est isomorphe à un $\GL(m_i,D)$. La représentation $V^{N_V \cap K}$ se décompose alors aussi en $V^{N_V \cap K} = V_1 \otimes \cdots \otimes V_s$ où chaque $V_i$ est un caractère de $M_i \cap K$. On écrit aussi \[ \delta_0 t = (t_1 , \dots, t_s) \in M_1 \times \dots \times M_s .\] On a alors 
\[ \Hh_{\Fpbar}\big( M_V, M_V \cap K, V^{N_V \cap K} \big) = \bigotimes\nolimits_i \Hh_{\Fpbar}(M_i,M_i \cap K, V_i) \] et on cherche à calculer 
\[ \Sss_{M_V} \big( T_{\delta_0 \lambda}^{M_V} \big) = \Sss_{M_1} \big( T_{t_1}^{M_1} \big) \otimes \Sss_{M_2} \big( T_{t_2}^{M_2} \big) \otimes \dots \otimes \Sss_{M_s} \big( T_{t_s}^{M_s} \big) .\]
Parce que l'on a $\Delta_\lambda = \Delta \smallsetminus \{ \alpha \}$, il existe un unique entier $j \in [1,s]$ tel que $t_j$ ne soit pas une homothétie. Par le lemme \ref{Satakehomothetie}, pour tout $i \neq j$, on a donc $\Sss_{M_i} \big( T_{t_i}^{M_i} \big) = \tau_{t_i}^{M_i}$. \\
Pour déterminer $\Sss_{M_j} \big( T_{t_j}^{M_j} \big)$, on va utiliser la proposition \ref{Prop5.1}. Si on a $d_0>1$, comme $V_j$ est un caractère de $M_j \cap K$, on a directement $\Sss_{M_j} \big( T_{t_j}^{M_j} \big) = \tau_{t_j}^{M_j}$ par la proposition \ref{Prop5.1}.(ii). Dans le cas $d_0 = 1$, on utilise la proposition \ref{Prop5.1}.(i) et le lemme \ref{SublemTo6.8} pour conclure. \hfill$\Box$

\begin{lem} \label{SublemTo6.8}
Gardons les hypothèses du lemme \ref{Prop6.8}.(ii). Soit $\mu \in X_*(T)^-$. Alors on a $\mu \geq_{M_V} 2 d\lambda$ si et seulement si on a $\mu = 2 d\lambda$ ou $\mu \geq_{M_V} 2d\lambda + \alpha^\vee$.
\end{lem}
\textsf{Remarque :} 
Ceci se trouve déjà dans le Claim de la preuve de la Proposition 6.7 de \cite{Her11b} dans le cas déployé; l'auteur ne le reproduit ici que par souci de lisibilité. \\
\textsf{Preuve :} \\
Parce que l'on a $\alpha \in \Delta_V$, il n'est pas dommageable de passer au Levi et de supposer $M_V = G$; c'est donc ce que l'on fera par la suite. L'implication $(\Leftarrow)$ est immédiate. Prouvons sa réciproque. Soit $\alpha_0 \in X^*(T)$ la somme des plus grandes racines de chaque composante irréductible du système de racines. On écrit $\alpha_0 = \sum_{\beta \in \Delta} n_\beta \beta$ et on a $n_\beta \geq 1$ pour tout $\beta \in \Delta$. Parce que $-\lambda$ est minuscule et que l'on a $n_\alpha \geq 1$, on en déduit $\langle \lambda, \alpha_0 \rangle = - n_\alpha = -1$. Comme $\mu$ est antidominant et vérifie $\mu -2d\lambda \geq 0$, on a 
\begin{equation} \label{mulambda} \langle \mu, \alpha \rangle \geq \langle \mu, \alpha_0 \rangle \geq \langle 2 d \lambda , \alpha_0 \rangle = - 2d .\end{equation}
Ecrivons $\mu - 2d\lambda = \sum_i \beta_i^\vee$ où les $\beta_i$ sont des racines simples, éventuellement avec répétition. Deux cas distincts se présentent alors à nous, selon si (\ref{mulambda}) est ou non une égalité (entre les deux termes extrémaux). \\ Si c'est une égalité, alors on a $\langle \mu, \alpha \rangle = \langle \mu, \alpha_0 \rangle$, et donc $\langle \mu, \beta \rangle = 0$ pour tout $\beta \in \Delta \smallsetminus \{\alpha\}$. De ce fait, on a $\langle \mu - 2d\lambda, \beta \rangle = 0$ pour toute racine simple $\beta$, et donc $\mu = 2d\lambda$. \\
Si (\ref{mulambda}) est une inégalité stricte, alors on a $\sum_i \langle \beta_i^\vee, \alpha \rangle > 0$. Comme on a $\langle \beta^\vee, \alpha \rangle \leq 0$ pour toute racine simple $\beta \neq \alpha$, il existe un indice $i$ tel que $\beta_i$ soit $\alpha$. Il en découle $\mu \geq 2d\lambda + \alpha^\vee$. Le lemme est prouvé. \hfill$\Box$ \\

Soient $V$ et $V'$ deux représentations irréductibles de $K$ qui sont conjuguées. Alors, par le lemme \ref{droiteconj} et la proposition \ref{HVsatake}, les algèbres $\Hh_{\Fpbar}(G,K,V)$ et $\Hh_{\Fpbar}(G,K,V')$ sont canoniquement isomorphes via l'isomorphisme de Satake:
\[ \Hh_{\Fpbar}(G,K,V) \mathop{\xrightarrow{\sim}}\limits_{\Sss_G} \Hh_{\Fpbar}^- \big( A, A \cap K, V^{U \cap K} \big) \xrightarrow{\sim} \Hh_{\Fpbar}^- \big( A, A \cap K, (V')^{U \cap K} \big) \mathop{\xrightarrow{\sim}}\limits_{\Sss_G^{-1}} \Hh_{\Fpbar}(G,K,V') .\] 
Tout caractère de $\Fpbar$-algèbres $\chi : \Hh_{\Fpbar}(G,K,V) \to \Fpbar$ induit alors un caractère \[ \chi' : \Hh_{\Fpbar}(G,K,V') \xrightarrow{\sim} \Hh_{\Fpbar}(G,K,V) \xrightarrow{\chi} \Fpbar .\]

\begin{prop} \label{Cor6.11}
Soit $\chi : \Hh_{\Fpbar}(G,K,V) \to \Fpbar$ un caractère de $\Fpbar$-algèbres de paramètre de Hecke-Satake $(M,\chi_M)$. On suppose $M \neq G$ et on prend $\alpha \in \Delta \smallsetminus \Delta_M$. On suppose l'hypothèse \ref{VtoV'} et on définit $V'$ comme en (\ref{defV'}), $d_0$ comme en (\ref{defd0}). On suppose de plus l'hypothèse \ref{V'toV}. Supposons enfin l'une des deux conditions suivantes satisfaite:
\begin{itemize}
\item[(a)] $\widetilde{\alpha}^\vee(\varpi) \notin \widetilde{Z}_M(V)$;
\item[(b)] $\chi_M(\widetilde{\alpha}^\vee(\varpi)) \neq 1$.
\end{itemize}   
Alors $\varphi_\lambda^+$ induit l'isomorphisme \[ \ind_K^G \, V \otimes_{\Hh,\chi} \Fpbar \xrightarrow{\sim} \ind_K^G \, V' \otimes_{\Hh',\chi'} \Fpbar \] de $G$-représentations.
\end{prop}
\textsf{Remarque :} 
L'hypothèse \ref{V'toV} est en particulier satisfaite dans le cas déployé $d=1$, de sorte que l'on recouvre l'énoncé de changement de poids de \cite{Her11b}. \\
\textsf{Remarque :} 
Combinatoirement, l'hypothèse \og non (a) \fg\ signifie que $\alpha$ est coincée entre deux blocs de $M$ de taille $1$ et que l'on a $d_0=1$. \\
\textsf{Preuve :} \\
On va montrer que $\varphi_\lambda^+$ réalise cet isomorphisme, d'inverse $\varphi_\lambda^?$ avec \[ ? = \begin{cases} \sim & \textrm{si } \widetilde{\alpha}^\vee(\varpi) \notin \widetilde{Z}_M(V); \\ - & \textrm{si } \widetilde{\alpha}^\vee(\varpi) \in \widetilde{Z}_M(V) \textrm{ et } \chi_M(\widetilde{\alpha}^\vee(\varpi)) \neq 1 .\end{cases}\] 
On remarque ensuite que $\varphi_\lambda^+$ et $\varphi_\lambda^?$ induisent apr\` es tensorisation par $\chi$ des op\' erateurs 
\[\xymatrix{ \ind_K^G \, V \otimes_{\Hh,\chi} \Fpbar \ar@<2pt>[r]^{\varphi_\lambda^? \otimes 1} & \ind_K^G \, V' \otimes_{\Hh', \chi'} \Fpbar \ar@<2pt>[l]^{\varphi_\lambda^+ \otimes 1} .}\]
En effet, on peut voir $\chi$ et $\chi'$ comme des caract\` eres de \[ \Hh_{\Fpbar}^- \big( A, A \cap K, V^{U \cap K} \big) \xrightarrow{\sim} \Hh_{\Fpbar}^- \big( A, A \cap K, (V')^{U \cap K} \big) \] gr\^ ace \` a la transform\' ee de Satake. Et une fois identifi\' ees $\Hh_{\Fpbar}^- \big( A, A \cap K, V^{U \cap K} \big)$ et $\Hh_{\Fpbar}^- \big( A, A \cap K, (V')^{U \cap K} \big)$, les structures de module \` a droite et \` a gauche sont la m\^ eme sur $\Hh_{\Fpbar} \big( A, A \cap K, V^{U \cap K}, (V')^{U \cap K} \big)$, ce qui justifie que $\varphi_\lambda^+$ et $\varphi_\lambda^?$ induisent bien $\varphi_\lambda^+ \otimes 1$ et $\varphi_\lambda^? \otimes 1$ comme voulu. \\
Pour $f \in \ind_K^G \, V$, l'image de $f \otimes 1$ par la composée \[ \ind_K^G \, V \otimes_{\Hh, \chi} \Fpbar \xrightarrow{\varphi_\lambda^+ \otimes 1} \ind_K^G \, V' \otimes_{\Hh',\chi'} \Fpbar \xrightarrow{\varphi_\lambda^? \otimes 1} \ind_K^G \, V \otimes_{\Hh, \chi} \Fpbar \] est \[ ( \varphi_\lambda^? * \varphi_\lambda^+).f \otimes 1 = f \otimes \chi \big( \varphi_\lambda^? * \varphi_\lambda^+ \big) ;\] et de même pour la composée associée à $\varphi_\lambda^+ * \varphi_\lambda^?$. Il s'agit donc de montrer que $\chi(\varphi_\lambda^+ * \varphi_\lambda^?)$ et $\chi(\varphi_\lambda^? * \varphi_\lambda^+)$ sont tous deux non nuls. \\
Supposons dans un premier temps que l'on a $\chi(\varphi_\lambda^? * \varphi_\lambda^+) \neq 0$. Alors $(\varphi_\lambda^? \otimes 1) \circ (\varphi_\lambda^+ \otimes 1)$ est un isomorphisme, et en particulier l'image de $\varphi_\lambda^+ \otimes 1$ est non réduite à $0$. Prenons un élément non nul $(\varphi_\lambda^+ \otimes 1). f \in \ind_K^G \, V' \otimes_{\Hh', \chi'} \Fpbar$; son image par \[ \ind_K^G \, V' \otimes_{\Hh', \chi'} \Fpbar \xrightarrow{\varphi_\lambda^? \otimes 1} \ind_K^G \, V \otimes_{\Hh,\chi} \Fpbar \xrightarrow{\varphi_\lambda^+ \otimes 1} \ind_K^G \, V' \otimes_{\Hh', \chi'} \Fpbar \] est alors \begin{align*} \varphi_\lambda^+(f) \otimes \chi(\varphi_\lambda^? * \varphi_\lambda^+) & = \varphi_\lambda^+ * (\varphi_\lambda^? * \varphi_\lambda^+).f \otimes 1 \\ & = (\varphi_\lambda^+ * \varphi_\lambda^?) * \varphi_\lambda^+ . f \otimes 1 = \varphi_\lambda^+(f) \otimes \chi(\varphi_\lambda^+ * \varphi_\lambda^?) .\end{align*} On en déduit alors \[ \chi(\varphi_\lambda^+ * \varphi_\lambda^?) = \chi(\varphi_\lambda^? * \varphi_\lambda^+) \neq 0 \] et $\varphi_\lambda^+ * \varphi_\lambda^?$ est aussi un isomorphisme. Cela nous assure l'énoncé voulu. \\
Il reste donc à montrer $\chi(\varphi_\lambda^? * \varphi_\lambda^+) \neq 0$. Deux cas se présentent à nous: d'abord si on a $\widetilde{\alpha}^\vee(\varpi) \notin \widetilde{Z}_M(V)$, c'est que ou bien les blocs de $M$ autour de $\alpha$ ne sont pas tous deux de taille $1$ ou bien $d_0$ est strictement supérieur à $1$. Lorsque $d_0 > 1$, par le lemme \ref{Prop6.8}.(i), $\Sss_G(\varphi_\lambda^\sim * \varphi_\lambda^+)$ est $\tau_{d_0 \lambda}$. De plus, l'hypothèse \ref{V'toV} nous assure que $d_0 \lambda$ est dans $X_*(T)_V^-$. Par le lemme \ref{Lem4.3}, on a alors \[ \chi( \varphi_\lambda^\sim * \varphi_\lambda^+) = \chi^{(A)} (\tau_{d_0 \lambda}) \neq 0 .\] 
Si on a $d_0 = 1$, mais que l'un des blocs de $M$ autour de $\alpha$ est de taille $>1$, on utilise le lemme \ref{Prop6.8}.(ii): $\Sss_G(\varphi_\lambda^- * \varphi_\lambda^+) = \tau_{2d \lambda} - \tau_{2d \lambda + \alpha^\vee}$. Par le lemme \ref{Lem4.3}, on a de nouveau \[ \chi(\varphi_\lambda^- * \varphi_\lambda^+) = \chi^{(A)}(\tau_{2d \lambda}) - \chi^{(A)}(\tau_{2d \lambda + \alpha^\vee}) = \chi^{(A)}(\tau_{2d \lambda}) \neq 0 .\]
Plaçons-nous maintenant dans le second cas de figure: $\widetilde{\alpha}^\vee(\varpi) \in \widetilde{Z}_M(V)$ et $\chi_M(\widetilde{\alpha}^\vee(\varpi)) \neq 1$. Dans ce cas-là, on a $d_0=1$ et $\Sss_G(\varphi_\lambda^- *\varphi_\lambda^+) = \tau_{2d \lambda} - \tau_{2d \lambda+ \alpha^\vee}$ par le lemme \ref{Prop6.8}.(ii). On utilise alors la proposition \ref{Prop4.1} pour affirmer 
\begin{align*} \chi(\varphi_\lambda^- * \varphi_\lambda^+) & = \chi_M \big( \widetilde{\lambda}(\varpi)^{2d} \big)^{-1} - \chi_M \big( \widetilde{\lambda}(\varpi)^{2d} \widetilde{\alpha}^\vee(\varpi) \big)^{-1} \\ \nonumber & = \chi_M \big( \widetilde{\lambda}(\varpi)^{2d}\big)^{-1} \big( 1- \chi_M \big( \widetilde{\alpha}^\vee(\varpi)\big)^{-1} \big) \neq 0 .\end{align*} Le résultat est donc prouvé. \hfill$\Box$ \\

On va appliquer ce qui précède à $B^-$ au lieu de $B$, de sorte que $V^*$ et le paramètre HS-dual de $\chi$ vont apparaître. \\ En effet, notons $(\chi, - \Delta_V^-)$ la paire de paramètres $\Par_{\overline{B}^-}(V)$: par la proposition \ref{paramV*}, on a alors $\Par_{\overline{B}}(V^*) = (\chi^{-1}, \Delta_V^-)$. On suppose l'hypothèse \ref{VtoV'} pour $B^-$, c'est-à-dire que $\alpha$ est un élément de $\Delta_V^-$. On prend $\lambda$ le copoids fondamental minuscule pour $\alpha$ tel que le dernier coefficient diagonal de $\lambda(x)$ est $1$ pour tout $x$; en particulier, $t := \widetilde{\lambda}(\varpi) \in A_\Lambda$ est maintenant dominant (pour $B$). On définit dès lors $V'$ comme la $K$-représentation irréductible satisfaisant 
\begin{equation} \label{defV'2} \Par_{\overline{B}^-}(V') = \big( \chi(t^{-1} \cdot t), -\Delta_V^- \smallsetminus \{ -\alpha \} \big) .\end{equation}
De cette manière, le diagramme (\ref{diagHSinduit}) nous assure que c'est le paramètre HS-dual de $\chi$ qui intervient. Aussi, on fait remarquer \[ \widetilde{Z}_M(V^*) = \widetilde{Z}_M \big( (V^*)^{U \cap K} \big) = \widetilde{Z}_M(V_{U \cap K}) .\] Tout cela se réécrit de la façon suivante. 

\begin{cor} \label{chgtpdsdual}
Soit $\chi : \Hh_{\Fpbar}(G,K,V) \to \Fpbar$ un caractère de $\Fpbar$-algèbres de paramètre HS-dual $(M,\chi_M)$. On suppose $M \neq G$ et on prend $\alpha \in \Delta \smallsetminus \Delta_M$. On suppose l'hypothèse \ref{VtoV'} pour $B^-$ et on définit $V'$ comme en (\ref{defV'2}), $d_0$ associé pour $V_{N_{(\alpha)} \cap K}$. On suppose de plus l'hypothèse \ref{V'toV}. Supposons enfin l'une des deux conditions suivantes satisfaite: 
\begin{itemize}
\item[(a)] $\widetilde{\alpha}^\vee(\varpi) \notin \widetilde{Z}_M(V^*)$;
\item[(b)] $\chi_M(\widetilde{\alpha}^\vee(\varpi)) \neq 1$.
\end{itemize}   
Alors $\varphi_\lambda^+$ induit l'isomorphisme \[ \ind_K^G \, V \otimes_{\Hh,\chi} \Fpbar \xrightarrow{\sim} \ind_K^G \, V' \otimes_{\Hh',\chi'} \Fpbar \] de $G$-représentations.
\end{cor}
\textsf{Preuve :} \\
Comme pour la proposition \ref{Cor6.11}, on cherche à évaluer $\chi(\varphi_\lambda^? * \varphi_\lambda^+)$. Cette fois-ci cependant, on écrit cela $\big( \chi \circ ({}' \Sss_G)^{-1} \big) \big( ({}' \Sss_G)(\varphi_\lambda^? * \varphi_\lambda^+) \big)$, où on a inversé \[ {}' \Sss_G : \Hh_{\Fpbar}(G,K,V) \hookrightarrow \Hh_{\Fpbar}(A, A \cap K, V_{U \cap K}) \] sur son image. Enfin, on calcule ${}' \Sss_G (\varphi_\lambda^? * \varphi_\lambda^+)$ en se ramenant à $\Sss_G^{(V^*)}$ via la proposition \ref{SatakeDual} : \[ {}' \Sss_G^{(V)} (\varphi_\lambda^? * \varphi_\lambda^+) = \iota_A \big( \Sss_G^{(V^*)} \big( \iota_G^{-1} (\varphi_\lambda^? * \varphi_\lambda^+) \big) \big). \] Le reste est identique à la preuve de la proposition \ref{Cor6.11}. \hfill$\Box$\\

Cet énoncé de changement de poids va nous permettre d'effectuer le travail nécessaire pour $\GL(m,D)$ avec certaines hypothèses sur le degré de $D$ sur $F$. Pour $m \leq 3$, ces hypothèses ne sont pas nécessaires mais il faut faire appel à un autre énoncé de changement de poids, dans le cas où $V$ n'est pas \og trop irrégulier \fg. C'est ce que l'on propose de présenter maintenant. 

\begin{hypo} \label{Vnontropirreg}
La $K$-représentation irréductible $V$ n'est pas un caractère.
\end{hypo}
On va utiliser la forme équivalente \og le stabilisateur $(P_V \cap K) K(1)$ de la droite $V^{U \cap K}$ n'est pas $K$ \fg, ou encore la condition équivalente $\Delta_V \neq \Delta$.

\begin{lem} \label{chgtpdsm=3}
Supposons $m = 3$ et $d>1$. Soient $V, \alpha, \lambda, t$ comme pour le lemme \ref{Prop6.8} et supposons les hypothèses \ref{VtoV'} et \ref{Vnontropirreg}. Soit $V'$ comme en (\ref{defV'}). L'opérateur $\varphi_\lambda^= * \varphi_\lambda^+$ de $\Hh_{\Fpbar}(G,K,V)$ est $T_{d \lambda}$. De plus, on a \[ \Sss_G(\varphi_\lambda^= * \varphi_\lambda^+) = \begin{cases} \tau_{d \lambda} & \textrm{si } d_0 >1 ; \\ \tau_{d \lambda} - \tau_{d \lambda + \alpha^\vee} & \textrm{si } d_0=1 .\end{cases} \]
\end{lem}
\textsf{Preuve :} \\
Le début de la démonstration du lemme \ref{Prop6.8} est encore valable et $\varphi_\lambda^= * \varphi_\lambda^+$ est simplement l'opérateur $T_{d \lambda}$. A cause des hypothèses \ref{VtoV'} et \ref{Vnontropirreg}, $\Delta_V$ est exactement le singleton $\{ \alpha \}$; si $M_V$ désigne le Levi standard de $P_V$, le lemme \ref{deGàMv} nous donne $\Sss_G(T_{d\lambda}) = \Sss_{M_V} \big( T_{d \lambda}^{M_V} \big)$. Ecrivons $M_V = M_1 \times M_2$ avec $\Delta_{M_2} = \{ \alpha \}$. La représentation $V^{N_V \cap K}$ est un caractère de $M_V \cap K$ et on peut la décomposer en $V^{N_V \cap K}  = V_1 \otimes V_2$ où $V_1 $ et $V_2$ sont respectivement des caractères de $M_1 \cap K$ et de $M_2 \cap K$. On a alors 
\[ \Hh_{\Fpbar} \big( M_V, M_V \cap K, V^{N_V \cap K} \big) = \Hh_{\Fpbar}(M_1, M_1 \cap K, V_1) \otimes \Hh_{\Fpbar}(M_2 , M_2 \cap K, V_2) ,\] et il s'agit de calculer 
\[ \Sss_{M_V} \big( T_{d \lambda}^{M_V} \big) = \Sss_{M_1} \big( T_{t_1}^{M_1} \big) \otimes \Sss_{M_2} \big( T_{t_2}^{M_2} \big) \] où $t_2$ est l'élément $\left( \begin{array}{cc} 1 & 0 \\ 0 & \varpi \end{array} \right)$ et $t_1$ est $(1)$ ou $(\varpi)$ selon les positions de $M_1$ et de $M_2$. Dans tous les cas, on a $\Sss_{M_1} \big( T_{t_1}^{M_1} \big) = \tau_{t_1}^{M_1}$ (car $\Sss_{M_1}$ est l'identité) et le calcul de $\Sss_{M_2} \big( T_{t_2}^{M_2} \big)$ fait appel à la proposition \ref{Prop5.1}.(iii) si $d_0>1$, et à la proposition \ref{Prop5.1}.(i) et à un analogue du lemme \ref{SublemTo6.8} pour $\mu \geq d \lambda$ si $d_0=1$. Cela conclut la preuve. \hfill$\Box$

\begin{prop}
Supposons $m=3$ et $d>1$. Soit $\chi : \Hh_{\Fpbar}(G,K,V) \to \Fpbar$ un caractère de $\Fpbar$-algèbres de paramètre de Hecke-Satake $(M,\chi_M)$. On suppose $M \neq G$ et on prend $\alpha \in \Delta \smallsetminus \Delta_M$. On suppose les hypothèses \ref{VtoV'} et \ref{Vnontropirreg} et on définit $V'$ comme en (\ref{defV'}). Supposons enfin l'une des deux conditions suivantes satisfaite:
\begin{itemize}
\item[(a)] $\widetilde{\alpha}^\vee(\varpi) \notin \widetilde{Z}_M(V)$;
\item[(b)] $\chi_M \big( \widetilde{\alpha}^\vee(\varpi) \big) \neq 1$.
\end{itemize}  
Alors $\varphi_\lambda^+$ induit l'isomorphisme \[ \ind_K^G \, V \otimes_{\Hh,\chi} \Fpbar \xrightarrow{\sim} \ind_K^G \, V' \otimes_{\Hh',\chi'} \Fpbar \] de $G$-représentations.
\end{prop}
\textsf{Preuve :} \\
Le même raisonnement que dans la première partie de la preuve de la proposition \ref{Cor6.11} nous indique qu'il s'agit de montrer $\chi( \varphi_\lambda^= * \varphi_\lambda^+) \neq 0$. On définit $d_0$ comme en (\ref{defd0}). Deux cas se présentent à nous: si on a $d_0 >1$, alors $\Sss_G(\varphi_\lambda^= * \varphi_\lambda^+)$ est $\tau_{d \lambda}$ par le lemme \ref{chgtpdsm=3}. Puisque $\widetilde{\lambda}(\varpi)^d$ est dans $\widetilde{Z}_M(V)$ et que $\chi$ se factorise par $\Sss_G^M$, on utilise par suite le lemme \ref{Lem4.3} pour affirmer \[ \chi(\varphi_\lambda^= * \varphi_\lambda^+) = \chi^{(A)}(\tau_{d \lambda}) \neq 0 .\] Si on a $d_0 =1$, $\Sss_G(\varphi_\lambda^= * \varphi_\lambda^+)$ est $\tau_{d \lambda} - \tau_{d \lambda + \alpha^\vee}$. On a de nouveau une subdivision de cas. D'abord, si $\widetilde{\alpha}^\vee(\varpi)$ n'appartient pas à $\widetilde{Z}_M(V)$, le lemme \ref{Lem4.3} nous dit: 
\[ \chi(\varphi_\lambda^= * \varphi_\lambda^+) = \chi^{(A)} (\tau_{d \lambda}) - \chi^{(A)} (\tau_{d \lambda + \alpha^\vee}) = \chi^{(A)} (\tau_{d \lambda}) \neq 0 .\] 
Enfin, si on a $\widetilde{\alpha}^\vee(\varpi) \in \widetilde{Z}_M(V)$ et $\chi_M \big( \widetilde{\alpha}^\vee(\varpi) \big) \neq 1$, on utilise la proposition \ref{Prop4.1}: \[ \chi( \varphi_\lambda^= * \varphi_\lambda^+ ) = \chi_M \big( \widetilde{\lambda}(\varpi)^d \big)^{-1} \big( 1 - \chi_M \big( \widetilde{\alpha}^\vee(\varpi) \big)^{-1} \big) \neq 0 .\] 
La preuve est terminée. \hfill$\Box$ \\

On reprend le formalisme précédent le corollaire \ref{chgtpdsdual}, et de la même manière on obtient:

\begin{cor} \label{chgtpdsdualm=3}
Supposons $m=3$ et $d>1$. Soit $\chi : \Hh_{\Fpbar}(G,K,V) \to \Fpbar$ un caractère de $\Fpbar$-algèbres de paramètre HS-dual $(M,\chi_M)$. On suppose $M \neq G$ et on prend $\alpha \in \Delta \smallsetminus \Delta_M$. On suppose les hypothèses \ref{VtoV'} et \ref{Vnontropirreg} pour $B^-$ et on définit $V'$ comme en (\ref{defV'2}). Supposons enfin l'une des deux conditions suivantes satisfaite:
\begin{itemize}
\item[(a)] $\widetilde{\alpha}^\vee(\varpi) \notin \widetilde{Z}_M(V^*)$;
\item[(b)] $\chi_M \big( \widetilde{\alpha}^\vee(\varpi) \big) \neq 1$.
\end{itemize} 
Alors $\varphi_\lambda^+$ induit l'isomorphisme \[ \ind_K^G \, V \otimes_{\Hh,\chi} \Fpbar \xrightarrow{\sim} \ind_K^G \, V' \otimes_{\Hh',\chi'} \Fpbar \] de $G$-représentations.
\end{cor}

\section{Irréductibilité des induites paraboliques} 

Dans tout ce paragraphe, on supposera la formule de comparaison de tranformées de Satake (conjecture \ref{conjSatake}) satisfaite. \\

Soit $P$ un parabolique standard de $G$ de radical unipotent $N$ et de Levi standard $M$. Ecrivons $M = \prod_i G_i$ avec $G_i = \GL(m_i,D)$ pour $1 \leq i \leq r$ et $\sum_i m_i = m$. Pour $1 \leq i \leq r$, soit $K_i$ le compact maximal $\GL(m_i, \Oo_D)$ de $G_i$.

\begin{lem} \label{decompGi}
$\phantom{,}$
\begin{itemize}
\item[(i)] Les représentations irréductibles admissibles de $M$ sont les $\sigma = \bigotimes_i \sigma_i$ pour $\sigma_i$ une représentation irréductible admissible de $G_i$.
\item[(ii)] Les $(M \cap K)$-poids de $\sigma$ sont de la forme $\bigotimes_i V_i$ où $V_i$ est un $K_i$-poids pour $\sigma_i$.
\item[(iii)] Les paramètres de Hecke-Satake correspondant aux vecteurs propres de $\Hom_K(V,\sigma)$ sont de la forme $(\prod_i M_i, \prod_i \chi_{M_i})$ où $(M_i,\chi_{M_i})$ est un paramètre de Hecke-Satake pour un vecteur propre de $\Hom_{K_i}(V_i,\sigma_i)$.
\end{itemize}
\end{lem}
\textsf{Preuve :} \\
C'est identique à \cite{Her11b}, Lemma 8.2. Revenons tout de m\^ eme un instant sur l'\' enonc\' e de factorisation en produit tensoriel du (i). Prouvons le cas $M = G_1 \times G_2$, le cas général s'en déduisant par récurrence. \\
Soient $\sigma_1$ et $\sigma_2$ deux représentations irréductibles admissibles de $G_1$ et $G_2$ respectivement. Montrons que $\sigma_1 \otimes \sigma_2$ est une représentation irréductible admissible de $M$. L'admissibilité de $\sigma_1 \otimes \sigma_2$ vient de ce que l'on a $(\sigma_1 \otimes \sigma_2)^{H_1 \times H_2} = \sigma_1^{H_1} \otimes \sigma_2^{H_2}$ pour tout sous-groupe compact ouvert $H_1 \times H_2 \leq G_1 \times G_2$. Voyons que $\sigma_1 \otimes \sigma_2$ est irréductible. Commençons par remarquer que, par lemme de Schur admissible, pour tout vecteur non nul $v_1$ de $\sigma_1$ et tout sous-$\Fpbar$-espace vectoriel $\pi_2$ de $\sigma_2$, \[ \begin{array}{ccc} \pi_2 & \to & \Hom_{G_1} (\sigma_1, \sigma_1 \otimes \pi_2) \\ v_2 & \mapsto & (v_1 \mapsto v_1 \otimes v_2) \end{array} \]
est un isomorphisme de $\Fpbar$-espaces vectoriels. Aussi, parce que toute sous-$G_1$-représentation $\pi$ de $\sigma_1 \otimes \sigma_2$ est isomorphe à une somme directe de $\sigma_1$, à nouveau grâce au lemme de Schur, pour toute telle $\pi$ on a l'isomorphisme de $G_1$-représentations
\[ \begin{array}{ccc} \sigma_1 \otimes \Hom_{G_1}(\sigma_1,\pi) & \xrightarrow{\sim} & \pi \\ v \otimes f & \mapsto & f(v) \end{array} .\]
On en déduit la bijection 
\[ \begin{array}{ccc} \{ \textrm{sous-} \Fpbar \textrm{-espaces vectoriels de } \sigma_2 \} & \xrightarrow{\sim} & \{ \textrm{sous-} G_1 \textrm{-représentations de } \sigma_1 \otimes \sigma_2 \} \\ \pi_2 & \to & \sigma_1 \otimes \pi_2 \\ \Hom_{G_1}(\sigma_1, \pi) & \mapsfrom & \pi \end{array} .\]
De ce fait, toute sous-$G_1$-représentation non nulle de $\sigma_1 \otimes \sigma_2$ est de la forme $\sigma_1 \otimes \pi_2$ où $\pi_2$ est une sous-$G_2$-représentation non nulle de $\sigma_2$. Parce que $\sigma_2$ est irréductible, on a $\pi_2 = \sigma_2$ et l'irréductibilité de $\sigma_1 \otimes \sigma_2$ voulue. \\
Montrons la réciproque. Soit $\sigma$ une représentation irréductible admissible de $M$. On choisit un vecteur non nul $v \in \sigma$ et on regarde la sous-$G_1$-représentation $\pi_1$ de $\sigma$ engendrée par $v$. Parce que $v$ est fixe par un sous-groupe compact ouvert $H_1 \times H_2$ de $M$, pour tout sous-groupe compact ouvert $K_1$ de $G_1$, $\pi_1^{K_1}$ est inclus dans $\sigma^{K_1 \times H_2}$ et est donc de dimension finie: $\pi_1$ est admissible. Par \cite{HenVig11b}, Lemma 7.10, il existe alors une sous-$G_1$-représentation irréductible admissible $\sigma_1$ dans $\pi_1$. Formons maintenant $\pi_2 = \Hom_{G_1}(\sigma_1,\sigma)$. C'est une $G_2$-représentation admissible: si on fixe $v_1 \in \sigma_1^{H_1}$ non nul, $f \in \pi_2$ est alors déterminée par $f(v_1)$, et si $f$ est dans $\pi_2^{H_2}$ pour un sous-groupe ouvert compact $H_2$ de $G_2$ alors $f(v_1)$ est dans $\sigma^{H_1 \times H_2}$, qui est de dimension finie. Une nouvelle application de \cite{HenVig11b}, Lemma 7.10 fournit une sous-$G_2$-représentation irréductible admissible $\sigma_2$ de $\pi_2$. Le morphisme non nul $\begin{array}{ccc} \sigma_1 \otimes \sigma_2 & \to & \sigma \\ v_1 \otimes f_2 & \mapsto & f_2(v_1) \end{array}$ est un isomorphisme par irréductibilité de $\sigma$ et de $\sigma_1 \otimes \sigma_2$. \hfill$\Box$ \\

On est maintenant en mesure de prouver le premier résultat important, qui est un critère d'irréductibilité d'une induite parabolique. On a besoin de rajouter soit une hypothèse sur le groupe $G$ lui-même, soit sur le parabolique $P$ à partir duquel on induit.

\begin{hypo} \label{Dpremier}
L'une des deux conditions suivantes est satisfaite: 
\begin{itemize}
\item[(a)] le degré $d$ est un nombre premier, ou égal à $1$;
\item[(b)] l'entier $m$ est inférieur ou égal à $3$.
\end{itemize}
\end{hypo}

\begin{hypo} \label{Pminimal}
Le parabolique $P$ est minimal, égal à $B$.
\end{hypo}

\begin{theo} \label{irreductibilite}
Supposons l'hypothèse \ref{Dpremier} ou l'hypothèse \ref{Pminimal} satisfaite. Soit, pour tout $1 \leq i \leq r$, $\sigma_i$ une représentation irréductible admissible de $G_i$ qui soit dans l'un des cas suivants:
\begin{itemize}
\item[(a)] $\sigma_i$ est supersingulière avec $m_i > 1$;
\item[(b)] $\sigma_i$ est une représentation de $D^\times$ (i.e. $m_i = 1$) avec $\dim \, \sigma_i > 1$;
\item[(c)] $\sigma_i \simeq \St_{Q_i} \rho_i^0$ pour un parabolique standard $Q_i \leq G_i$ et un caractère $\rho_i : D^\times \xrightarrow{\Nrm} F^\times \xrightarrow{\rho_i^0} \Fpbar^\times$.
\end{itemize}
Supposons $\rho_i \neq \rho_{i+1}$ s'il y a deux blocs adjacents dans le cas (c). Alors $\Ind_{P}^G(\sigma_1 \otimes \dots \otimes \sigma_r)$ est irréductible admissible.
\end{theo}
\textsf{Preuve :} \\
On pose $\sigma = \sigma_1 \otimes \cdots \otimes \sigma_r$. Commençons par l'admissibilité\footnote{cet argument est en fait présenté en toute généralité au (2) du chapitre I, paragraphe 5.6 de \cite{Vig96}}. Grâce à la décomposition d'Iwasawa (voir \cite{HenVig11}, Proposition 6.5), l'injection naturelle \[ P \cap K \backslash K \hookrightarrow P \backslash G \] est un isomorphisme. Alors, par restriction des fonctions à $K$, $(\Ind_P^G \, \sigma)^{K(h)}$ est égal à $(\Ind_{P \cap K}^K \, \sigma)^{K(h)}$ pour $h \geq 1$. En réduisant modulo $K(h)$, on obtient alors que ce dernier est isomorphe à $\Ind_{\overline{P}^{(h)}}^{\overline{G}^{(h)}}(\sigma^{M(h)}) $ où $M(h)$ désigne $M \cap K(h)$; cet espace est de dimension finie par admissibilité de $\sigma$. \\
Venons-en à l'irréductibilité. Soient $\pi$ une sous-$G$-représentation irréductible non nulle de $\Ind_P^G \, \sigma$ et $V$ un $K$-poids de $\pi$. Par admissibilité de $\pi$, $\Hom_K(V,\pi)$ est de dimension finie sur $\Fpbar$. La $\Fpbar$-algèbre commutative $\Hh_{\Fpbar}(G,K,V)$ agit sur $\Hom_K(V,\pi)$, et on choisit alors un vecteur propre $f \neq 0$ pour cette action. Il lui est ainsi associé un caractère $\chi_f : \Hh_{\Fpbar}(G,K,V) \to \Fpbar$ de $\Fpbar$-algèbres et un morphisme surjectif de $G$-représentations
\[ \Phi_f : \ind_K^G \, V \otimes_{\Hh, \chi_f} \Fpbar \twoheadrightarrow \pi .\]
Par le lemme \ref{paramHSinduite}, $\chi_f$ se factorise en \[ \Hh_{\Fpbar}(G,K,V) \mathop{\hookrightarrow}\limits_{{}' \Sss_G^M} \Hh_{\Fpbar} \big( M, M \cap K , V_{N \cap K} \big) \mathop{\to}\limits_{\chi_0} \Fpbar .\] De plus, $\chi_0$ et $\chi_f$ ont même paramètre HS-dual $(L,\chi_L)$ avec $L \leq M$. Soit $P'$ le parabolique standard contenu dans $P$ de radical unipotent $N'$ et de Levi standard $M' = \prod_i M_i'$ défini par: 
\begin{itemize}
\item $M_i' = G_i$ si $\sigma_i$ est dans le cas $(a)$ ou $(b)$;
\item $M_i'$ est le Levi standard de $Q_i$ si $\sigma_i$ est dans le cas $(c)$.
\end{itemize}
En particulier, on remarque que $P'$ est $P$ si et seulement si toutes les Steinberg généralisées \og dans $\sigma$ \fg\ sont des caractères. Par la proposition \ref{paramSt}, $\chi_0$ se factorise en \[ \Hh_{\Fpbar}(M,M \cap K, V_{N \cap K}) \mathop{\hookrightarrow}_{{}' \Sss_M^{M'}} \Hh_{\Fpbar}(M',M' \cap K, V_{N' \cap K}) \mathop{\to}_{\chi} \Fpbar ,\] et $\chi$ a encore $(L,\chi_L)$ comme paramètre HS-dual, donnant donc $L \leq M'$. On définit de plus la $M'$-représentation $\sigma' = \bigotimes_i \sigma_i'$ par 
\begin{itemize}
\item $\sigma_i' = \sigma_i$ si dans le cas $(a)$ ou $(b)$;
\item $\sigma_i' = \rho_i^0 \circ \det_{M_i'}$ si dans le cas $(c)$.
\end{itemize}
Par r\' eciprocit\' e de Frobenius, on a 
\[ \Hom_K(V,\pi) \subseteq \Hom_K \big( V, \Ind_P^G \, \sigma \big) \simeq \Hom_{P \cap K}(V,\sigma) = \Hom_{M \cap K} \big( V_{N \cap K}, \sigma \big) ,\]
et on notera $f_1: V_{N \cap K} \to \sigma$ le morphisme de $(M \cap K)$-repr\' esentations qui est l'image de $f \in \Hom_K(V,\pi)$ par cette injection. Par le lemme \ref{decompGi}, on peut \' ecrire $V_{N \cap K} = \bigotimes_i V_i$ et on a des fl\` eches de $K_i$-repr\' esentations $V_i \to \sigma_i$. Lorsque $\sigma_i$ est une repr\' esentation de Steinberg g\' en\' eralis\' ee, par la proposition \ref{Stsocle}, on a une fl\` eche de $(K_i \cap M_i')$-repr\' esentations $(V_i)_{N_i' \cap K_i} \to \sigma_i'$. Lorsque $\sigma_i$ est d'un autre type, on a $M_i' = M_i$ et il n'y a rien \` a dire pour pouvoir affirmer l'existence d'un morphisme non nul de $(M' \cap K)$-repr\' esentations $V_{N' \cap K} \to \sigma'$. Ainsi $f_1$ se factorise par 
\[ V_{N \cap K} \twoheadrightarrow V_{N' \cap K} \xrightarrow{f_2} \sigma' \subseteq \sigma \] avec $f_2$ un morphisme de $(M' \cap K)$-représentations. \\
Puisque $\sigma'$ est une $M'$-représentation irréductible, on a la surjection \[ \ind_{M' \cap K}^{M'} \, V_{N' \cap K} \otimes_{\Hh_{M'} , \chi} \Fpbar \twoheadrightarrow \sigma' ,\]
que l'on notera encore $f_2$. Par exactitude du foncteur $\Ind_{P' \cap M}^M$ (voir \cite{Eme10}, Proposition 4.1.5, et \cite{Vig12}, Proposition 1.1), on a la surjection de $M$-représentations \[ \Ind \, f_2 : \Ind_{P' \cap M}^M \big( \ind_{M' \cap K}^{M'} \, V_{N' \cap K} \otimes_{\Hh_{M'} , \chi} \Fpbar \big) \twoheadrightarrow \Ind_{P' \cap M}^M \, \sigma' .\]
Par l'irréductibilité de $\sigma$ et la définition des Steinberg généralisées, on en déduit la surjection de $M$-représentations \[ f_3 : \Ind_{P' \cap M}^M \big( \ind_{M' \cap K}^{M'} \, V_{N' \cap K} \otimes_{\Hh_{M'} , \chi} \Fpbar \big) \twoheadrightarrow \Ind_{P' \cap M}^M \, \sigma' \twoheadrightarrow \sigma .\] Par transitivité de $\Ind$ et exactitude de $\Ind_P^G$, on a finalement le morphisme $G$-équivariant \[ \Ind \, f_3 : \Ind_{P'}^G \big( \ind_{M' \cap K}^{M'} \, V_{N' \cap K} \otimes_{\Hh_{M'} , \chi} \Fpbar \big) \twoheadrightarrow \Ind_P^G \, \sigma .\]
Si $V$ est $(P')^-$-régulière, par la proposition \ref{Thm3.1}, on obtient \[ (\Ind \, f_3 ) \circ \iota_V : \ind_K^G \, V \otimes_{\Hh, \chi_f} \Fpbar \xrightarrow{\sim} \Ind_{P'}^G \big( \ind_{M' \cap K}^{M'} \, V_{N' \cap K} \otimes_{\Hh_{M'} , \chi} \Fpbar \big) \twoheadrightarrow \Ind_P^G \, \sigma .\] 
Vérifions que $(\Ind \, f_3) \circ \iota_V$ s'identifie à 
\[ \Phi_f : \ind_K^G \, V \otimes_{\Hh, \chi_f} \Fpbar \twoheadrightarrow \pi \subseteq \Ind_P^G \, \sigma .\]
Parce que l'on a 
\[ \iota_V ([1,v] \otimes 1)(1) = [1,p_{N'}(v)] \otimes 1 ,\] et donc 
\[ \big( (\Ind f_3) \circ \iota_V \big) ([1,v] \otimes 1)(1) = f_2 \circ p_{N'}(v) = f(v)(1) \]
pour tout $v \in V$, on a bien $(\Ind f_3) \circ \iota_V = \Phi_f$. Il en résulte $\pi = \Ind_P^G \, \sigma$ et l'irréductibilité de $\Ind_P^G \, \sigma$ comme voulu. \\
Si $V$ n'est pas $(P')^-$-régulière, par la proposition \ref{chgtpdsreg} ci-dessous que nous admettons pour l'instant, il existe une $K$-représentation irréductible $(P')^-$-régulière $V'$ conjuguée à $V$ et vérifiant \[ \ind_K^G \, V \otimes_{\Hh, \chi_f} \Fpbar \xrightarrow{\sim} \ind_K^G \, V' \otimes_{\Hh, \chi_f} \Fpbar .\] Ainsi on a une surjection $\ind_K^G \, V' \otimes_{\Hh, \chi_f} \Fpbar \twoheadrightarrow \pi$, et $\pi|_K$ contient $V'$. De cette manière, on peut réitérer le raisonnement précédent avec $V'$ au lieu de $V$ et conclure quant à l'irréductibilité de $\Ind_P^G \, \sigma$. \hfill$\Box$\\

Par la remarque suivant le lemme \ref{paramHSinduite}, on va pouvoir se placer sous l'hypothèse suivante pour les $K$-poids de $\pi$ où $\pi$ est une sous-représentation de $\Ind_B^G \, \sigma$; en d'autres termes, \og l'hypothèse \ref{Pminimal} implique l'hypothèse \ref{chifactorbas} \fg.

\begin{hypo} \label{chifactorbas}
Le paramètre HS-dual du caractère $\chi : \Hh_{\Fpbar}(G,K,V) \to \Fpbar$ est de la forme $(A, \chi_A)$ où $\chi_A$ est un caractère $\widetilde{Z}_A(V^*) \to \Fpbar^\times$.
\end{hypo}

\begin{prop} \label{chgtpdsreg}
Soient $P$ et $\sigma$ comme dans le théorème \ref{irreductibilite}, $P'$ comme dans sa preuve. Soient $V$ un $K$-poids de $\Ind_P^G \, \sigma$ et $\chi: \Hh_{\Fpbar}(G,K,V) \to \Fpbar$ un caractère propre sur $\Hom_K(V,\Ind_P^G \, \sigma)$. Supposons l'hypothèse \ref{Dpremier} ou l'hypothèse \ref{chifactorbas} satisfaite. Alors il existe une $K$-représentation irréductible $V'$ satisfaisant aux trois conditions suivantes:
\begin{itemize}
\item[(i)] $V_{U \cap K}$ et $V'_{U \cap K}$ sont conjuguées en tant que $(A \cap K)$-représentations (de sorte que l'on peut identifier $\Hh_{\Fpbar}(G,K,V)$ et $\Hh_{\Fpbar}(G,K,V')$; on notera $\Hh$ chacune de ces deux algèbres);
\item[(ii)] on a un isomorphisme de $G$-représentations \[ \ind_K^G \, V \otimes_{\Hh,\chi} \Fpbar \xrightarrow{\sim} \ind_K^G \, V' \otimes_{\Hh,\chi} \Fpbar ;\]
\item[(iii)] $V'$ est $(P')^-$-régulière.
\end{itemize}
\end{prop}

On repousse la preuve de cette proposition aux trois sous-sections suivantes. \\ Pour l'instant, récoltons les fruits du théorème \ref{irreductibilite} à travers un premier corollaire. On veut supprimer l'hypothèse $\rho_i \neq \rho_{i+1}$ des blocs adjacents de type \textit{(c)} et voir ce qui en découle. Pour cela, si on a la même donnée que dans l'énoncé du théorème \ref{irreductibilite}, lorsque deux blocs de type \textit{(c)} sont adjacents, on les regroupe. Ainsi si $G_{i_k+1}, G_{i_k+2}, \dots , G_{i_{k+1}}$ sont des blocs adjacents de type \textit{(c)} avec $\rho_{i_k+1} = \rho_{i_k+2} = \dots = \rho_{i_{k+1}}$ (et que les blocs alentours ne sont soit pas du type \textit{(c)}, soit ont un $\rho_i$ différent), alors on note $G'_k$ le facteur $\GL(m_{i_k+1} + \cdots + m_{i_{k+1}},D)$ dont $M_k := G_{i_k+1} \times \cdots \times G_{i_{k+1}}$ est un Levi. On laisse les autres $G_i$ inchangés (sauf qu'on les a renumérotés). On note $Q$ le parabolique standard de Levi $G'_1 \times \cdots \times G'_s$.

\begin{cor} \label{JHtotal}
Soit, pour tout $1 \leq i \leq r$, $\sigma_i$ une représentation irréductible admissible de $G_i$ qui soit dans l'un des cas\footnote{on ne suppose cependant pas $\rho_i \neq \rho_{i+1}$ pour les blocs adjacents de type \textit{(c)}} (a), (b), (c) du théorème \ref{irreductibilite}. Alors $\Ind_P^G(\sigma_1 \otimes \cdots \otimes \sigma_r)$ est de longueur finie, de constituants de Jordan-Hölder les $\Ind_Q^G(\sigma'_1 \otimes \cdots \otimes \sigma'_s)$ où $\sigma'_i$ est la représentation $\sigma_{j_i}$ associée si on n'a pas changé le bloc et $\sigma'_i$ parcourt l'ensemble des $\St_{R_i} \rho_i^0$ avec \[ \{\Delta_{R_i} \supseteq \Delta_{M_i} \ | \ \Delta_{M'_i} \smallsetminus \Delta_{M_i} \subseteq \Delta \smallsetminus \Delta_{R_i} \} \] pour des blocs adjacents de type (c). Chacune de ces représentations irréductibles admissibles y apparaît avec multiplicité $1$.
\end{cor}
\textsf{Preuve :} \\
C'est une conséquence du corollaire \ref{GLmJHsteinberg}, du théorème \ref{irreductibilite} et de l'exactitude de $\Ind_Q^G$ (voir \cite{Eme10}, Proposition 4.1.5, et \cite{Vig12}, Proposition 1.1). \hfill$\Box$\\

Soit $P''$ le parabolique standard de radical unipotent $N''$ et de Levi standard $M'' = \prod_i M_i''$ défini par:
\begin{itemize}
\item $M_i'' = G_i$ si $\sigma_i$ est dans le cas $(a)$ ou $(b)$;
\item $M_i'' = A \cap G_i$ si $\sigma_i$ est dans le cas $(c)$.
\end{itemize}
Avant que d'attaquer la preuve de la proposition \ref{chgtpdsreg}, qui diffère selon l'hypothèse vérifiée, on va préciser un peu le paramètre HS-dual de $\chi$. 

\begin{lem} \label{paramHSpoids}
Soient $P$ et $\sigma$ comme dans le théorème \ref{irreductibilite}, $M''$ comme précédemment. Soient $V$ un $K$-poids de $\Ind_P^G \, \sigma$ et $\chi : \Hh_{\Fpbar}(G,K,V) \to \Fpbar$ un caractère propre sur $\Hom_K(V,\Ind_P^G \, \sigma)$ de paramètre HS-dual $(L,\chi_L)$. Alors on a $L = M''$.
\end{lem}
\textsf{Preuve :} \\
On a déjà remarqué (grâce au lemme \ref{paramHSinduite} et à la proposition \ref{paramSt}) l'inclusion $L \leq M''$. On a la factorisation suivante de $\iota_G^*(\chi)$: 
\[\xymatrix{ \Hh_{\Fpbar}(G,K,V^*) \ar@{^{(}->}[r]_{\Sss_G^{M'} \phantom{xxxxxxxx}} \ar@/^2pc/[rr]^{\iota_G^*(\chi)} & \Hh_{\Fpbar} \big( M'', M'' \cap K , (V^*)^{N'' \cap K} \big) \ar[r] \ar[d]_{\iota_{M''}} & \Fpbar \\ & \Hh_{\Fpbar} \big( M'', M'' \cap K, V_{N'' \cap K} \big) \ar[ur]_{\chi'} & }.\]
On écrit \[ V_{N'' \cap K} = V_1 \otimes \cdots \otimes V_r \] où chaque $V_i$ est une représentation irréductible de $M_i'' \cap K$. Alors, par le lemme \ref{decompGi}, le paramètre HS-dual de $\chi'$ est $(L,\chi_L) =  \big( \prod_i L_i , \prod_i \chi_{L_i} \big)$ où on a $L_i \leq M_i''$ pour tout $i$. \\
Lorsque $\sigma_i$ est une représentation de $D^\times$ ou une représentation de Steinberg généralisée, $M_i''$ ne possède pas de sous-groupe de Levi propre et on a $L_i = M_i''$. Lorsque $\sigma_i$ est une représentation supersingulière, on a $L_i = M_i = M_i''$ par définition. Le résultat en découle. \hfill$\Box$\\

\textsf{Remarque :}
Le premier cas que l'on n'arrive pas à résoudre (proposition \ref{chgtpdsreg}, théorème \ref{irreductibilite}) est celui de $d>1$ entier composé et $m=4$ avec $m_1 = m_2 =2$, $\sigma_1$ et $\sigma_2$ deux représentations supersingulières.

\subsection{Preuve de la proposition \ref{chgtpdsreg} sous l'hypothèse \ref{Dpremier}.(a)} \label{GLmpreuve1}

On rappelle que l'hypothèse \ref{Dpremier}.(a) nous place dans les conditions \og $d$ premier ou égal à $1$ \fg. On note $\Par_{\overline{B}}(V^*) = (\chi_{V^*}, \Delta_V^-)$, de sorte que l'on a $\Par_{\overline{B}^-}(V) = (\chi_{V^*}^{-1},-\Delta_V^-)$ par la proposition \ref{paramV*}. On procède par récurrence descendante sur l'ensemble fini $\Delta_V^- \smallsetminus \Delta_V^- \cap \Delta_{M'}$. Si ce dernier est vide, alors on a $\Delta_V^- \subseteq \Delta_{M'}$ et $V$ est $(P')^-$-régulière: on peut prendre $V = V'$. \\
Supposons donc $\Delta_V^- \smallsetminus \Delta_V^- \cap \Delta_{M'} \neq \varnothing$, et prenons $\alpha$ une racine de cet ensemble. Par définition, l'hypothèse \ref{VtoV'} pour $B^-$ est satisfaite. On définit $\lambda$, $t$ comme avant le corollaire \ref{chgtpdsdual}, $d_0$ associé pour $V_{N_{(\alpha)} \cap K}$ et $V'$ comme en (\ref{defV'2}). A cause de l'hypothèse \ref{Dpremier}.(a), l'hypothèse \ref{V'toV} est satisfaite (puisque $d_0>1$ implique $d_0 = d$). Il reste à montrer $\widetilde{\alpha}^\vee(\varpi) \notin \widetilde{Z}_L(V^*)$ ou $\chi_L \big( \widetilde{\alpha}^\vee(\varpi) \big) \neq 1$ pour pouvoir appliquer le corollaire \ref{chgtpdsdual}. Grâce au lemme \ref{paramHSpoids}, on a $L = M''$. Par hypothèse, on a $\alpha \notin \Delta_{M'}$. De plus, $V_{N \cap K}$ se décompose en $V_{N \cap K} = V_1 \otimes \cdots \otimes V_r$ où $V_i$ est une représentation irréductible de $K_i$. Dans le cas où $\sigma_i$ est une représentation de Steinberg généralisée, par la proposition \ref{Stsocle}, $V_i$ est uniquement déterminée et est $Q_i^-$-régulière. De ce fait, $\alpha$ n'appartient pas à $\Delta_M$. Notons $i \in [1,r-1]$ l'entier tel que $\alpha \notin \Delta_M$ est situé entre les blocs $G_i$ et $G_{i+1}$. \\
Si $\sigma_i$ ou  $\sigma_{i+1}$ est une représentation de dimension $>1$ de $D^\times$, alors on a $\widetilde{\alpha}^\vee(\varpi) \notin \widetilde{Z}_A(V^*) \supseteq \widetilde{Z}_L(V^*)$. Si $\sigma_i$ ou $\sigma_{i+1}$ est une représentation supersingulière, on a $\widetilde{\alpha}^\vee(\varpi) \notin \widetilde{Z}_L \supseteq \widetilde{Z}_L(V^*)$. Si enfin $\sigma_i$ et $\sigma_{i+1}$ sont toutes deux des représentations de Steinberg généralisées, on a $\widetilde{\alpha}^\vee(\varpi) \in \widetilde{Z}_L(V^*)$. Mais alors, comme on a $\rho_i |_{\Oo_D^\times} = \rho_{i+1} |_{\Oo_D^\times}$ (puisque $\alpha$ est dans $\Delta_V^-$) et $\rho_i \neq \rho_{i+1}$, on obtient \[ \chi_L \big( \widetilde{\alpha}^\vee(\varpi) \big) = \rho_{i+1}(\varpi_D) \rho_i(\varpi_D)^{-1} \neq 1.\]
Dans tous les cas, on peut appliquer le corollaire \ref{chgtpdsdual} et obtenir \[ \ind_K^G \, V \otimes_{\Hh,\chi} \Fpbar \xrightarrow{\sim} \ind_K^G \, V' \otimes_{\Hh,\chi} \Fpbar .\] 
De ce fait, $V'$ vérifie $(i)$, $(ii)$ et \[ \Delta_{V'}^- \smallsetminus \Delta_{V'}^- \cap \Delta_{M'} \subsetneq \Delta_V^- \smallsetminus \Delta_V^- \cap \Delta_{M'} .\] Par récurrence descendante, on obtient le résultat voulu.

\subsection{Preuve de la proposition \ref{chgtpdsreg} sous l'hypothèse \ref{chifactorbas}} \label{GLmpreuve2}

L'hypothèse \ref{chifactorbas} nous dit $L=A$, et on cherche $V'$ $(P')^-$-régulière vérifiant $(i)$ et $(ii)$. On note $\Par_{\overline{B}}(V^*) = (\chi_{V^*}, \Delta_V^-)$. On procède par récurrence descendante sur l'ensemble fini $\Delta_V^- \smallsetminus \Delta_V^- \cap \Delta_{M'}$. Si ce dernier est vide, alors $V$ est $(P')^-$-régulière et on prend $V'=V$. \\
Supposons donc $\Delta_V^- \smallsetminus \Delta_V^- \cap \Delta_{M'} \neq \varnothing$, et prenons $\alpha$ une racine de cet ensemble. On note $P_V^{(-)}$ le parabolique antistandard\footnote{attention, ce n'est pas le parabolique opposé à $P_V$} tel que $\big( P_V^{(-)} \cap K \big) K(1)$ est le stabilisateur de la droite $V^{U^- \cap K} \simeq V_{U \cap K}$, et $P_V^{(-)} = M_V^{(-)} N_V^{(-)}$ sa décomposition de Levi antistandard. On note $P_V^{(+)} = M_V^{(-)} N_V^{(+)}$ le parabolique opposé à $P_V^{(-)}$, et qui est donc standard. Parce que $V$ est $P_V^{(-)}$-régulière, par la proposition \ref{Thm3.1}, on a l'isomorphisme de $G$-représentations \[ \ind_K^G \, V \otimes_{\Hh,\chi} \Fpbar \xrightarrow{\sim} \Ind_{P_V^{(+)}}^G \big( \ind_{M_V^{(-)} \cap K}^{M_V^{(-)}} \, V_{N_V^{(+)} \cap K} \otimes_{\Hh_{M_V^{(-)}}, \chi'} \Fpbar \big) ,\] où on a factorisé $\chi$ en (ce qui est loisible par le lemme \ref{Lem4.3} et la proposition \ref{Prop4.1} au vu de $L=A$): \[ \Hh_{\Fpbar}(G,K,V) \mathop{\hookrightarrow}\limits_{{}' \Sss_G^{M_V^{(-)}}} \Hh_{\Fpbar} \big( M_V^{(-)} , M_V^{(-)} \cap K , V_{N_V^{(+)} \cap K} \big) \mathop{\to}\limits_{\chi'} \Fpbar .\] 
On décompose ensuite $M_V^{(-)} = M^{(1)} \times \cdots \times M^{(s)}$ où chaque $M^{(i)}$ est isomorphe à un $\GL(m_i^{(-)},D)$. La représentation $V_{N_V^{(+)} \cap K}$ se décompose alors en \[ V_{N_V^{(+)} \cap K} = V_1 \otimes \cdots \otimes V_s\] où chaque $V_i$ est un caractère de $M^{(i)} \cap K$. On écrit aussi 
\[ \Hh_{\Fpbar} \big( M_V^{(-)} , M_V^{(-)} \cap K , V_{N_V^{(+)} \cap K} \big) = \bigotimes\nolimits_i \Hh_{\Fpbar}(M^{(i)}, M^{(i)} \cap K, V_i) \] et $\chi' = \prod_i \chi_i'$ avec $\chi_i' : \Hh_{\Fpbar}(M^{(i)}, M^{(i)} \cap K, V_i) \to \Fpbar$ un caractère de $\Fpbar$-algèbres pour tout $i$. \\
Il existe un unique entier $j \in [1,s]$ tel que la racine $\alpha$ appartient à $\Delta_{M^{(j)}}$. On ne travaille alors plus qu'à l'intérieur du groupe réductif $M^{(j)}$. On a $\alpha \in \Delta_{V_j}^- \subseteq \Delta_V^-$ et l'hypothèse \ref{VtoV'} pour $B^- \cap M^{(j)}$ est vérifiée; de plus, comme $V_j$ est un caractère de $M^{(j)} \cap K$, l'hypothèse \ref{V'toV} est vérifiée. On prend $\lambda, t, V_j'$ associés pour $(V_j)_{N_{(\alpha)} \cap M^{(j)} \cap K}$ dans $M^{(j)}$. Par le même raisonnement qu'au début du paragraphe \ref{GLmpreuve1}, on a en fait $\alpha \notin \Delta_M \supseteq \Delta_{M'}$. On note alors $i \in [1,r-1]$ l'entier tel que $\alpha$ est situé entre les blocs $G_i$ et $G_{i+1}$. \\
A cause de l'hypothèse \ref{chifactorbas}, ni $\sigma_i$ ni $\sigma_{i+1}$ ne peut être supersingulière. Si $\sigma_i$ ou $\sigma_{i+1}$ est une représentation de $D^\times$ de dimension $>1$, alors on $\widetilde{\alpha}^\vee(\varpi) \notin \widetilde{Z}_A(V^*) = \widetilde{Z}_L(V^*)$. Si $\sigma_i$ et $\sigma_{i+1}$ sont toutes deux des représentations de Steinberg généralisées, on a $\widetilde{\alpha}^\vee(\varpi) \in \widetilde{Z}_L(V^*)$ mais \[ (\chi_j')_A \big( \widetilde{\alpha}^\vee(\varpi) \big) = \rho_{i+1}(\varpi_D) \rho_i(\varpi_D)^{-1} \neq 1 ,\] où $\big( A \cap M^{(j)}, (\chi_j')_A \big)$ désigne le paramètre HS-dual de $\chi_j'$. Dans tous les cas on peut appliquer la proposition \ref{chgtpdsdual} et obtenir 
\begin{equation} \label{isomMj} \ind_{M^{(j)} \cap K}^{M^{(j)}} \, V_j \otimes_{\Hh_{M^{(j)}}, \chi_j'} \Fpbar \xrightarrow{\sim} \ind_{M^{(j)} \cap K}^{M^{(j)}} \, V_j' \otimes_{\Hh_{M^{(j)}}, \chi_j'} \Fpbar .\end{equation}
Notons 
\[ \Par_{\overline{B}}(V^*) = \big( \chi_1 \otimes \cdots \otimes \chi_{j-1} \otimes \chi_j \otimes \chi_{j+1} \otimes \cdots \otimes \chi_s, \Delta_V^- \big) \]
où chaque $\chi_i$ est un caractère de $\overline{A \cap M^{(i)}}$. On définit le caractère $\chi'_j$ de $\overline{A \cap M^{(j)}}$ égal à $\chi_j(t^{-1} \cdot t)$ (on rappelle que $t$ est associé à $\alpha$), puis on note $V'$ l'unique représentation irréductible de $K$ (le lien explicite étant donné par la proposition \ref{paramV*}) vérifiant 
\begin{equation} \label{defV'3} \Par_{\overline{B}} \big( (V')^* \big) = \big( \chi_1 \otimes \cdots \otimes \chi_{j-1} \otimes \chi'_j \otimes \chi_{j+1} \otimes \cdots \otimes \chi_s , \Delta_V^- \smallsetminus \{ \alpha \}\big) ,\end{equation} 
ce qui est bien autorisé car on a (avec la notation de (\ref{defDeltachi}))
\[ \Delta_{\chi_1 \otimes \cdots \otimes \chi_{j-1} \otimes \chi'_j \otimes \chi_{j+1} \otimes \cdots \otimes \chi_s} = \Delta_{\chi_1 \otimes \cdots \otimes \chi_{j-1} \otimes \chi_j \otimes \chi_{j+1} \otimes \cdots \otimes \chi_s} \smallsetminus \{ \alpha \}  .\] 
Grâce à (\ref{isomMj}), on a \[ \ind_{M_V^{(-)} \cap K}^{M_V^{(-)}} \, V_{N_V^{(+)} \cap K} \otimes_{\Hh_{M_V^{(-)}}, \chi'} \Fpbar \xrightarrow{\sim} \ind_{M_V^{(-)} \cap K}^{M_V^{(-)}} \, V_{N_V^{(+)} \cap K}' \otimes_{\Hh_{M_V^{(-)}}, \chi'} \Fpbar .\] De plus, $V'$ est a fortiori $P_V^{(-)}$-régulière, et en appliquant la proposition \ref{Thm3.1}, on a :
\[\xymatrix{ \ind_K^G \, V \otimes_{\Hh,\chi} \Fpbar \ar[r]^{\sim \phantom{xxxxxxxxxxxxxx}} & \Ind_{P_V^{(+)}}^G \big( \ind_{M_V^{(-)} \cap K}^{M_V^{(-)}} \, V_{N_V^{(+)} \cap K} \otimes_{\Hh_{M_V^{(-)}}, \chi'} \Fpbar \big) \ar[d]^\sim \\ \ind_K^G \, V' \otimes_{\Hh,\chi} \Fpbar \ar[r]^{\sim \phantom{xxxxxxxxxxxxxx}} & \Ind_{P_V^{(+)}}^G \big( \ind_{M_V^{(-)} \cap K}^{M_V^{(-)}} \, V_{N_V^{(+)} \cap K}' \otimes_{\Hh_{M_V^{(-)}}, \chi'} \Fpbar \big) }.\] Enfin, $V_{U \cap K}'$ et $V_{U \cap K}$ sont conjugués (voir (\ref{defV'3})) et on a \[ \Delta_{V'}^- \smallsetminus \Delta_{V'}^- \cap \Delta_{M'} \subsetneq \Delta_V^- \smallsetminus \Delta_V^- \cap \Delta_{M'} .\] Par récurrence descendante, on a fini.

\subsection{Preuve de la proposition \ref{chgtpdsreg} sous l'hypothèse \ref{Dpremier}.(b)}

On rappelle que l'hypothèse \ref{Dpremier}.(b) est \og $m \leq 3$ \fg. Dans le cas $m=1$, il n'y a rien à faire. Lorsque $m$ est égal à $2$, de deux choses l'une: ou bien $P$ est égal $G$ et il n'y a rien à prouver, ou bien $P$ est minimal et l'hypothèse \ref{chifactorbas} est satisfaite, de sorte que l'on peut utiliser le paragraphe \ref{GLmpreuve2} (en une version simplifiée dans laquelle la récurrence se termine au bout d'une étape au maximum). Supposons donc $m=3$ et examinons les différents cas possibles. \\
Si $P$ est égal à $G$, $\sigma$ est une représentation supersingulière de $G$ et on a $P'=P=G$: il n'y a rien à prouver. Si $P$ est minimal (c'est-à-dire égal à $B$), on utilise à nouveau le paragraphe \ref{GLmpreuve2} pour conclure. On note $\Par_{\overline{B}}(V^*) = (\chi_{V^*}, \Delta_V^-)$ et $\alpha$ la racine $\left( \begin{array}{ccc} x & & \\ & y & \\ & & z \end{array} \right) \mapsto xy^{-1}$; on appelle $\beta$ l'autre racine de $\Delta$. Il nous reste à examiner les cas $\Delta_P = \{ \alpha \}$ ou $\Delta_P = \{ \beta \}$. \\

Supposons dans un premier temps $\Delta_P = \{ \beta \}$. On suppose $\Delta_V^- \smallsetminus \Delta_V^- \cap \Delta_{M'}$ non vide (sinon il n'y a rien à faire) et on prend $\delta$ dans cet ensemble. Par le même raisonnement qu'au paragraphe \ref{GLmpreuve1}, on a $\delta \notin \Delta_M$, et donc $\delta = \alpha$. Prenons alors $\lambda$ le copoids fondamental minuscule associé à $\alpha$ de dernier coefficient diagonal \' egal \` a $1$, et $t = \widetilde{\lambda}(\varpi) \in A_\Lambda^+$. On prend $d_0$ associé pour $V_{N_{(\alpha)} \cap K}$ et $V'$ comme en (\ref{defV'2}). On remarque que l'expression de $t$ ($t$ associé à $\alpha \neq \beta$ et $\Delta_\chi \supseteq \Delta_P = \{\beta\}$) implique que l'hypothèse \ref{V'toV} est automatiquement satisfaite (puisque l'on a $\chi(t^{-d_0} \cdot t^{d_0}) = \chi$). Récapitulons la situation: $\sigma_1$ est une représentation de $D^\times$ de dimension $>1$ ou bien un caractère de $D^\times$; et $\sigma_2$ est une représentation supersingulière ou bien une Steinberg généralisée. Si $\sigma_1$ est de dimension $>1$ ou si $\sigma_2$ est une supersingulière, alors on a $\widetilde{\alpha}^\vee(\varpi) \notin \widetilde{Z}_L(V^*)$. Si $\sigma_1$ est un caractère et $\sigma_2$ une Steinberg généralisée, par le lemme \ref{paramHSpoids}, on a $L=A$ et $\widetilde{\alpha}^\vee(\varpi) \in \widetilde{Z}_L(V^*)$. Mais alors on calcule \[ \chi_L \big( \widetilde{\alpha}^\vee(\varpi) \big) = \rho_2(\varpi_D) \rho_1 (\varpi_D)^{-1} \neq 1 .\] 
Dans tous les cas, on peut appliquer le corollaire \ref{chgtpdsdual} et obtenir: \[ \ind_K^G \, V \otimes_{\Hh,\chi} \Fpbar \xrightarrow{\sim} \ind_K^G \, V' \otimes_{\Hh,\chi} \Fpbar .\] La représentation $V'$ est $(P')^-$-régulière (à cause éventuellement de la proposition \ref{Stsocle}) et convient. \\

Supposons maintenant $\Delta_P = \{\alpha\}$. On suppose encore $\Delta_V^- \smallsetminus \Delta_V^- \cap \Delta_{M'}$ non vide et on prend $\delta$ dans cet ensemble: on a cette fois-ci $\delta = \beta$, par la même argumentation. Soient $\lambda$ le copoids fondamental minuscule associé à $\beta$ de dernier coefficient diagonal \' egal \` a $1$, et $t = \widetilde{\lambda}(\varpi) \in A_\Lambda^+$. On prend $d_0$ associé pour $V_{N_{\beta} \cap K}$ et $V'$ comme en (\ref{defV'2}). Si on est dans le cas déployé $d=1$, on peut utiliser le corollaire \ref{chgtpdsdual} en raisonnant comme dans le cas $\Delta_P = \{ \beta \}$. On suppose donc à présent $d>1$. Supposons dans un premier temps $\Delta_{\chi_{V^*}} = \Delta$; alors l'hypothèse \ref{V'toV} est satisfaite. La même discussion de cas que pour $\Delta_P = \{\beta\}$ nous permet de voir $\widetilde{\beta}^\vee(\varpi) \notin \widetilde{Z}_L(V^*)$ ou $\chi_L \big( \widetilde{\beta}^\vee(\varpi) \big) \neq 1$. On applique alors encore le corollaire \ref{chgtpdsdual} pour conclure. \\
Il reste donc à examiner le cas $\Delta_{\chi_{V^*}} \neq \Delta$. On a alors $\Delta_{\chi_{V^*}} = \Delta_V^-$, et en particulier $\sigma_2$ ne peut pas être une représentation de Steinberg généralisée (à cause de la proposition \ref{Stsocle}). On a ainsi $L=M$ et $\widetilde{\beta}^\vee(\varpi) \notin \widetilde{Z}_L \supseteq \widetilde{Z}_L(V^*)$. On peut donc appliquer le corollaire \ref{chgtpdsdualm=3} et obtenir \[ \ind_K^G \, V \otimes_{\Hh,\chi} \Fpbar \xrightarrow{\sim} \ind_K^G \, V' \otimes_{\Hh,\chi} \Fpbar .\] Aussi, la représentation $V'$ est $(P')^-$-régulière et la preuve est terminée. 

\section{Utilisation du foncteur parties ordinaires} \label{prgordinaire}

Soit $P$ un parabolique standard de Levi standard $M$ et de radical unipotent $N$. Dans \cite{Eme10}, Emerton définit un foncteur $\Ord_P$ de la catégorie des $G$-représentations lisses admissibles vers celle des $M$-représentations lisses admissibles. En fait, \cite{Eme10} traite uniquement le cas où $F$ est de caractéristique $0$ et il faut se reporter à \cite{Vig12} pour le cas d'égale caractéristique. C'est un adjoint à droite du foncteur $\Ind_{P^-}^G$, et en ce sens le comprendre nous sera bien utile. \\

Soit $V$ une représentation irréductible de $K$. Le point important du résultat suivant est de comprendre comment le groupe abélien $\widetilde{Z}_M(V)$ est engendré par $Z(M) \cap A_\Lambda$ et $\widetilde{Z}_M^+(V)$. Aussi, on note $M^+$ le sous-monoïde de $M$ qui contracte $N \cap K$ par conjugaison (ce qui généralise la définition sur $A$).

\begin{lem} \label{pasdeZM?}
Soient $M$ un Levi standard de $G$ et $V$ une représentation irréductible de $K$. Il existe un élément $z' \in Z(M)^+ \cap A_\Lambda$ tel que l'on a:
\begin{itemize}
\item[(i)] \[ Z(M) = Z(M)^+ (Z(M) \cap A_\Lambda) = Z(M)^+ (z')^{- \N} ;\]
\item[(ii)] \begin{equation} \label{ZMdecomp} \widetilde{Z}_M(V) = \widetilde{Z}_M^+(V) (Z(M) \cap A_\Lambda) = \widetilde{Z}_M^+(V) (z')^{- \N};\end{equation}
\item[(iii)] \[ M = M^+ (Z(M) \cap A_\Lambda) = M^+ (z')^{-\N} .\]
De plus, dans chacune des égalités précédentes, l'intersection des monoïdes des termes du milieu est égale à $Z(M)^+ \cap A_\Lambda$.
\end{itemize}
\end{lem}
\textsf{Preuve :} \\
Le \textit{(i)} vient des égalités
\[ Z(M) = (Z(M) \cap K)(Z(M) \cap A_\Lambda) = (Z(M)^+ \cap K)(Z(M) \cap A_\Lambda) .\]
Pour \textit{(ii)} et \textit{(iii)}, les inclusions $\supseteq$ sont claires. Soit $z' \in Z(M) \cap A_\Lambda$ un élément \og strictement dominant hors de $\Delta_M$ \fg: par exemple, on peut prendre 
\[ z' = \Diag(\varpi_D^{d a_1}, \varpi_D^{d a_2}, \dots, \varpi_D^{d a_m}) \in Z(M) \cap A_\Lambda \] 
avec $a_i = a_{i+1} + 1$ si la racine entre $i$ et $i+1$ est dans $\Delta \smallsetminus \Delta_M$, et $a_i = a_{i+1}$ sinon. L'inclusion $\subseteq$ dans (\ref{ZMdecomp}) vient de ce que, pour tout \' el\' ement $z \in \widetilde{Z}_M(V)$, $z (z')^n$ est dominant pour tout $n \geq 0$ assez grand. \\
Le \textit{(iii)} se prouve de la même manière que le \textit{(ii)}: on multiplie tout élément $m \in M$ par une puissance assez grande de $z'$. \hfill$\Box$\\

Lorsque $\pi$ est une représentation lisse de $G$, on définit une autre action de $\widetilde{Z}_M^+$ sur $\pi^{N \cap K}$ comme suit. Si $t$ est un élément de $M^+$, alors $t(N \cap K)t^{-1}$ est un sous-groupe compact ouvert de $N \cap K$, donc d'indice fini, et $t$ agit de la manière suivante: 
\begin{equation} \label{actioncontract} \begin{array}{ccccc} \pi^{N \cap K} & \to & \pi^{t(N \cap K)t^{-1}} & \to & \pi^{N \cap K} \\ v & \mapsto & tv & \mapsto & \sum\limits_{n \in (N \cap K)/t(N \cap K) t^{-1}} ntv \end{array} .\end{equation}
On note $\Hom_{Z(M)^+}(Z(M),\pi^{N \cap K})$ le $\Fpbar$-espace vectoriel des fonctions $f : Z(M) \to \pi^{N \cap K}$ qui vérifient $f(z^+z) = z^+. f(z)$ pour tous $z^+ \in Z(M)^+$ et $z \in Z(M)$, où l'action de $Z(M)^+$ sur $\pi^{N \cap K}$ est celle que l'on vient de définir.
L'espace $\Ord_P \pi$ des $P$-parties ordinaires de $\pi$ est alors le sous-espace des vecteurs $Z(M)$-finis dans l'espace $\Hom_{Z(M)^+}(Z(M),\pi^{N \cap K})$. \\ Et $\Ord_P \pi$ hérite d'une structure de $M$-représentation lisse: si $f$ est un vecteur $Z(M)$-fini de $\Hom_{Z(M)^+}(Z(M),\pi^{N \cap K})$ et $m \in M$ s'écrit $m = m^+ z$ avec $m^+ \in M^+$ et $z \in Z(M) \cap A_\Lambda$ par le lemme \ref{pasdeZM?}.(iii), on définit $m.f: z' \mapsto m^+ f(zz')$ (qui est bien défini et est $Z(M)$-fini). De plus, par les \textit{(i)} et \textit{(ii)} du lemme \ref{pasdeZM?}, on a le résultat suivant.

\begin{lem}
Soient $V$ une représentation irréductible de $K$ et $\pi$ une représentation lisse de $G$. Soit $P = MN$ un parabolique standard de $G$. On a 
\begin{equation} \label{defOrdbis} \Ord_P \pi \simeq \Hom_{\widetilde{Z}_M^+(V)}(\widetilde{Z}_M(V),\pi^{N \cap K})_{Z(M) \textrm{-fin}} ,\end{equation}
où $\Hom_{\widetilde{Z}_M^+(V)}(\widetilde{Z}_M(V),\pi^{N \cap K})$ désigne l'espace des fonctions $f : \widetilde{Z}_M(V) \to \pi^{N \cap K}$ qui vérifient $f(z^+z) = z^+. f(z)$ pour tous $z^+ \in \widetilde{Z}_M^+(V)$ et $z \in \widetilde{Z}_M(V)$ avec l'action (\ref{actioncontract}) de $\widetilde{Z}_M^+(V)$ sur $\pi^{N \cap K}$. 
\end{lem}
\textsf{Preuve :} \\
Il s'agit de voir que 
\[ \mathrm{Hom}_{Z(M)^+}(Z(M), \pi^{N \cap K}) \xrightarrow{\alpha} \mathrm{Hom}_{Z(M)^+ \cap A_\Lambda}(Z(M) \cap A_\Lambda, \pi^{N \cap K}) \mathop{\leftarrow}\limits^{\beta} \mathrm{Hom}_{\widetilde{Z}^+_M(V)}(\widetilde{Z}_M(V), \pi^{N \cap K}) \]
sont bien définis et sont des isomorphismes: on prendra ensuite le sous-espace des vecteurs $Z(M)$-finis sur chacun de ces trois espaces. Tout d'abord, les flèches $\alpha$ et $\beta$ sont toutes deux des flèches de restriction, qui arrivent bien dans l'espace des fonctions $Z(M)^+ \cap A_\Lambda$-équivariantes par la dernière assertion du lemme \ref{pasdeZM?}. \\ 
Pour définir leur inverse, on se sert des égalités du lemme \ref{pasdeZM?}: voyons le cas de $\beta^{-1}$. Soit $f \in \mathrm{Hom}_{Z(M)^+ \cap A_\Lambda}(Z(M) \cap A_\Lambda, \pi^{N \cap K})$. Tout élément $z$ de $\widetilde{Z}_M(V)$ s'écrit $z = z_+ (z')^{-n}$ avec $z_+ \in \widetilde{Z}_M^+(V)$ et $(z')^{-n} \in Z(M) \cap A_\Lambda$ où $z'$ est l'élément donné par le lemme \ref{pasdeZM?}.(ii) et $n \geq 0$ assez grand. On pose alors $\beta^{-1}(f)(z) = z_+.f((z')^{-n})$ à l'aide de l'action en (\ref{actioncontract}). Il s'agit ensuite de voir que cela définit bien $\beta^{-1}(f)$: si $z = z'_+ (z')^{-n'}$ est une autre écriture comme précédemment, alors disons que l'on a $n' \geq n$. On a alors $(z')^{n'-n} \in Z(M)^+ \cap A_\Lambda$ et donc 
\[ z'_+ . f((z')^{-n'}) = z_+ (z')^{n'-n}.f((z')^{-n'}) = z_+ f((z')^{-n}) .\] 
Aussi, $\beta^{-1}(f)$ est bien $\widetilde{Z}_M^+(V)$-équivariante: si $z_1$ est dans $\widetilde{Z}_M^+(V)$ alors, $z = (z_1 z_+) (z')^{-n}$ est encore une écriture de la forme voulue et on a
\[ \beta^{-1}(f)(z_1z) = (z_1 z_+).f((z')^{-n}) = z_1. f(z) .\]
De plus, si $z $ est dans $Z(M) \cap A_\Lambda$, une écriture $z = z_+ (z')^{-n}$ comme précédemment implique $z_+ \in Z(M)^+ \cap A_\Lambda$; pour $f \in \mathrm{Hom}_{Z(M)^+ \cap A_\Lambda}(Z(M) \cap A_\Lambda, \pi^{N \cap K})$, on en déduit
\[ \beta \circ \beta^{-1}(f)(z) = z_+.f((z')^{-n}) = f(z) .\]
Enfin, si $z = z_+ (z')^{-n} \in \widetilde{Z}_M(V)$ est une écriture avec $z_+ \in \widetilde{Z}_M^+(V)$ et $(z')^{-n} \in Z(M) \cap A_\Lambda$ et que $g$ est un élément de $\mathrm{Hom}_{\widetilde{Z}^+_M(V)}(\widetilde{Z}_M(V), \pi^{N \cap K})$, on a 
\[ \beta^{-1} \circ \beta (f)(z) = z_+.g((z')^{-n}) = g(z) .\]
Le cas de $\alpha^{-1}$ se traite de manière similaire. \hfill$\Box$ \\

Supposons que $V$ est un caractère de $M \cap K$, que l'on notera $\psi$. Si $h$ est dans $\widetilde{Z}_M(V)$, on a alors $\psi(h \cdot h^{-1}) = \psi$. Ceci permet, pour tout caractère $\chi_M : \widetilde{Z}_M(V) \to \Fpbar^\times$, de faire agir $\widetilde{Z}_M(V)$ par $\chi_M$ sur $V$, ce qui fait de $V$ un caractère de $(M \cap K) \widetilde{Z}_M(V)$.

\begin{lem} \label{isomZM?}
Soit $\pi$ une représentation lisse de $G$. Soient $V$ un caractère de $(M \cap K)$ et $\chi_M : \widetilde{Z}_M(V) \to \Fpbar^\times$ un caractère, conférant à $V$ une structure de caractère de $(M \cap K) \widetilde{Z}_M(V)$. Alors on a un isomorphisme \[ \Hom_{(M \cap K) \widetilde{Z}_M(V)}  \big( V, \Ord_P \pi \big) \simeq \Hom_{(M \cap K) \widetilde{Z}_M^+(V)} \big( V, \pi^{N \cap K} \big) ,\] où l'action de $\widetilde{Z}_M^+(V)$ sur $\pi^{N \cap K}$ est donnée par (\ref{actioncontract}).
\end{lem}
\textsf{Preuve :} \\
Le membre de gauche est par définition de $\Ord_P \pi$ (voir \cite{Eme10}, Definition 3.1.9 et \cite{Vig12}, Remark 1.9.(ii) car $V$ est de dimension finie) \[ H_{V,\pi} := \Hom_{(M \cap K) \widetilde{Z}_M(V)} \big( V , \Hom_{Z(M)^+} (Z(M), \pi^{N \cap K}) \big) .\] 
Grâce à (\ref{defOrdbis}), on a alors 
\[ H_{V,\pi} \simeq \Hom_{(M \cap K) \widetilde{Z}_M(V)} \big( V, \Hom_{\widetilde{Z}_M^+(V)}(\widetilde{Z}_M(V), \pi^{N \cap K}) \big) .\]
Puis la flèche $f \mapsto (v \mapsto f(v)(1))$ induit l'isomorphisme \[ H_{V,\pi} \simeq \Hom_{(M \cap K) \widetilde{Z}_M^+(V)} \big( V, \pi^{N \cap K} \big) \] et le résultat. \hfill$\Box$

\begin{lem} \label{VhV}
Soient $V$ une repr\' esentation irr\' eductible de $M \cap K$ et $h$ un \' el\' ement de $\widetilde{Z}_M(V)$. Alors l'espace $h V = \{ hv \ | \ v \in V \}$ est une $(M \cap K)$-repr\' esentation \` a travers l'action $m. hv = h (h^{-1}mh) v$, qui est irr\' eductible et isomorphe \` a $V$.
\end{lem}
\textsf{Preuve :} \\
Parce que l'on a $h^{-1}(M \cap K) h = M \cap K$, $hV$ est bien une $(M \cap K)$-repr\' esentation irr\' eductible. Notons $\Par_{\overline{M} \cap \overline{B}}(V) = (\chi,\Delta_V)$ suivant la proposition \ref{irredGbar}. Parce que l'on a $h \in \widetilde{Z}_M(V)$, que l'on a \[ h^{-1}(M \cap U \cap K) h = M \cap U \cap K \] et que $\Delta_V$ est inclus dans $\Delta_M$, les param\` etres de $hV$ sont 
\[ \Par_{\overline{M} \cap \overline{B}}(hV) = (\chi,\Delta_V) = \Par_{\overline{M} \cap \overline{B}}(V) ,\] 
et $hV$ est alors bien isomorphe \` a $V$. \hfill$\Box$ \\

On note $I_P$ le sous-groupe parahorique de $G$ correspondant à $P$. Soit $V$ une représentation irréductible de $I_P$; comme $V$ est triviale sur le pro-$p$-radical $I_P(1)$ de $I_P$, $V$ se factorise en fait par $M \cap K$. Et inversement, on fera de toute représentation de $M \cap K$ une représentation de $I_P$ par inflation. On va ensuite noter $\Hh^G_{I_P}$ la sous-algèbre de $\Hh_{\Fpbar}(G,I_P,V)$ engendrée par les éléments de support contenu dans $I_P h I_P$ pour $h \in \widetilde{Z}_M^-(V)$.

\begin{lem}
Soient $V$ une $(M \cap K)$-représentation irréductible et $h \in \widetilde{Z}_M^-(V)$. L'algèbre $\Hh_{I_P}^G$ possède un opérateur $E_h$ de support $I_P h I_P$. \\ De plus, si $V$ est un caractère, on peut choisir un tel $E_h$ dont la valeur en $h$ est l'identité de $V$.
\end{lem}
\textsf{Preuve :} \\
On cherche un opérateur $E_h$ de support $I_P h I_P$. Il s'agit de vérifier que l'on a bien $p_1 E_h(h) = E_h(h) p_2$ pour tous $p_1, p_2 \in I_P$ vérifiant $p_1 h = h p_2$. Ecrivons $p_1 = p_1^- p_1^\circ p_1^+$ et $p_2 = p_2^- p_2^\circ p_2^+$ selon la décomposition d'Iwahori (voir (\ref{decompIwahori})) \[ I_P = (I_P \cap N^-)(I_P \cap M)(I_P \cap N) .\] On obtient par suite \[ p_1^- (p_1^\circ h)(h^{-1} p_1^+ h) = (h p_2^- h^{-1})(h p_2^\circ) p_2^+ ,\] et comme (\ref{decompIwahori}) est une décomposition, on a $p_1^\circ h = h p_2^\circ$. Or on a $p_1 E_h(h) = p_1^\circ E_h(h)$ et $E_h(h) p_2 = E_h(h) p_2^\circ$, vus dans $\End_{\Fpbar}(V)$. On a alors le diagramme commutatif de $(M \cap K)$-représentations \[ \xymatrix{ V \ar[r]^{E_h(h)} \ar[d]_{p_1^\circ} & V \ar[d]^{h^{-1} p_1^\circ h} \\ V \ar[r]^{E_h(h)} & V } .\] Réécrit autrement, $h E_h(h)$ est un élément de $\Hom_{M \cap K}(V,hV)$. 
De plus, $hV$ est aussi une repr\' esentation irr\' eductible de $M \cap K$, isomorphe \` a $V$ par le lemme \ref{VhV}. De ce fait, le lemme de Schur nous assure que $\Hom_{M \cap K}(V,hV)$ est de dimension $1$. Donc $E_h(h)$ est un isomorphisme linéaire $V \to V$ et on a bien \[ p_1 E_h(h) = E_h(h) p_2 \in \End_{\Fpbar}(V) \] comme voulu. \\ 
Lorsque $V$ est un caractère, on a $E_h(h) = k \, \Id_V$ avec $k \in \Fpbar^\times$ et on peut choisir $k=1$. \hfill$\Box$\\

Pour toute $(M \cap K)$-représentation $V$, on a une application de restriction $r_{\Hh}: \Hh_{\Fpbar}(G,I_P,V) \to \Hh_{\Fpbar}(M, M \cap K,V)$, qui fait que tout caractère de $\Fpbar$-algèbres $\Hh_{\Fpbar}(M,M \cap K, V) \to \Fpbar$ induit\footnote{on n'affirme pas ici que $r_{\Hh}$ est un morphisme d'alg\` ebres, mais sa restriction \` a $\Hh_{I_P}^G$ l'est. Pour voir cela, prenons $a,b \in \widetilde{Z}_M^-(V)$ et calculons $\id_{I_P a I_P} * \id_{I_P b I_P}$ selon (\ref{GLmdefconvol}). Il s'agit de d\' eterminer $I_P a I_P b I_P$: on va voir que cet ensemble est inclus dans $I_P \widetilde{Z}_M(V) I_P$ et est \' egal \` a $I_P ab I_P$, donnant alors la compatibilit\' e de $r_{\Hh}|_{\Hh_{I_P}^G}$ avec $*$ comme voulu. Par (\ref{decompIwahori}), on a la d\' ecomposition d'Iwahori
\[ I_P = (I_P \cap N^-) (I_P \cap M) (I_P \cap N) .\]
On a alors l'inclusion 
\[ I_P a I_P b I_P = I_P a (I_P \cap N^-) (I_P \cap M) (I_P \cap N) b I_P \subseteq I_P (I_P \cap N^-)a (I_P \cap M) b (I_P \cap N) I_P = I_P ab I_P .\]
Parce que l'inclusion inverse est \' evidente, on a l'\' egalit\' e.
} un caractère $\Hh_{I_P}^G \to \Fpbar$, que l'on notera encore $\chi$. On remarque que, pour tout $h \in \widetilde{Z}_M^-(V)$, cette application de restriction $r_{\Hh}$ envoie $E_h$ sur $T_h^M$.

\begin{prop} \label{Prop9.14}
Soient $V$ un caractère de $M \cap K$ et $\chi_M : \widetilde{Z}_M(V) \to \Fpbar^\times$ un caractère, faisant de $V$ un caractère de $(M \cap K) \widetilde{Z}_M(V)$. Il existe un unique caract\` ere d'alg\` ebres \[ \chi : \Hh_{\Fpbar}(M,M \cap K,V) \to \Fpbar \] vérifiant $\chi(T^M_h) = \chi_M(h)^{-1}$ pour tout $h \in \widetilde{Z}_M^-(V)$. De plus, pour toute représentation lisse $\pi$ de $G$, on a un isomorphisme \[ \Hom_{(M \cap K) \widetilde{Z}_M(V)}  \big( V, \Ord_P \pi \big) \simeq \Hom_G \big( \ind_{I_P}^G \, V \otimes_{\Hh_{I_P}^G , \chi} \Fpbar, \pi \big) .\]
\end{prop}
\textsf{Preuve :} \\ 
Par le lemme \ref{isomZM?}, on cherche à obtenir un isomorphisme \[ \Hom_{(M \cap K) \widetilde{Z}_M^+(V)} \big( V, \pi^{N \cap K} \big) \simeq \Hom_G \big( \ind_{I_P}^G \, V \otimes_{\Hh_{I_P}^G , \chi} \Fpbar, \pi \big) .\] 
Commençons par montrer que toute application linéaire $(M \cap K) \widetilde{Z}_M^+(V)$-équivariante $f : V \to \pi^{N \cap K}$ se factorise à travers l'espace d'invariants $\pi^{I_P(1)}$ du pro-$p$-radical $I_P(1)$ de $I_P$. Parce que $\pi$ est lisse, il existe un $h \in \widetilde{Z}_M^{--}(V)$ tel que $f(V)$ est fixé point par point par $h(I_P \cap N^-)h^{-1}$. Alors $h^{-1}.f(V)$ est fixé par le groupe engendré par $I_P \cap N^-$, $M \cap K(1)$ et $h^{-1}(I_P \cap N) h$, c'est-à-dire par le groupe $I_P(1) \cap h^{-1} I_P(1) h$. Comme $f$ est $\widetilde{Z}_M^+(V)$-équivariante, on réécrit l'action (voir (\ref{actioncontract})) de $h^{-1}$ sur $f(v)$, pour tout $v \in V$: 
\begin{equation} \label{Ordaction} \chi_M(h^{-1}) f(v) = \sum_{n \in (N \cap K)/h^{-1}(N \cap K)h} n h^{-1} f(v) .\end{equation} 
Comme l'application \[ (N \cap K) / h^{-1} (N \cap K) h \to I_P(1) / (I_P(1) \cap h^{-1} I_P(1) h) \] est une bijection, (\ref{Ordaction}) se réécrit \[ \chi_M(h^{-1}) f(v) = \sum_{n \in I_P(1)/(I_P(1) \cap h^{-1} I_P(1) h)} n h^{-1} f(v) .\] Cela implique l'égalité 
\begin{equation} \label{chiMf} \chi_M(h^{-1}) f(v) = \sum_{n \in I_P(1)/(I_P(1) \cap h^{-1} I_P(1) h)} \chi_M(h^{-1}) n f(v) \end{equation} puisque $f$ est $\widetilde{Z}_M^+(V)$-équivariante et que $h^{-1}$ agit à travers $\chi_M$ sur $V$. De ce fait, $f(v)$ est $I_P(1)$-invariante et $f : V \to \pi^{I_P(1)} \subseteq \pi$ est donc $I_P$-équivariante. \\ 
On note que l'existence et l'unicit\' e du caractère d'algèbres \[ \chi : \Hh_{\Fpbar}(M,M \cap K,V) \to \Fpbar \] de paramètres de Hecke-Satake $(M,\chi_M)$ et vérifiant $\chi(T^M_h) = \chi_M(h)^{-1}$ pour tout $h \in \widetilde{Z}_M^-(V)$ (directement grâce à la formule (\ref{Prop3.6})) suivent de la proposition \ref{Prop4.1}. Par réciprocité de Frobenius, on voit maintenant $f \in \Hom_{I_P}(V,\pi)$ comme un élément de $\Hom_G(\ind_{I_P}^G \, V, \pi)$ et on veut vérifier que $f$ est un vecteur propre pour $\Hh^G_{I_P}$ de caractère propre $\chi$. Pour $h \in \widetilde{Z}_M^-(V)$ et $v \in V$, on a: 
\begin{equation} \label{f*Eh} (f * E_h)(v) = \sum_{i \in I_P/(I_P \cap h^{-1} I_P h)} i h^{-1} f(E_h(h i^{-1}) v) .\end{equation} 
Comme l'application \[ (N \cap K) / h^{-1} (N \cap K) h \to I_P / (I_P \cap h^{-1} I_P h) \] est bijective, (\ref{f*Eh}) vaut, pour la même raison qu'en (\ref{chiMf}), \[ \sum_{i \in (N \cap K)/h^{-1}(N \cap K)h} ih^{-1} f(v) = \chi_M(h^{-1}) f(v) .\] On vient donc de construire une application \[ \Hom_{(M \cap K) \widetilde{Z}_M^+(V)} \big( V, \pi^{N \cap K} \big) \to \Hom_G \big( \ind_{I_P}^G \, V \otimes_\chi \Fpbar , \pi \big) .\] Dans l'autre sens, un élément $\phi_0 \in \Hom_G \big( \ind_{I_P}^G  \, V \otimes_\chi \Fpbar , \pi \big)$ donne par réciprocité de Frobenius un élément $\phi \in \Hom_{I_P}(V,\pi)$. Par la condition sur $\chi$ et $\chi_M$, $\phi$ est de fait un élément de $\Hom_{I_P \widetilde{Z}_M^+(V)}(V,\pi)$; et comme $\phi$ est d'image $I_P(1)$-invariante (car $I_M(1)$ agit trivialement sur $V$), en relâchant les conditions on a $\phi \in \Hom_{(M \cap K) \widetilde{Z}_M^+(V)} \big( V, \pi^{N \cap K} \big)$. Cela nous fournit une application \[ \Hom_G \big( \ind_{I_P}^G \, V \otimes_\chi \Fpbar , \pi \big) \to \Hom_{(M \cap K) \widetilde{Z}_M^+(V)} \big( V, \pi^{N \cap K} \big) \] qui est inverse de la précédente. \hfill$\Box$ \\

On en profite pour énoncer tout de suite un lemme qui utilise aussi (dans sa preuve) ce foncteur des parties ordinaires; son utilité ne se ressentira qu'au cours du paragraphe \ref{classifadmiss}.

\begin{lem} \label{findadmiss}
Soit $P = MN$ un parabolique standard. Soient $\pi$ une représentation irréductible admissible de $G$ et $\tau$ une représentation lisse de $M$ à caractère central et induisant une surjection $\Ind_P^G \, \tau \twoheadrightarrow \pi$ de $G$-représentations. Alors il existe une représentation irréductible admissible $\sigma$ de $M$ induisant une surjection $\Ind_P^G \, \sigma \twoheadrightarrow \pi$ de $G$-représentations.
\end{lem}
\textsf{Preuve :} 
Voir \cite{HenVig11b}, Proposition 7.9 et Lemma 7.10. \hfill$\Box$

\section{Vecteurs $K$-invariants et algèbres de Hecke} \label{KinvHecke}

Soit $\pi$ une représentation lisse de $G$. Lorsque l'espace des $K$-invariants $\pi^K$ est non réduit à $0$, c'est un objet intéressant à considérer. En effet, d'après la réciprocité de Frobenius, \[ \pi^K \simeq \Hom_K(\id,\pi) \simeq \Hom_G(\ind_K^G \, \id, \pi) \] est naturellement muni d'une structure de module à droite sur $\Hh_{\Fpbar}(G,K)$.
Explicitement, l'opérateur $T_g$ de $\Hh_{\Fpbar}(G,K)$ correspondant à la fonction caractéristique de $K g K$ pour $g \in G$ agit sur $\pi^K$ par 
\begin{equation} \label{Kagit} v \mapsto v T_g = \sum_{\gamma \in K \backslash K g K} \gamma^{-1} v = \sum_{u \in (K \cap g^{-1} K g) \backslash K} u^{-1} g^{-1} v .\end{equation}

On va aussi s'intéresser à l'algèbre de Hecke-Iwahori $\Hh_{\Fpbar}(G,I)$ où $I$ désigne le parahorique correspondant au parabolique déployé minimal $B$. Donnons tout d'abord quelques notations pour des éléments particuliers de $\Hh_{\Fpbar}(G,I)$. \\
Soit $t_i = \Diag(\varpi_D, \dots, \varpi_D, 1, \dots, 1)$ l'élément de $A$ avec $i$ copies de $\varpi_D$ (où $0 \leq i \leq m$). Soit aussi, pour $1 \leq i \leq m-1$, $s_i$ l'élément du Levi standard par blocs $\GL(1)^{i-1} \times \GL(2) \times \GL(1)^{m-i-1}$ où les blocs de dimension $1$ sont $(1)$ et le bloc de dimension $2$ est $\left( \begin{array}{cc} 0 & 1 \\ 1 & 0 \end{array} \right)$. Ce sont des relèvements des générateurs $S_\Delta$ de $W$ dans $G$. Soit enfin $s_\varpi$ l'élément $\left( \begin{array}{cccc} 0 & 1 & & \\ & \ddots & \ddots & \\ & & \ddots & 1 \\ \varpi_D & & & 0 \end{array} \right)$. La classe de $s_\varpi$ dans $W_e$ et les éléments de $S_\Delta$ constituent un système de générateurs minimal du groupe de Weyl étendu $W_e$. \\

Comme au paragraphe \ref{LKglm} avec l'algèbre de Hecke-Satake, on considère ici l'algèbre de Hecke-Iwahori d'une représentation triviale de $A \cap K$ et on va le comparer à celle de $\GE = \GL(m,E)$. On va noter $\IE$ l'Iwahori standard de $\GE$ correspondant au Borel $\GE \cap B$, et $s^{(E)}_\varpi$ l'élément $\left( \begin{array}{cccc} 0 & 1 & & \\ & \ddots & \ddots & \\ & & \ddots & 1 \\ \varpi & & & 0 \end{array} \right)$ qui le normalise. On note alors $U_i, S_i, \Pi$ les éléments de $\Hh_{\Fpbar}(G,I)$ de support respectif $I t_i^{-1} I$, $I s_i I$ et $I s_\varpi^{-1} I$. De manière analogue, $S_i^{(E)}$ et $\Pi^{(E)}$ désigneront les opérateurs de $\Hh_{\Fpbar}(\GE,\IE)$ de support respectif $\IE s_i \IE$ et $\IE (s_\varpi^{(E)})^{-1} \IE$.

\begin{prop} \label{IwahoripourE}
Les algèbres de Hecke $\Hh_{\Fpbar}(\GE,\IE)$ et $\Hh_{\Fpbar}(G,I)$ sont isomorphes via $S_i^{(E)} \mapsto S_i$ et $\Pi^{(E)} \mapsto \Pi$.
\end{prop}
\textsf{Preuve :} \\
La décomposition de Bruhat nous dit que l'algèbre $\Hh_{\Fpbar}(G,I)$ est engendrée vectoriellement par les éléments de support $I g I$ et valant $1$ en $g$, où $g$ parcourt l'ensemble de représentants des éléments du groupe de Weyl étendu $W_e$ engendré par les $s_i$ et $s_\varpi$. Le fait analogue est bien entendu aussi valable pour $\Hh_{\Fpbar}(\GE,\IE)$, qui a pour $\Fpbar$-base les fonctions caractéristiques de $I_{(E)} g I_{(E)}$ pour $g$ parcourant un ensemble de représentants de $W_e^{(E)}$. En vertu du lemme \ref{Weylextend}, on a donc un isomorphisme de $\Fpbar$-espaces vectoriels \[ f: \Hh_{\Fpbar}(G_{(E)},I_{(E)}) \xrightarrow{\sim} \Hh_{\Fpbar}(G,I) ,\] qui envoie $\Pi^{(E)}$ sur $\Pi$ et $S_i^{(E)}$ sur $S_i$ pour tout $1 \leq i \leq m-1$. On veut établir que $f$ est un morphisme de $\Fpbar$-algèbres. Comme le paragraphe 5.4 de \cite{BusKut93} nous dit que\footnote{de fait, le résultat cité est pour $\Hh_\Z(\GE,\IE)$ mais ce qui est fait se réduit bien modulo $p$} l'algèbre $\Hh_{\Fpbar}(\GE,\IE)$ est de type fini, engendrée par les $S_i^{(E)}$ et $\Pi^{(E)}$, il reste à voir que les opérateurs $\Pi$ et $S_i$ vérifient exactement les mêmes relations que les $\Pi^{(E)}$ et $S_i^{(E)}$. Mais ceci r\' esulte du lemme \ref{Weylextend} et de \cite{Bou68VII4}, paragraphe 2, exercices 23 et 24. \hfill$\Box$

\begin{cor} \label{braid}
On a les relations suivantes: 
\begin{itemize} 
\item[a)] $S_i^2 = - S_i$ pour tout $i$;
\item[b)] $S_i S_j = S_j S_i$ lorsque $|i-j|>1$;
\item[c)] $S_i S_{i+1} S_i = S_{i+1} S_i S_{i+1}$ pour $1 \leq i \leq m-2$;
\item[d)] $\Pi S_i = S_{i+1} \Pi$ pour $1 \leq i \leq m-2$.
\end{itemize}
\end{cor} 
\textsf{Preuve :} \\
Les mêmes relations pour $\GE$ sont l'objet du paragraphe 5.4 de \cite{BusKut93}. On remarquera que n'y est pas écrite la toute dernière relation mais celle-ci n'est pas difficile puisque $\Pi^{(E)}$ est en fait porté par la classe simple $\IE (s_\varpi^{(E)})^{-1} = (s_\varpi^{(E)})^{-1} \IE$ ($s_\varpi^{(E)}$ normalise $\IE$). On utilise ensuite l'isomorphisme d'algèbres de la proposition \ref{IwahoripourE}. \hfill$\Box$

\begin{cor} \label{combinSi}
$\phantom{,}$
\begin{itemize}
\item[(i)] Soient $i \leq k \leq l \leq j$. On a\footnote{avec les conventions $S_{k-1} S_{k-2} \cdots S_k = S_k$ et $S_k S_{k-1} \cdots S_{k+1} = S_{k+1}$} \[ S_{l-1} S_{l-2} \cdots S_k S_j S_{j-1} \cdots S_i = S_j S_{j-1} \cdots S_i S_l S_{l-1} \cdots S_{k+1} .\]
\item[(ii)] Pour tout $1 \leq i \leq m-1$, on a \[ U_i = (\Pi S_{m-1} S_{m-2} \cdots S_i)^i .\]
\item[(iii)] Soient $\pi$ une $G$-représentation et $v \in \pi^I$. Si $v S_1$ est nul, alors on a $v \Pi^2 S_{m-1} = 0$.
\end{itemize}
\end{cor}
\textsf{Preuve :} \\
Ce sont les \cite{Her11b}, Sublemma 9.3, 9.4 et 9.5, qui sont purement combinatoires en utilisant le corollaire \ref{braid}. On notera que dans \cite{Her11b}, les actions sont à gauche alors que l'on a pris la convention de faire agir à droite. \hfill$\Box$\\

On va maintenant s'intéresser à l'algèbre de Hecke-Satake $\Hh_{\Fpbar}(G,K)$, dont on va essayer de relier l'action à celle de l'algèbre de Hecke-Iwahori. Soit, pour $1 \leq i \leq m$, $T_i$ l'opérateur de $\Hh_{\Fpbar}(G,K)$ de support $K t_i^{-1} K$.

\begin{lem}
Soient $\pi$ une $G$-représentation et $v \in \pi^K$. Alors on a\footnote{avec la convention $\Pi S_{m-1} S_{m-2} \cdots S_m = \Pi$} \begin{equation} \label{IwahoriSatake} v T_1 = \sum_{i=1}^m v \Pi S_{m-1} S_{m-2} \cdots S_i .\end{equation}
\end{lem}
\textsf{Preuve :} \\
Notons, pour $1 \leq i \leq m$, $d_i = \Diag(1, \dots, 1, \varpi_D, 1, \dots, 1)$ où $\varpi_D$ est à la $i$-ème place, et $u_i(a_{i+1}, a_{i+2}, \dots, a_m)$ l'élément de $U$ dont tous les coefficients surdiagonaux sont nuls sauf sur la i-ème ligne pour laquelle ils sont égaux à $(a_{i+1}, \dots, a_m)$. Lorsque les $a_i$ sont pris dans $\Rr_{(1)}$ (voir (\ref{defR})), on a d'abord $u_i(a_{i+1}, \dots, a_m) d_i \in K t_1 K$, et $u_i(a_{i+1}, \dots, a_m) d_i K = u_j(b_{j+1}, \dots, b_m) d_j K$ si et seulement si $i=j$ et $a_{i+1} = b_{i+1}, \dots, a_m=b_m$. Enfin, $K t_1 K/K$ est de cardinal $(q^{d m}-1)(q^d -1)^{-1}$, de sorte que les $u_i(a_{i+1}, \dots, a_m) d_i$ pour $1 \leq i \leq m$ et $a_j \in \Rr_{(1)}$ forment un système de représentants de $K t_1 K/K$. La formule d'action (\ref{Kagit}) se réécrit alors, pour $g= t_1^{-1}$: \[ v T_1 = \sum_{i=1}^m \sum_{a_{i+1}, \dots, a_m} u_i(a_{i+1}, \dots, a_m) d_i v ,\] où les $a_j$ parcourent l'ensemble $\Rr_{(1)}$ de représentants de $k_D$ dans $\Oo_D$. Pour tout $1 \leq i \leq m$, on note $\sigma_i$ l'élément de la diagonale par blocs $\GL(i) \times \GL(m-i)$ de premier bloc $\left( \begin{array}{cccc} 0 & 1 & & \\ & \ddots & \ddots & \\ & & \ddots & 1 \\ 1 & & & 0 \end{array} \right)$ et de second bloc l'identité. Parce que $v$ est invariant par $K$, on a encore 
\begin{equation} \label{vT1expr} v T_1 = \sum_{i=1}^m \sum_{a_{i+1},\dots,a_m} u_i(a_{i+1}, \dots, a_m) (d_i \sigma_i) v .\end{equation} 
Comme $d_i \sigma = s_i s_{i+1} \cdots s_{m-1} s_\varpi$ est une écriture réduite, on a l'égalité \[I d_i \sigma_i I = I s_i I s_{i+1} I \cdots I s_{m-1} I s_\varpi I \] (voir \cite{Ric69}). Parce que les $u_i(a_{i+1}, \dots, a_m) d_i \sigma_i$ sont des représentants de $I d_i \sigma_i I/I$, on en tire que les termes \[ \sum\limits_{a_{i+1}, \dots, a_m} u_i(a_{i+1}, \dots, a_m) d_i \sigma_i v \] dans (\ref{vT1expr}) sont respectivement les $v \Pi S_{m-1} S_{m-2} \cdots S_i$. \hfill$\Box$

\section{Vers une classification des représentations lisses irréductibles admissibles} \label{classifadmiss}

Dans tout ce paragraphe, on supposera la formule de comparaison de tranformées de Satake (conjecture \ref{conjSatake}) satisfaite. \\

On conserve les notations du paragraphe \ref{KinvHecke}. Comme auparavant, on abrégera $\Hh_{\Fpbar}(G,K,V)$ en $\Hh$ étant donnée une $K$-représentation irréductible $V$. Par exemple, dans l'énoncé suivant, $\Hh$ désigne $\Hh_{\Fpbar}(G,K)$ et $\chi$ un caractère $\Hh \to \Fpbar$ de $\Fpbar$-algèbres.

\begin{prop} \label{alternateOrd}
Soit $\pi$ une $G$-représentation lisse telle qu'il existe une surjection \[ \ind_K^G \, \id \otimes_{\Hh,\chi} \Fpbar \twoheadrightarrow \pi \] où $\chi$ est le caractère ayant $(A, \id : \widetilde{Z}_A \to \Fpbar^\times)$ pour paramètre de Hecke-Satake. Alors on a l'une des possibilités suivantes:
\begin{itemize}
\item[(i)] ou bien $\pi$ contient la représentation triviale $\id$;
\item[(ii)] ou bien il existe un parabolique standard propre $P$ avec $\Ord_P \pi \neq 0$.
\end{itemize}
\end{prop}
\textsf{Preuve :} \\
Par hypothèse, on a $\pi^K \neq 0$. Choisissons donc un vecteur $v \in \pi^K$ propre pour tous les opérateurs de Hecke (de caractère propre paramétrisé par $(A,1)$): par la proposition \ref{Prop4.1}, on a alors \[ v T_1 = \chi(T_1) v = v .\] Supposons $\Ord_P \pi = 0$ pour tout parabolique standard propre $P$. On va commencer par montrer que l'action des $U_i$ est nilpotente pour tout $1 \leq i \leq m-1$. \\
Soient $P$ le parabolique standard de $G$ de Levi standard $M = \GL(i) \times \GL(m-i)$ et $I_P$ le parahorique correspondant. Alors $\Hh_{I_P}^G$ (voir définition avant la proposition \ref{Prop9.14}) désigne la sous-algèbre de $\Hh_{\Fpbar}(G,I_P)$ engendrée par $E_i$ et $E_m$, les opérateurs de support respectif $I_P t_i^{-1} I_P$ et $I_P t_m^{-1} I_P = I_P \varpi_D^{-1}$. Si $\pi^{I_P}$ possédait un vecteur propre associé à une valeur propre non nulle de $E_i$ et de $E_m$, la proposition \ref{Prop9.14} appliquée à $V = \id$ et $\chi_M : \widetilde{Z}_M \to \Fpbar^\times$ un caractère encodant ces valeurs propres nous fournirait un caractère $\chi : \Hh_{\Fpbar}(M, M \cap K) \to \Fpbar$ et un élément non nul de $\Hom_{(M \cap K) \widetilde{Z}_M}\big( \id, \Ord_P \pi \big)$. En particulier, l'espace $\Ord_P \pi$ ne serait pas réduit à $0$, ce qui est absurde. De ce fait, l'action de tous les $E_i$ (pour $1 \leq i \leq m-1$) est nilpotente sur $\pi^{I_P}$. Comme la flèche naturelle $I \backslash I t_i^{-1} I \to I_P \backslash I_P t_i^{-1} I_P$ est une bijection (des deux côtés, un ensemble de représentants est donné par $\overline{N}$), les actions de $E_i$ et de $U_i$ coïncident sur $\pi^{I_P} \subseteq \pi^I$. En particulier, les $U_i$ sont nilpotents sur $v \in \pi^K \subseteq \pi^{I_P}$ comme annoncé plus haut. \\
On veut maintenant montrer que l'on a $v U_i = 0$ pour tout $i$. Commençons par examiner $v U_1$. Par le corollaire \ref{combinSi}.(ii), on a $v U_1 = v \Pi S_{m-1} S_{m-2} \cdots S_1$. En nous servant de (\ref{IwahoriSatake}), on écrit ensuite, puisque l'on a $v T_1 = v$, \[ v (\Pi S_{m-1} S_{m-2} \cdots S_1)^2 = \Big( v - \sum_{i=2}^m v \Pi S_{m-1} S_{m-2} \cdots S_i \Big) \Pi S_{m-1} S_{m-2} \cdots S_1 .\] Grâce au corollaire \ref{braid}, on peut ramener $\Pi$ tout à gauche et obtenir 
\[ v (\Pi S_{m-1} S_{m-2} \cdots S_1)^2 = v \Pi S_{m-1} S_{m-2} \cdots S_1 - \sum_{i=2}^m v \Pi^2 S_{m-2} S_{m-3} \cdots S_{i-1} S_{m-1} S_{m-2} \cdots S_1 .\] Par le corollaire \ref{combinSi}.(i), cela se réécrit encore 
\begin{equation} \label{actioncarree} v \Pi S_{m-1} S_{m-2} \cdots S_1 - \sum_{i=2}^m v \Pi^2 S_{m-1} S_{m-2} \cdots S_1 S_{m-1} S_{m-2} \cdots S_i .\end{equation} 
Parce que les $K s_1 u_1([a],0, \dots, 0)$, pour $a \in k_D$, sont les classes de $K \backslash K s_1 K$, par (\ref{Kagit}), on a \[ v S_1 = \sum_{a \in k_D} u_1([a], 0, \dots, 0) s_1 v = 0 \] en caractéristique $p$. Par le corollaire \ref{combinSi}.(iii), les termes de droite de (\ref{actioncarree}) sont donc tous nuls. Il en résulte $v U_1^2 = v U_1$, et comme $U_1$ est nilpotent sur $v$, on a $v U_1 = 0$. \\
Pour $v U_2$ maintenant, on se sert de nouveau du corollaire \ref{combinSi}.(ii) pour écrire \[ v U_2 = v(\Pi S_{m-1} S_{m-2} \cdots S_2)^2 .\] Puis, ensuite, on sait que l'on a un terme de moins dans (\ref{IwahoriSatake}): \[ v = vT_1 = \sum_{i=2}^m v \Pi S_{m-1} S_{m-2} \cdots S_i .\] On répète le calcul précédent en utilisant cette remarque et on obtient de même \[ v(\Pi S_{m-1} S_{m-2} \cdots S_2)^2 = v \Pi S_{m-1} S_{m-2} \cdots S_2 .\] Parce que $U_2$ est nilpotent sur $v$, on a aussi $v U_2 = 0$. On répète l'opération et il s'ensuit $v U_3=0, \dots, v U_{m-1} = 0$. \\ 
Cela nous permet de conclure que $v$ est un élément de $\pi^G$. En effet, l'égalité (\ref{IwahoriSatake}) est maintenant réduite à $v = v T_1 = v \Pi$. Parce que $s_\varpi$ normalise $I$, l'action de $\Pi$ sur $\pi^I$ est celle de $s_\varpi$. Et comme $G$ est engendré par $K$ et $s_\varpi$ (cela résulte de (\ref{decompCartan}) et de ce que $K$ contient un ensemble de relèvements de $W$), le résultat suit. \hfill$\Box$ \\

On est maintenant en mesure de faire un premier pas vers une classification des représentations irréductibles admissibles de $G$. Après cela, il restera à déterminer \og l'unicité \fg, c'est-à-dire à étudier les entrelacements entre de telles représentations.

\begin{theo} \label{classifadmis}
Supposons l'hypothèse \ref{Dpremier} satisfaite. Soit $\pi$ une représentation irréductible admissible de $G$. Alors il existe un parabolique standard $P$ de Levi standard $\prod\limits_{i=1}^m \GL(m_i,D)$ avec $\sum_i m_i = m$ et des représentations irréductibles admissibles $\sigma_i$ de $\GL(m_i,D)$ tels que l'on ait $\pi \simeq \Ind_P^G (\sigma_1 \otimes \cdots \otimes \sigma_r)$. De plus, pour tout $i$, on est dans l'un des cas suivants: 
\begin{itemize}
\item[(a)] $\sigma_i$ est supersingulière avec $m_i > 1$;
\item[(b)] $\sigma_i$ est une représentation de $D^\times$ (i.e. $m_i = 1$) avec $\dim \, \sigma_i > 1$;
\item[(c)] $\sigma_i \simeq \St_{Q_i} \rho_i^0$ pour un parabolique standard $Q_i \leq \GL(m_i,D)$ et un caractère $\rho_i : D^\times \xrightarrow{\Nrm} F^\times \xrightarrow{\rho_i^0} \Fpbar^\times$.
\end{itemize}
Aussi, s'il y a deux blocs adjacents dans le cas (c), alors on a $\rho_i \neq \rho_{i+1}$.
\end{theo}
\textsf{Preuve :} \\
On va effectuer une récurrence sur $m \geq 1$. Pour l'étape d'initiation $m=1$, il n'y a rien à prouver: $\pi$ est une représentation irréductible admissible de $D^\times$, qui est ou bien un caractère (cas $(c)$), ou bien de dimension strictement supérieure à $1$ (cas $(b)$). Supposons donc $m>1$ et le résultat vrai pour tout $n<m$. Soit $\pi$ une représentation irréductible admissible de $G$. Si $\pi$ est supersingulière, il n'y a rien à faire. Supposons donc qu'il existe un $K$-poids $V$ de $\pi$ et un caractère propre $\chi$ sur $\Hom_K(V,\pi)$ de paramètre HS-dual une paire $(L,\chi_L)$ avec $L \neq G$. On prend les notations suivantes pour la paramétrisation de $V$: \[ \Par_{\overline{B}^-}(V) = (\chi_{V^*}, - \Delta_V^-) .\] Plusieurs cas distincts de présentent à nous, que l'on va discuter séparément. \\
\underline{Cas 1} : $\Delta_V^- \cup \Delta_L \neq \Delta$. \\ On choisit $\alpha \in \Delta \smallsetminus (\Delta_V^- \cup \Delta_L)$ et $\lambda$ le copoids fondamental associé à $\alpha$ de dernier coefficient diagonal \' egal \` a $1$. Soit $P_\lambda = M_\lambda N_\lambda$ le parabolique standard associé, avec $\Delta_{M_\lambda} = \Delta \smallsetminus \{ \alpha \} \supseteq \Delta_L \cup \Delta_V^-$. Alors $V$ est $P_\lambda^-$-régulier et $\chi$ se factorise à travers \[ \chi^{(M)} : \Hh_{\Fpbar} \big( M_\lambda , M_\lambda \cap K , V_{N_\lambda \cap K}\big) \to \Fpbar .\] Par la proposition \ref{Thm3.1}, on a alors l'isomorphisme de comparaison \[ \ind_K^G \, V  \otimes_{\Hh,\chi} \Fpbar \xrightarrow{\sim} \Ind_{P_\lambda}^G \big( \ind_{M_\lambda \cap K}^{M_\lambda} \, V_{N_\lambda \cap K} \otimes_{\Hh_{M_\lambda}, \chi^{(M)}} \Fpbar \big) ;\] et par irréductibilité de $\pi$, on a alors une surjection \[ \Ind_{P_\lambda}^G \big( \ind_{M_\lambda \cap K}^{M_\lambda} \, V_{N_\lambda \cap K} \otimes_{\Hh_{M_\lambda}, \chi^{(M)}} \Fpbar \big) \twoheadrightarrow \pi .\] Comme $\ind_{M_\lambda \cap K}^{M_\lambda} \, V_{N_\lambda \cap K} \otimes_{\Hh_{M_\lambda}, \chi^{(M)}} \Fpbar$ est lisse à caractère central, on lui applique le lemme \ref{findadmiss} et on obtient une représentation $\sigma$ irréductible admissible de $M_\lambda$ avec une surjection $\Ind_{P_\lambda}^G \, \sigma \twoheadrightarrow \pi$. Ecrivons alors $M_\lambda = \GL(m_1,D) \times \GL(m_2,D)$ avec $0<m_1 <m$ et $\sigma = \sigma^{(1)} \otimes \sigma^{(2)}$; on applique l'hypothèse de récurrence à $\sigma^{(1)}$ et $\sigma^{(2)}$. Ces représentations $\sigma^{(1)}$ et $\sigma^{(2)}$ sont alors de la forme voulue et par transitivité de l'induction parabolique on a une surjection \[ \Ind_P^G ( \sigma_1 \otimes \cdots \otimes \sigma_r) \twoheadrightarrow \pi \] avec $P \leq P_\lambda$ et les $\sigma_i$ comme voulus. Si deux quelconques blocs adjacents de type \textit{(c)} vérifient de plus la condition $\rho_i \neq \rho_{i+1}$, alors le théorème \ref{irreductibilite} nous affirme que $\Ind_P^G ( \sigma_1 \otimes \cdots \otimes \sigma_r)$ est irréductible, et que l'on a par suite \[ \Ind_P^G ( \sigma_1 \otimes \cdots \otimes \sigma_r) \xrightarrow{\sim} \pi .\] Sinon, le corollaire \ref{JHtotal} nous donne exactement les constituants de Jordan-Hölder de $\Ind_P^G ( \sigma_1 \otimes \cdots \otimes \sigma_r)$, qui sont tous de la forme désirée avec la condition supplémentaire $\rho_i \neq \rho_{i+1}$; et l'un d'eux se surjectant sur $\pi$, c'est l'isomorphisme qu'on attendait. \\
\underline{Cas 2} : $\Delta_V^- \cup \Delta_L = \Delta$, et $\widetilde{\alpha}^\vee(\varpi) \in \widetilde{Z}_L(V^*)$, $\chi_L \big( \widetilde{\alpha}^\vee (\varpi) \big) =1$ pour tout $\alpha \in \Delta \smallsetminus \Delta_L$. \\ 
La condition \og $\widetilde{\alpha}^\vee(\varpi) \in \widetilde{Z}_L(V^*)$ pour tout $\alpha \in \Delta \smallsetminus \Delta_L$ \fg\ implique $L = A$ et $\widetilde{Z}_A(V^*) = \widetilde{Z}_A$. En particulier, on a $\Delta_L = \varnothing$ et donc $\Delta_V^-$ est tout $\Delta$: $V$ est un caractère de $K$. Combiné à $\widetilde{Z}_A(V^*) = \widetilde{Z}_A$, cela nous assure que $V$ se factorise par $\det_G$. Quitte à tensoriser $\pi$ par un $\xi \circ \det_G$ pour un certain caractère $\xi$ de $F^\times$, on peut ainsi supposer que $\chi_L$ est le caractère trivial de $\widetilde{Z}_A$ (voir proposition \ref{paramtorsion}) et que $V_{U \cap K}$ est le caractère trivial de $A \cap K$; c'est ce que l'on va faire par la suite. De plus, parce que $\chi$ est un caractère propre sur $\Hom_K(V,\pi)$, la proposition \ref{alternateOrd} (appliquée à $B^-$) nous apporte l'alternative suivante: ou bien $\pi$ est $\id$ et la preuve est terminée, ou bien il existe un parabolique standard $Q$ propre de $G$ avec $\Ord_{Q^-} \pi \neq 0$. Mais, dans ce dernier cas, $\Ord_{Q^-} \pi$ étant admissible par \cite{Vig12}, Theorem 1.10.(i), il contient une sous-représentation irréductible $\pi'$ pour laquelle on a une surjection $\Ind_Q^G \, \pi' \twoheadrightarrow \pi$. On peut appliquer l'hypothèse de récurrence à $\pi'$ et la transitivité de l'induction permet de conclure. \\
\underline{Cas 3} : $\Delta_V^- \cup \Delta_L = \Delta$, et il existe $\alpha \in \Delta_V^- \smallsetminus \Delta_V^- \cap \Delta_L$ avec $\widetilde{\alpha}^\vee(\varpi) \notin \widetilde{Z}_L(V^*)$ ou $\chi_L \big( \widetilde{\alpha}^\vee (\varpi) \big) \neq 1$. \\
Supposons dans un premier temps l'hypothèse \ref{Dpremier}.(a) satisfaite. On prend $\alpha$ comme dans la condition du Cas 3, $V'$ comme en (\ref{defV'2}), $d_0$ associé pour $V_{N_{(\alpha)} \cap K}$. A cause de l'hypothèse \ref{Dpremier}.(a), l'hypothèse \ref{V'toV} est automatique. On peut donc appliquer le corollaire \ref{chgtpdsdual}, et obtenir l'isomorphisme \[ \ind_K^G \, V \otimes_{\Hh, \chi} \Fpbar \xrightarrow{\sim} \ind_K^G \, V' \otimes_{\Hh, \chi} \Fpbar \] de $G$-représentations. On peut ensuite appliquer le Cas 1 à $V'$. \\
Supposons maintenant l'hypothèse \ref{Dpremier}.(b). Si on a $L=A$, on prend $\alpha$ comme dans la condition du Cas 3, $V'$ comme en (\ref{defV'2}), $d_0$ associé pour $V_{N_{(\alpha)} \cap K}$. La condition $\Delta_V^- \cup \Delta_L = \Delta$ implique $\Delta_V^- = \Delta$: $V$ est un caractère de $K$. Mais alors, l'hypothèse \ref{V'toV} est automatique, et on peut appliquer de nouveau le corollaire \ref{chgtpdsdual}, puis appliquer le Cas 1 à $V'$. Enfin, si $L$ est différent de $A$, pour $m=2$, $L$ est $G$ et il n'y a rien à faire. On suppose donc $m=3$ et $L \neq A$: il existe une racine $\alpha \in \Delta_V^- \smallsetminus \Delta_V^- \cap \Delta_L$ telle que $\widetilde{\alpha}^\vee(\varpi)^{d_0}$ ne soit pas dans $\widetilde{Z}_L(V^*)$ (où $d_0$ est associé à $V_{N_{(\alpha)} \cap K}$). On prend $V'$ comme en (\ref{defV'2}) on peut appliquer le corollaire \ref{chgtpdsdualm=3}: on a \[ \ind_K^G \, V \otimes_{\Hh, \chi} \Fpbar \xrightarrow{\sim} \ind_K^G \, V' \otimes_{\Hh, \chi} \Fpbar ,\] et on conclut en appliquant le Cas 1 à $V'$. \hfill$\Box$

\section{Représentations supersingulières, représentations supercuspidales, entrelacements} \label{scuspssing}

Dans tout ce paragraphe, on supposera la formule de comparaison de tranformées de Satake (conjecture \ref{conjSatake}) satisfaite. \\

Une représentation irréductible $\pi$ de $G$ est dite supercuspidale si elle n'est isomorphe à aucun sous-quotient d'un $\Ind_P^G \, \sigma$ pour un parabolique standard $P = MN$ propre de $G$ et une représentation irréductible admissible $\sigma$ de $M$. On renvoit au paragraphe \ref{prgssing} pour la définition de la notion de supersingularité. 

\begin{theo}
Supposons l'hypothèse \ref{Dpremier} satisfaite. Soient $m \geq 2$ et $\pi$ une représentation irréductible admissible de $G$. Alors $\pi$ est supersingulière si et seulement si elle est supercuspidale.
\end{theo}
\textsf{Preuve :} \\
Si $\pi$ est non supersinguli\` ere, par le th\' eor\` eme \ref{classifadmis}, il existe un parabolique standard $P = MN$ et une repr\' esentation irr\' eductible admissible $\sigma$ de $M$ tels que $\pi$ est isomorphe \` a $\Ind_P^G \, \sigma$. Ainsi $\pi$ n'est pas supercuspidale. \\
Avant que de prouver la r\' eciproque, faisons la remarque que le corollaire \ref{JHtotal} affirme que toute induite parabolique d'une repr\' esentation irr\' eductible admissible est de longueur finie, de constituants de Jordan-H\" older bien d\' etermin\' es, tous isomorphes \` a des induites paraboliques de repr\' esentations irr\' eductibles admissibles: si $Q$ est \' egal \` a $G$, c'est une repr\' esentation irr\' eductible admissible (a), (b) ou (c); et le cas (a) ne se produit que lorsque l'on a $P=Q=G$. De ce fait, une repr\' esentation irr\' eductible admissible de $G$ est supercuspidale si et seulement si elle n'est pas isomorphe \` a une induite parabolique propre ou une repr\' esentation de Steinberg g\' en\' eralis\' ee. \\
Supposons que $\pi$ n'est pas supercuspidale. Si $\pi$ est isomorphe \` a une repr\' esentation de Steinberg g\' en\' eralis\' ee, la proposition \ref{paramSt} implique que $\pi$ n'est pas supersinguli\` ere. On suppose donc $\pi$ isomorphe \` a une induite parabolique propre $\Ind_P^G \, \sigma$ avec $P=MN$.
Soient $V$ un $K$-poids de $\pi$ et $\chi : \Hh_{\Fpbar}(G,K,V) \to \Fpbar$ un caractère propre sur $\Hom_K(V,\pi)$. Soit $(L,\chi_L)$ le paramètre HS-dual de $\chi$. Par le lemme \ref{paramHSpoids}, on a $L \leq M \lneq G$ et $\pi$ n'est pas supersingulière. Le résultat en découle. \hfill$\Box$

\begin{cor} \label{ABmain2scusp}
Supposons l'hypothèse \ref{Dpremier} satisfaite. 
\begin{itemize}
\item[(i)] Soit $P$ un parabolique standard de $G$, de Levi standard $M = \prod_{i=1}^r G_i$. Soit, pour tout $1 \leq i \leq r$, $\sigma_i$ une représentation irréductible admissible de $G_i = \GL(m_i,D)$ qui est dans l'un des cas suivants:
\begin{itemize}
\item[(a)] $\sigma_i$ est supercuspidale avec $m_i > 1$;
\item[(b)] $\sigma_i$ est une représentation de $D^\times$ (i.e. $m_i = 1$) avec $\dim \, \sigma_i > 1$;
\item[(c)] $\sigma_i \simeq \St_{Q_i} \rho_i^0$ pour un parabolique standard $Q_i \leq G_i$ et un caractère $\rho_i : D^\times \xrightarrow{\Nrm} F^\times \xrightarrow{\rho_i^0} \Fpbar^\times$.
\end{itemize}
Supposons $\rho_i \neq \rho_{i+1}$ s'il y a deux blocs adjacents dans le cas (c). Alors $\Ind_{P}^G(\sigma_1 \otimes \dots \otimes \sigma_r)$ est irréductible admissible.
\item[(ii)] Soit $\pi$ une représentation irréductible de $G$. Alors il existe un parabolique standard $P$ de Levi standard $\prod_i G_i$ et des représentations irréductibles admissibles $\sigma_i$ de $G_i$ comme en $(i)$ tels que $\pi$ soit isomorphe à $\Ind_P^G(\sigma_1 \otimes \cdots \otimes \sigma_r)$.
\end{itemize}
\end{cor}

On va finir par regarder les entrelacements entre les représentations irréductibles admissibles que l'on a obtenues pour conclure quant à une classification.

\begin{lem} \label{onlyparam}
Supposons l'hypothèse \ref{Dpremier} satisfaite. Soit $\pi$ une représentation irréductible admissible de $G$. Alors il existe une paire $(L, \omega_L : Z(L) \to \Fpbar^\times)$ telle que, pour tout $K$-poids $V$ de $\pi$ et tout caractère $\chi : \Hh_{\Fpbar}(G,K,V) \to \Fpbar$ propre sur $\Hom_K(V,\pi)$, on a:
\begin{itemize}
\item[(i)] le paramètre HS-dual de $\chi$ est une paire $(L,\chi_L)$ où $\chi_L$ vérifie \[ \chi_L |_{Z(L) \cap A_\Lambda} = \omega_L|_{Z(L) \cap A_\Lambda} ;\]
\item[(ii)] la restriction de $V_{U \cap K}$ à $Z(L) \cap K$ est égale à $\omega_L^{-1}|_{Z(L) \cap K}$.
\end{itemize}
\end{lem}
\textsf{Remarque :} 
En particulier, on ne sait pas s'il existe un sous-groupe $\widetilde{Z}_L^\pi$ de $\widetilde{Z}_L$ tel que pour tout $K$-poids $V$ de $\pi$ on a $\widetilde{Z}_L(V^*) = \widetilde{Z}_L^\pi$ (à cause du cas supersingulier); c'est propre au cas non déployé. \\ 
\textsf{Preuve :} \\
Par le théorème \ref{classifadmis}, il existe un parabolique standard $P$ de Levi standard $M$ et une représentation irréductible admissible $\sigma$ de $M$ telle que $\pi$ soit isomorphe à $\Ind_P^G \, \sigma$; de plus $\sigma$ vérifie les conditions énoncées dans le théorème \ref{classifadmis}. Par suite, le lemme \ref{paramHSpoids} nous affirme que, pour tout $K$-poids $V$ de $\pi$, tout caractère propre $\chi$ sur $\Hom_K(V,\pi)$ est de paramètre HS-dual $(L,\chi_L)$ où $L$ est uniquement déterminé par $M$ et $\sigma$, et est donc indépendant de $\chi$. \\ Le fait que $\chi_L |_{Z(L) \cap A_\Lambda}$ est uniquement déterminé provient des lemmes \ref{Lem4.3} et \ref{Lem4.6}. \\
Prouvons enfin (ii): par r\' eciprocit\' e de Frobenius, on a 
\[ \Hom_K \big( V, \Ind_P^G \, \sigma \big) \simeq \Hom_{P \cap K}(V,\sigma) = \Hom_{M \cap K} \big( V_{N \cap K} , \sigma \big) .\]
Comme $V_{U \cap K}$ est un quotient de $V_{N \cap K}$ et que $\sigma$ et $V_{N \cap K}$ se factorisent par le lemme \ref{decompGi}, on est ramen\' e au cas $M=G$, ce que l'on suppose \` a pr\' esent. D\` es lors, ou bien $\sigma$ n'est pas une repr\' esentation de Steinberg g\' en\' eralis\' ee et on a $L=M=G$. Parce que $\sigma$ est irr\' eductible admissible, le lemme de Schur pr\' edit que l'action du centre $Z(G)=Z(L)$ est donn\' ee par un caract\` ere $\psi$. Tout $K$-poids $V$ de $\sigma$ voit alors son caract\` ere central \^ etre \' egal \` a $\psi|_{Z(G) \cap K}$ et, $Z(G) \cap K$ agit via $\psi|_{Z(G) \cap K}$ sur $V_{U \cap K}$ aussi. Ou bien $\sigma$ est une repr\' esentation de Steinberg g\' en\' eralis\' ee $\St_Q \xi$ et on a $L = A$. Dans ce cas-ci, (ii) r\' esulte de ce que le $K$-socle de $\St_Q \xi$ est irr\' eductible par la proposition \ref{Stsocle} et l'unicit\' e est alors \' evidente. La preuve du lemme est termin\' ee. \hfill$\Box$\\

Ainsi, en général, seule la restriction à $Z(L) \cap A_\Lambda$ de $\chi_L^\chi$ est connue. Cependant, lorsque l'on a affaire à des représentations de Steinberg généralisées, on connaît entièrement ce paramètre comme ce sera transparent au cours de la preuve du résultat suivant. 

\begin{lem} \label{paramStrec}
Soit $\pi$ une représentation irréductible admissible de $G$. On suppose qu'il existe un $K$-poids $V$ de $\pi$ et un caractère $\chi : \Hh_{\Fpbar}(G,K,V) \to \Fpbar$ propre sur $\Hom_K(V,\pi)$ vérifiant:
\begin{itemize}
\item[(a)] $\widetilde{Z}_A(V^*) = \widetilde{Z}_A$ et $V_{U \cap K}$ est un caractère de $A \cap K$ qui se factorise par $\det_G$;
\item[(b)] le paramètre HS-dual de $\chi$ est $(A,\xi_\varpi^{-1} \circ \det_G : \widetilde{Z}_A \to \Fpbar)$ pour un certain caractère $\xi_\varpi$ de $((-1)^{d-1} \varpi)^\Z$.
\end{itemize}
Alors il existe un parabolique standard $Q$ de $G$ et un caractère $\xi$ de $F^\times$ prolongeant $\xi_\varpi$ tels que $\pi$ soit isomorphe à $\St_Q \xi$.
\end{lem}
\textsf{Remarque :} 
A posteriori, $(a)$ et $(b)$ sont vérifiées par tout $K$-poids $V$ de $\pi$ et tout caractère propre $\chi$ sur $\Hom_K(V,\pi)$ puisque le poids $V$ est unique et qu'il n'y a qu'un seul caractère propre (propositions \ref{Stsocle} et \ref{paramSt}). \\
\textsf{Preuve :} \\
Par le théorème \ref{classifadmis}, il existe un parabolique standard $P$ de Levi standard $M = \prod\limits_{i=1, \dots,r} M_i$ et des représentations $\sigma_i$ de $M_i$ (qui sont dans le cas $(a)$, $(b)$ ou $(c)$) telles que $\pi$ soit isomorphe à $\Ind_P^G(\sigma_1 \otimes \dots \otimes \sigma_r)$. Par les lemmes \ref{onlyparam} et \ref{paramHSpoids}, chaque $\sigma_i$ est une représentation de Steinberg (éventuellement un caractère) $\St_{Q_i} \xi_i$ pour un parabolique standard $Q_i$ de $M_i$  et un caractère $\xi_i$ de $F^\times$. Soient $V$ un $K$-poids de $\pi$ et $\chi : \Hh_{\Fpbar}(G,K,V) \to \Fpbar$ un caractère propre sur $\Hom_K(V,\pi)$ satisfaisant $(a)$ et $(b)$. Le paramètre HS-dual de $\chi$ est $(L,\chi_L^\chi) = (A,\xi_\varpi^{-1} \circ \det_G)$ par hypothèse. Si $\alpha$ est la racine entre les blocs $M_i$ et $M_{i+1}$, en calculant $\chi_L^\chi \big( \widetilde{\alpha}^\vee(\varpi) \big)$, par les lemmes \ref{decompGi} et \ref{paramSt}, et l'hypothèse $(b)$, on a 
\[ \xi_{i+1}((-1)^{d-1} \varpi) = \xi_i((-1)^{d-1} \varpi) = \xi_\varpi((-1)^{d-1} \varpi) .\] 
Enfin, par l'hypothèse $(a)$, on a $\xi_{i+1}|_{\Oo_F^\times} = \xi_i|_{\Oo_F^\times}$. Comme le théorème \ref{classifadmis} impose $\xi_i \neq \xi_{i+1}$, cela veut dire que l'on a $r=1$ et $P=G$. Le résultat est prouvé en prenant $\xi = \xi_1$. \hfill$\Box$

\begin{prop}
Supposons l'hypothèse \ref{Dpremier} satisfaite. Il n'y a pas d'entrelacement entre les induites paraboliques propres irréductibles $\Ind_P^G(\sigma_1 \otimes \cdots \otimes \sigma_r)$ comme dans le théorème \ref{irreductibilite}. 
\end{prop}
\textsf{Preuve :} \\
Soient $P$ et $\sigma_1, \dots, \sigma_r$ comme dans le théorème \ref{irreductibilite}. Soit $\pi = \Ind_P^G(\sigma_1 \otimes \dots \otimes \sigma_r)$, qui est donc irréductible admissible. Pour tout $K$-poids $V$ de $\pi$ et tout caractère $\chi : \Hh_{\Fpbar}(G,K,V) \to \Fpbar$ propre sur $\Hom_K(V,\pi)$, on note $\chi^V$ le caractère $V_{U \cap K}$ de $A \cap K$ et $(L,\chi_L^\chi)$ le paramètre HS-dual de $\chi$.  Voyons que l'on a :  
\[ \Delta_P = \Delta_L \cup \ \bigcap\nolimits_{V,\chi} \big\{ \alpha \in \Delta \ \big| \ \widetilde{\alpha}^\vee(\varpi) \in \widetilde{Z}_L(V^*) , \, \chi_L^\chi \big( \widetilde{\alpha}^\vee(\varpi) \big) =1, \, \chi^V(s_\alpha \cdot s_{\alpha}^{-1})|_{A \cap K} = \chi^V|_{A \cap K} \big\} .\] 
En effet, par le lemme \ref{paramHSpoids}, $P$ et $\sigma_1 \otimes \dots \otimes \sigma_r$ d\' eterminent uniquement $L$, et le Levi standard $M$ de $P$ co\" incide avec $L$ sur les blocs o\` u il n'y a pas de repr\' esentation de Steinberg. Comme $L$ est uniquement d\' etermin\' e par $\pi$ (voir lemme \ref{onlyparam}), il s'agit de voir que l'on peut d\' etecter intrins\` equement les blocs o\` u se trouvent des repr\' esentations de Steinberg g\' en\' eralis\' ees. La condition $\widetilde{\alpha}^\vee(\varpi) \in \widetilde{Z}_L(V^*)$ signifie que $\alpha$ est une racine \` a l'int\' erieur d'un bloc avec une Steinberg ou entre deux blocs adjacents portant des Steinberg. Les deux autres conditions viennent du lemme \ref{paramStrec}; on remarquera que, en vertu de la proposition \ref{Stsocle} sur l'irr\' eductibilit\' e du $K$-socle et l'unicit\' e des param\` etres HS pour les repr\' esentations de Steinberg, la formule donnant $\Delta_P$ peut se r\' e\' ecrire pour un choix de $(V,\chi)$ quelconque, sans en prendre l'intersection. Mais on a pr\' ef\' er\' e cette forme, sur laquelle on lit facilement l'ind\' ependance en $V$ et en $\chi$. \\
Ainsi $P$ est uniquement déterminé, et il nous reste à examiner les entrelacements entre deux induites paraboliques irréductibles $\Ind_P^G \, \sigma$ et $\Ind_P^G \, \tau$ pour $\sigma$ et $\tau$ des représentations irréductibles admissibles du Levi standard $M$ de $P$. Mais alors, par \cite{Vig12}, Theorem 1.10.(ii), ou \cite{Eme10}, Corollary 4.3.5 (en caractéristique $0$), on a: \[ \sigma \simeq \Ord_{P^-} \big( \Ind_P^G \, \sigma \big) \simeq \Ord_{P^-} \big( \Ind_P^G \, \tau \big) \simeq \tau .\] Ceci termine la preuve. \hfill $\Box$

\begin{cor}
Supposons l'hypothèse \ref{Dpremier} satisfaite. Soit $\pi$ une représentation irréductible de $G$. Alors il existe un unique parabolique standard $P$ de Levi standard $\prod_i G_i$ et des représentations irréductibles admissibles $\sigma_i$ de $G_i$ comme en $(i)$, uniques à isomorphisme près, tels que $\pi$ soit isomorphe à $\Ind_P^G(\sigma_1 \otimes \cdots \otimes \sigma_r)$.
\end{cor}

Terminons en disant que tous les résultats des paragraphes \ref{classifadmiss} et \ref{scuspssing} nécessitant l'hypothèse \ref{Dpremier} ne l'utilisent qu'à travers la proposition \ref{chgtpdsreg}. Si on arrive à montrer cette dernière sous une hypothèse plus faible, le reste suivra.

\section{Appendice: diviseurs élémentaires pour $\GL(2,D)$ et $\GL(3,D)$} \label{appdivelem}

On garde les notations du lemme \ref{egalcardinaux} et on veut prouver que la caractérisation des diviseurs élémentaires par les mineurs est encore valable pour $\GL(2,D)$ et $\GL(3,D)$ (dans le cadre utile au lemme \ref{egalcardinaux} uniquement). \\
Tout d'abord, on remarque que $g = \left( \begin{array}{cc} r & s \\ t & u \end{array} \right) \in \GL(2,D)$ est dans la double classe de Cartan $K \left( \begin{array}{cc} \varpi_D^\alpha & 0 \\ 0 & \varpi_D^\beta \end{array} \right) K$ (avec $\alpha \leq \beta$) si et seulement si on a \[ \begin{cases} \min(v_D(r),v_D(s),v_D(t),v_D(u)) = \alpha \\ d^{-1} v_D \circ \det_G (g) = \alpha + \beta \end{cases} .\] Cela vient de la preuve de l'existence des diviseurs élémentaires et de ce que le $\det_G$ d'un élément de $K$ a valuation égale à $0$. Ce fait sera sans cesse utilisé pour l'étude de cas dans $\GL(3,D)$. \\

Enonçons le cas de $\GL(3,D)$ sous forme d'un lemme.

\begin{lem}
Supposons $m=3$. Soit \[ g = \left( \begin{array}{ccc} \varpi_D^x & a & b \\ 0 & \varpi_D^y & c \\ 0 & 0 & \varpi_D^z \end{array} \right) \in \GL(3,D) \] 
avec $x \leq y \leq z$. Alors $g$ appartient à la double classe $K \left( \begin{array}{ccc} \varpi_D^\alpha & 0 & 0 \\ 0 & \varpi_D^\beta & 0 \\ 0 & 0 & \varpi_D^\gamma \end{array} \right) K$ avec $\alpha \leq \beta \leq \gamma$ si et seulement si on a\footnote{à cause des inégalités $x \leq y \leq z$, certains des termes peuvent être supprimés; on a choisi de les garder dans l'énoncé pour rappeler que l'on considère bien tous les mineurs}
\begin{equation} \label{divelemD2} \left\{ \begin{array}{l} \min(x,y,z,
      \overline{a}, \overline{b}, \overline{c}) = \alpha \\ x+y+z = \alpha
      + \beta + \gamma \\ \min(x+y,y+z,x+z,z+\overline{a},x+ \overline{c},
      d^{-1} v_D \Big( \det_G \left( \begin{array}{cc} a & b \\ \varpi_D^y
          & c \end{array} \right) \Big) ) = \alpha + \beta \end{array}
  \right. .\end{equation}
\end{lem} 
\textsf{Preuve :} \\
On va voir que l'on a bien les conditions (\ref{divelemD2}) en distinguant suivant le coefficient de valuation minimale de $g$. On remarque que l'on établit uniquement que les conditions (\ref{divelemD2}) sont nécessaires, mais elles sont alors aussi suffisantes car $G$ est union disjointe de ses doubles classes de Cartan. \\ 
De la même manière que pour $\GL(2,D)$, les deux premières conditions sont faciles. Il reste à établir la troisième condition, qui devient
\begin{equation} \label{cond3soft} \min \big( x+y,z+\overline{a},x+ \overline{c}, d^{-1} v_D \circ \det\nolimits_G(g') \big) = \alpha + \beta \end{equation}
à cause des inégalités $x \leq y \leq z$, où on a noté $g' = \left( \begin{array}{cc} a & b \\ \varpi_D^y & c \end{array} \right)$. \\ On s'attache à présent à établir (\ref{cond3soft}) pour $g$: pour cela, on effectue des pivots avec des matrices de transvection dans $K$. \\

\textsf{Cas 1: $x$ minimal} \\
On écrit 
\[ \left( \begin{array}{ccc} \varpi_D^x & a & b \\ 0 & \varpi_D^y & c \\ 0 & 0 & \varpi_D^z \end{array} \right) \left( \begin{array}{ccc} 1 & - \varpi_D^{-x} a & - \varpi_D^{-x} b \\ 0 & 1 & 0 \\ 0 & 0 & 1 \end{array} \right) = \left( \begin{array}{ccc} \varpi_D^x & 0 & 0 \\ 0 & \varpi_D^y & c \\ 0 & 0 & \varpi_D^z \end{array} \right) .\] 
Les valuations des facteurs invariants de $g$ sont donc \[ \big( x, \min(\overline{c},y), y+z- \min(\overline{c},y) \big) \] et on veut montrer 
\begin{equation} \label{cond3pourx} x+ \min(\overline{c},y) = \min \big( x+y,z+\overline{a},x+ \overline{c}, d^{-1} v_D \circ \det\nolimits_G(g') \big) .\end{equation}
Parce que $x$ est minimal, on a $\overline{b} \geq x$ et, grâce à (\ref{calculdet2}) pour $c \neq 0$, il s'ensuit
\[ d^{-1} v_D \circ \det\nolimits_G(g') \geq \overline{b}+y \geq x+y \geq x+ \min(\overline{c},y) ;\]
de même, on a $\overline{a} \geq x$, et donc
\[ z + \overline{a} \geq z+x \geq y+x \geq x+ \min(\overline{c},y) .\]
Ainsi, l'égalité (\ref{cond3pourx}) est vérifiée, ce qui donne bien la condition (\ref{cond3soft}) lorsque $x$ est minimal. \\

\textsf{Cas 2: $\overline{a}$ minimal} \\
On écrit
\[\hspace{-0.9cm} \left( \begin{array}{ccc} 1 & 0 & 0 \\ - \varpi_D^y a^{-1} & 1 & 0 \\ 0 & 0 & 1 \end{array} \right) \left( \begin{array}{ccc} \varpi_D^x & a & b \\ 0 & \varpi_D^y & c \\ 0 & 0 & \varpi_D^z \end{array} \right) \left( \begin{array}{ccc} 1 & 0 & 0 \\ - a^{-1} \varpi_D^x & 1 & -a^{-1} b \\ 0 & 0 & 1 \end{array} \right) = \left( \begin{array}{ccc} 0 & a & 0 \\ - \varpi_D^y a^{-1} \varpi_D^x & 0 & c- \varpi_D^y a^{-1} b \\ 0 & 0 & \varpi_D^z \end{array} \right) .\]
Les valuations des facteurs invariants de $g$ sont donc 
\[\hspace{-0.3cm} \big( \overline{a} , \min \big( z,x+y-\overline{a}, v_D(c - \varpi_D^y a^{-1} b) \big) , x+y+z - \overline{a} - \min \big( z,x+y-\overline{a}, v_D(c - \varpi_D^y a^{-1} b) \big) \big) \] 
et on veut établir
\begin{equation} \label{cond3poura} \overline{a} + \min \big( z,x+y-\overline{a}, v_D(c - \varpi_D^y a^{-1} b) \big) = \min \big( x+y,z+\overline{a},x+ \overline{c}, d^{-1} v_D \circ \det\nolimits_G(g') \big) .\end{equation}
Parce que $a$ est non nul, on écrit 
\[ \left( \begin{array}{cc} 1 & 0 \\ - \varpi_D^y a^{-1} & 1 \end{array} \right) \left( \begin{array}{cc} a & b \\ \varpi_D^y & c \end{array} \right) = \left( \begin{array}{cc} a & b \\ 0 & c - \varpi_D^y a^{-1} b \end{array} \right) \]
et donc 
\begin{equation} \label{calculdet2b} d^{-1} v_D \circ \det\nolimits_G(g') = v_D(ac - a \varpi_D^y a^{-1} b) .\end{equation}
On sépare ensuite deux sous-cas: ou bien on a $\overline{a}+ \overline{c} \neq \overline{b}+y$, et alors par minimalité de $\overline{a}$, on a
\[x + \overline{c} \geq \overline{a} + \overline{c} \geq v_D(ac - a \varpi_D^y a^{-1} b) ;\]
ou bien, on a $\overline{a}+ \overline{c} = \overline{b}+y$, et toujours par minimalité de $\overline{a}$, il s'ensuit
\[ x + \overline{c} = x+y + \overline{b} - \overline{a} \geq x+y .\]
Dans les deux cas, (\ref{cond3poura}) est vérifiée, ce qui donne bien la condition (\ref{cond3soft}) lorsque $\overline{a}$ est minimal. \\

\textsf{Cas 3: $\overline{c}$ minimal} \\
On écrit
\[ \left( \begin{array}{ccc} 1 & -bc^{-1} & 0 \\ 0 & 1 & 0 \\ 0 & - \varpi_D^z c^{-1} & 1 \end{array} \right) \left( \begin{array}{ccc} \varpi_D^x & a & b \\ 0 & \varpi_D^y & c \\ 0 & 0 & \varpi_D^z \end{array} \right) \left( \begin{array}{ccc} 1 & 0 & 0 \\ 0 & 1 & 0 \\ 0 & - c^{-1} \varpi_D^y & 1 \end{array} \right) = \left( \begin{array}{ccc} \varpi_D^x & a - bc^{-1} \varpi_D^y & 0 \\ 0 & 0 & c \\ 0 & - \varpi_D^z c^{-1} \varpi_D^y & 0 \end{array} \right) .\]
Parce que $\overline{c}$ est minimal, on a $x \leq y+z-\overline{c}$, et donc les valuations des facteurs invariants sont
\[ \big( \overline{c}, \min \big( x, v_D(a-bc^{-1} \varpi_D^y) \big),x+y+z-\overline{c}-\min \big( x, v_D(a-bc^{-1} \varpi_D^y) \big) \big) ;\]
on cherche à montrer
\begin{equation} \label{cond3pourc} \overline{c} + \min \big( x, v_D(a-bc^{-1} \varpi_D^y) \big) = \min \big( x+y,z+\overline{a},x+ \overline{c}, d^{-1} v_D \circ \det\nolimits_G(g') \big) .\end{equation}
Comme $\overline{c}$ est minimal et que l'on a $z \geq x$, on a 
\[ x+y \geq x+ \overline{c} ; \quad z+ \overline{a} \geq z+\overline{c} \geq x + \overline{c} .\]
Alors (\ref{cond3pourc}) est vérifiée, ce qui donne bien la condition (\ref{cond3soft}) lorsque $\overline{c}$ est minimal. \\

\textsf{Cas 4: $\overline{b}$ minimal} \\
On écrit
\[\hspace{-0.9cm} \left( \begin{array}{ccc} 1 & 0 & 0 \\ -cb^{-1} & 1 & 0 \\ - \varpi_D^z b^{-1} & 0 & 1 \end{array} \right) \left( \begin{array}{ccc} \varpi_D^x & a & b \\ 0 & \varpi_D^y & c \\ 0 & 0 & \varpi_D^z \end{array} \right) \left( \begin{array}{ccc} 1 & 0 & 0 \\ 0 & 1 & 0 \\ - b^{-1} \varpi_D^x & - b^{-1} a & 1 \end{array} \right) = \left( \begin{array}{ccc} 0 & 0 & b \\ -c b^{-1} \varpi_D^x & \varpi_D^y - cb^{-1} a & 0 \\ - \varpi_D^z b^{-1} \varpi_D^x & -\varpi_D^z b^{-1} a & 0 \end{array} \right) .\]
Pour le cas $\overline{a} \geq x$, on continue: 
\[ \left( \begin{array}{cc} -cb^{-1} \varpi_D^x & \varpi_D^y - cb^{-1} a \\ - \varpi_D^z b^{-1} \varpi_D^x & - \varpi_D^z b^{-1} a \end{array} \right) \left( \begin{array}{cc} 1 & - \varpi_D^{-x} a \\ 0 & 1 \end{array} \right) = \left( \begin{array}{cc} - cb^{-1} \varpi_D^x & \varpi_D^y \\ - \varpi_D^z b^{-1} \varpi_D^x & 0 \end{array} \right) .\]
Parce que l'on a $z+x- \overline{b} \geq y$, les valuations des facteurs invariants de $g$ sont 
\[ \big( \overline{b} , \min(y,x+ \overline{c}-\overline{b}), x+y+z- \overline{b} - \min(y,x+ \overline{c}-\overline{b}) \big).\]
On cherche à montrer 
\begin{equation} \label{cond3pourb1} \overline{b} + \min(y,x+ \overline{c}-\overline{b}) = \min \big( x+y,z+\overline{a},x+ \overline{c}, d^{-1} v_D \circ \det\nolimits_G(g') \big).\end{equation}
On a d'ores et déjà, par minimalité de $\overline{b}$
\[ x+y \geq \overline{b} + y ; \quad z+ \overline{a} \geq z +x \geq y+x \geq y+ \overline{b} .\]
Ensuite, il faut distinguer: 
\begin{itemize}
\item si $\overline{b}+y<\overline{a}+\overline{c}$, alors on a $d^{-1} v_D \circ \det_G(g') = \overline{b}+y$;
\item si $\overline{b}+y=\overline{a}+\overline{c}$, alors on a $d^{-1} v_D \circ \det_G(g') \geq \overline{b}+y = \overline{a}+\overline{c} \geq x + \overline{c}$;
\item si $\overline{b}+y>\overline{a}+\overline{c}$, alors on a $\overline{b}+y > d^{-1} v_D \circ \det_G(g') = \overline{a}+\overline{c} \geq x + \overline{c}$.
\end{itemize}
Dans chacun de ces cas, (\ref{cond3pourb1}) est donc vérifiée, ce qui donne bien la condition (\ref{cond3soft}) lorsque $\overline{b}$ est minimal et que l'on a $\overline{a} \geq x$. \\
Pour le cas $\overline{a} < x$ (en particulier $a \neq 0$), on écrit 
\[ \left( \begin{array}{cc} -cb^{-1} \varpi_D^x & \varpi_D^y - cb^{-1} a \\ - \varpi_D^z b^{-1} \varpi_D^x & - \varpi_D^z b^{-1} a \end{array} \right) \left( \begin{array}{cc} 1 & 0 \\ - a^{-1} \varpi_D^x & 1 \end{array} \right) = \left( \begin{array}{cc} - \varpi_D^y a^{-1} \varpi_D^x & \varpi_D^y -cb^{-1}a \\ 0 & - \varpi_D^z b^{-1} a \end{array} \right) .\]
Les valuations des facteurs invariants de $g$ sont
\[ \big( \overline{b} , \min \big( x+y- \overline{a} , z+ \overline{a} - \overline{b} , v_D(\varpi_D^y - c b^{-1}a) \big) , x+y+z- \overline{b} - \min \big( x+y- \overline{a} , z+ \overline{a} - \overline{b} , v_D(\varpi_D^y - c b^{-1}a) \big) \big) \]
et on cherche à établir 
\begin{equation} \label{cond3pourb2} \overline{b} + \min \big( x+y- \overline{a} , z+ \overline{a} - \overline{b} , v_D(\varpi_D^y - c b^{-1}a) \big) = \min \big( x+y,z+\overline{a},x+ \overline{c}, d^{-1} v_D \circ \det\nolimits_G(g') \big) .\end{equation}
On doit encore distinguer trois sous-cas:
\begin{itemize}
\item si $\overline{b}+y<\overline{a}+\overline{c}$, alors on a  
\[ d^{-1} v_D \circ \det\nolimits_G(g') = \overline{b}+y = v_D(\varpi_D^y - c b^{-1}a) + \overline{b} ;\]
aussi, on a, à l'aide de $x>\overline{a}$ et $x \geq \overline{b}$, 
\[ x+y \geq \overline{b} + y ; \quad x+ \overline{c}> x - \overline{a} +y + \overline{b} > \overline{b}+ y ;\]
\item si $\overline{b}+y>\overline{a}+\overline{c}$, alors on a 
\[ d^{-1} v_D \circ \det\nolimits_G(g') = \overline{a}+ \overline{c} = v_D(\varpi_D^y - c b^{-1}a) + \overline{b} ;\] 
de plus, par minimalité de $\overline{b}$, on a 
\[ x+y \geq x+y - \overline{a} + \overline{b} > x + \overline{c} ;\]
on finit en observant $x + \overline{c} > \overline{a}+\overline{c}$ puisque l'on a $x>\overline{a}$;
\item si $\overline{b}+y=\overline{a}+\overline{c}$, alors on a 
\[ x+y \geq x+y+\overline{b}-\overline{a} = x+ \overline{c} ;\]
enfin on remarque, grâce à (\ref{calculdet2b}): 
\[ \overline{b} + v_D(\varpi_D^y - cb^{-1}a) = \overline{b} + v_D(a \varpi_D^y a^{-1} - ac b^{-1}) = d^{-1} v_D \circ \det\nolimits_G(g') .\]
\end{itemize}
Dans chacun de ces sous-cas, (\ref{cond3pourb2}) est donc vérifiée, ce qui donne bien la condition (\ref{cond3soft}) lorsque $\overline{b}$ est minimal et que l'on a $\overline{a} < x$. \\

Parce que l'on a $x \leq y \leq z$, on a fini d'examiner tous les cas possibles et on a bien les conditions caractéristiques des mineurs (\ref{divelemD2}) pour $g$. La preuve est terminée.

\section{Appendice: défaut de multiplicativité du prolongement $\widetilde{\chi}_\rho$} \label{appendixmult}

On garde les notations du paragraphe \ref{prgcharacter} et on va présenter une réciproque au lemme \ref{multsiNrm}. Autrement dit, on se demande à quelle condition un caractère de $K$ se prolonge à $G$.

\begin{prop} \label{conversemult}
Soit $\chi$ un caractère de $K$. Alors $\chi$ se prolonge en un caractère de $G$ si et seulement si $\chi$ se factorise par $\det_G$.
\end{prop}
\textsf{Remarque :}
Cela suit aussi du Lemme 6.1 du chapitre \ref{Ly11b}. \\
\textsf{Preuve :} \\
Le sens $(\Leftarrow)$ est l'objet du lemme \ref{multsiNrm}. Supposons que $\chi$ se prolonge à $G$ et montrons que $\chi$ se factorise par $\det_G$. Parce que $\chi$ est un caractère de $K$, il existe un morphisme de groupes $\chi_0 : k_D^\times \to \Fpbar^\times$ tel que l'on ait $\chi = \chi_0 \circ \overline{\det}$. Pour tout $x \in D^\times$, notons \[ a_1(x) = \left( \begin{array}{cccc} x & & & \\ & 1 & & \\ & & \ddots & \\ & & & 1 \end{array} \right) \in A .\] Parce que $\chi$ se prolonge à $G$, si $x$ est un élément de $\Oo_D^\times$, on a \[ \chi \big( a_1(\varpi_D x \varpi_D^{-1}) \big) = \chi \big( a_1(\varpi_D) \big) \, \chi \big( a_1(x) \big) \, \chi \big( a_1(\varpi_D^{-1}) \big) = \chi \big( a_1(x) \big) ,\] c'est-à-dire \[ \chi_0 \big( \overline{\varpi_D x \varpi_D^{-1}} \big) = \chi_0 \big( \overline{x} \big) .\] De ce fait, $\chi_0$ se factorise par la norme $N_{k_D/k_F}$ de l'extension galoisienne $k_D/k_F$. On peut alors écrire \[ \chi = \chi_1 \circ N_{k_D/k_F} \circ \overline{\det} \] sur $K$, pour un certain caractère $\chi_1 : k_F^\times \to \Fpbar^\times$. Soit $b$ un élément de $B \cap K$, de coefficients diagonaux $a_1, \dots, a_m \in \Oo_D^\times$. On a alors \[ \chi(b) = \prod\nolimits_i \chi_1 \circ N_{k_D/k_F} \big( \overline{a}_i\big) .\] Par (\ref{GLmdetNrm}) et la preuve de \cite{MatNak43}, Satz 1, on a ensuite \[ \chi(b) = \prod\nolimits_i \chi_1 \circ \overline{\Nrm(a_i)} = \chi_1 \circ \overline{\det\nolimits_G(b)} .\] Soit $s \in S_\Delta$, que l'on relève à $K$ en une matrice à coefficients dans $\{0,1\}$ (que l'on notera encore $s$). On a \[ \chi(s) = \chi_1((-1)^d) = \chi_1 \circ \overline{\det\nolimits_G(s)}.\] Parce que $S_\Delta$ engendre le groupe de Weyl $W$, $\chi$ se factorise par $\det_G$ sur un ensemble $W_K$ de relèvements de $W$ dans $K$. Comme $\chi$ se factorise par $K/K(1)$, par décomposition de Bruhat (pour $\overline{K}$), le fait que $\chi$ se factorise par $\det_G$ sur $B \cap K$ et sur $W_K$ suffit à dire que $\chi$ se factorise par $\det_G$ sur $K$. \hfill$\Box$\\

Enfin, remarquons que si pour $\GL(2,D)$ il n'y a jamais aucun problème à prolonger $\chi$ en une application pseudo-multiplicative, pour $\GL(3,D)$ déjà ce n'est plus le cas. \\
Soient $d_0 \geq 1$ un diviseur de $d$, et $a \leq b \leq c$, $x \leq y \leq z$ deux suites croissantes d'entiers. Alors, si on prend $t_1 = \Diag(\varpi_D^a,\varpi_D^b, \varpi_D^c)$, $t_2 = \Diag(\varpi_D^x, \varpi_D^y,\varpi_D^z)$ et $\alpha, \beta, \gamma \in \Oo_D$ vérifiant \[ v_D(\alpha) \in d_0 \Z, \ v_D(\beta)>0 , \ v_D(\gamma) >0 , \ v_D(\beta) + v_D(\gamma) \in d_0 \Z ,\] et $a+x+v_D(\alpha)$ la valuation minimale des coefficients de $t_1 \left( \begin{array}{ccc} \alpha & \gamma & 1 \\ \beta & 1 & 0 \\ 1 & 0 & 0 \end{array} \right) t_2$, alors on obtient (après de fastidieux calculs): \[ t_1 \left( \begin{array}{ccc} \alpha & \gamma & 1 \\ \beta & 1 & 0 \\ 1 & 0 & 0 \end{array} \right) t_2 \in H_K A_\Lambda t_3^0 H_K\] avec 
\[ t_3^0 = \Diag \big( \varpi_D^{-v_D(\alpha)} \alpha, - \varpi_D^{-y+v_D(\alpha)-v_D(\beta)-v_D(\gamma)} \beta \alpha^{-1} \gamma \varpi_D^y , \varpi_D^{-z+v_D(\beta)+v_D(\gamma)} \beta^{-1} \gamma^{-1} \varpi_D^z \big) .\] 
En particulier, on trouve alors \[ \widetilde{\chi}_\rho \Bigg( t_1 \left( \begin{array}{ccc} \alpha & \gamma & 1 \\ \beta & 1 & 0 \\ 1 & 0 & 0 \end{array} \right) t_2 \Bigg) = \widetilde{\chi}_\rho(t_1 t_2) \, \chi( \overline{- \alpha \beta \alpha^{-1} \gamma \beta^{-1} \gamma^{-1}} ) \] et le second facteur du membre de droite peut ne pas valoir $\widetilde{\chi}_\rho(t_1 t_2) \, \chi(-1)$ (lorsque $d$ est strictement supérieur à $1$).

\newpage
$\phantom{,}$
\newpage
\bibliographystyle{amsalpha}
\bibliography{PhD_Main_References}

\providecommand{\bysame}{\leavevmode\hbox to3em{\hrulefill}\thinspace}
\providecommand{\MR}{\relax\ifhmode\unskip\space\fi MR }
% \MRhref is called by the amsart/book/proc definition of \MR.
\providecommand{\MRhref}[2]{%
  \href{http://www.ams.org/mathscinet-getitem?mr=#1}{#2}
}
\providecommand{\href}[2]{#2}
\begin{thebibliography}{Dem11b}

\bibitem[Abd11]{Abd11}
Ramla Abdellatif, \emph{Autour des représentations modulo $p$ des groupes
  réductifs $p$-adiques de rang $1$}, Ph.D. thesis, Université Paris-Sud 11,
  2011.

\bibitem[Abe11]{Abe11}
Noriyuki Abe, \emph{On a classification of irreducible admissible modulo $p$
  representations of a $p$-adic split reductive group}, preprint (2011),
  \url{http://arxiv.org/pdf/1103.2525v3.pdf}.

\bibitem[BB10]{BerBre10}
Laurent Berger and Christophe Breuil, \emph{Sur quelques représentations
  potentiellement cristallines de $\mathrm{GL}_2(\mathbb{Q}_p)$}, Astérisque
  \textbf{330} (2010), 155--211.

\bibitem[Ber10]{Ber10}
Laurent Berger, \emph{On some modular representations of the {B}orel subgroup
  of $\mathrm{GL}_2(\mathbb{Q}_p)$}, Composition Mathematica \textbf{146}
  (2010), 58--80.

\bibitem[BH06]{BusHen06}
Colin~J. Bushnell and Guy Henniart, \emph{The local {L}anglands conjecture for
  $\mathrm{GL}(2)$}, 2006.

\bibitem[BK93]{BusKut93}
Colin~J. Bushnell and Philip~C. Kutzko, \emph{The admissible dual of
  $\mathrm{GL}(n)$ via compact open subgroups}, Annals of Mathematics Studies,
  no. 129, 1993.

\bibitem[BL94]{BarLiv94}
Laure Barthel and Ron Livne, \emph{Irreducible modular representations of
  $\mathrm{GL}_2$ of a local field}, Duke Math. J. \textbf{75} (1994),
  261--292.

\bibitem[BL95]{BarLiv95}
\bysame, \emph{Modular representations of $\mathrm{GL}_2$ of a local field: the
  ordinary, unramified case}, Journal of Number Theory \textbf{55} (1995),
  1--27.

\bibitem[Bou58]{Bou58II8}
Nicolas Bourbaki, \emph{Eléments de mathématique}, vol.~II, ch.~8, 1958.

\bibitem[Bou68]{Bou68VII4}
\bysame, \emph{Eléments de mathématique: {G}roupes et algèbres de {L}ie}, vol.
  VII, ch.~4, 1968.

\bibitem[BP12]{BrePas12}
Christophe Breuil and Vytautas Paskunas, \emph{Towards a modulo $p$ {L}anglands
  correspondence for $\mathrm{GL}_2$}, Memoirs of American Math. Soc.
  \textbf{216} (2012).

\bibitem[Bre03a]{Bre03}
Christophe Breuil, \emph{Sur quelques représentations modulaires et $p$-adiques
  de $\mathrm{GL}_2(\mathbb{Q}_p)$ {I}}, Compositio Math. \textbf{138} (2003),
  165--188.

\bibitem[Bre03b]{Bre03b}
\bysame, \emph{Sur quelques représentations modulaires et $p$-adiques de
  $\mathrm{GL}_2(\mathbb{Q}_p)$ {II}}, J. Inst. Math. Jussieu \textbf{2}
  (2003), 1--36.

\bibitem[Bre04]{Bre04}
\bysame, \emph{Invariant $\mathcal{L}$ et série spéciale $p$-adique}, Ann. Sci.
  Ecole Norm. Sup. \textbf{37} (2004), 559--610.

\bibitem[BT65]{BorTit65}
Armand Borel and Jacques Tits, \emph{Groupes réductifs}, Publications
  mathématiques de l'IHES \textbf{27} (1965).

\bibitem[BT72]{BruTit72}
François Bruhat and Jacques Tits, \emph{Groupes réductifs sur un corps local
  {I}, {D}onnées radicielles valuées}, Publications mathématiques de l'IHES
  \textbf{41} (1972).

\bibitem[BT84]{BruTit84}
\bysame, \emph{Groupes réductifs sur un corps local {II}, {S}chémas en groupes,
  existence d'une donnée radicielle valuée}, Publications mathématiques de
  l'IHES \textbf{60} (1984).

\bibitem[Bum97]{Bum97}
Daniel Bump, \emph{Automorphic forms and representations}, 1997.

\bibitem[Car79]{Car79}
Pierre Cartier, \emph{Representations of $p$-adic groups: a survey},
  Automorphic forms, representations and {L}-functions, 1979.

\bibitem[Car85]{Car85}
Roger~W. Carter, \emph{Finite groups of {L}ie type: conjugacy classes and
  complex characters}, 1985.

\bibitem[CE04]{CabEng04}
Marc Cabanes and Michel Enguehard, \emph{Representation theory of finite
  reductive groups}, 2004.

\bibitem[Col10a]{Col10}
Pierre Colmez, \emph{La série principale unitaire de
  $\mathrm{GL_2(\mathbb{Q}_p)}$}, Astérisque \textbf{330} (2010), 213--262.

\bibitem[Col10b]{Col10c}
\bysame, \emph{Représentations de $\mathrm{GL_2(\mathbb{Q}_p)}$ et
  $(\varphi,{\Gamma})$-modules}, Astérisque \textbf{330} (2010), 281--509.

\bibitem[Col10c]{Col10b}
\bysame, \emph{$(\varphi,{\Gamma})$-modules et représentations du mirabolique
  de $\mathrm{GL_2(\mathbb{Q}_p)}$}, Astérisque \textbf{330} (2010), 61--153.

\bibitem[Cur66]{Cur66}
C.W. Curtis, \emph{The {S}teinberg character of a finite group with a
  $({B},{N})$-pair}, J. Algebra (1966), no.~4, 433--441.

\bibitem[Cur70]{Cur70}
\bysame, \emph{Modular representations of finite groups with split
  $({B},{N})$-pairs}, Seminar on algebraic groups and related finite groups
  (Springer-Verlag, ed.), Lecture Notes in Math., no. 131, 1970, pp.~57--95.

\bibitem[Dem11a]{Dem11}
Michel Demazure, \emph{Groupes réductifs, déploiements, sous-groupes, groupes
  quotients ({E}xposé {XXII})}, Schémas en groupes (SGA3) (Société~Mathématique
  de~France, ed.), Documents Mathématiques, no.~8, 2011, pp.~109--176.

\bibitem[Dem11b]{Dem11b}
\bysame, \emph{Sous-groupes paraboliques des groupes r\' eductifs ({E}xposé
  {XXVI})}, Schémas en groupes (SGA3) (Société~Mathématique de~France, ed.),
  Documents Mathématiques, no.~8, 2011.

\bibitem[Dia07]{Dia07}
Fred Diamond, \emph{A correspondence between representations of local {G}alois
  groups and {L}ie-type groups}, L-functions and Galois representations
  (London~Mathematical Society, ed.), 2007, pp.~187--206.

\bibitem[Eme08]{Eme08}
Matthew Emerton, \emph{On a class of coherent rings, with applications to the
  smooth representation theory of $\mathrm{GL}_2(\mathbb{Q}_p)$ in
  characteristic $p$}, preprint (2008).

\bibitem[Eme10]{Eme10}
\bysame, \emph{Ordinary parts of admissible representations of $p$-adic
  reductive groups {I}, {D}efinition and first properties}, Astérisque
  \textbf{331} (2010), 335--381.

\bibitem[Fon91]{Fon91}
Jean-Marc Fontaine, \emph{Représentations $p$-adiques des corps locaux}, The
  Grothendieck Festschrift \textbf{2} (1991), 249--309.

\bibitem[GK09]{Grk09}
Elmar Grosse-Klönne, \emph{On special representations of $p$-adic reductive
  groups}, preprint (2009),
  \url{http://www.math.hu-berlin.de/~zyska/Grosse-Kloenne/specsubm.pdf}.

\bibitem[Hai09]{Hai09}
Thomas Haines, \emph{Corrigendum : {T}he base change fundamental lemma for
  central elements in parahoric {H}ecke algebras}, webpage (2009),
  \url{http://www2.math.umd.edu/~tjh//fl_corr3.pdf}.

\bibitem[Hen00]{Hen00}
Guy Henniart, \emph{Une preuve simple des conjectures de {L}anglands pour
  $\mathrm{GL}(n)$ sur un corps $p$-adique}, Inventiones Mathematicae
  \textbf{139} (2000), 439--455.

\bibitem[Her11a]{Her11b}
Florian Herzig, \emph{The classification of irreducible admissible mod $p$
  representations of a $p$-adic $\mathrm{GL}_n$}, Inventiones Mathematicae
  \textbf{186} (2011), no.~2, 373--434.

\bibitem[Her11b]{Her11}
\bysame, \emph{A {S}atake isomorphism modulo $p$}, Compositio Math.
  \textbf{147} (2011), 263--283.

\bibitem[HKP10]{HaiKotPra09}
Thomas Haines, Robert Kottwitz, and Amrithanshu Prasad, \emph{Iwahori-{H}ecke
  algebras}, Journal of the Ramanujan Math. Soc. \textbf{25} (2010), no.~2,
  113--145.

\bibitem[HR08]{HaiRap08}
Thomas Haines and Michael Rapoport, \emph{Appendix: {O}n parahoric subgroups},
  Advances in Mathematics \textbf{219} (2008), 188--198.

\bibitem[HR10]{HaiRos10}
Thomas Haines and Sean Rostami, \emph{The {S}atake isomorphism for special
  maximal parahoric {H}ecke algebras}, Electronic Journal of Representation
  Theory \textbf{14} (2010), 264--284.

\bibitem[HT01]{HarTay01}
Michael Harris and Richard Taylor, \emph{On the geometry and cohomology of some
  simple {S}himura varieties}, no. 151, 2001.

\bibitem[Hu10]{Hu10}
Yongquan Hu, \emph{Sur quelques représentations supersingulières de
  $\mathrm{GL}_2(\mathbb{Q}_{p^f})$}, Journal of Algebra \textbf{324} (2010),
  1577--1615.

\bibitem[Hum92]{Hum92}
James~E. Humphreys, \emph{Reflection groups and {C}oxeter groups}, 1992.

\bibitem[Hum06]{Hum06}
\bysame, \emph{Modular representations of finite groups of {L}ie type}, 2006.

\bibitem[HV11]{HenVig11}
Guy Henniart and Marie-France Vignéras, \emph{A {S}atake isomorphism for
  representations modulo $p$ of reductive groups over local fields}, preprint
  (2011),
  \url{http://www.math.jussieu.fr/~vigneras/satake_isomorphism-07032012.pdf}.

\bibitem[HV12]{HenVig11b}
\bysame, \emph{Comparison of compact induction with parabolic induction},
  Pacific Journal of Mathematics \textbf{260} (2012), no.~2, 457--495.

\bibitem[Jan87]{Jan87}
Jens~Carsten Jantzen, \emph{Representations of algebraic groups}, 1987.

\bibitem[Kaz86]{Kaz86}
David Kazhdan, \emph{Representations of groups over close local fields},
  Journal d'analyse mathématique \textbf{47} (1986), 175--179.

\bibitem[Kis10]{Kis10}
Mark Kisin, \emph{Deformations of ${G}_{{\mathbb{Q}}_p}$ and
  $\mathrm{GL}_2(\mathbb{Q}_p)$ representations}, Astérisque \textbf{330}
  (2010), 529--542.

\bibitem[KL79]{KazLus79}
David Kazhdan and George Lusztig, \emph{Representations of {C}oxeter groups and
  {H}ecke algebras}, Inventiones Math. (1979), no.~53, 165--184.

\bibitem[Kot97]{Kot97}
Robert~E. Kottwitz, \emph{Isocrystals with additional structure {II}},
  Compositio Math. \textbf{109} (1997), 255--339.

\bibitem[Koz12]{Koz12}
Karol Koziol, \emph{A classification of the irreducible mod $p$ representations
  of ${U}(1,1)(\mathbb{Q}_{p^2}/\mathbb{Q}_p)$}, preprint (2012).

\bibitem[KX12]{KozXu12}
Karol Koziol and Peng Xu, \emph{Hecke modules and supersingular representations
  of ${U}(2,1)$}, preprint (2012).

\bibitem[Lan]{Lan67}
Robert Langlands, \emph{Letter to {A}ndré {W}eil}, \newline
  www.sunsite.ubc.ca/DigitalMathArchive/Langlands/pdf/langlands-weil-ps.pdf,
  accessed 14.11.2012.

\bibitem[Lan00]{Lan00}
Erasmus Landvogt, \emph{Some functorial properties of the {B}ruhat-{T}its
  building}, J. reine. angew. Math. \textbf{518} (2000), 213--241.

\bibitem[MN43]{MatNak43}
Yozô Matsushima and Tadasi Nakayama, \emph{Über die multiplikative {G}ruppe
  einer $p$-adischen {D}ivisionalgebra}, Proc. Imp. Acad. \textbf{19} (1943),
  no.~10, 622--628.

\bibitem[NR03]{NelRam03}
K.~Nelsen and A.~Ram, \emph{Kostka-{F}oulkes polynomials and {M}acdonald
  spherical functions}, Surveys in Combinatorics 2003 (Cambridge~University
  Press, ed.), 2003, pp.~325--370.

\bibitem[Pas04]{Pas04}
Vytautas Paskunas, \emph{Coefficient systems and supersingular representations
  of $\mathrm{GL}_2({F})$}, Bulletin de la S.M.F., vol.~99, 2004.

\bibitem[Pas13]{Pas10}
\bysame, \emph{The image of {C}olmez's {M}ontreal functor}, Publications
  mathématiques de l'IHES (2013).

\bibitem[Pie82]{Pie82}
Richard~Scott Pierce, \emph{Associative algebras}, Graduate Texts in
  Mathematics, no.~88, 1982.

\bibitem[Rei75]{Rei75}
Irving Reiner, \emph{Maximal orders}, 1975.

\bibitem[Ric69]{Ric69}
Forrest Richen, \emph{Modular representations of split {BN} pairs}, Trans. AMS
  \textbf{140} (1969), 435--460.

\bibitem[Sch12]{Sch12}
Benjamin Schraen, \emph{Sur la présentation des représentations
  supersingulières de $\mathrm{GL}_2({F})$}, preprint (2012).

\bibitem[Ser67]{Ser67}
Jean-Pierre Serre, \emph{Local class field theory}, Algebraic number theory
  (Academic~Press London, ed.), 1967.

\bibitem[Ser77]{Ser77}
\bysame, \emph{Arbres, amalgames, $\mathrm{SL}_2$}, Astérisque, no.~46, 1977.

\bibitem[Ser98]{Ser98}
\bysame, \emph{Représentations linéaires des groupes finis}, 1998.

\bibitem[Ste51]{Ste51}
Robert Steinberg, \emph{A geometric approach to the representations of the full
  linear group over a {G}alois field}, Transactions of the A.M.S. \textbf{71}
  (1951), no.~2, 274--282.

\bibitem[Séc04]{Sec04}
Vincent Sécherre, \emph{Représentations lisses de $\mathrm{GL}(m,{D})$ {I}:
  caractères simples}, Bull. Soc. Math. France \textbf{132} (2004), 327--396.

\bibitem[Vig96]{Vig96}
Marie-France Vignéras, \emph{Représentations $l$-modulaires d'un groupe
  réductif $p$-adique avec $l \neq p$}, Progr. Math., no. 137, 1996.

\bibitem[Vig01]{Vig01}
\bysame, \emph{Correspondance de {L}anglands semi-simple pour
  $\mathrm{GL}(n,{F})$ modulo $\ell \neq p$}, Inventiones Mathematicae
  \textbf{144} (2001), 197--223.

\bibitem[Vig04]{Vig04}
\bysame, \emph{Representations modulo $p$ of the $p$-adic group
  $\mathrm{GL}(2,{F})$}, Comp. Math. (2004), no.~140, 333--358.

\bibitem[Vig08]{Vig08}
\bysame, \emph{Série principale modulo $p$ de groupes réductifs $p$-adiques},
  GAFA \textbf{17} (2008).

\bibitem[Vig11]{Vig07}
\bysame, \emph{Représentations $p$-adiques de torsion admissibles}, Number
  Theory, Analysis and Geometry: In memory of Serge Lang (Springer, ed.), 2011.

\bibitem[Vig12]{Vig12}
\bysame, \emph{Emerton's ordinary parts in positive characteristic},
  communication personnelle (2012).

\end{thebibliography}

\end{document}